\documentclass{article}
\usepackage[utf8]{inputenc}
\usepackage{enumerate}
\usepackage[hang,flushmargin]{footmisc}
\usepackage{amsmath}
\usepackage{amsfonts}
\usepackage{amssymb}
\usepackage{amsthm}
\usepackage{xcolor} 
\usepackage{graphicx} 
\usepackage{subfigure}
\usepackage{cite}
\usepackage{paralist}
\usepackage[acronym,toc]{glossaries}

\usepackage{algorithm}
\usepackage{algpseudocode}

\usepackage[colorlinks=false]{hyperref}

\title{Geodesic complexity of the octahedron, and an algorithm for cut loci on convex polyhedra}
\author{
Florian Frick\footnote{supported by NSF CAREER Grant DMS 2042428}
\\\href{mailto:frick@cmu.edu}{frick@cmu.edu}
\and
Pranav Rajbhandari\footnote{primary author}
\\\href{mailto:prajbhan@cs.cmu.edu}{prajbhan@cs.cmu.edu}
}
\date{September 10, 2025}

\begin{document}

\def\AFSOC{Assume for the sake of contradiction }
\def\aFSOC{assume for the sake of contradiction }
\def\WLOG{Without loss of generality }
\def\wLOG{without loss of generality }

\newcommand{\proofclose}[1]{$\blacksquare$}

\newcommand{\lrq}[1]{\lq#1\rq}

\newcommand{\meminisection}[1]{\textit{\textbf{#1}}}
\newcommand{\meminisubsection}[1]{\textit{#1}}

\newcommand{\relint}[1]{\text{relint}^*\left(#1\right)}

\newcommand{\img}[1]{\text{img}\left(#1\right)}
\def\Sone{S^1}

\def\funcprodsymb{\boxtimes}

\newacronym{GMPR}{GMPR}{Geodesic Motion Planning Rule}

\newtheorem{define}{Definition}
\numberwithin{define}{section}
\newcommand{\defineautorefname}{Definition}

\newtheorem{lemma}{Lemma}
\numberwithin{lemma}{section}
\newcommand{\lemmaautorefname}{Lemma}

\newtheorem{theorem}{Theorem}
\numberwithin{theorem}{section}

\newtheorem{corr}{Corollary}
\numberwithin{corr}{section}
\newcommand{\corrautorefname}{Corollary}

\newtheorem{conj}{Conjecture}
\numberwithin{conj}{section}
\newcommand{\conjautorefname}{Conjecture}

\newcommand{\algorithmautorefname}{Algorithm}
\renewcommand{\sectionautorefname}{Section}
\renewcommand{\subsectionautorefname}{Section}
\renewcommand{\subsubsectionautorefname}{Section}

\newcommand{\imgwidth}[0]{0 pt}
\newcommand{\imgheight}[0]{0 pt}

\maketitle

\begin{abstract}
    The geodesic complexity of a length space $X$ quantifies the required number of case distinctions to continuously choose a shortest path connecting any given start and end point.
    We prove a local lower bound for the geodesic complexity of $X$ obtained by embedding simplices into $X\times X$.
    We additionally create and prove correctness of an algorithm to find cut loci on surfaces of convex polyhedra, as the structure of a space's cut loci is related to its geodesic complexity.
    We use these techniques to prove the geodesic complexity of the octahedron is four.
    Our method is inspired by earlier work of Recio-Mitter and Davis, and thus recovers their results on the geodesic complexity of the $n$-torus and the tetrahedron, respectively.
\end{abstract}

\section{Introduction}\label{section:intro}

In our work, we inspect the geodesic complexity of surfaces of polyhedra, taking inspiration from Davis's results on the tetrahedron \cite{tetra} and cube \cite{cube}.
We prove \autoref{thm:quite_simplex}, which gives a general criterion for lower bounding the geodesic complexity of a length space $X$, witnessed by embedding simplices into $X\times X$.
We additionally give an algorithm to calculate the cut locus of a point on a polyhedral surface.
The methods refined and developed here in particular yield:
\begin{theorem}
    \label{thm:octahedron_GC}
    The geodesic complexity of the octahedron is four.
\end{theorem}

\subsection{Topological Complexity}
Topological complexity was introduced by Farber \cite{top_complex} as a way to describe the difficulty of creating a continuous motion plan over a space $X$.

\begin{define}
    The \textbf{free path space} $PX$ is the space of paths $C([0,1], X)$ imbued with the compact-open topology.
    This space comes with a \textbf{projection function} $\pi\colon PX\to X\times X$ via $\pi(\gamma)=
    (\gamma(0),\gamma(1))$.
\end{define}
A natural question to ask is if there is a local right inverse to $\pi$.
Such a function $\sigma$ would describe a continuous mapping from $X\times X$ (the space of start and end points in $X$, equipped with the product topology) into $PX$ (the space of paths in $X$) such that $\sigma(x,y)$ is a path from $x$ to $y$.
\begin{define}
    For $U\subseteq X\times X$, a \textbf{motion planning rule} is a map $\sigma\colon U\to PX$ such that $\pi\circ\sigma=\text{id}$.
\end{define}

A global motion planning rule on a space imposes strong topological constraints: The space must be contractible. 
Thus, topological complexity is defined through covering $X\times X$ with open sets that each have a local motion planning rule:
\begin{define}
    The \textbf{topological complexity} $TC(X)$ is the minimal $k$ such that there exists open cover $U_1\cup \dots\cup U_k=X\times X$ where each $U_i$ has a motion planning rule $\sigma_i\colon U_i\to PX$.
\end{define}

\subsection{Geodesic Complexity}

For a metric space $(X,d)$, the length of a path $\gamma\colon [0,1]\to X$ is
\[L(\gamma)=
\sup\limits_{0=t_0\leq \dots<t_N=1}
\left[\sum\limits_{i=1}^Nd(\gamma(t_{i-1}),\gamma(t_i))\right].\]

Since we have a notion of length, a natural extension of topological complexity is to require the paths in each motion planning rule to be shortest paths.
This was introduced by Recio-Mitter~\cite{geo_complex}, and is constructed analogously to topological complexity.

\begin{define}
A path $\gamma$ is a \textbf{geodesic} if there
exists $\lambda\in\mathbb R$ such that $\forall t,t'$ with $0\leq t<t'\leq 1$, 
\begin{align*}
d(\gamma(t),\gamma(t'))=\lambda|t-t'|.
\end{align*}

\end{define}

This definition of a geodesic is \lrq{the shortest path with constant speed}.
If a path has minimal length, it can always be reparameterized into a geodesic.
Additionally, each path with minimal length has a unique reparameterization into a geodesic.

We will only consider \textbf{length spaces} \cite{metic_non_pos}, where the metric $d(x,y)$ is equal to the infinimal path length of all paths between $x$ and $y$:
$d(x,y)=\inf\limits_{\gamma\in \pi^{-1}(x,y)}L(\gamma)$.

\begin{define}
    Let $GX\subseteq PX$ be the subset containing only geodesics. 
    We define a projection function $\pi_{GX}:=\pi|_{GX}$ via the restriction of the projection function for $PX$.
\end{define}
\begin{define}
\label{def:path_distances}
    Since $[0,1]$ is compact and $(X,d)$ is a metric space, $GX$ is metrizable \cite{cont_func}.
    The metric we will use is 
    $d_{GX}(f,g)=
    \sup\limits_{t\in[0,1]}d(f(t),g(t))$ 
    for $f,g\in GX$.
    
\end{define}

\begin{define}
For $E\subseteq X\times X$, a \textbf{\gls{GMPR}} is a map $\sigma\colon E\to GX$ such that $\pi_{GX}\circ\sigma=\text{id}$.
We say a set of geodesics $\mathcal G$ \textbf{extends} to a \gls{GMPR} $\sigma$ if $\forall G\in \mathcal G$, $\sigma(\pi_{GX}(G))=G$.
\end{define}

\begin{define}
The \textbf{geodesic complexity} $GC(X)$ of a length space $X$ is the smallest $k$ for which there is a partition into $k+1$ locally compact sets 
$\bigsqcup\limits_{i=0}^kE_i=X\times X$ with \gls{GMPR}s $\sigma_i\colon E_i\to GX$.
We call the collection of $\sigma_i$ a \textbf{geodesic motion planner}.
\end{define}

We switch to covering $X\times X$ with locally compact sets out of necessity, as continuing to use an open cover would result in simple spaces such as $\Sone$ having an undefined geodesic complexity (see \cite[Remark 3.17]{geo_complex}).
Additionally we add one to the number of sets due to convention.
To reduce confusion, we will clarify the number of sets along with any statement about geodesic complexity.

\subsection{Document Structure}
In \autoref{section:intro}, we introduce the concepts of topological and geodesic complexity.
In \autoref{section:methods}, we develop methods we use to bound geodesic complexity.
In \autoref{section:examples}, we apply our methods to produce novel results, as well as to reprove existing bounds in the literature.

\section{Methods}\label{section:methods}
\subsection{Geodesic complexity? It's actually quite simplex}

We will prove a local obstruction that bounds geodesic complexity of a space $X$ from below by the dimension of certain simplices we can embed into $X\times X$.
This is inspired by lower bound arguments of Davis~\cite[Thm.~4.1]{tetra,cube} and Recio-Mitter~\cite[Thm.~4.4]{geo_complex}, whose proofs utilize sequences, sequences of sequences, and so on of elements in $X\times X$.
Their proofs construct each successive iteration so that it cannot be in a \gls{GMPR} with previous points.
This shows that the geodesic complexity of $X$ must be at least the number of iterations in this procedure. We adapt this by replacing $n$-dimensional sequences with embeddings of~$\Delta_n$ into $X\times X$.
Our hypotheses will let us use the same recursive sequence argument on these simplices, packaging all the machinery to show a lower bound from a collection of embeddings.
Moving forward, we utilize this to lower bound geodesic complexity without the need to reference iterations of sequences.

\begin{define}
The \textbf{affine hull} $\displaystyle \operatorname {aff} (S)$ of a set $S\subseteq \mathbb R^n$ is the space of all affine combinations of $S$:
\begin{align*}
{\displaystyle \operatorname {aff} (S)=\left\{\sum _{i=1}^{k}\alpha _{i}x_{i}\,{\Bigg |}\,k>0,\,x_{i}\in S,\,\alpha _{i}\in \mathbb {R} ,\,\sum _{i=1}^{k}\alpha _{i}=1\right\}}.
\end{align*}
\end{define}

\begin{define}
The \textbf{relative interior} $\text{relint}(S)$ for $S\subseteq \mathbb R^n$ is the interior of $S$ in the affine hull of $S$:
\begin{align*}
    {\displaystyle \operatorname {relint} (S)=
    \{x\in S:{\text{ there exists }}\epsilon >0{\text{ such that }}N_{\epsilon }(x)\cap \operatorname {aff} (S)\subseteq S\}}.
\end{align*}
If we embed set $S\subseteq \mathbb R^n$ via $f\colon S\hookrightarrow X$, we denote $\relint{f}=f(\text{relint}(S))$.

\end{define}

\begin{theorem}
\label{thm:quite_simplex}

For integers $0\leq k_0<k$, if for each $i\in \{k_0, k_0+1, \dots, k\}$ there exists a non-empty collection $\mathcal F_{i}$ of embeddings into $X\times X$ such that:

\begin{compactenum}[(a)]
    \item
    Each $\mathcal F_i$ consists of embeddings of $i$-dimensional simplices $\Delta_i\hookrightarrow X\times X$.

    \item
    For $F\in \mathcal F_i$, any geodesic in $\pi_{GX}^{-1}(\relint{F})$ extends to a \gls{GMPR} on all of $\img{F}$.

        \item
    With $k_0\leq i<k$, for any $F\in \mathcal F_i$ and \gls{GMPR} $\Gamma$ on $\img{F}$, there is an $F'\in\mathcal F_{i+1}$ such that:
    \begin{compactenum}[(i)]
        \item
        The map $F$ is the restriction of $F'$ to a face of $\Delta_{i+1}$. 
        We describe this relation as $F$ being a face of $F'$.
        \item
        Each geodesic in $\Gamma(\img{F})$ is separated from any geodesics in $\pi_{GX}^{-1}(\relint{F'})$ by an open set, that is,
        $\Gamma(\img{F})$ does not intersect $\overline{\pi_{GX}^{-1}(\relint{F'})}$, the closure of all geodesics on $\relint{F'}$.
        
                \item
        There are a finite number of \gls{GMPR}s on $\img{F'}$.
    \end{compactenum}
\end{compactenum}

Then $GC(X)\geq k-k_0$. 
(This space requires at least $k-k_0+1$ sets in any geodesic motion planner).

\end{theorem}

We consider this a local obstruction since most of these hypotheses will hold when decreasing the size of the simplex embeddings.
This is trivially true for hypothesis (a) and (b).
For hypothesis (c), if we take care in the way we decrease the simplex embedding size, hypotheses (i) and (ii) will still hold.
We display an example on the 2-torus in \autoref{subsec:ntorus}, where the hypotheses hold at arbitrarily small scale.

\textbf{Proof:} 
\AFSOC that $GC(X)<k-k_0$.
Then partition $X\times X$ into $k-k_0$ locally compact sets $E_{k_0},\dots,E_{k-1}$ such that each $E_i$ has a \gls{GMPR} $\sigma_i$.
We will proceed by inductively constructing $k-k_0+1$ collections of points, where each collection belongs to a unique $E_i$.

\meminisubsection{Base Case (Step $k_0$):}
Pick an embedding $F_{k_0}\in\mathcal F_{k_0}$, and let $(p_*,q_*)\in \relint{F_{k_0}}$.
\WLOG let $E_{k_0}$ be the set that contains $(p_*,q_*)$.
Let $\Gamma_{k_0}$ be a \gls{GMPR} on $\img{F_{k_0}}$ that extends $\sigma_{k_0}(p_*,q_*)$.
This is possible from assumption (b).

\meminisubsection{Induction Hypothesis (Step $i$ with $k_0< i\leq k$):}
Assume for each $j\in \{k_0+1,k_0+2,\dots, i-1\}$, we have chosen an embedding $F_j\in\mathcal F_j$ and sequences $(p_{\mathbf{s'}},q_{\mathbf{s'}})_{\mathbf{s'}\in \mathbb N^{j-k_0}}\subseteq\relint{F_j}$ such that:
\begin{compactenum}[$(1)_{\text{emb}}^{(j)}$:]
    \item
    $F_{j-1}$ is a face of $F_j$.
    \item
    For all $j'\in \{k_0,\dots,j-1\}$, any geodesic in $\Gamma_{j'}(\img{F_{j'}})$ is some distance $\delta$ from all geodesics in $\pi_{GX}^{-1}(\relint{F_j})$.
    Equivalently, $\Gamma_{j'}(\img{F_{j'}})\cap\overline{\pi_{GX}^{-1}(\relint{F_j})}$ is empty.
        \item
    Any geodesic in $\pi_{GX}^{-1}(\relint{F_j})$ extends to a \gls{GMPR} on $\img{F_j}$.
    \item
    There are a finite number of \gls{GMPR}s on $\img{F_j}$.
\end{compactenum}

\begin{compactenum}[$(1)_{\text{seq}}^{(j)}$]
    \item
    $d((p_{\mathbf{s}},q_{\mathbf{s}}),(p_{\mathbf{s},n},q_{\mathbf{s},n}))\leq \frac1n$ for $\mathbf{s}\in \mathbb N^{(j-k_0)-1}$.    
    In the 1-dimensional case, $d((p_*,q_*),(p_{n},q_{n}))\leq \frac1n$.
    \item
    The points $\left\{(p_{\mathbf{s}'},q_{\mathbf{s}'}):\mathbf{s}'\in \mathbb N^{j-k_0}\right\}\subseteq E_j$.
        \item
    There is a \gls{GMPR} $\Gamma_{j}$ on $\img{F_j}$ that extends $\left\{\sigma_j(p_{\mathbf{s}'},q_{\mathbf{s}'}): \mathbf{s}'\in \mathbb N^{j-k_0}\right\}$.
    \item
    All of these properties still hold under taking subsequences.
\end{compactenum}

This hypothesis is trivially true for $i=k_0+1$, as there is no integer $j$ such that $k_0<j<k_0+1$.

\meminisubsection{Induction Step (Step $i$ with $k_0< i\leq k$):}
From assumptions (c) and (b), we can find an embedding $F_i\in \mathcal F_i$ such that:
\begin{compactenum}[$(1)_{\text{emb}}^{(i)*}$:]
    \item
    $F_{i-1}$ is a face of $F_i$.
    \item
    For $j=i-1$, the intersection $\Gamma_{j}(\img{F_{j}})\cap\overline{\pi_{GX}^{-1}(\relint{F_i})}$ is empty.
    \item
    Any geodesic in $\pi_{GX}^{-1}(\relint{F_i})$ extends to a \gls{GMPR} on $\img{F_i}$.
    \item
    There are a finite number of \gls{GMPR}s on $\img{F_i}$.
\end{compactenum}

To show $(1)_{\text{emb}}^{(i)},\dots,(4)_{\text{emb}}^{(i)}$, it is enough to show that 
$(2)_{\text{emb}}^{(i)*}$ holds for all $j\in \{k_0,\dots,i-1\}$.

\AFSOC that this were not true for some $k_0\leq j<i-1$ and geodesic $G\in \Gamma_{j}(\img{F_{j}})\cap\overline{\pi_{GX}^{-1}(\relint{F_i})}$.
Then $G=\Gamma_{j}(p,q)$ for $(p,q)\in \img{F_{j}}$, and we can find for each $\delta\in(\frac1n)_n$, a geodesic in $\pi_{GX}^{-1}(\relint{F_i})$ that is $\delta$-close to$~G$.
By $(3)_{\text{emb}}^{(i)*}$, each of these geodesics extend to a \gls{GMPR} across $\img{F_i}$.
By $(4)_{\text{emb}}^{(i)*}$ and pigeonhole principle, we can find a \gls{GMPR} $\Gamma$ on $\img{F_i}$ that agrees with an infinite subsequence of these geodesics.
Then for any distance $\delta$, there exists $(p',q')\in \relint{F_i}$ such that $\Gamma(p',q')$ is $\delta$-close to$~G$.
By \autoref{def:path_distances}, this implies $(p,q)$ and $(p',q')$ are $(2\delta)$-close.
Thus, by choosing $\delta\in (\frac1n)_n$, we find a sequence $(p'_n,q'_n)_n$ that approaches $(p,q)$ while $\Gamma(p'_n,q'_n)_n$ approaches $G$.
Since $GX$ is Hausdorff, this shows $\Gamma(p,q)=G$.
However, since $\Gamma$ is a \gls{GMPR} restricted to $\img{F_{j+1}}$, we contradict $(2)_{\text{emb}}^{(j+1)}$, which asserts $G$ is some distance from all geodesics in $\pi_{GX}^{-1}(\relint{F_{j+1}})$, in particular, those in $\Gamma(\relint{F_{j+1}})$.

Thus, we establish $(1)_{\text{emb}}^{(i)},\dots,(4)_{\text{emb}}^{(i)}$.
Next, for each $\mathbf{s}\in\mathbb N^{(i-k_0)-1}$, we can create a sequence\footnote{In the 1-dimensional case, this is $(p_n,q_n)_n$} 
\\$(p_{\mathbf{s},n},q_{\mathbf{s},n})_n\subseteq \relint{F_i}$ such that:
\begin{compactenum}[$(1)_{\text{seq}}^{(i)*}$:]
    \item
    $d((p_{\mathbf{s}},q_{\mathbf{s}}),(p_{\mathbf{s},n},q_{\mathbf{s},n}))\leq \frac1n$.
    
    This is possible, as $(p_{\mathbf{s}},q_{\mathbf{s}})\in \relint{F_{i-1}}$, and $\relint{F_{i-1}}\subseteq \overline{\relint{F_i}}$ by $(1)^{(i)}_{\text{emb}}$.
    
    \item
    The points $\{(p_{\mathbf{s}'},q_{\mathbf{s}'}):\mathbf{s}'\in \mathbb N^{i-k_0}\}\subseteq E_\ell$ for some $E_\ell$.
        We show this is possible through induction:

    First, each point is trivially in a single set, showing our base case. 
    Fix some $\mathbf{t}\in \mathbb N^{j-k_0}$ for $k_0\leq j<i$. 
    We inductively assume for each $n\in\mathbb N$, $\{(p_{\mathbf{t},n,\mathbf{t}'},q_{\mathbf{t},n,\mathbf{t}'}):(\mathbf{t},n,\mathbf{t}')\in \mathbb N^{i-k_0}\}\subseteq E_{\ell_n}$.
    However, since there are finitely many choices for $E_{\ell_n}$, one of these sets (say $E_\ell$) occurs infinitely many times in the sequence $(E_{\ell_n})_{n}$.
    This corresponds to an infinite set $M\subseteq \mathbb N$ such that
    $\{(p_{\mathbf{t},m,\mathbf{t}'},q_{\mathbf{t},m,\mathbf{t}'}):(\mathbf{t},m,\mathbf{t}')\in \mathbb N^{i-k_0}, m\in M\}\subseteq E_{\ell}$. 
    By $(4)_{\text{seq}}^{(j+1)}$, we can assume \wLOG that this subsequence $(p_{\mathbf{t},m},q_{\mathbf{t},m})_{m\in M}$ is the whole sequence $(p_{\mathbf{t},n},q_{\mathbf{t},n})_n$.
    
    \item
    There is a \gls{GMPR} $\Gamma_{i}$ on $\img{F_i}$ that extends $\left\{\sigma_{\ell}(p_{\mathbf{s}'},q_{\mathbf{s}'}):\mathbf{s}'\in \mathbb N^{i-k_0}\right\}$.
    
    This can be proved possible similarly to $(2)_{\text{seq}}^{(i)*}$. 
    The key observation is that using $(3)_{\text{emb}}^{(i)}$, we may associate each $\sigma_{\ell}(p_{\mathbf{s}'},q_{\mathbf{s}'})$ (for $\mathbf{s}'\in \mathbb N^{i-k_0}$) with a \gls{GMPR} across $\img{F_i}$.
    By $(4)_{\text{emb}}^{(i)}$, there are finitely many of these, which allows us to use the same machinery as before.
    
    \item
    All of these properties still hold under taking subsequences.
\end{compactenum}
To show $(1)_{\text{seq}}^{(i)},\dots,(4)_{\text{seq}}^{(i)}$, it is enough to show that $E_\ell\neq E_j$ for $j<i$.

\AFSOC that $E_\ell=E_j$ for some $j<i$.
Then take some $(p_{\mathbf{r}},q_{\mathbf{r}})$ for $\mathbf{r}\in\mathbb N^{j-k_0}$ (take $(p_{\mathbf{r}},q_{\mathbf{r}})=(p_*,q_*)$ if $j=k_0$).
Consider the sequence $(p_{\mathbf{r},n,\dots,n},q_{\mathbf{r},n,\dots,n})_n$ (where $(\mathbf{r},n,\dots,n)\in \mathbb N^{i-k_0}$).
By $(1)^{(i)*}_{\text{seq}}, 
(1)^{(i-1)}_{\text{seq}},\dots,
(1)^{(j+1)}_{\text{seq}}$, we have that $(p_{\mathbf{r},n,\dots,n},q_{\mathbf{r},n,\dots,n})_n$ approaches $(p_\mathbf{r},q_\mathbf{r})$.
Through continuity of $\sigma_j=\sigma_\ell$, we construct a sequence $\sigma_j(p_{\mathbf{r},n,\dots,n},q_{\mathbf{r},n,\dots,n})_n$ of geodesics approaching $\sigma_j(p_{\mathbf{r}},q_{\mathbf{r}})$.
However, this is a contradiction with $(2)_{\text{emb}}^{(i)}$, which claims since $(p_{\mathbf{r},n,\dots,n},q_{\mathbf{r},n,\dots,n})_n\subseteq \relint{F_i}$, there is an open set separating $\sigma_j(p_{\mathbf{r}},q_{\mathbf{r}})=\Gamma_j(p_{\mathbf{r}},q_{\mathbf{r}})$ from all $\sigma_j(p_{\mathbf{r},n,\dots,n},q_{\mathbf{r},n,\dots,n})$.
Thus, $E_\ell\neq E_j$ for all $j<i$, and \wLOG we can set $E_\ell=E_i$.

We can continue until $i=k$, creating distinct sets $E_{k_0},\dots,E_{k}$ and contradicting the number of sets assumed.
Therefore, $GC(X)\geq k-k_0$, needing at least $k-k_0+1$ sets in any geodesic motion planner.
\proofclose{\autoref{thm:quite_simplex}}

We illustrate the mechanism behind the proof with an example where $X$ is the 2-torus, showing $GC(X)\geq 2$.
For collections of embeddings, let $\mathcal F_0=\{\Delta_0\mapsto (p_*,q_*)\}$ for some point $p_*$ and its antipode $q_*$.
Let $\varepsilon > 0$ be small. 
There are four axis-aligned lines on $X$ of length $\varepsilon$ that end in $q_*$ (displayed as solid gray lines in \autoref{fig:proof-walkthrough}(a)).
Let $\mathcal F_1$ consist of four maps $\Delta_1\hookrightarrow X\times X$ with their first factor fixed at $p_*$ and their second factor one of those four lines.
Similarly, there are four axis-aligned squares on $X$ with side length $\varepsilon$ with a vertex as $q_*$ (displayed as gray squares in \autoref{fig:proof-walkthrough}(a)).
Let $\mathcal F_2$ consist of four maps $\Delta_2\hookrightarrow X\times X$ with their first factor fixed at $p_*$ and their second factor one of these four squares.
With this construction, for $F_1\in \mathcal F_1$, $\img{F_1}$ can be written $\{p_*\}\times F_1'$ for a line $F_1'$.
Similarly, for $F_2\in \mathcal F_2$, $\img{F_2}=\{p_*\}\times F_2'$ for a square $F_2'$.
We display $F_1'$ and $F_2'$ for particular $F_1\in \mathcal F_1$ and $F_2\in \mathcal F_2$ in \autoref{fig:proof-walkthrough}(b) and (c) respectively. 
In \autoref{subsec:ntorus}, we describe this construction in more detail and show it satisfies the hypotheses of \autoref{thm:quite_simplex}.

\renewcommand{\imgwidth}{.3\linewidth}
\begin{figure}[ht!]
\centering
\subfigure[Base case $(p_*,q_*)$ and possible $\sigma_0$]{
\includegraphics[width=\imgwidth]{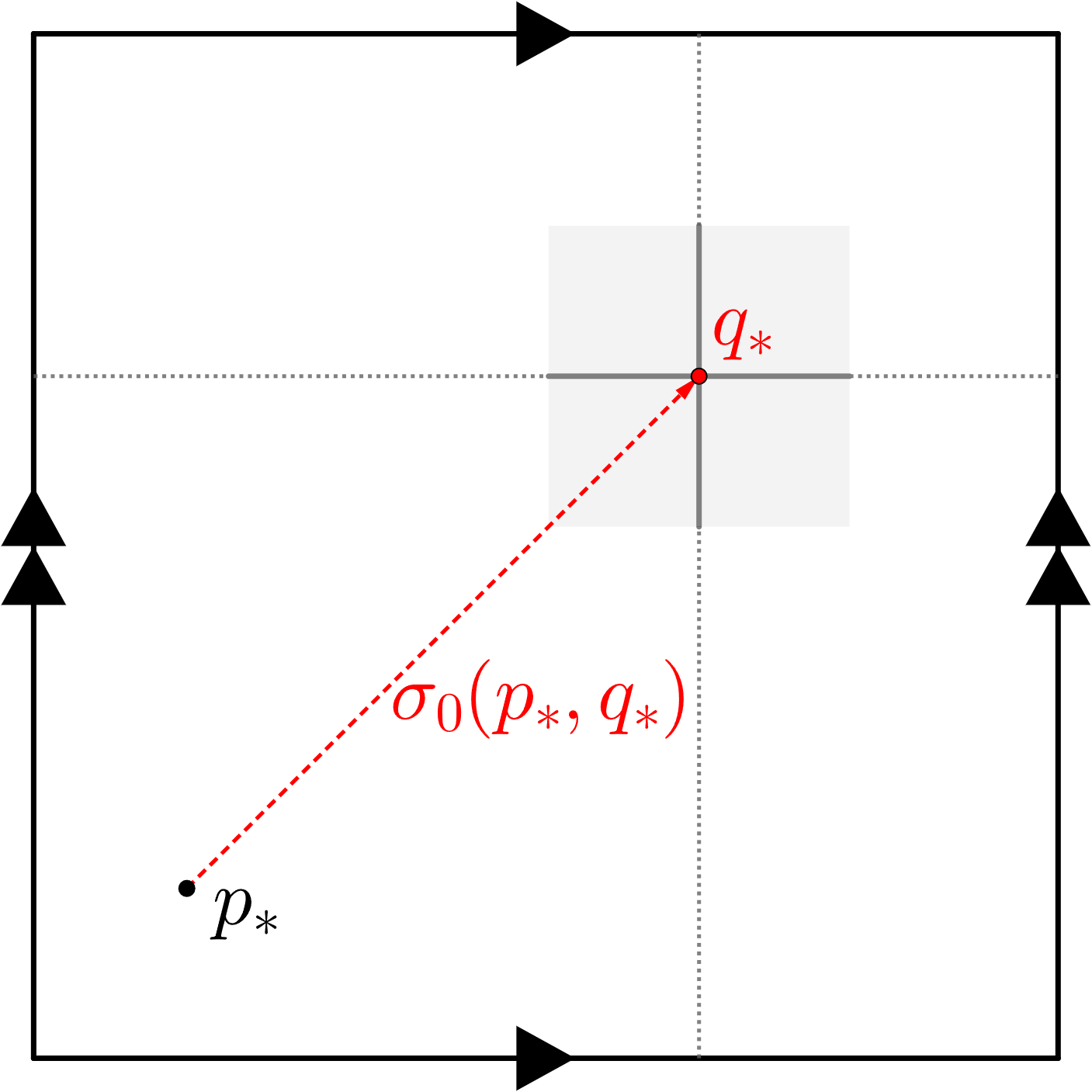}
}
\quad
\subfigure[Choice of $F_1$, sequence $(p_n,q_n)_n$, and possible $\sigma_1$]{
\includegraphics[width=\imgwidth]{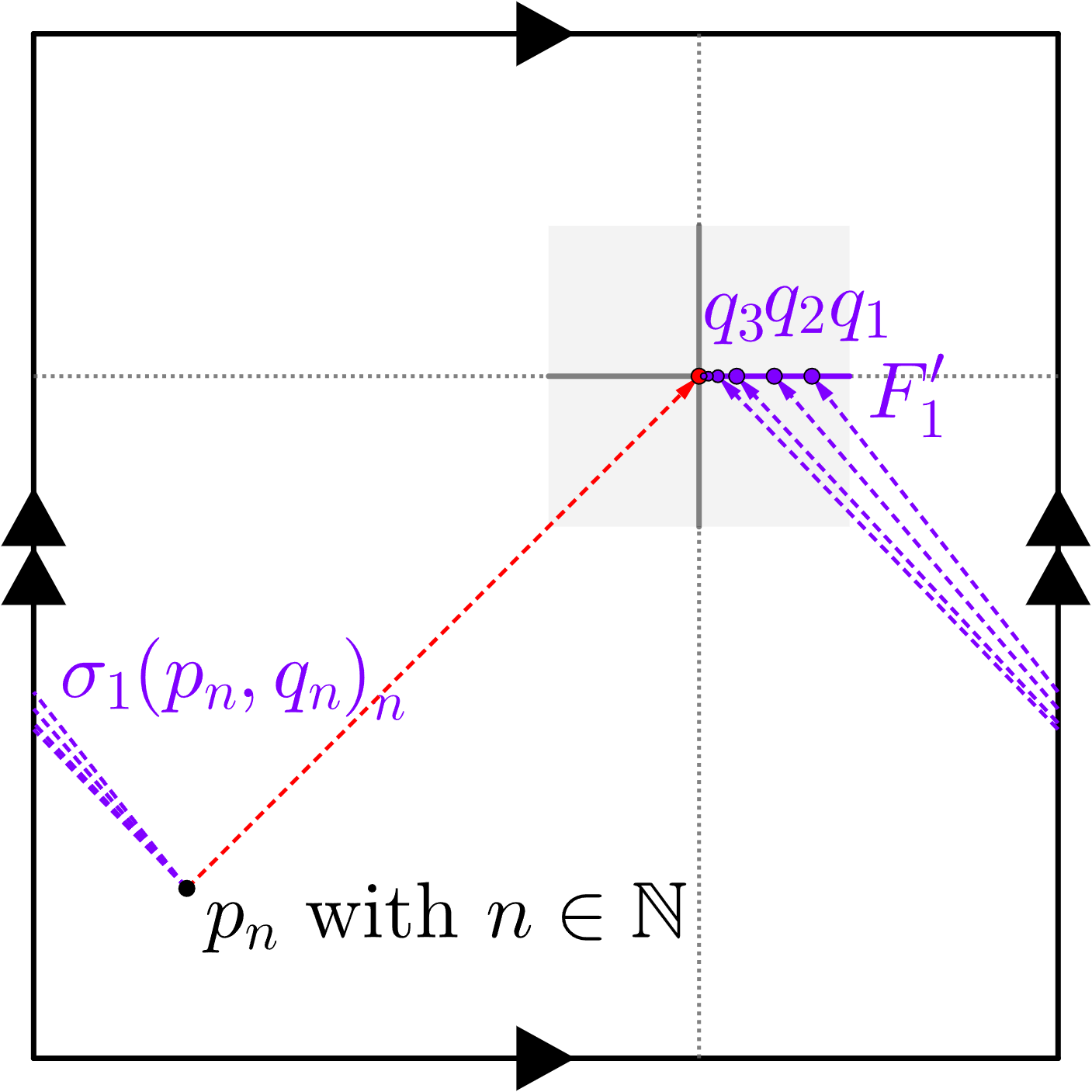}
}
\quad
\subfigure[Choice of $F_2$, sequences $(p_{n,m},q_{n,m})_{n,m}$, and the only possible $\sigma_2$]{
\includegraphics[width=\imgwidth]{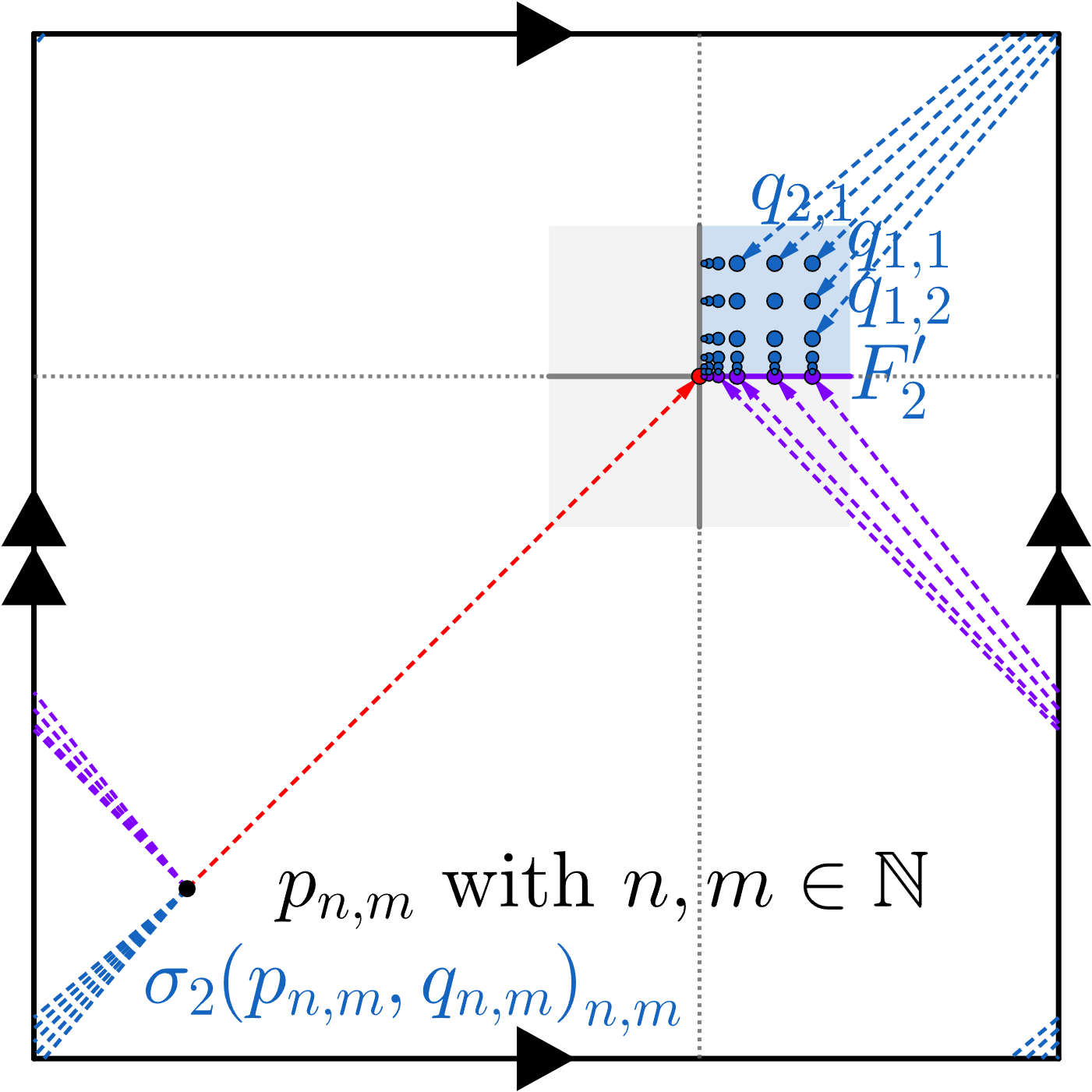}
}
\caption{Mechanism of \autoref{thm:quite_simplex} when used on the 2-torus}
\label{fig:proof-walkthrough}
\end{figure}

We will follow along the construction in \autoref{thm:quite_simplex} for this explicit example.
We begin by assuming for sake of contradiction that $GC(X)<2$. 
We partition $X\times X$ into sets $E_0$ and $E_1$ with \gls{GMPR}s $\sigma_0$ and $\sigma_1$.

\meminisubsection{Base Case:}
We select $F_0\in\mathcal F_0$, with $\img{F_0}=(p_*,q_*)$.
\WLOG we find $(p_*,q_*)\in E_0$, and find $\sigma_0(p_*,q_*)$ is as displayed in \autoref{fig:proof-walkthrough}(a).
The geodesic $\sigma_0(p_*,q_*)$ extends to $\Gamma_0$, a \gls{GMPR} over $\img{F_0}$.

\meminisubsection{Step 1:}
We select $F_1$, an embedding in $\mathcal F_1$ where $F_0$ is a face of $F_1$ and $\Gamma_0(\img{F_0})\cap \overline{\pi_{GX}^{-1}(\relint{F_1})}$ is empty.
We choose sequence $(p_n,q_n)_n$ within $\relint{F_1}$, where $\lim\limits_{n\to\infty}(p_n,q_n)=(p_*,q_*)$.
\WLOG we assume that $(p_n,q_n)_n\subseteq E_\ell$ and that $\sigma_\ell(p_n,q_n)_n$ is as displayed in \autoref{fig:proof-walkthrough}(b).
Since $\sigma_\ell(p_n,q_n)_n$ does not approach $\sigma_0(p_*,q_*)$, this would yield a discontinuity if $\ell=0$.
Thus, $\ell=1$.
The geodesics $\{\sigma_1(p_n,q_n):n\in\mathbb N\}$ extend to a \gls{GMPR} $\Gamma_1$ over $\img{F_1}$.

\meminisubsection{Step 2:}
We select $F_2$, the only embedding in $\mathcal F_2$ where $F_1$ is a face of $F_2$ and $\Gamma_1(\img{F_1})\cap \overline{\pi_{GX}^{-1}(\relint{F_2})}$ is empty.
For each $n$, we choose sequence $(p_{n,m},q_{n,m})_m$ within $\relint{F_2}$, where $\lim\limits_{m\to\infty}(p_{n,m},q_{n,m})=(p_n,q_n)$.
\WLOG we assume that $(p_{m,n},q_{m,n})_{m,n}\subseteq E_\ell$ and that $\sigma_\ell(p_{m,n},q_{m,n})_{m,n}$ is as displayed in \autoref{fig:proof-walkthrough}(c).
Since for each $n$, $\sigma_\ell(p_{m,n},q_{m,n})_m$ does not approach $\sigma_1(p_n,q_n)$, this would yield a discontinuity if $\ell=1$.
Similarly, $\sigma_\ell(p_{k,k},q_{k,k})_k$ does not approach $\sigma_0(p_*,q_*)$, and this would yield a discontinuity if $\ell=0$.
This is a contradiction, as $\sigma_\ell$ can only be continuous if $\ell$ is a unique value (say $\ell=2$).

\subsection{Embeddings of Simplices}

Since \autoref{thm:quite_simplex} is concerned with embeddings of simplices into $X\times X$, we will describe a simple way to construct such embeddings in subsets of Euclidean spaces.
The construction is motivated by the fact that after embedding the start point $p$ into the first factor of $X\times X$, the embedding of the second factor can be chosen in a convex hull of points in $X$.

\begin{define}
    \label{def:simplex_shenanigans}
    Let $T_k:=\left\{r\in \mathbb R^k:\sum\limits_{i=1}^kr_i\leq1;\;
     \forall i,\;r_i\geq 0
    \right\}$, and let $S_{k}:=\left\{r\in T_k:\sum\limits_{i=1}^kr_i=1\right\}$, which are $k$-dimensional and $(k-1)$-dimensional simplices respectively.
    For $\textbf{t}\in T_k$, $1\leq k'\in\mathbb N$, let
    \\$D_{\textbf{t},k'}:=\{\textbf{x}\in\Delta_{k+k'-1}:(\textbf{x}_{1},\dots,\textbf{x}_k)=\textbf{t}\}$ be the subset of simplex $\Delta_{k+k'-1}$ whose first $k$ coordinates are $\textbf{t}$.
\end{define}
\begin{define}
    A finite set of points $A$ in a Euclidean space is \textbf{affine independent} if no $a\in A$ can be written as $\sum\limits_{a'\in A\setminus \{a\}}\lambda_{a'}a'$, where $\lambda_{a'}\in\mathbb R$ and $\sum\limits_{a'\in A\setminus \{a\}}\lambda_{a'}=1$. 
    It is easily seen that if $A$ is affine independent, its convex hull is homeomorphic to $\Delta_{|A|-1}$.
\end{define}

\begin{lemma}
    For any space $Z$, subset $Y$ of a Euclidean space, and integers $k,k'>0$, let there be an embedding $p\colon T_k\hookrightarrow Z$ and a map $q\colon T_k\to Y^{k'}$.
    Assume that for $\textbf{t}\in T_k\setminus S_k$, $q(\textbf{t})$ maps to $k'$ affinely independent points whose convex hull is contained within $Y$, and for $\textbf{t}\in S_k$, $q(\textbf{t})=(y_{\textbf{t}},\dots,y_{\textbf{t}})$ for some $y_{\textbf{t}}\in Y$.
    Then there exists an embedding
    $f\colon \Delta_{k+k'-1}\hookrightarrow Z\times Y$ such that:
    \begin{compactitem}
        \item
        For all $\textbf{t}\in T_k$ and $\textbf{x}\in D_{\textbf{t},k'}$, the first factor $\pi_1(f(\textbf{x}))=p(\textbf{t})$.

        \item
        For all $\textbf{t}\in T_k$, the image $\pi_2(f(D_{\textbf{t},k'}))$ is the convex hull of the points defined by $q(\textbf{t})$.
    \end{compactitem}
    \label{lem:nice_embeddings}
\end{lemma}

\textbf{Proof:}
The map $q'\colon T_k\to C(\Delta_{k'-1}, Y)$ can be constructed by mapping the $i$th vertex of $\Delta_{k'-1}$ to $q(\textbf{t})_i$ for each $\textbf{t}\in T_k$, and constructing the rest of $q'(\textbf{t})$ linearly.
The map $q'(\textbf{t})$ is an embedding for $\textbf{t}\in T_k\setminus S_k$ and maps all of $\Delta_{k'-1}$ to a point for $\textbf{t}\in S_k$.

All elements of $\Delta_{k+k'-1}$ can be written $(\textbf{t},\textbf{t}')$ for $\textbf{t}\in T_k$.
Using this, we define $f\colon \Delta_{k+k'-1}\to Z\times Y$:

$f(\textbf{t},\textbf{t}')=\begin{cases}
    \left(p(\textbf{t}),q'(\textbf{t})\left(\frac{\textbf{t}'}{\sum\limits_{i}\textbf{t}'_i}\right)\right)&\textbf{t}\in T_k\setminus S_k\\
    \left(p(\textbf{t}),y\right)&\textbf{t}\in S_k \text{, and for } q'(\textbf{t})(\Delta_{k'-1}) =\{(y,\dots,y)\}
\end{cases}$

Continuity follows from continuity of $p$ and $q$.
Since $p$ is injective, and $\pi_1\circ f(\textbf{t},\textbf{t}')=p(\textbf{t})$, it is sufficient to verify injectivity independently for each value of $\textbf{t}$.
For $\textbf{t}\notin S_k$, injectivity of $f$ follows from injectivity of $q'(\textbf{t})$.
For $\textbf{t}\in S_k$, we have that $\textbf{t}'$ can only be $(0,\dots,0)$, so injectivity is trivial in this case.
Thus, $f$ is a well defined embedding $\Delta_{k+k'-1}\hookrightarrow Z\times Y$.
The desired properties follow from construction.
\proofclose{\autoref{lem:nice_embeddings}}

To create a map $q\colon T_k\to Y^{k'}$, we will often combine $q_1,\dots,q_{k'}\colon T_k\to Y$ to create a \lrq{product,} a map $T_k\to Y^{k'}$ via $\textbf{t}\mapsto (q_1(\textbf{t}),\dots,q_{k'}(\textbf{t}))$.
We will denote this map $q_1\funcprodsymb\dots\funcprodsymb q_{k'}$.

\subsection{Polyhedra}
We inspect \textbf{convex polyhedra}, defined as the convex hull of finitely many points in $\mathbb R^3$. 
In particular, the spaces we are interested in are the 2-D surfaces of convex polyhedra, equipped with the \textbf{flat metric} that measures the length of the shortest path along the surface between points.
These spaces may be equivalently viewed as the union of their polygonal faces, glued at appropriate edges.
We will prove various lemmas for these spaces, and use them to prove correctness of an algorithm to calculate all geodesics from a point.

\subsubsection{General Lemmas}

\begin{define}
    \label{def:facewalkfacepath}
    Let $X$ be a polyhedron, and consider the graph whose vertices are the polygon faces of $X$.
    Connect two faces if they share an edge.
    We will define a \textbf{face walk} as a walk in this graph, and a \textbf{face path} as a path (a walk with no repeat faces).
    The graph considered is equivalent to the 1-skeleton of the dual polyhedron \cite{duality} of $X$.
\end{define}

\begin{define}
    Let $F_0,\dots,F_n$ be a face walk.
    A \textbf{walk unfolding} is a map $f\colon \bigcup F_i\to \mathbb R^2$ such that $f|_{F_i}$ is an isometric embedding for each $F_i$, and $f|_{F_i\cup F_{i+1}}$ is a homeomorphism for each $i<n$. 
    We similarly define a \textbf{path unfolding} for face paths.
\end{define}

This is similar to the concept of an unfolding discussed in Shephard's Conjecture \cite{shephard}, with the main distinction being that in walk/path unfoldings, all faces of the polyhedron need not be represented.
We display an example of a path unfolding obtained from an icosahedron in \autoref{fig:lemonstarconvex}(b).

\begin{lemma}
    \label{lem:unfold_unique}
    Given a face walk $F_0,\dots,F_n$ on a convex polyhedron and isometric embedding $g\colon F_i\hookrightarrow \mathbb R^2$, there is a unique walk unfolding $f$ such that $f|_{F_i}=g$.
    We call this the walk unfolding with respect to $g$.
\end{lemma}
This can be shown through induction by considering the restriction of $f$ to each face.

\begin{lemma}
    \label{lem:nofacerevisit}
    Geodesics on a convex polyhedron cannot exit a face then revisit it.
\end{lemma}

\textbf{Proof:}
\AFSOC a segment of geodesic $G$ starts and ends on face $F$, but $\img{G}\not\subseteq F$.
We can \wLOG assume $G$ is this segment, with end points $a$ and $b$.
Consider a natural\footnote{
Polyhedra are defined as objects in $\mathbb R^3$, so there is a natural embedding, unique up to isometry.
}
embedding $i$ of the whole polyhedron into $\mathbb R^3$ such that each face is isometrically embedded.
It is clear that both $i$ and its inverse preserve the lengths of paths, as the restriction to each face is an isometry.
Thus, it is sufficient to find a path between $i(a)$ and $i(b)$ that is shorter than $i\circ G$, and is contained in the image of $i$.
The \lrq{straight line} path $G'$ between $i(a)$ and $i(b)$ is contained in the image of $F$ by convexity of faces, and thus is a valid path.
We have that $i\circ G$ is distinct from $G'$, as $i\circ G([0,1])\not\subseteq \img{F}$.
This implies $i\circ G$ is longer than $G'$ as in $\mathbb R^3$, any two points have a unique shortest path.
Thus, $i^{-1}\circ G'$ is shorter than $G$, and we have a contradiction.
\proofclose{\autoref{lem:nofacerevisit}}

\begin{lemma}
    \label{lem:geodesicdecomp}
    Let $G$ be a geodesic on convex polyhedron $X$. 
    There exists a decomposition $G=G_0\cup\dots\cup G_n$ where each $G_i$ is a line segment or point that intersects $G_{i+1}$ at a point, and $G_i\subseteq F_i$ for some face $F_i$ of $X$. 
    Additionally, $F_0,\dots,F_n$ is a face path, and the image of $G$ is a line on the path unfolding of $F_0,\dots,F_n$.
\end{lemma}
\textbf{Proof:} 
To decompose $G$ we first claim that the intersection of $G$ with each edge of $X$ must be connected, and thus either a point or a line segment.
This may be shown similar to \autoref{lem:nofacerevisit}.
We may then consider a particular edge $e$ and split $G$ into one or two segments based on a point $G\cap e$, or one to three segments based on the endpoints of line segment $G\cap e$.
Continuing the decomposition on each segment for the remaining edges will create a finite ordered decomposition of non-point elements $G''_0,\dots,G''_{n''}$, where for any $G''_i$ and any edge $e$ of $X$, $G''_i\cap e$ is either empty, all $G''_i$, or a single endpoint of $G''_i$.
Thus, each element must be contained in a face, as leaving a face would cause an intersection with an edge.
Let $F''_0,...,F''_{n''}$ be these faces, where $G''_i\subseteq F''_i$. 
If the endpoint of $G''_i$ (i.e. $G''_i\cap G''_{i+1}$) is on the interior of an edge of $X$, then $F''_i$ and $F''_{i+1}$ must share that edge.
Otherwise, $G''_i\cap G''_{i+1}$ must be a vertex of $X$ adjacent to $F''_i$ and $F''_{i+1}$.
In this case, all faces that share this vertex are joined in a cycle, and we may insert the appropriate faces in between $F''_i$ and $F''_{i+1}$ such that each adjacent pair share an edge. 
We will insert copies of the vertex $G''_i\cap G''_{i+1}$ to ensure the decomposition matches the list of faces.
After doing this for all the original $F''_i$, we obtain a new decomposition $G'_0,\dots,G'_{n'}$, where $G'_i\subseteq F'_i$ and $F'_0,\dots,F'_{n'}$ is a face walk. 
Finally, if $F'_i=F'_j$ for $i<j$, $G'_i\cup\dots\cup G'_j$ is a continuous segment of a geodesic that begins and ends on the same face.
Thus, by \autoref{lem:nofacerevisit}, this segment never leaves $F'_i$, and we may replace $F'_i,\dots,F'_j$ with just $F'_i$, and replace $G'_i,\dots,G'_j$ with the union $G'_i\cup\dots\cup G'_j$.
Repeating this until no duplicates remain will create a face path $F_0,\dots,F_n$.

\renewcommand{\imgheight}{100 pt}
\begin{figure}[ht!]
    \centering
\subfigure[Corner within an edge]{
\includegraphics[height=\imgheight]{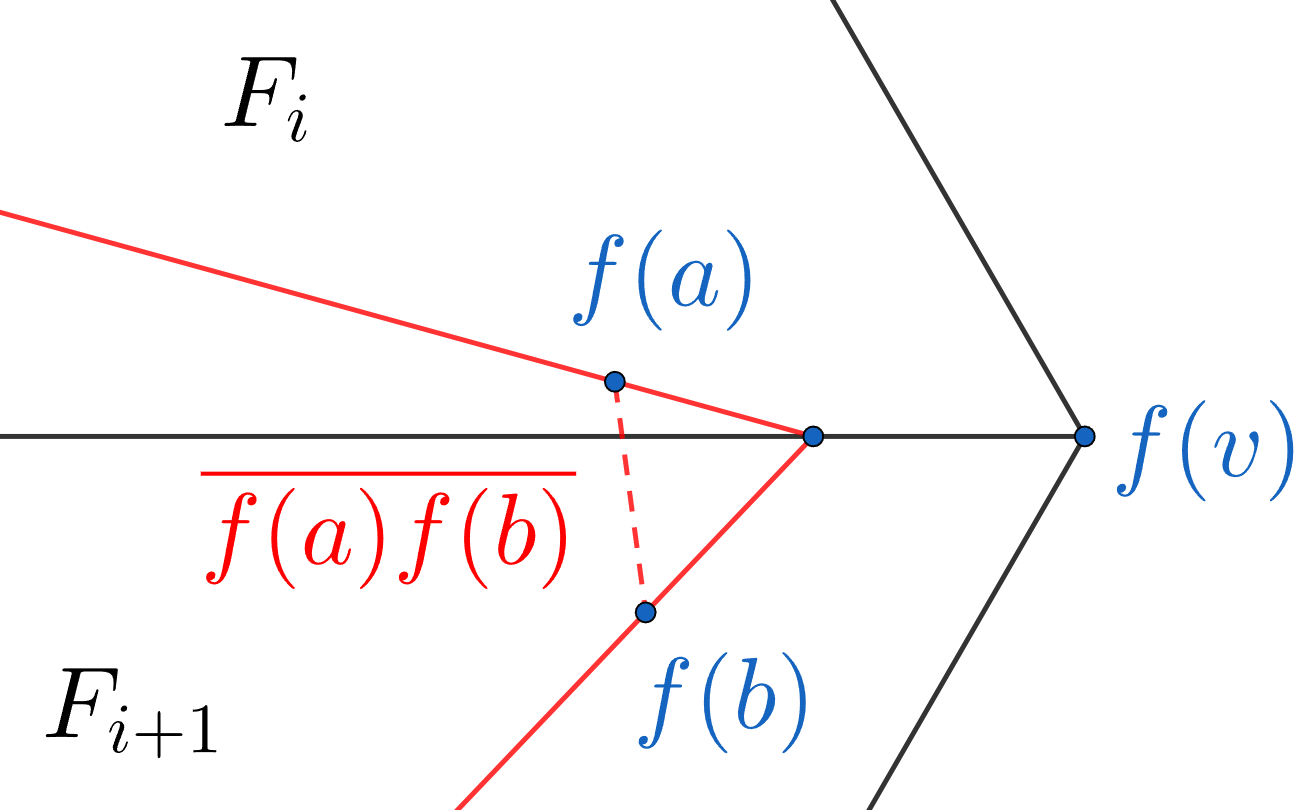}
}
\subfigure[Corner on a vertex with $\theta<\pi$]{
\includegraphics[height=\imgheight]{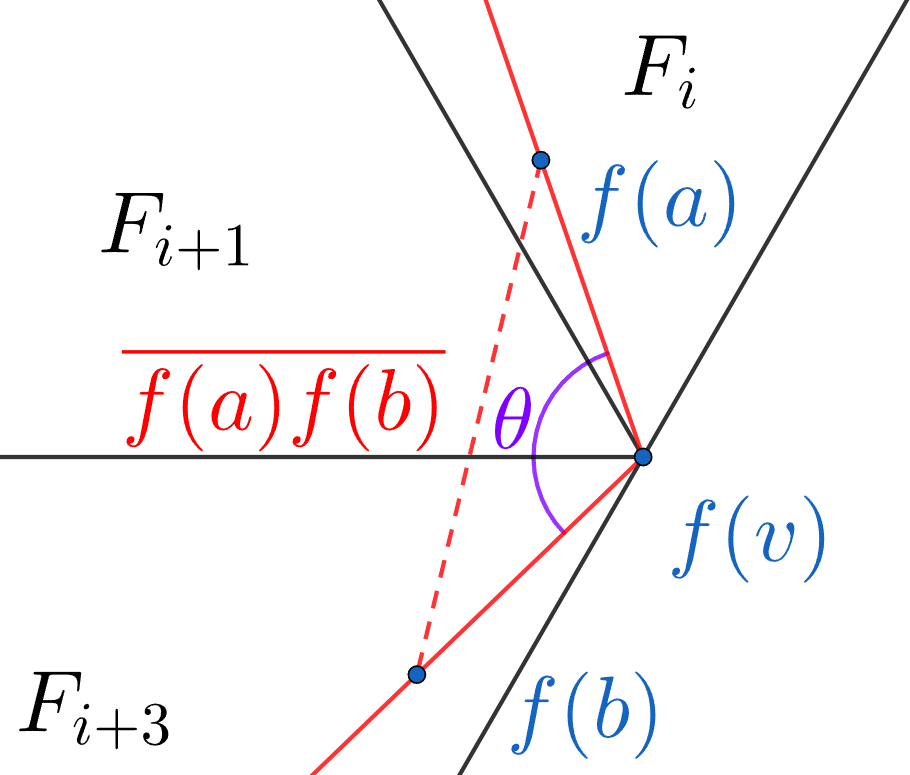}
}
\subfigure[Corner on a vertex with $\theta>\pi$]{
\includegraphics[height=\imgheight]{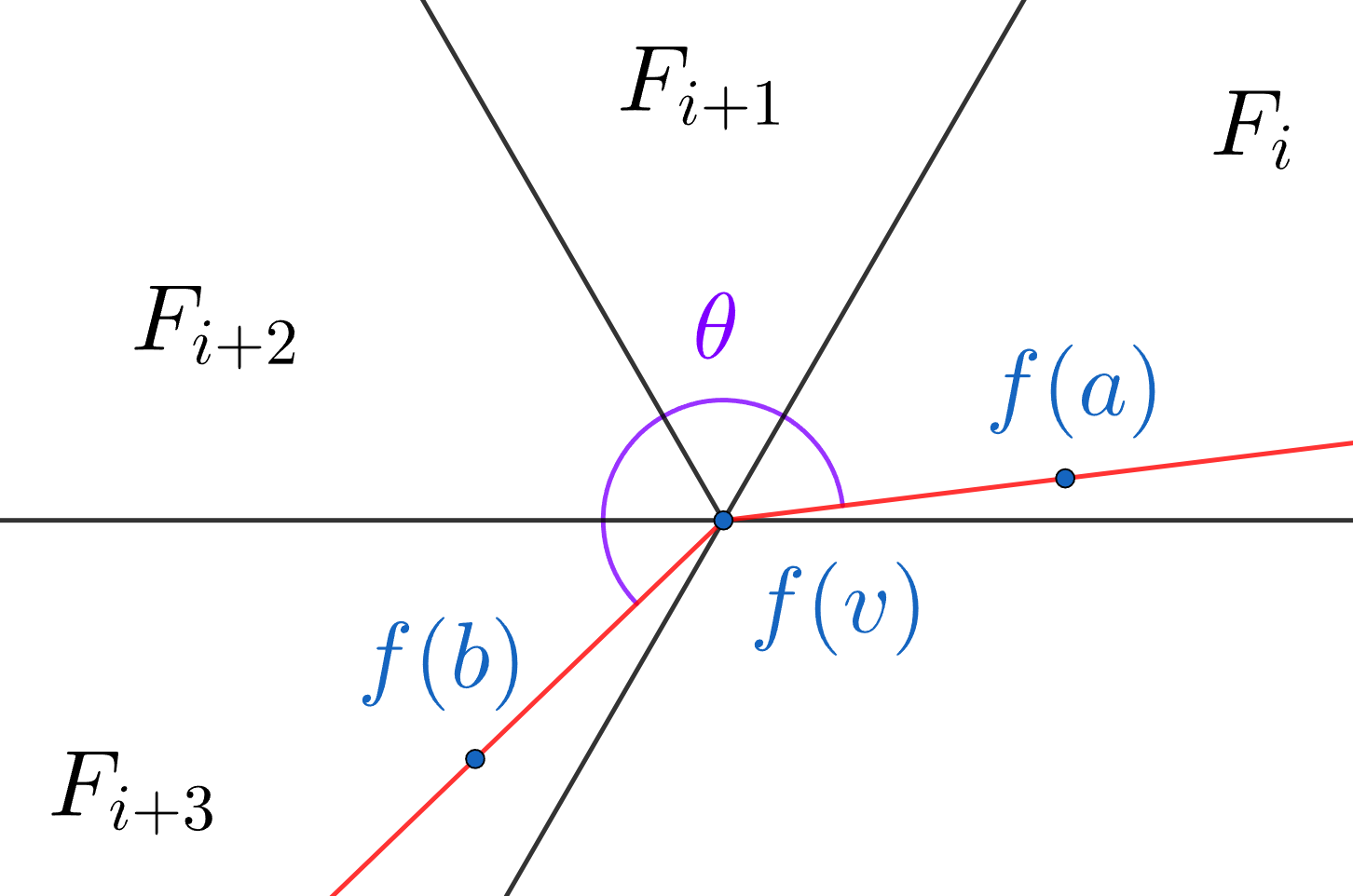}
}
\caption{Contradictions arising from a corner within a geodesic}
\label{fig:shortcuts}
\end{figure}

Consider a path unfolding $f$ of $F_0,\dots,F_n$ onto $\mathbb R^2$, and consider the image of $G$ wrt. this map.
It is clear that each $f(G_i)$ must be a line segment or point, so $f(G)$ must be a union of line segments.
We must show that the segments are collinear. 
If a \lrq{corner} exists, it must either be on the edge or vertex of a face.
In either case, we will reach a contradiction by creating a shorter path on a small neighborhood of the corner.

\AFSOC there is a corner on some edge.
The image of $f|_{F_i\cup F_{i+1}}$ is two polygons that share an edge, and the assumed corner is on the interior of this edge.
Thus, the corner is on the interior of the image of $f|_{F_i\cup F_{i+1}}$.
Thus we can find a shorter path within the image between points $f(a)$ and $f(b)$ sufficiently close to the corner, visualized in \autoref{fig:shortcuts}(a).
Since both $f|_{F_i\cup F_{i+1}}$ and its inverse preserve the length of paths, we have a contradiction.

Now \aFSOC there is a corner on a vertex $v$.
Let $\epsilon$ be much smaller than the length of any edge containing $v$.
Consider the section of the path that is $\epsilon$-close to $v$, and let the endpoints be $a$ and $b$.
It is enough to show that the segment of $G$ from $a$ to $b$ is not a geodesic.
Let $\theta$ be the angle between $f(a), f(v)$, and $f(b)$, oriented as in \autoref{fig:shortcuts}(b, c).
If $\theta<\pi$, then we can create a shorter path by simply connecting the line $\overline{f(a)f(b)}$.
We visualize this in \autoref{fig:shortcuts}(b).
If $\theta>\pi$ as in \autoref{fig:shortcuts}(c), we will show we can choose a different face walk that reduces this to the case where $\theta<\pi$.
All faces on the polyhedra connected to vertex $v$ create a cycle.
We will go the other way around the cycle (with associated walk unfolding $g$).
Since polyhedra have positive curvature, the angle between $g(a), g(v)$, and $g(b)$ in this direction must be less than $2\pi-\theta$, and thus less than $\pi$.
Then this reduces to the first case, and we can find a shorter path between $a$ and $b$.
\proofclose{\autoref{lem:geodesicdecomp}}

\begin{define}
    For a convex polyhedron $X$, face path $F_0,\dots,F_n$, and path unfolding $f$, let $p\in F_0$.
    We will say $q\in F_n$ has property $\star$ with respect to $p$ if the line $P:=\overline{f(p)f(q)}$ has the following properties:
    \begin{compactenum}[(a)]
        \item
        $P\subseteq \img{f}$.
        \item
        The intersection $P\cap f(F_0)$ is a path from $f(p)$ to edge $f(F_0\cap F_{1})$, and the intersection $P\cap f(F_n)$ is a path from edge $f(F_{n-1}\cap F_{n})$ to $f(q)$.
        \item
        For $i\in \{1,\dots,n-1\}$, the intersection $P\cap f(F_i)$ is a path between edges $f(F_i\cap F_{i-1})$ and $f(F_i\cap F_{i+1})$.
            \end{compactenum}
\end{define}

\begin{lemma}
    \label{lem:starconvex}
    For convex polyhedron $X$, face path $F_0,\dots,F_n$, path unfolding $f$, and some $p\in F_0$, the set of $q\in F_n$ that satisfy property $\star$ (with respect to $p$) is convex.
\end{lemma}

\textbf{Proof:}
Pick $q_1,q_2\in F_n$ such that $q_1$ and $q_1$ have $\star$.
Let $P_1=\overline{f(p)f(q_1)}$ and $P_2=\overline{f(p)f(q_2)}$.
For some $\lambda\in(0,1)$, let $q=\lambda q_1+(1-\lambda)q_2$, and let $P=\overline{f(p)f(q)}$.
We can assume $P_1\cap P_2=\{p\}$, as if this were not true, then \wLOG $P_1\subseteq P_2$.
Then $\star$ is trivial for $q$, as $P_1\subseteq P\subseteq P_2$.
We abuse notation and represent $P_1$ and $P_2$ as lines $[0,1]\to \mathbb R^2$ via $P_i(t)=t\cdot f(q_i)+(1-t)\cdot f(p)$.
Then consider $P':[0,1]\to \mathbb R^2$ via $P'(t)=\lambda P_1(t)+(1-\lambda)P_2(t)$.
This map is a line with end points $p$ and $q$.
Thus, $P'=P$, since there a unique line between endpoints in $\mathbb R^2$.

\renewcommand{\imgheight}{165 pt}
\begin{figure}[ht!]
    \centering
\subfigure[A line between $P_1$ and $P_2$ must intersect $P$]{
\includegraphics[height=\imgheight]{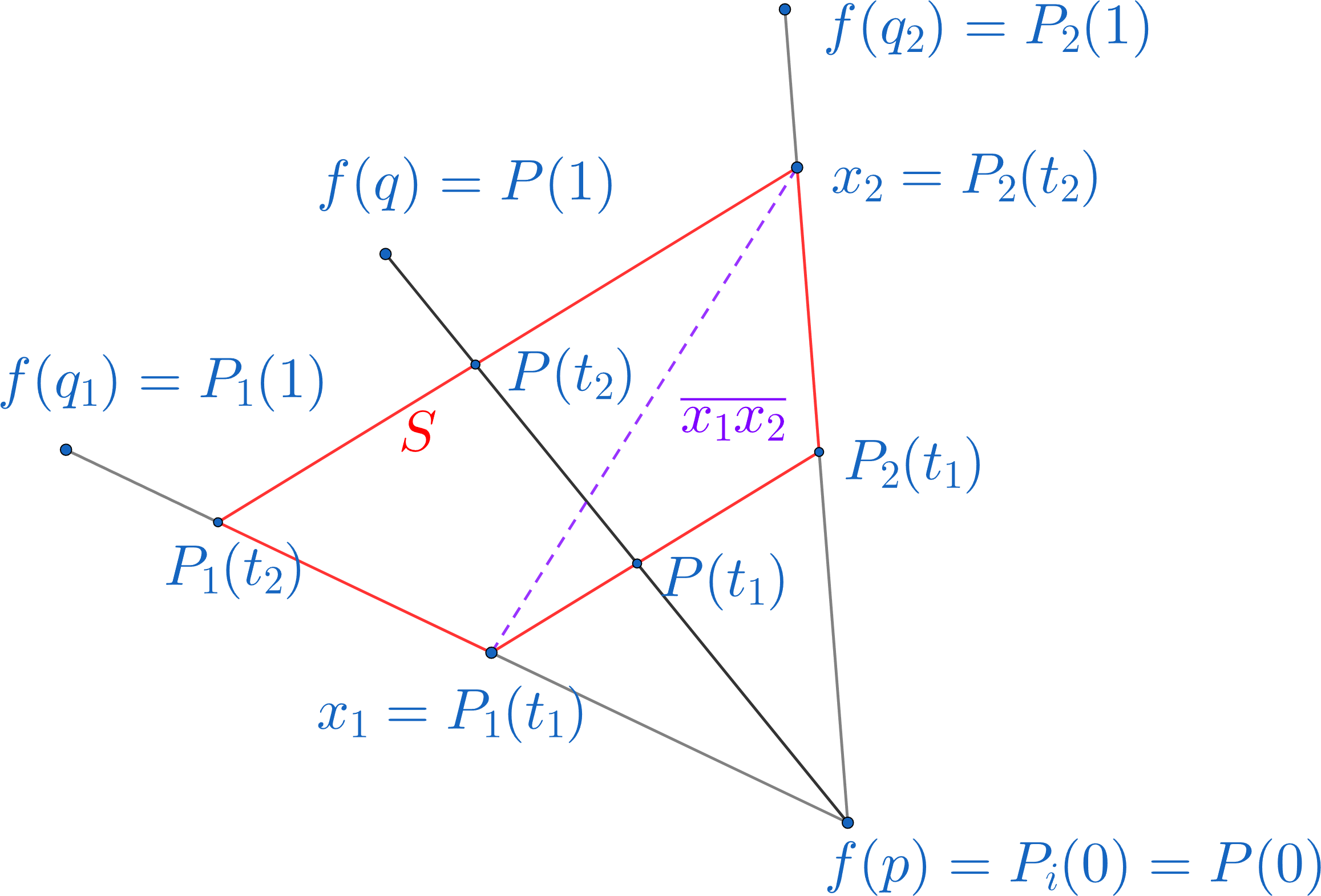}
}
\subfigure[Three points with $\star$]{
\includegraphics[height=\imgheight]{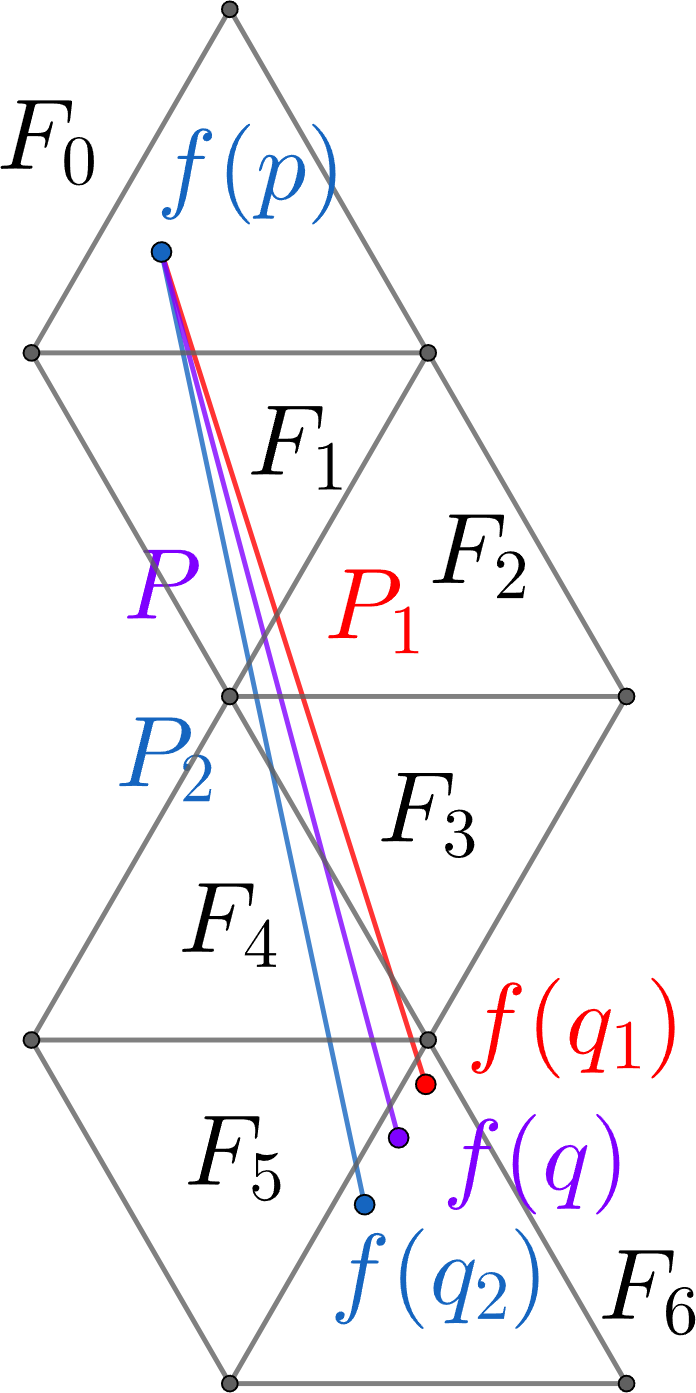}
}\quad
\subfigure[Point without $\star$, as $\overline{f(p)f(q)}$ exits $\img{f}$]{
\includegraphics[height=\imgheight]{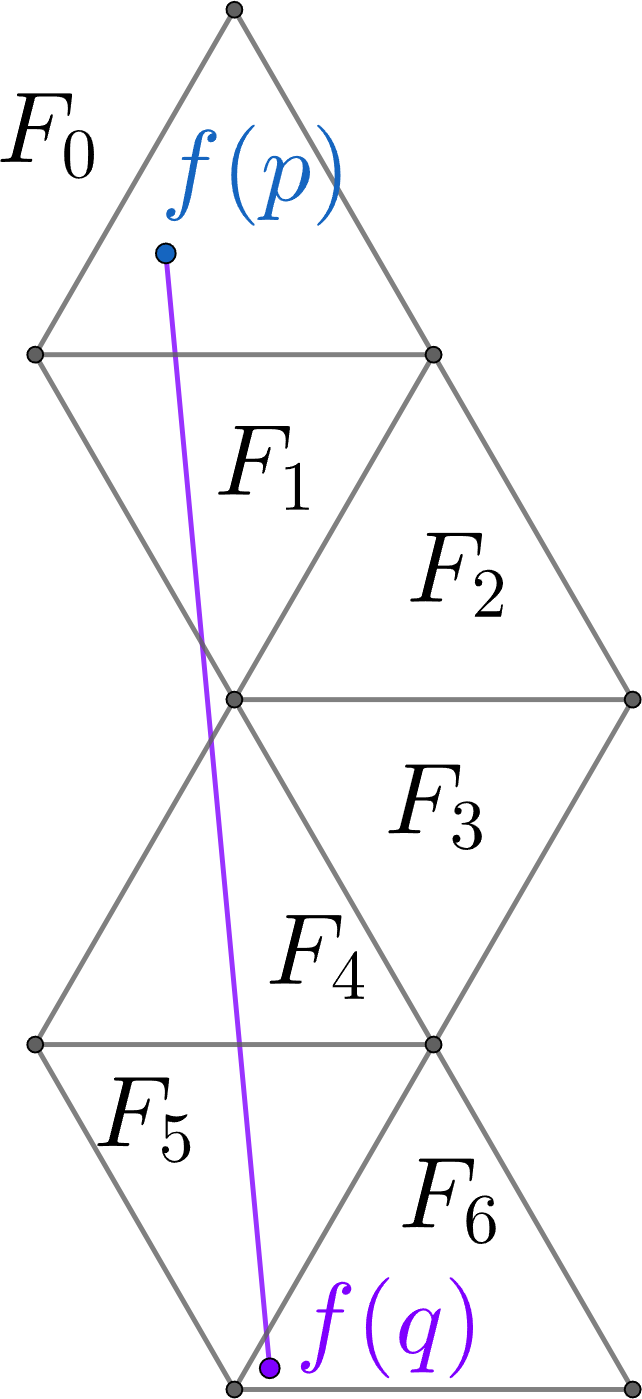}
}
\caption{Visualizations for \autoref{lem:starconvex}. Face paths in (b) and (c) are obtained from an icosahedron.}
\label{fig:lemonstarconvex}
\end{figure}

We will show for any $x_1\in \img{P_1}$ and $x_2\in\img{P_2}$, $\overline{x_1x_2}$ intersects $P$.
Let $x_1=P_1(t_1)$ and $x_2=P_2(t_2)$ for $t_1,t_2\in[0,1]$.
This is trivial for $t_1=t_2$ as we can set $t=t_1=t_2$ and obtain $P(t)=\lambda P_1(t_1)+(1-\lambda)P_2(t_2)$.
Additionally, if $t_i=0$ for either $i$, then $x_i=P_i(0)=f(p)=P(0)$.
Then \wLOG we will assume $0<t_1<t_2$.
Since $0<t_1,t_2$, $P_1(t_1)\neq P_2(t_1)$ and $P_1(t_2)\neq P_2(t_2)$, so $P_1(t_1),P_1(t_2),P_2(t_2),P_2(t_1)$ are four distinct points, and form a cycle $S$ in $\mathbb R^2$.
Note that $S$ is formed from the triangle $f(p),P_1(t_2),P_2(t_2)$ by adding the edge $\overline{P_1(t_1)P_2(t_1)}$ and removing the segments $\overline{P_1(0)P_1(t_1)}$ and $\overline{P_1(0)P_2(t_1)}$.
Thus, $S$ is a polygon, and since any internal angle of $S$ must be at most $\pi$, $S$ is convex.
We consider the line $\overline{x_1x_2}$, and the line $\overline{P(t_1)P(t_2)}$.
By convexity, both of these lines are contained in $S$.
The order of their endpoints on the boundary of $S$ is $x_1,P(t_1),x_2,P(t_2)$, so these lines must intersect.
We visualize this in \autoref{fig:lemonstarconvex}(a).

We will inspect $P\cap f(F_0)$.
This path certainly begins at $f(p)$.
From assumptions, we know that $P_1\cap f(F_0)$ and $P_2\cap f(F_0)$ have endpoints on edge $f(F_0\cap F_1)$.
Thus, a convex combination of these endpoints must be on $P$, implying $P$ intersects that edge.
This intersection must be the first time $P$ intersects an edge since $P$ begins within $f(F_0)$.
We will proceed with induction on $i\in\{1,\dots,n-1\}$.
We use the hypothesis that $P\cap f(F_{i-1})$ has an endpoint on $f(F_{i-1}\cap F_i)$.
Since both $P_1$ and $P_2$ have endpoints on edge $f(F_{i}\cap F_{i+1})$, we may repeat the argument from the $F_0$ case to find that $P\cap f(F_i)$ must have an endpoint on $f(F_{i}\cap F_{i+1})$.
Finally, for $P\cap f(F_n)$, we have from induction that the startpoint is on edge $f(F_{n-1}\cap F_n)$.
Since $P$ ends at point $f(q)$, we have from \autoref{lem:nofacerevisit} that $P\cap f(F_n)$ must be a line from the startpoint to $f(q)$.
This shows properties (b) and (c).
Property (a) follows from noticing that during the proof for (b) and (c), we have accounted for all of line $P$ as a sequence of line segments whose union is a line from $f(p)$ to $f(q)$.
\proofclose{\autoref{lem:starconvex}}

We visualize an example in \autoref{fig:lemonstarconvex}(b), for a path unfolding from an icosahedron.

\begin{define}
    For metric space $X$, let $P\subseteq X$ be a finite collection of points.
    For each $p\in P$, the associated \textbf{Voronoi cell} is the set of all points in $X$ closest to $p$, that is, $\left\{x\in X:d(x,p)=\min\limits_{p'\in P}d(x,p')\right\}$.
    The collection of all Voronoi cells is called the \textbf{Voronoi diagram}, and covers $X$ \cite{voronoibook}.
    \end{define}

\begin{lemma}
    \label{lem:voronoisubcell}
    For some metric space $X$, let $P\subseteq P'\subseteq X$ be collections of points with Voronoi diagrams $V$ and $V'$ respectively.
    Let $p\in P$ have Voronoi cell $C_p$ in $V$, and cell $C'_p$ in $V'$. Then $C_p'\subseteq C_p$.
\end{lemma}

\textbf{Proof:}
By definition, if a point $x\in C_p'$, it is closest to $p$ among all points in $P'$.
Thus, $x$ must be the closest to $p$ among all points in $P$, since $P\subseteq P'$.
Then $x\in C_p$.
\proofclose{\autoref{lem:voronoisubcell}}

\begin{define}
    Let $X$ be a metric space, and take some $p\in X$.
    The \textbf{cut locus} of $p$ is the set of all $q\in X$ such that there are multiple distinct geodesics from $p$ to $q$.
    The \textbf{total cut locus} is $\{(p,q)\in X\times X:q\text{ is on cut locus of }p\}$.
\end{define}

\begin{lemma}
    \label{lem:cutlocvoronoi}
    For a convex polyhedron, the cut locus on face $F^*$ of point $p$ (on face $F_0$) has the form of the boundaries of a Voronoi diagram in $\mathbb R^2$, and the points that generate this Voronoi diagram are copies of $p$ arising from path unfoldings of face paths from $F_0$ to $F^*$.
    Additionally, for each copy $p'$ of $p$ with associated path unfolding $f$, and every point $q$ of its Voronoi cell in $f(F^*)$, $\overline{p'q}$ is the image of a geodesic between $p$ and $(f|_{F^*})^{-1}(q)$, and $q$ has property $\star$ with respect to $p'$.
\end{lemma}

The first part of this is shown by Agarwal in \cite{star}, and properties (b) and (c) of $\star$ arise from the face path of each $p'$ being constructed so that geodesics from $p'$ visit the faces in assigned order.
Given that $\overline{p'q}$ is the image of a geodesic, it is trivial to show (a), as $\overline{p'q}$ is certainly contained within the image of its associated path unfolding.
\proofclose{\autoref{lem:cutlocvoronoi}}

\subsubsection{Cut Locus Algorithm}

The total cut locus of a space is relevant to its geodesic complexity, as it represents areas where \gls{GMPR}s must make a choice, and thus could potentially be discontinuous.
Because of this, simplex embeddings we construct for \autoref{thm:quite_simplex} will be made with the total cut locus of each space in mind.
Due to the significance of the cut locus, we give an algorithm\footnote{
Our Python implementation is available at \cite{cut_locus_alg}.
} to calculate the cut locus of a point on a polyhedron surface, and will proceed to show its correctness.

\meminisection{Filter Algorithm:}
We will show correctness of \autoref{alg:filter}, designed to filter face paths $F_0,\dots,F_n$ based on whether certain points in face $F_n$ have property $\star$.
To show termination, it is enough to show that the loop on line \ref{line:filteralgloop} is over a finite set.
The set $f(F_n)$ is a polygon, and cells of a Voronoi diagram of a finite point set in $\mathbb R^2$ have boundaries composed of line segments, rays, and lines.
Thus, $C\cap f(F_n)$ must also be a polygon, so $V$ is well defined and finite.
If $C\cap f(F_n)$ is empty, we return true, which is consistent with our desired output.
Otherwise, since both $C$ and $f(F_n)$ are convex (cells of a Voronoi diagram in $\mathbb R^2$ are convex), $C\cap f(F_n)$ is a convex polygon.
Thus, by \autoref{lem:starconvex}, all points in $C\cap f(F_n)$ have property $\star$ if and only if all vertices of $C\cap f(F_n)$ have property $\star$.
Then our algorithm returns the correct output, since the loop on line \ref{line:filteralgloop} tests if $q$ has property $\star$.
\begin{algorithm}[ht!]
\algnotext{EndLoop}
\algnotext{EndIf}
\caption{Filter Algorithm (Pseudocode)}\label{alg:filter}
\hspace*{\algorithmicindent}\textbf{Input:}\\
\hspace*{\algorithmicindent}\hspace{27 pt}$F_0,\dots,F_n$: Face path in a convex polyhedron\\
\hspace*{\algorithmicindent}\hspace{27 pt}$f\colon \bigcup F_i\to \mathbb R^2$: Path unfolding of $F_0,\dots,F_n$\\
\hspace*{\algorithmicindent}\hspace{28 pt}$p^*\in f(F_0)$: Point to inspect\\
\hspace*{\algorithmicindent}\hspace{27 pt}$C\subseteq \mathbb R^2$: Voronoi cell with $p^*\in C$\\
\hspace*{\algorithmicindent} \textbf{Output:} \\
\hspace*{\algorithmicindent}\hspace{27 pt}Whether $\forall q\in C\cap f(F_n)$, $q$ has property $\star$ with respect to $p^*$
\begin{algorithmic}[1]
\State
$V\gets \text{vertices of polygon }C\cap f(F_n)$\label{line:filteralgvertexfinder}
\Loop{ for $q\in V$\label{line:filteralgloop}}
\If{$\overline{p^*q}\not\subseteq \text{img}(f)$\label{line:filteralgpropertyatest}}
\Return{False}
\Comment{Test property (a)}
\EndIf
\State{$s\gets p^*$}
\Loop{ for $i\in (0,\dots,n-1)$\label{line:filteralginnerloop}}
\Comment{Test property (c), part of (b)}
\If{$\overline{p^*q}\cap f(F_i)$ is not a path from $s$ to $f(F_i\cap F_{i+1})$\label{line:filteralgpropertybctest}}
\Return{False}
\EndIf
\State{$s\gets\text{endpoint of }\overline{p^*q}\cap f(F_i)$}
\EndLoop
\If{$\overline{p^*q}\cap f(F_{n})$ is not a path from $s$ to $f(q)$\label{line:filteralgpropertybtest}}
\Return{False}
\Comment{Finish test of property (b)}
\EndIf
\EndLoop
\\\Return{True}
\end{algorithmic}
\end{algorithm}

In practice, we compute line \ref{line:filteralgvertexfinder} by considering all vertices of $C$ and $f(F_n)$, as well as intersections of their boundaries.
We keep vertices that lie in both regions.
We check line \ref{line:filteralgpropertyatest} by splitting $\overline{p^*q}$ into segments based on intersections with face boundaries, then verifying each segment is contained in a face.
For every equivalence or incidence check in lines$~$\ref{line:filteralgpropertybctest} and \ref{line:filteralgpropertybtest}, we use a tolerance to mitigate rounding errors.

Checking if a line segment is within a face can be done in linear time wrt. the face's edge count, and checking if a point is on a line can be done in constant time. 
Repeated applications of these allow lines$~$\ref{line:filteralgpropertyatest}, \ref{line:filteralgpropertybctest}, and \ref{line:filteralgpropertybtest} to be implemented in polynomial time wrt. the number of edges in all of $F_0,\dots, F_n$.
Thus, considering the loops on lines$~$\ref{line:filteralgloop} and \ref{line:filteralginnerloop}, the overall complexity of \autoref{alg:filter} is polynomial with respect to the number of edges in all of $F_0,\dots, F_n, C$.

\meminisection{Cut Locus Algorithm:}
We will show correctness of \autoref{alg:cutlocus}, designed to obtain the cut locus on face $F_*$ of a point $p\in F_0$.
The algorithm terminates since the size of finite set $S$ decreases each loop.
From \autoref{lem:cutlocvoronoi}, the cut locus of point $p$ must have the structure of a Voronoi diagram in $\mathbb R^2$. 
The diagram is generated by images of $p$ from path unfoldings of face paths from $F_0$ to $F_*$.
Let $P^*$ be these images, and let $V^*$ be the resulting Voronoi diagram.
Before the loop on line \ref{line:cutlocalgloop}, let $P=\{f(p):(\mathbb P,f)\in S\}$, and observe $P^*\subseteq P$, since $P$ contains all possible images of $p$ from path unfoldings.
For $x\in P^*$, we will denote $C_x^*$ as the Voronoi cell of $x$ in $V^*$, restricted to face $F_*$.
Similarly, for $x\in P$, we will denote $C_x$ as the Voronoi cell of $x$ in $V$, restricted to face $F_*$.

\def\algcolornew{purple}
\begin{algorithm}[ht!]
\algnotext{EndLoop}
\caption{Cut Locus Algorithm (Pseudocode)}\label{alg:cutlocus}
\hspace*{\algorithmicindent}\textbf{Input:}\\
\hspace*{\algorithmicindent}\hspace{27 pt}$F_0\subseteq X$: Source face of convex polyhedron $X$\\
\hspace*{\algorithmicindent}\hspace{28 pt}$p\in F_0$: Point to inspect\\
\hspace*{\algorithmicindent}\hspace{27 pt}$F_*\subseteq X$: Sink face of $X$\\
\hspace*{\algorithmicindent}\hspace{27 pt}$i_*\colon \hookrightarrow \mathbb R^2$: Embedding of $F_*$\\
\hspace*{\algorithmicindent} \textbf{Output:} \\
\hspace*{\algorithmicindent}\hspace{27 pt}Cut locus of $p$ on face $F_*$\\
\hspace*{\algorithmicindent}\hspace{27 pt}Finite set $S$ with elements $(\mathbb P,f)$ such that:\\
\hspace*{\algorithmicindent}\hspace{34 pt}$\bullet$
$\mathbb P$ is a face path from $F_0$ to $F_*$\\
\hspace*{\algorithmicindent}\hspace{34 pt}$\bullet$
$f\colon \bigcup\limits_{F\in \mathbb P}F\to \mathbb R^2$ is the path unfolding of $\mathbb P$ wrt. $i_*$\\
\hspace*{\algorithmicindent}\hspace{34 pt}$\bullet$
The cut locus of $p$ on $F_*$ is the Voronoi diagram of $\{f(p):(\mathbb P,f)\in S\}$
\begin{algorithmic}[1]
\State
$S\gets\{(\mathbb P,f):\mathbb P \text{ is a face path from }F_0\text{ to } F_*,\, f \text{ is the path unfolding of }\mathbb P\text{ wrt. }i_*\}$\label{line:cutlocalgcreateS}
\State
$V\gets $ Voronoi cell diagram of $\{f(p):(\mathbb P,f)\in S\}$ \label{line:cutlocalgvoronoicreate}
\State For each $(\mathbb P,f)\in S$, $C_{f(p)}\gets$ Voronoi cell of $f(p)$ in $V$
\Loop{ until $\forall (\mathbb P,f)\in S$, (\autoref{alg:filter})$(\mathbb P,f,f(p),C_{f(p)})$ is true\label{line:cutlocalgloop}}
    \State $(\mathbb P',f')\gets$ element of $S$ such that (\autoref{alg:filter})$(\mathbb P',f',f'(p),C_{f'(p)})$ is false
    \State $S\gets S\setminus \{(\mathbb P',f')\}$\label{line:cutlocalgremoval}
    \State $V\gets $ Voronoi cell diagram of $\{f(p):(\mathbb P,f)\in S\}$\label{line:cutlocalgvoronoirecreate}
    \State For each $(\mathbb P,f)\in S$, $C_{f(p)}\gets$ Voronoi cell of $f(p)$ in $V$
\EndLoop
\\\Return{($i_*^{-1}(V)$, $S$)}
\end{algorithmic}
\end{algorithm}

\begin{lemma}
\label{lem:starcorrectcorr}
For $x\in P$ such that $C_x$ is nonempty, let $\mathbb P$ be the face path associated with $x$, and let $f$ be the path unfolding.
We claim $x\in P^*$ if and only if for all $q\in C_x$, $q$ has property $\star$ with respect to $\mathbb P$.
\end{lemma}
\textbf{Proof:}
One direction is simple, as for $x\in P^*$ and $q\in C_x$, \autoref{lem:voronoisubcell} shows $q\in C^*_x$, and thus by \autoref{lem:cutlocvoronoi}, $q$ has property $\star$.
For the other direction, let $x\in P\setminus P^*$ and $q\in C_x$, and \aFSOC $q$ has property $\star$.
Since the Voronoi diagram $V^*$ covers $F_*$, we can choose some $x^*\in P^*$ such that $q\in C_x^*$.
Now $x^*$ must also contribute to $V$, as $P^*\subseteq P$.
Thus, by definition of Voronoi diagram, the length of $\overline{xq}$ is at most the length of $\overline{x^*q}$.
Since $f$ and its inverse preserve the lengths of paths, the preimage of $\overline{xq}$ is at most the length of the preimage of $\overline{x^*q}$.
Then the preimage of $\overline{xq}$ is a geodesic, and $x$ should be included in $P^*$, which contradicts our assumption.
\proofclose{\autoref{lem:starcorrectcorr}}

The condition $\forall q\in C_x$, $q$ has property $\star$ is exactly the condition checked by \autoref{alg:filter}.
Thus, with $P=\{f(p):(\mathbb P,f)\in S\}$, if the first iteration of the loop (line \ref{line:cutlocalgloop}) does not run, every $x\in P$ either is in $P^*$ or has $C_x$ empty.
This implies the cut locus restricted to face $F_*$ is identical to $V^*$, and we return the correct result. 
On the other hand, assume $\exists x\in P$ with associated face path $\mathbb P$ and path unfolding $f$ such that (\autoref{alg:filter})($\mathbb P, f, x, C_x$) is false. 
Then by \autoref{lem:starcorrectcorr}, $x\notin P^*$, so when we remove this element in line$~$\ref{line:cutlocalgremoval}, we maintain that $P^*\subseteq P$.

Then it is always true that $P^*\subseteq P$, and we exit the loop only if our output is correct.
Since the algorithm must terminate, it will always return the correct output.
Before the final line, we may remove all $(\mathbb P,f)\in S$ such that $C_{f(p)}\cap V$ is empty, as this will simplify $S$ without changing the validity of the algorithm.

Line \ref{line:cutlocalgcreateS} can be done with an exhaustive search over the finite graph described in \autoref{def:facewalkfacepath}. For lines$~$\ref{line:cutlocalgvoronoicreate} and \ref{line:cutlocalgvoronoirecreate}, there are many existing algorithms \cite{voronoibook,voronoicalcparallel,voronoicalcrandinc,voronoicalcsweep,voronoicalchulls,voronoicalcenvelopes} to find Voronoi diagrams in $\mathbb R^2$.

The output of the algorithm is can be visualized as in \autoref{fig:tetra_locus}(b, c), where we have all points and path unfoldings contributing to the cut locus on a particular face.
We may use this algorithm on all faces of the polyhedron to produce a result as in \autoref{fig:tetra_locus}(a).
We may also plot relevant path unfoldings starting from the source face $F_0$, and transform the cut loci on each face appropriately, as in \autoref{fig:tetra_locus}(d).
After removing unnecessary edges, this becomes a \textbf{Voronoi star unfolding} of the polyhedron, described in \cite{star}.

Since line \ref{line:cutlocalgcreateS} enumerates all possible face paths between two faces of $X$, it has exponential complexity wrt. the number of faces in $X$. 
Thus, \autoref{alg:cutlocus} must also be at least exponential complexity.
This will limit its usability for polyhedra with many faces.
With our implementation, we find that it is feasible to calculate cut loci on the icosahedron and dodecahedron within a few seconds.

\section{Examples}\label{section:examples}
\subsection{\textit{n}-Torus}\label{subsec:ntorus}
We will use \autoref{thm:quite_simplex} to show the geodesic complexity of an $n$-torus equipped with the flat metric is at least $n$ (any geodesic motion planner needs at least $n+1$ sets).
This result was initially shown by Recio-Mitter \cite{geo_complex}, and we use this as the first example since the collections of simplex embeddings we use are instructive and easy to visualize.

We represent the $n$-torus $X:=(\Sone)^n$ with elements in $(\mathbb R/\mathbb Z)^n$, identifying $(\dots,x,\dots)\sim(\dots,x+m,\dots)$ for $m\in\mathbb Z$ in each dimension.
All geodesics between $p,q\in X$ are diagonal segments that project to a geodesic in each dimension.
We say a path from $p$ to $q$ goes \lrq{down} from $p$ in dimension $k$ if the image in this dimension contains interval $(p_k-\delta,p_k]$ for some $\delta>0$. 
We say the path goes \lrq{up} from $p$ otherwise.

Pick any $p\in X$, let $q:=p+(0.5,\dots,0.5)$ be the antipodal point, and take some $0<\varepsilon\ll 0.5$.
For each $i\in \{0,\dots,n\}$, we will construct $\mathcal F_i'$ as embeddings to all axis-oriented $\varepsilon$-boxes of dimension $i$ with $q$ as a vertex.
Formally, we build each $\mathcal F_i'$ as follows: 
For some set $D\subseteq [n]$ such that $|D|=i$, and \lrq{direction} assignments $M\colon D\to \{0,1\}$, we will construct an $i$-dimensional axis-oriented box with sides $B_{D,M,k}\subseteq \Sone$ for $k\in [n]$.
For $k\notin D$, $B_{D,M,k}=\{q_k\}$.
For $k\in D$, $B_{D,M,k}=\left[q_k,q_k+\varepsilon\right]$ if $M(k)=1$, and $\left[q_k-\varepsilon,q_k\right]$ otherwise.
We construct $\mathcal F_i'$ by choosing an embedding for each possible $B_{D,M}:=B_{D,M,1}\times\dots\times B_{D,M,n}$.
Whenever we choose an embedding from $\Delta_i$ to an $i$-dimensional box $B_{D,M}$, we map $i$ of the simplex faces as the $i$ box faces that contain $q$.
This ensures restrictions to these faces are represented in $\mathcal F'_{i-1}$.

Finally, for each $i\in\{0,\dots,n\}$, we construct $\mathcal F_i:=\{\textbf{t}\mapsto (p, F(\textbf{t})) : F\in \mathcal F_i'\}$.
We display the images of $F\in \mathcal F_i'$ for the 2-torus in \autoref{fig:torus_proof}(a).
We will verify that the hypotheses of \autoref{thm:quite_simplex} hold.

\renewcommand{\imgwidth}{.25\linewidth}
\begin{figure}[ht!]
\centering
\subfigure[Sets $\mathcal F_i'$]{
\includegraphics[width=\imgwidth]{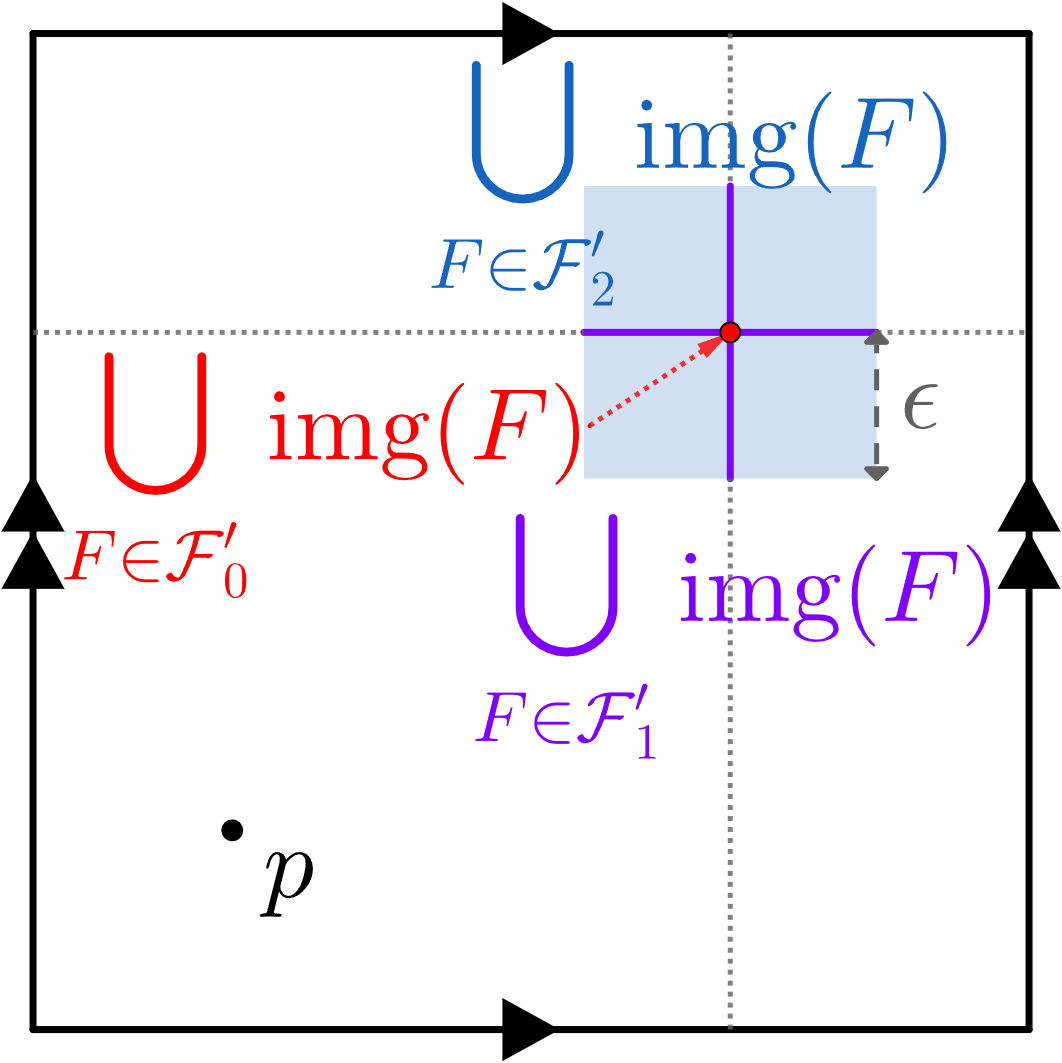}
}
\subfigure[Examples of \gls{GMPR}s on some sets in $\mathcal F_i$]{
\includegraphics[width=\imgwidth]{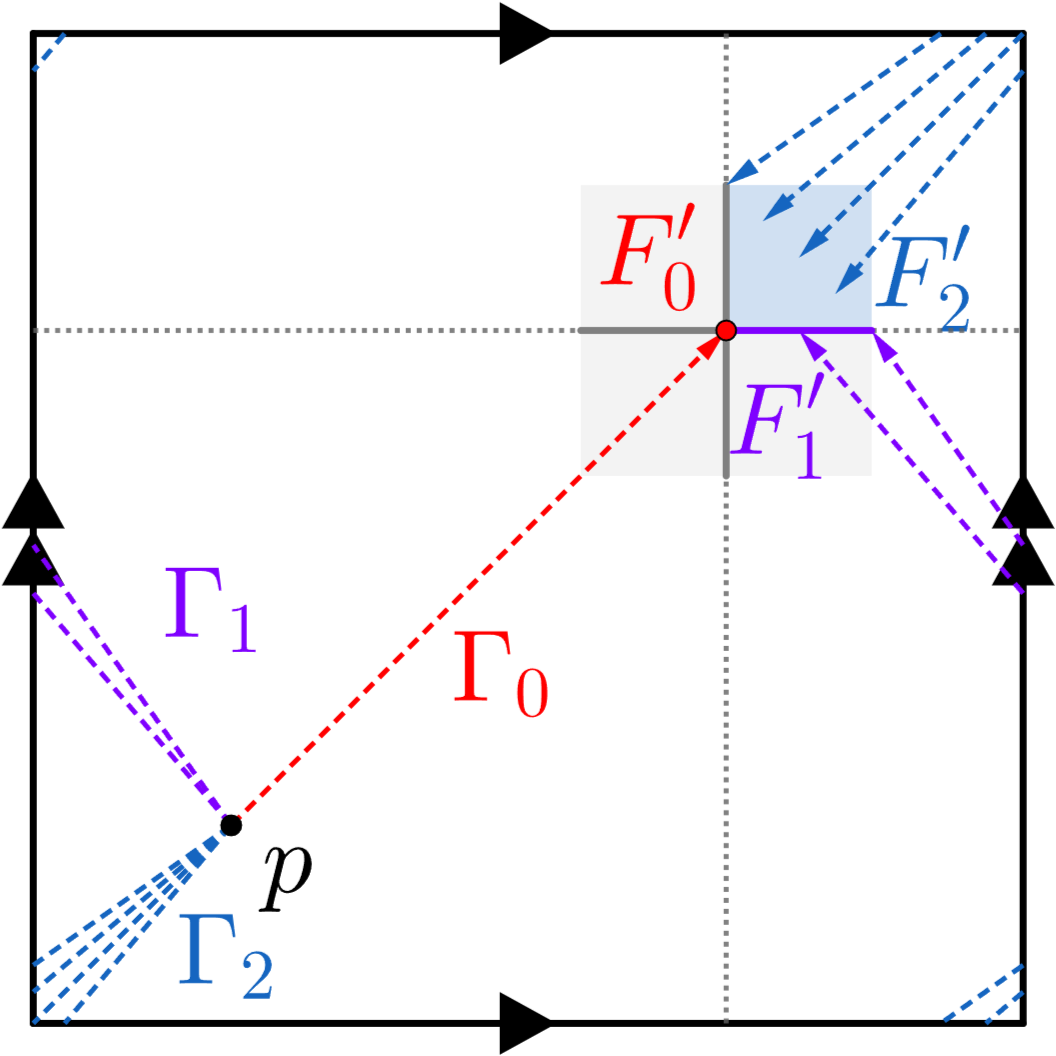}
}
\caption{Visualization of sets used in the 2-torus}
\label{fig:torus_proof}
\end{figure}

\begin{compactenum}[(a)]
    \item
    By construction, each $\mathcal F_i$ contains embeddings of $i$-dimensional simplices.
    \item
    Take some $F_i\in \mathcal F_i$.
    \WLOG $\img{F_i}=\{p\}\times \left(\prod\limits_{k=1}^i[q_k,q_k+\varepsilon]\times \prod\limits_{k=i+1}^n\{q_k\}\right)$, as the other cases are symmetric.
    We will characterize all \gls{GMPR}s on $\img{F_i}$.
    Since geodesics in $X$ are geodesics in each dimension, the \gls{GMPR} must go \lrq{down} from $p$ in dimensions $1,\dots,i$, and can choose a direction for the others.
    This creates $2^{n-i}$ possible \gls{GMPR}s, distinguished by these choices.
    Continuity follows from continuity of projections to each dimension.

    Let $(p,q')\in \relint{F_i}$, and take geodesic $G$ between them.
    We have $q'=((q_k+\delta_k)_{k\in [i]},(q_k)_{k>i})$ for $0<\delta_k<\varepsilon$, so $G$ must go \lrq{down} from $p$ in dimensions $1,\dots,i$.
    Thus, there exists a \gls{GMPR} on $\img{F_i}$ extending $G$.
    \item
    Take the same $F_i\in \mathcal F_i$, and take a \gls{GMPR} $\Gamma_i$ on $\img{F_i}$ that \wLOG goes \lrq{up} from $p$ in dimension $i+1$.
    We will show that $F_{i+1}\in\mathcal F_{i+1}$ where $\img{F_{i+1}}=\{p\}\times \left(\prod\limits_{k=1}^{i+1}[q_k,q_k+\varepsilon]\times \prod\limits_{k=i+2}^n\{q_k\}\right)$ satisfies our requirements.
    \begin{compactenum}[(i)]
        \item
        $F_{i+1}\in\mathcal F_{i+1}$ and has $F_i$ as a face by construction.
        \item
        Pick some geodesic $G$ from $\Gamma_i$.
        We will show there is an open set separating $G$ from any geodesic in $\pi_{GX}^{-1}(\relint{F_{i+1}})$, which all go \lrq{down} from $p$ in dimension $i+1$.
        By our choice of $q$ and $\varepsilon$, the image of any of these geodesics in dimension $i+1$ intersects $(q_{i+1},q_{i+1}+\varepsilon)$ and does not intersect $(q_{i+1}-\varepsilon,q_{i+1})$. 
        In contrast, the image of $G$ in dimension $i+1$ intersects only $(q_{i+1}-\varepsilon,q_{i+1})$.
        Let $t_0\in[0,1]$ such that $G(t_0)_{i+1}\in (q_{i+1}-\varepsilon,q_{i+1})$.
        We construct an open set $U\subseteq GX$ via $U=\{G'\in GX:G'(t_0)_{i+1}\in (q_{i+1}-\varepsilon,q_{i+1})\}$.
        From the discussion above, $G\in U$, while no geodesics in $\pi_{GX}^{-1}(\relint{F_{i+1}})$ are in $U$.
        \item
        There are $2^{n-(i+1)}$ \gls{GMPR}s on $\img{F_{i+1}}$, distinguished by choice of direction in dimensions $i+2,\dots,n$.
    \end{compactenum}
    These show hypothesis (c).
    We visualize an example on the 2-torus in \autoref{fig:torus_proof}(b), where for step $0\leq i<2$, we arbitrarily choose $\Gamma_i$, one of $2^{2-i}$ possible \gls{GMPR}s on the image of $F_i\in \mathcal F_i$, and choose an $F_{i+1}$ that \lrq{avoids} $\Gamma_i$.
    We display each $\Gamma_i$, as well as $\Gamma_2$, the only \gls{GMPR} on $\img{F_2}$.
    With our construction, $\img{F_{i}}=\{p\}\times F_{i}'$ for some axis-oriented box $F_{i}'$, so we also display each $F'_i$.
\end{compactenum}

Thus, $GC((\Sone)^n)\geq n$, so any geodesic motion planner needs at least $n+1$ sets.
We have $GC((\Sone)^n)= n$ from the explicit geodesic motion planner by Cohen and Pruidze \cite{torus}.

\subsection{Tetrahedron}
We show how the results of Davis \cite{tetra} fit into this framework, and show that the geodesic complexity of a tetrahedron is at least three (needing at least four sets in any geodesic motion planner).

\subsubsection{Cut Locus}
\label{subsubsec:tetra_cut_locus_proof}

First, we will calculate the cut locus of a point $p$ on the tetrahedron $X$.
We represent all faces as equilateral triangles with side length $2\sqrt3$, oriented as in \autoref{fig:tetra_locus}(a) with the face's center point at $(0,0)$.
We label the faces Face 0,\dots , Face 3 as in \autoref{fig:tetra_locus}(a).
When considering walk unfoldings to a particular face (as in \autoref{fig:tetra_locus}(b, c, d)), we use the walk unfolding that preserves the center point and orientation of that face.

\renewcommand{\imgheight}{110 pt}
\begin{figure}[ht!]
\centering
\begin{minipage}{.5\textwidth}
  \centering
    \includegraphics[height=\imgheight]{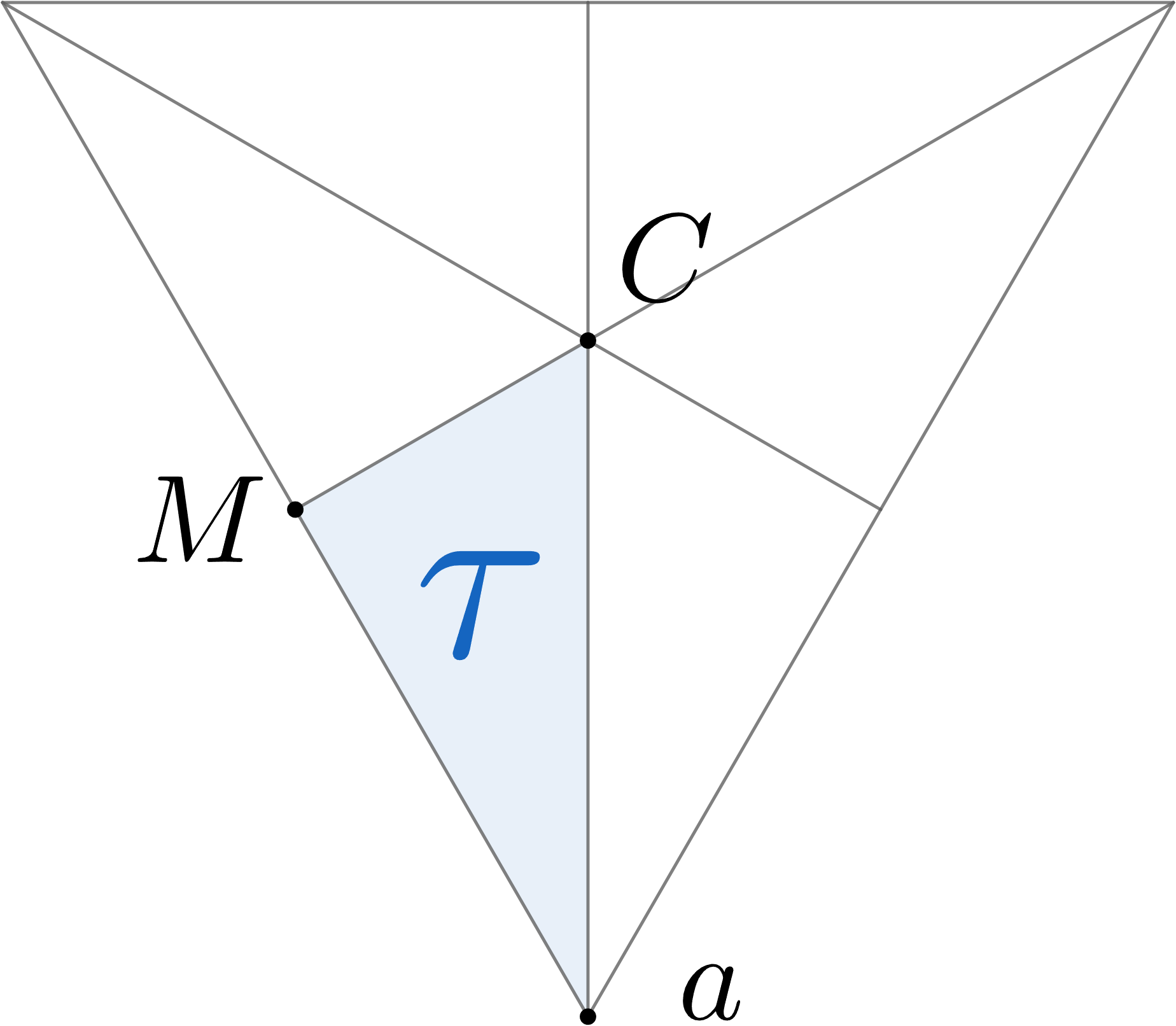}
    \caption{Representative region of the tetrahedron}
    \label{fig:dihedral_smmetry}
\end{minipage}\begin{minipage}{.5\textwidth}
  \centering
    \includegraphics[height=\imgheight]{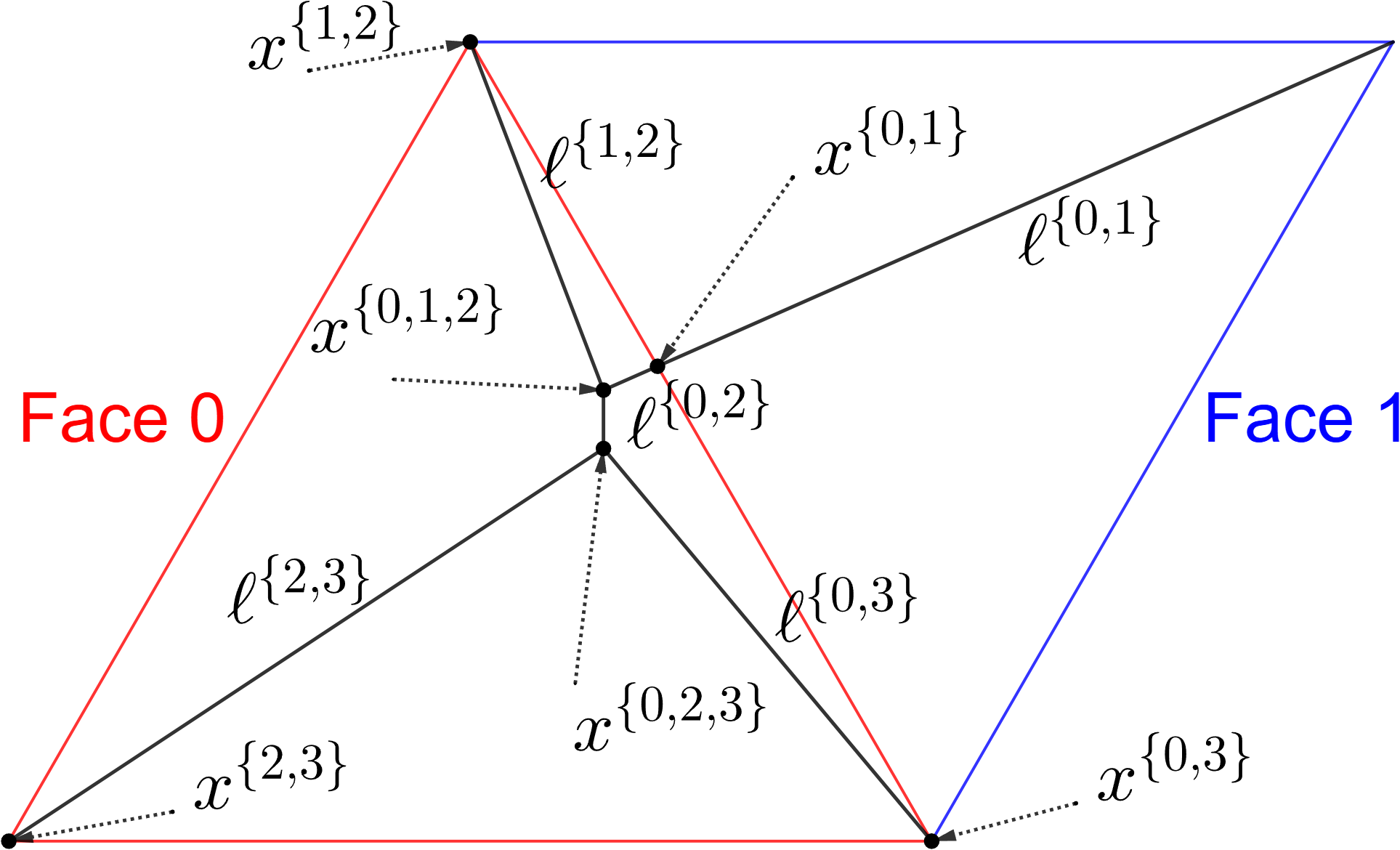}
    \caption{Labeled $\ell^{I}$ and $x^{I'}$ on cut locus of
    \\$p=(-\frac12,-\frac{\sqrt3}2)\in \text{int}(\tau)$}
    \label{fig:tetralabelexample}
\end{minipage}
\end{figure}

\renewcommand{\imgheight}{150 pt}
\begin{figure}[ht!]
\centering
\subfigure[Labeling of tetrahedron faces, boundaries between them, and the cut locus of $p$ on each of them]{
\includegraphics[height=\imgheight]{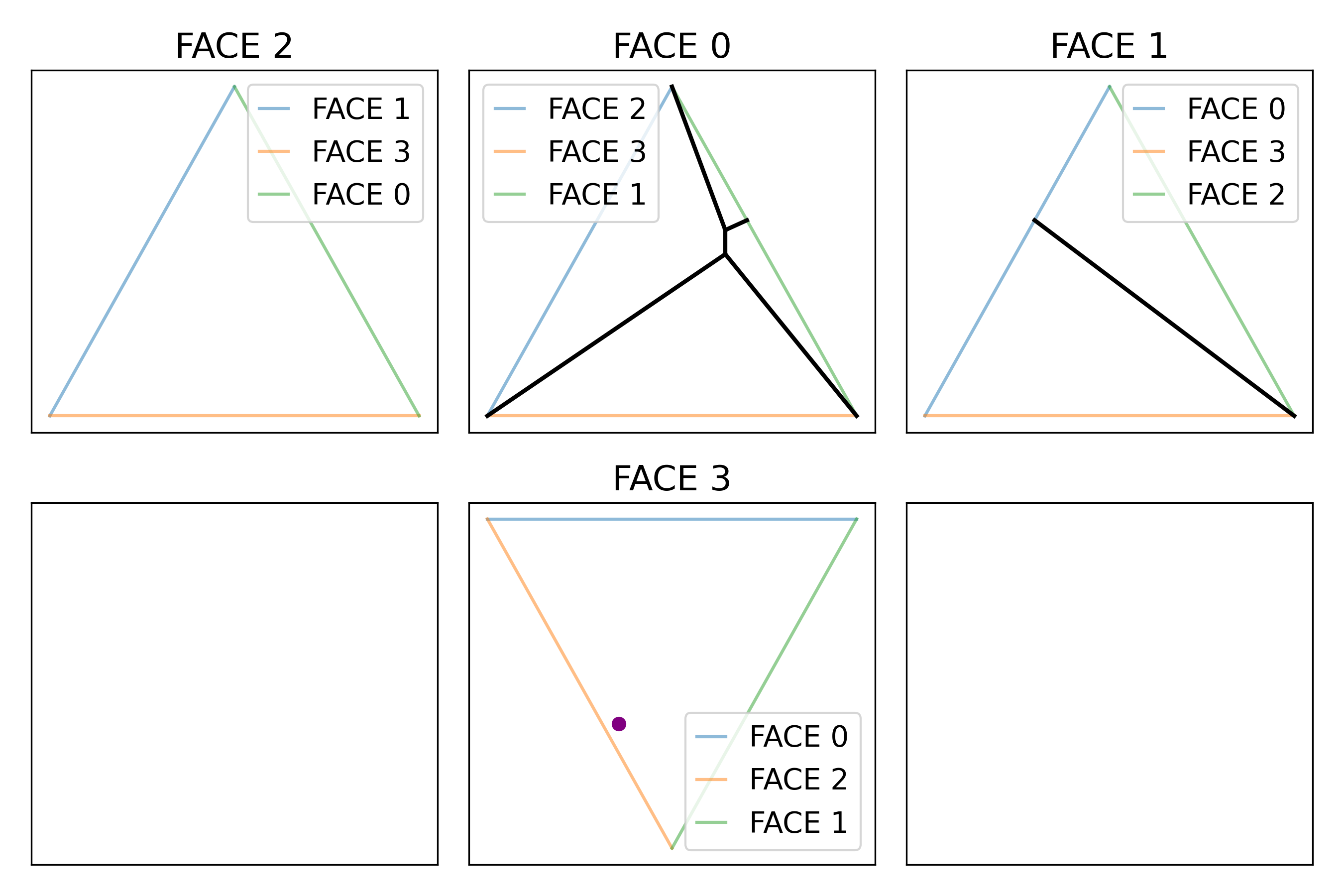}}
\quad
\subfigure[Copies of $p$ and path unfoldings that contribute to the cut locus on Face 0]{
\includegraphics[height=\imgheight]{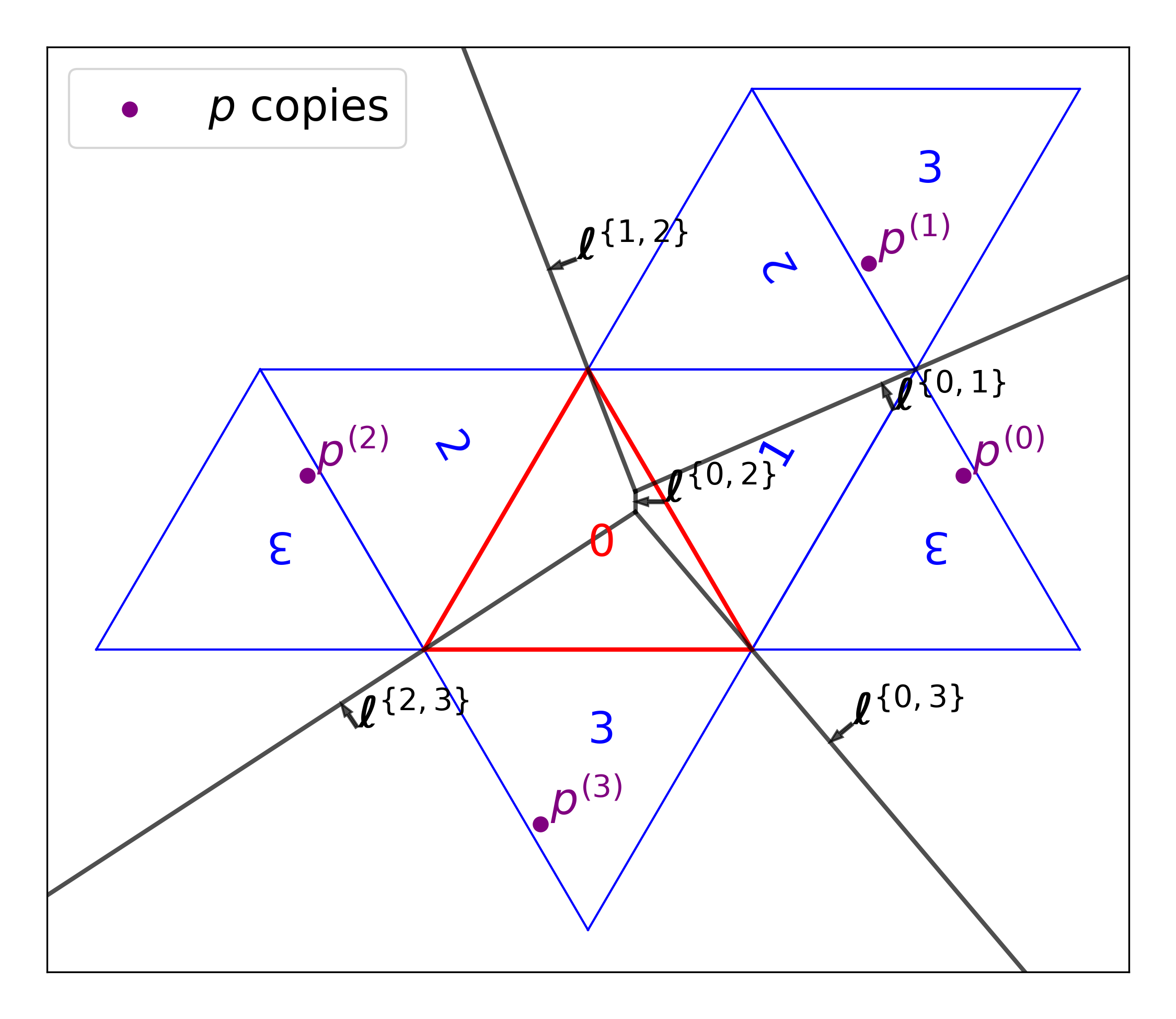}
}
\subfigure[Copies of $p$ and path unfoldings that contribute to the cut locus on Face 1]{
\includegraphics[height=\imgheight]{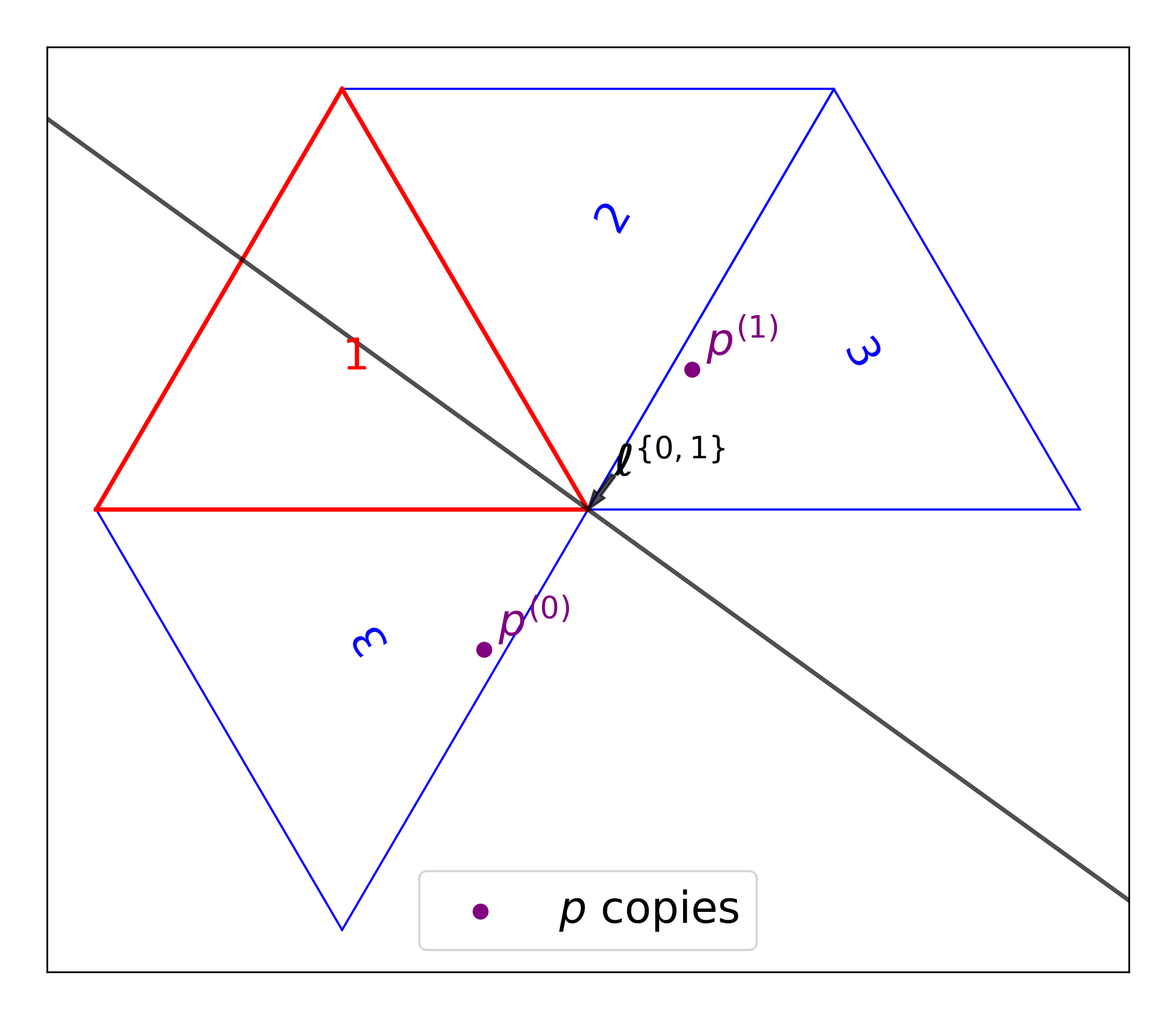}
}
\subfigure[Voronoi star unfolding from $p$]{
\includegraphics[height=\imgheight]{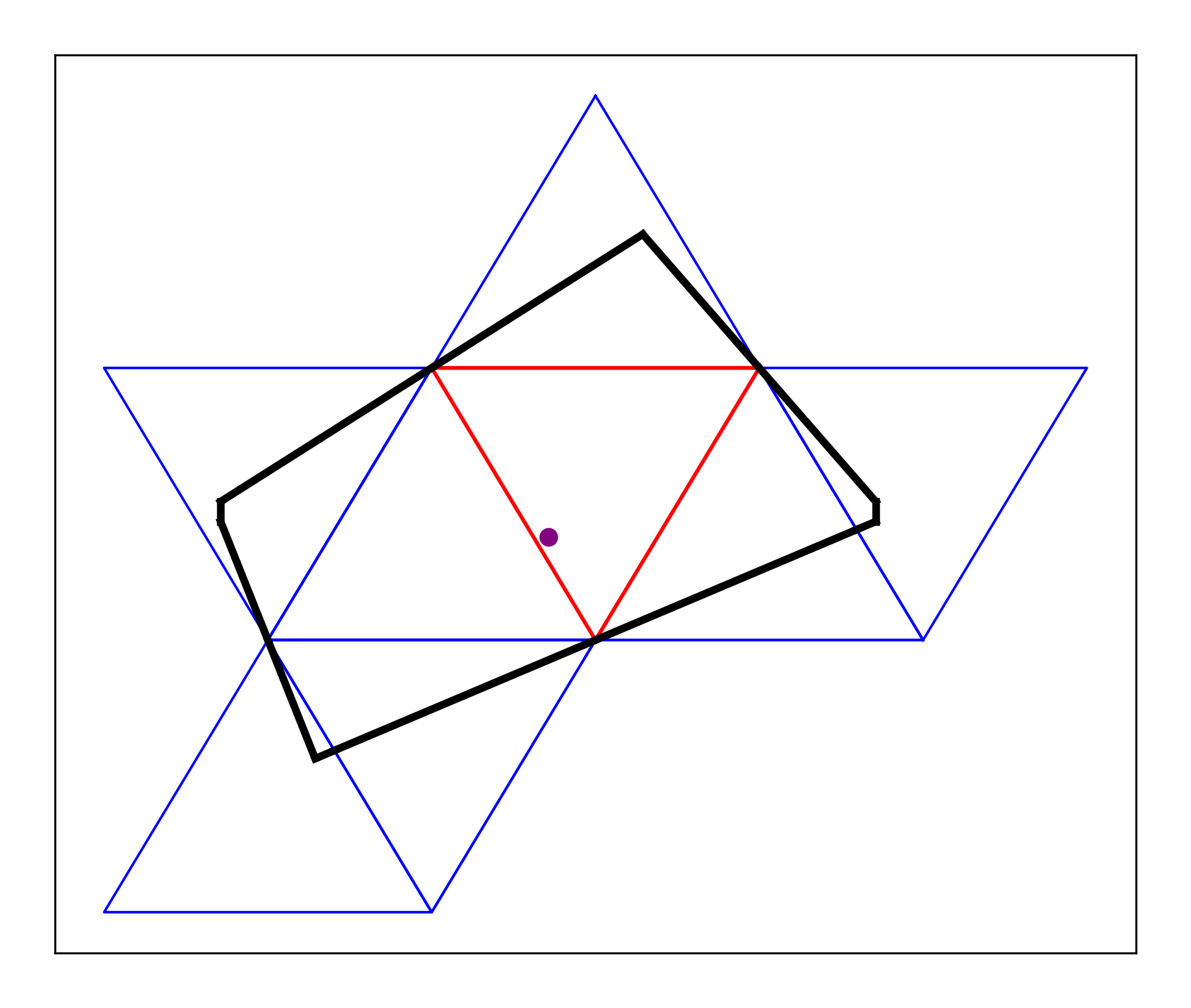}
}
\caption{Output of \autoref{alg:cutlocus} on a tetrahedron for $p=(-\frac12,-\frac{\sqrt3}2)\in \text{int}(\tau)$ on Face 3}
\label{fig:tetra_locus}
\end{figure}

Up to symmetry, a point on one of six congruent triangles on a face (displayed in \autoref{fig:dihedral_smmetry}) represents any point on the tetrahedron surface.
Thus, we need only calculate the cut locus on one triangle.
We choose closed set $\tau$ as the bottom left region of Face 3 with vertices $a$, $C$, and $M$, and first consider $p\in \text{int}(\tau)$.
From our implementation of \autoref{alg:cutlocus}, we obtain the cut locus on each face (displayed in \autoref{fig:tetra_locus}), as well as the path unfoldings that contribute to them.
We find the cut locus appears on Faces 0 and 1.

Let $R_\theta:=\begin{bmatrix}
    \cos\theta&-\sin\theta\\
    \sin\theta&\cos\theta
\end{bmatrix}$ be the rotation matrix, so $R_\theta p$ represents a rotation of $p$ by $\theta$.
From \autoref{alg:cutlocus}, we find the cut locus on Face 0 of $p$ arises from the Voronoi diagram of four points, displayed in \autoref{fig:tetra_locus}(b).
We enumerate these points $(p^{(0)},p^{(1)},p^{(2)},p^{(3)})=
\left(
\begin{bmatrix}
    2\sqrt3\\0
\end{bmatrix}+
    R_\pi p,
\begin{bmatrix}
    2\sqrt3\\4
\end{bmatrix}+
p,
\begin{bmatrix}
    -2\sqrt3\\0
\end{bmatrix}+
R_\pi p,
\begin{bmatrix}
    0\\-2
\end{bmatrix}+
p
\right)$.

We denote $\ell^{\{i,j\}}(p)$ as the line that bisects $p^{(i)}$ and $p^{(j)}$.
Note that which $\ell^{\{i,j\}}(p)$ exist on the cut locus depend on the choice of $p$.
This is displayed in \autoref{fig:tetra_first} (a) and (c), where $\ell^{\{0,2\}}$ is replaced by $\ell^{\{1,3\}}$.

We denote $x^{\{i,j\}}(p)$ as a specified intersection of $\ell^{\{i,j\}}(p)$ with an edge of Face 0, and $x^{\{i,j,k\}}(p)$ as the intersection of $\ell^{\{i,j\}}(p)$ and $\ell^{\{j,k\}}(p)$.
If two points coincide $x^I(p)=x^{I'}(p)$, we call this point $x^{I\cup I'}(p)$.
We choose our notation so that $x^I(p)$ is equidistant from points $(p^{(i)})_{i\in I}$.
We display $\ell^{I}$ and $x^{I'}$ for a point $p\in\text{int}(\tau)$ in \autoref{fig:tetralabelexample}.
We will calculate the locations of some of these lines and intersecting points.
\begin{compactitem}
    \item
    $\ell^{\{0,1\}}(p)$ has equation 
    $t\mapsto (2\sqrt3,2)+((-2,0)+R_{\pi/2}p)t$.

    Let $x^{\{0,1\}}(p)$ be the intersection of $\ell^{\{0,1\}}(p)$ with the right edge of Face 0.
    
    \item
    $\ell^{\{1,2\}}(p)$ has equation 
    $t\mapsto (0,2)+((-2,2\sqrt3)+R_{\pi/2}p)t$.

    $\ell^{\{1,2\}}(p)$ always intersects the top vertex of the triangle.
    Let this point be $x^{\{1,2\}}(p):=(0,2)$.

    \item
    $\ell^{\{2,3\}}(p)$ has equation 
    $t\mapsto (-\sqrt3,-1)+((1,\sqrt3)+R_{\pi/2}p)t$.

    $\ell^{\{2,3\}}(p)$ always intersects the left vertex of the triangle.
    Let this point be $x^{\{2,3\}}(p):=(-\sqrt3,-1)$.
    
    \item
    $\ell^{\{0,3\}}(p)$ has equation
    $t\mapsto (\sqrt3,-1)+((1,-\sqrt3)+R_{\pi/2}p)t$.

    $\ell^{\{0,3\}}(p)$ always intersects the right vertex of the triangle.
    Let this point be $x^{\{0,3\}}(p):=(\sqrt3,-1)$.

\end{compactitem}

Using SymPy, a symbolic calculator \cite{sympy}, we find the other $x^I$ by calculating intersections of lines:
\begin{compactitem}
    \item
    $x^{\{0,1\}}(p)=
    \left(\frac{2\sqrt{3}p_1}{p_1 - \sqrt{3}p_2 - 2\sqrt{3}}, 2\frac{2p_1 + \sqrt{3}p_2 + 2\sqrt{3}}{-p_1 + \sqrt{3}p_2 + 2\sqrt{3}}\right)$.
    
                            \item
    $x^{\{0,1,2\}}(p)=
    \left(-p_1, \frac{p_1^{2} + 2\sqrt{3}p_1 + 2p_2 + 4}{p_2 + 2}\right)$.
    \item
    $x^{\{0,2,3\}}(p)=
    \left(-p_1, \frac{p_1^{2} - p_2 - 2}{p_2 - 1}\right)
    $.
    \item
    $x^{\{0,1,3\}}(p)=
    \left(\frac{\sqrt{3}p_1p_2 - 4\sqrt{3}p_1 + 3p_2^{2} + 6p_2}{-3p_1 + \sqrt{3}p_2 + 2\sqrt{3}}, \frac{-\sqrt{3}p_1^{2} - 3p_1p_2 + 2\sqrt{3}p_2 + 4\sqrt{3}}{-3p_1 + \sqrt{3}p_2 + 2\sqrt{3}}\right)$.
    
    \item
    $x^{\{1,2,3\}}(p)=
    \left(\frac{-\sqrt{3}p_1p_2 - 2\sqrt{3}p_1 - 3p_2^{2} - 6p_2}{3p_1 - \sqrt{3}p_2 + 4\sqrt{3}}, \frac{\sqrt{3}p_1^{2} + 3p_1p_2 + 12p_1 + 4\sqrt{3}p_2 + 8\sqrt{3}}{3p_1 - \sqrt{3}p_2 + 4\sqrt{3}}\right)$.
\end{compactitem}

Using our calculated expressions, we will show that the cut locus structure in \autoref{fig:tetra_locus}(b) is the same across the interior of $\tau$.
We do this by verifying the following are true:
\begin{compactitem}
    \item
    Lines $\ell^{\{0,1\}}$, $\ell^{\{1,2\}}$, $\ell^{\{2,3\}}$, and $\ell^{\{0,3\}}$ exist on the cut locus and are rays in the corresponding Voronoi diagram.
    \item
    Lines $\ell^{\{0,1\}}$ and $\ell^{\{1,2\}}$ end at point $x^{\{0,1,2\}}$, and this point is distinct from other $x^{\{0,1,k\}}$ or $x^{\{1,2,k\}}$.
    \item
    Lines $\ell^{\{2,3\}}$ and $\ell^{\{0,3\}}$ end at point $x^{\{0,2,3\}}$, and this point is distinct from other $x^{\{2,3,k\}}$ or $x^{\{0,3,k\}}$.
    \item
    Points $x^{\{0,1,2\}}$ and $x^{\{0,2,3\}}$ are on the interior of Face 0.
\end{compactitem}

We claim this is sufficient to fully prove the structure of the cut locus on Face 0, as we can account for the boundaries of all four cells.
For example, the cell corresponding to $p^{(1)}$ must be partially bounded by $\ell^{\{0,1\}}$ and $\ell^{\{1,2\}}$.
If these are both rays that begin at the same point, this must account for the entire boundary of this cell.
Similar arguments account for the full boundary of cells corresponding to $p^{(0)}$, $p^{(2)}$, and $p^{(3)}$.

To show these hypotheses, we utilize the Wolfram Mathematica symbolic calculator \cite{mathematica}, which can find regions in $\mathbb R^2$ where inequalities hold.
\begin{compactitem}
\item 
The fact that a line $\ell^{\{i,j\}}$ exists on the cut locus may be verified by simply checking that a point $m\in \ell^{\{i,j\}}$ is closer to $p^{(i)}$ and $p^{(j)}$ than any other copy of $p$.
Specifically, we verify for all $p\in\text{int}(\tau)$ that $d(p^{(i)},m)<d(p^{(k)},m)$ for all $k\neq i,j$.
We pick $m$ using our calculation of $\ell^{\{i,j\}}$ and its endpoints.
\item 
To verify that $\ell^{\{i,j\}}$ is a ray on the cut locus, it is enough to show that it contains a ray, and has an endpoint. 
We pick a point $m\in\ell^{\{i,j\}}$ and a predicted direction vector $v\in \mathbb R^2$, where we expect $\ell^{\{i,j\}}$ to contain $m+tv$ for $t\geq 0$. 
As before, we verify $m$ is on the cut locus line corresponding to $\ell^{\{i,j\}}$. 
From this, it is enough to show that for $k\neq j$, $\ell^{\{i,k\}}$ does not intersect $\{m+tv:t\geq 0\}$, since if a line segment on a Voronoi diagram in $\mathbb R^2$ ends, it must be at a vertex.
Thus, we verify for all $p\in\text{int}(\tau)$ that for all $k\neq i,j$, we have $v\cdot (x^{\{i,j,k\}}-m)<0$.
To ensure we have a ray and not a line, we must verify that an endpoint $x^{\{i,j,k\}}$ of $\ell^{\{i,j\}}$ is on the cut locus, which can be done by showing for $p\in\text{int}(\tau)$, that $d(p^{(i)},x^{\{i,j,k\}})\leq d(p^{(k')},x^{\{i,j,k\}})$ for all $k'\neq i,j,k$.
To show that this point is distinct from any other $x^{\{i,j,k'\}}$, it is enough to show this inequality is strict.
\item 
Finally, checking if a point $x^{\{i,j,k\}}$ is on the interior of each edge on Face 0 can be done using a point on the edge and the edge's normal vector.
\end{compactitem}

We find that all hypotheses hold on the interior of $\tau$, showing the cut locus on Face 0 is as expected.
We find the cut locus on Face 1 arises from the Voronoi diagram of 
two points, as in \autoref{fig:tetra_locus}(c).
We also find that the line $\ell^{\{0,1\}}$ is equivalent to the line of the same name in  \autoref{fig:tetra_locus}(b).
Thus, we can extend the same analysis to compute the entire cut locus of any $p\in\text{int}(\tau)$.

Similar calculations for each of the edges and vertices of $\tau$ show that the cut loci from \autoref{alg:cutlocus} match the results in the literature.
We display examples of these cut loci in \autoref{fig:tetra_locus_boundaries}, and in \autoref{apx:tetrapendix}, we reproduce all figures from~\cite[Sec.~2]{tetra}.
We also verify that when choosing $p$ within $aM$, the cut locus on Face 0 is created from four copies of $p$.
These copies arise from the same path unfoldings as the copies that create the cut locus within $\tau$.
Thus, as $p$ varies from the interior of $\tau$ to $aM$, two vertices of the cut locus approach each other, causing line $\ell^{\{0,2\}}$ to disappear from the cut locus (see \autoref{fig:tetra_first}(a) and (b)).

\renewcommand{\imgwidth}{.45\linewidth}
\begin{figure}[ht!]
\centering
\subfigure[$p$ chosen as $C$]{
\includegraphics[width=\imgwidth]{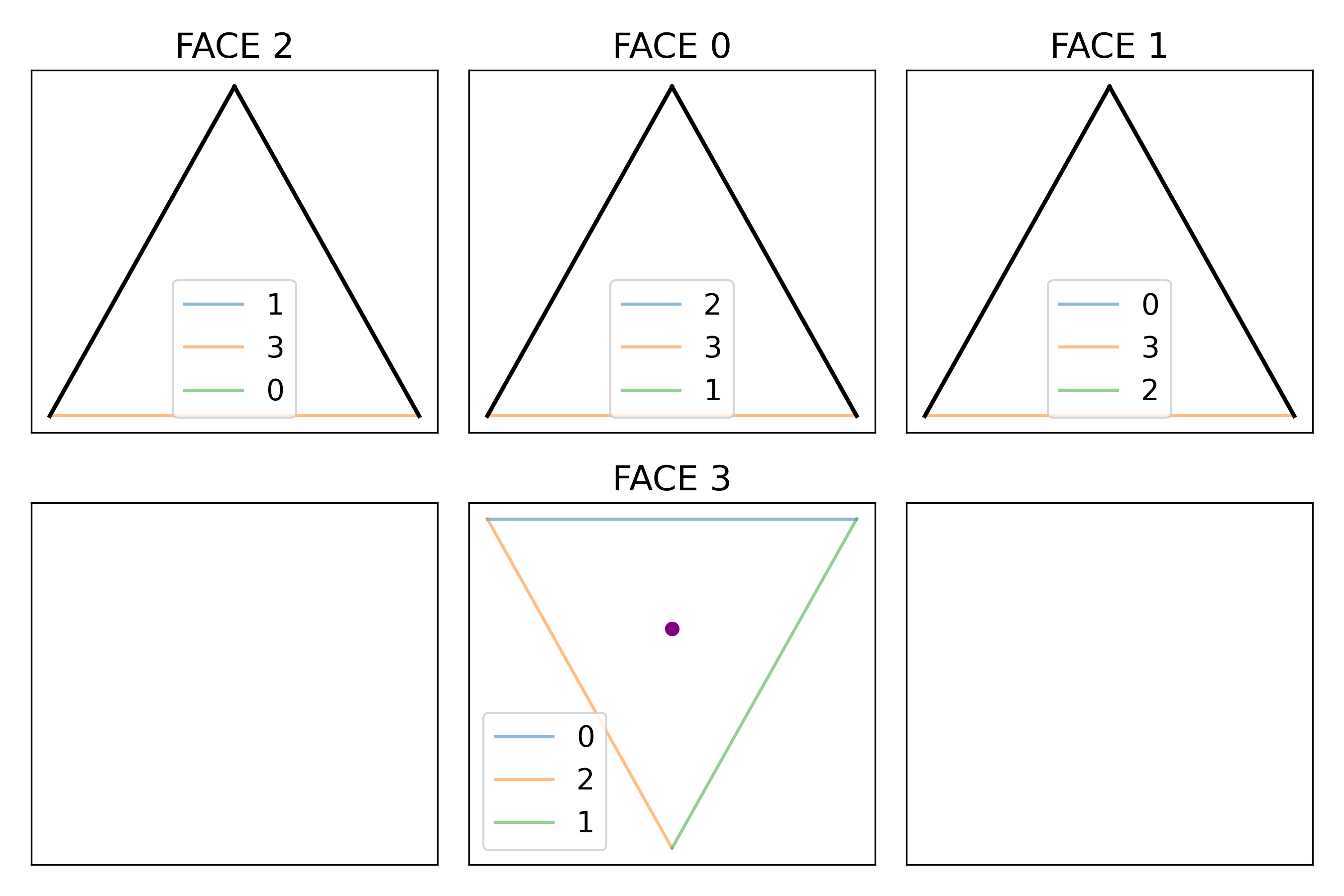}}
\subfigure[$p$ chosen within $aC$]{
\includegraphics[width=\imgwidth]{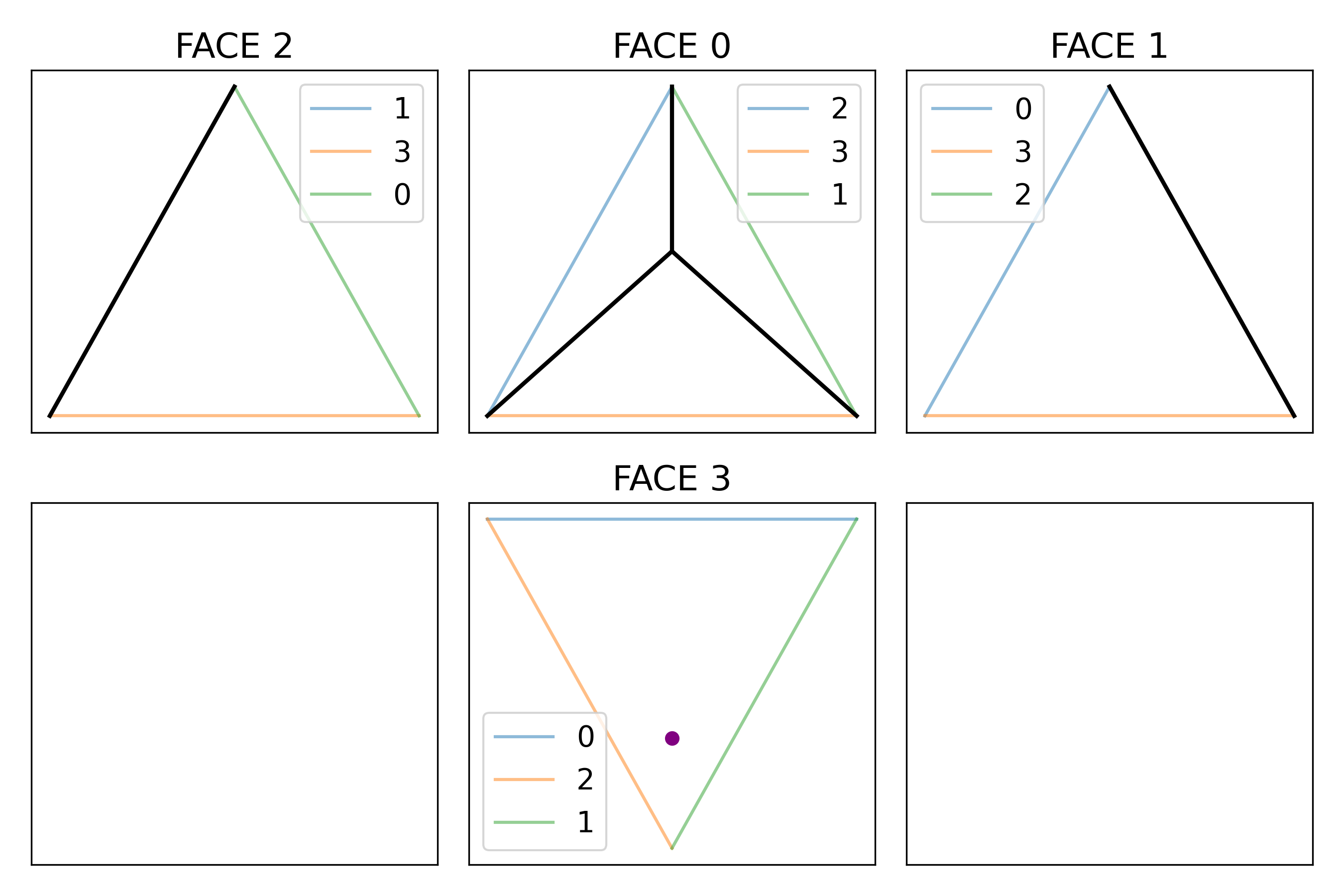}
}
\subfigure[$p$ chosen as $a$]{
\includegraphics[width=\imgwidth]{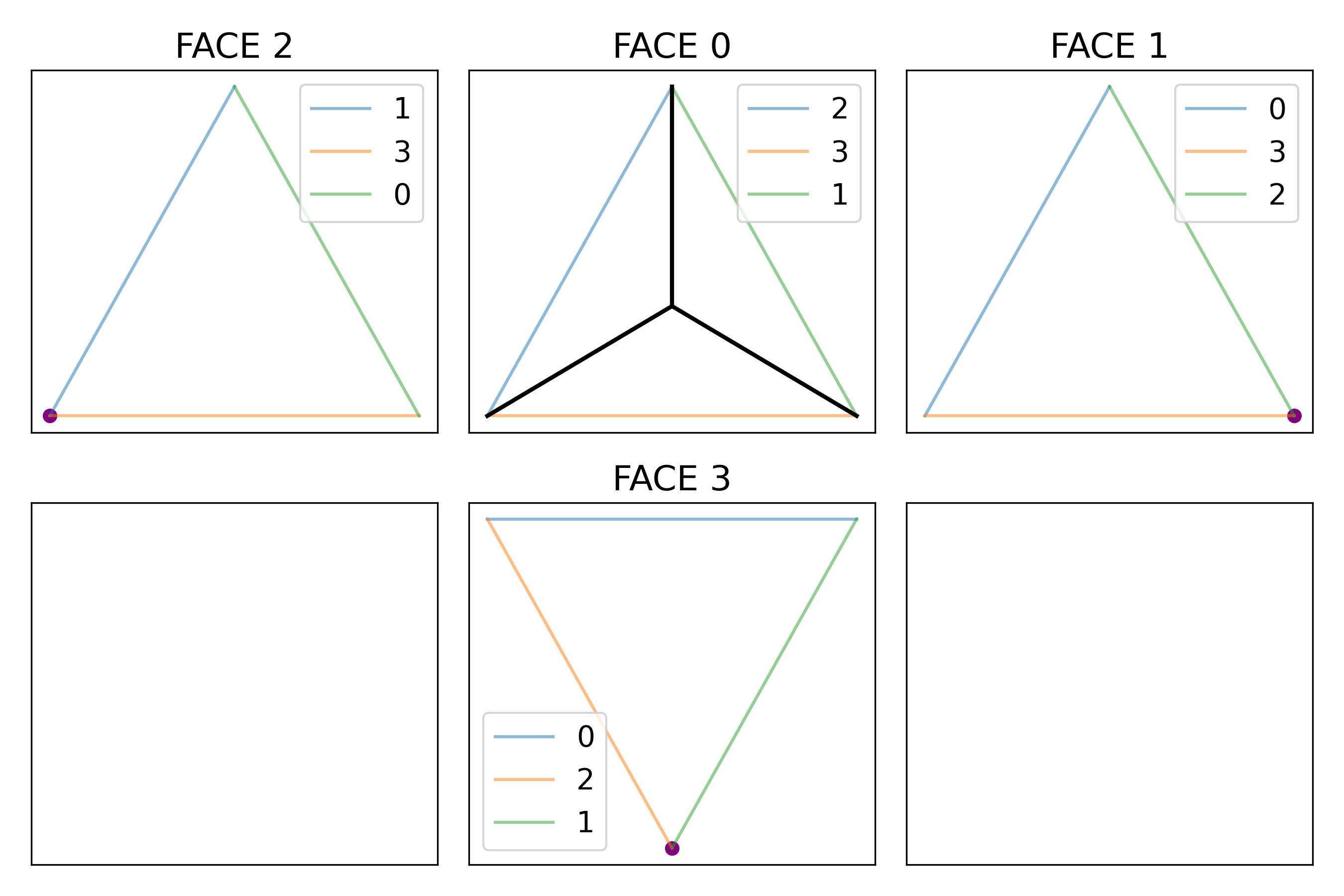}
}
\subfigure[$p$ chosen within $aM$]{
\includegraphics[width=\imgwidth]{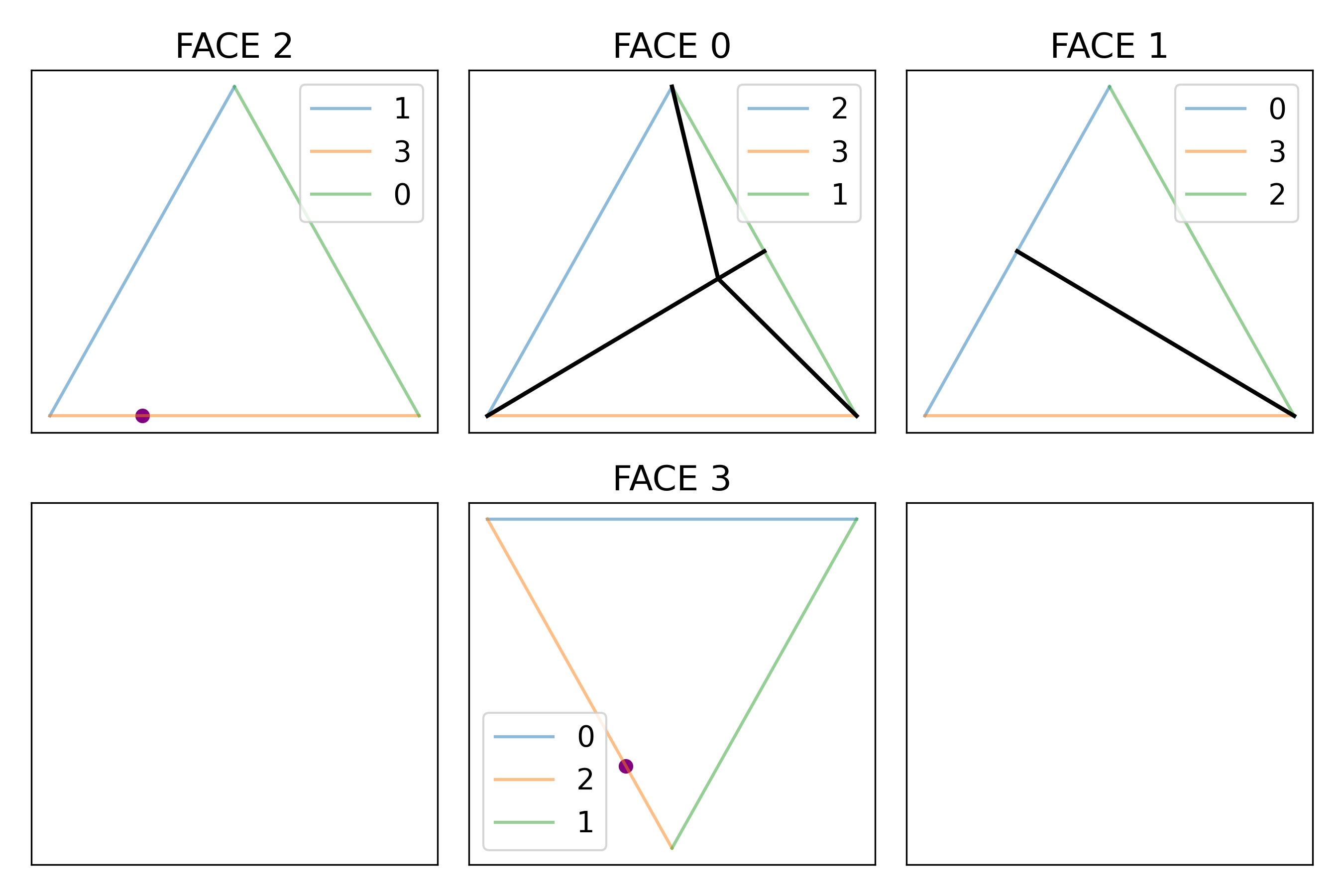}
}
\subfigure[$p$ chosen as $M$]{
\includegraphics[width=\imgwidth]{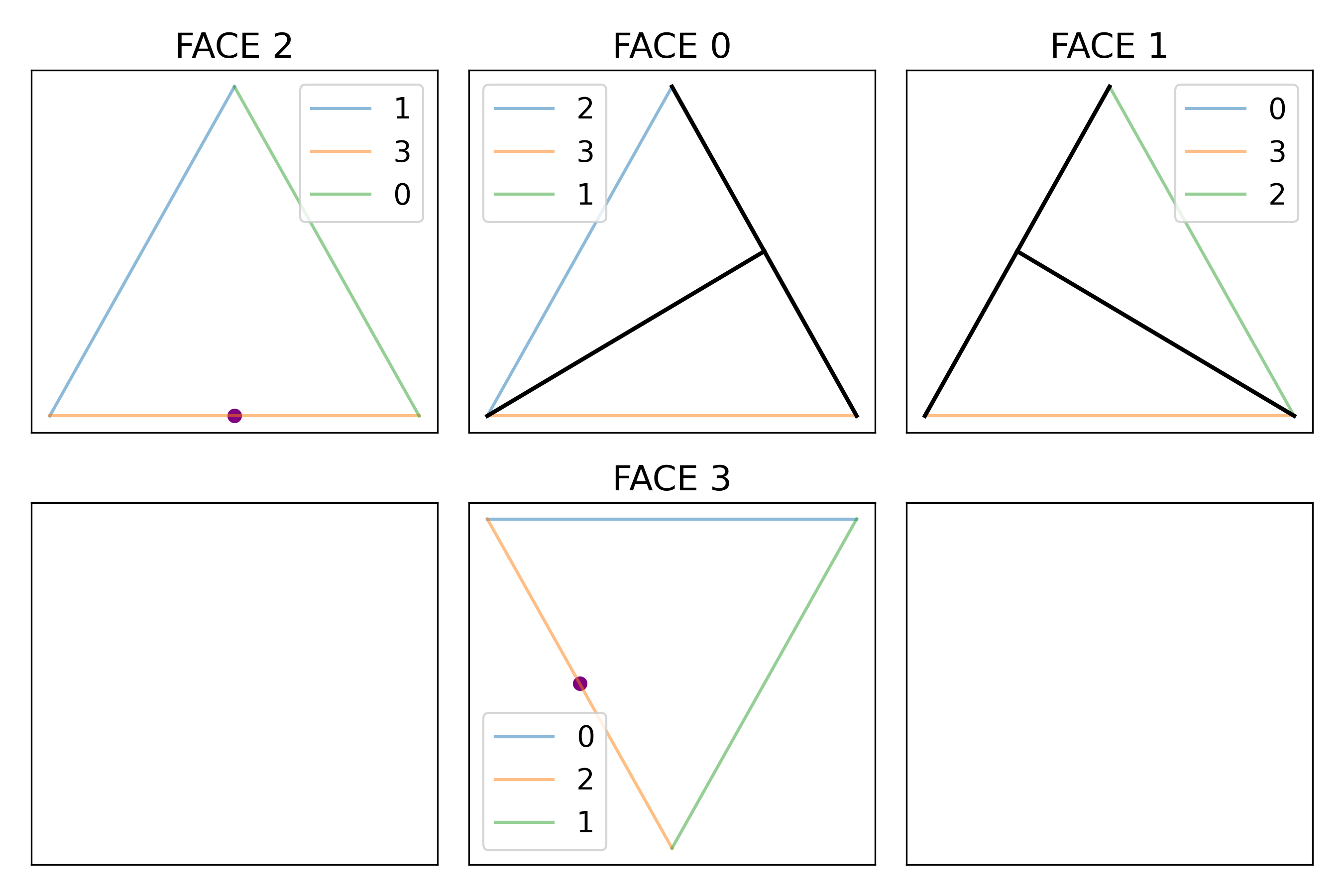}
}
\subfigure[$p$ chosen within $CM$]{
\includegraphics[width=\imgwidth]{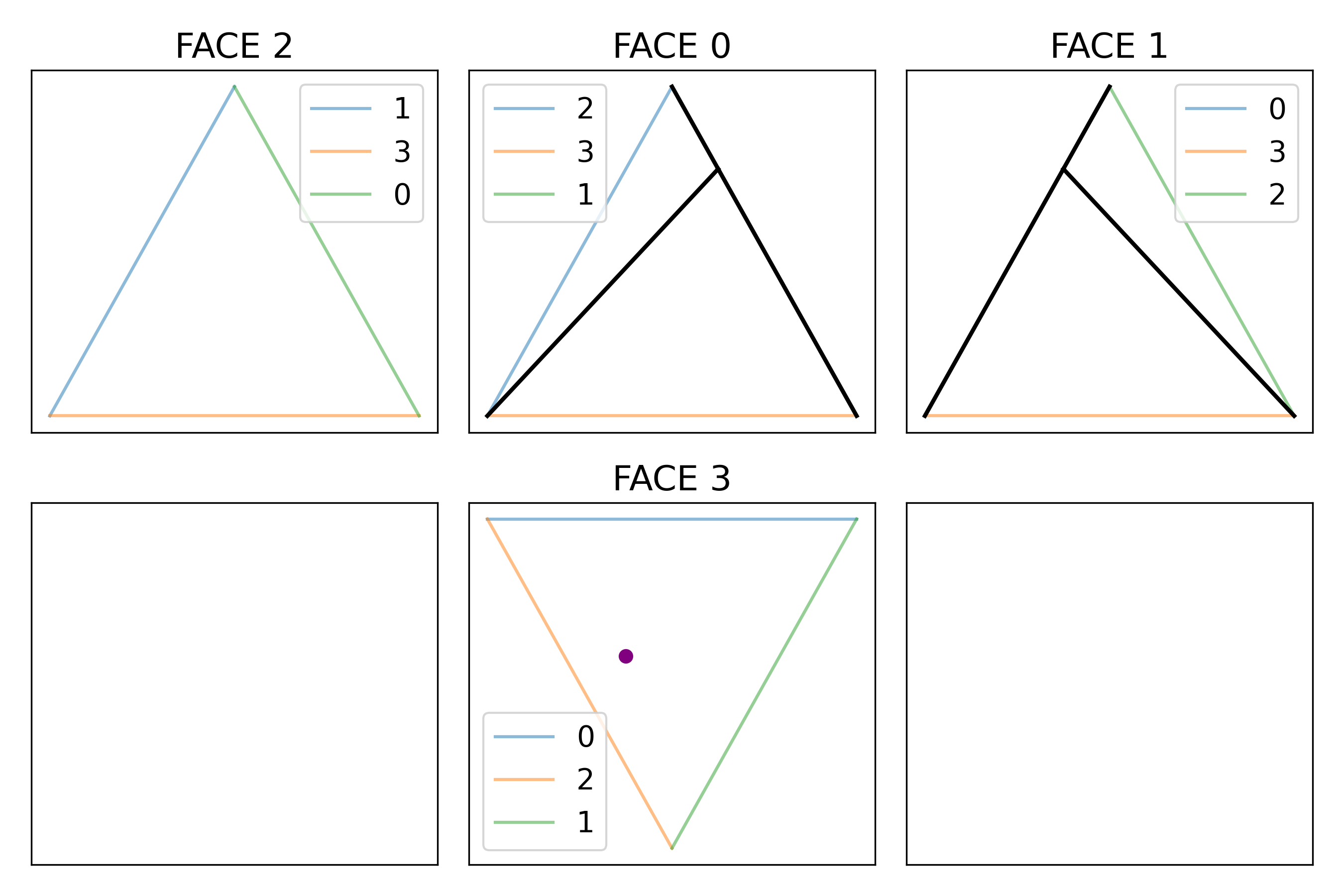}
}
\caption{Cut locus of $p$ chosen on the boundary of $\tau$}
\label{fig:tetra_locus_boundaries}
\end{figure}

\subsubsection{Geodesic Complexity Lower Bound}
\label{subsubsec:tetra_lower_bound_proof}

\renewcommand{\imgwidth}{150 pt}
\begin{figure}[htp!]
    \centering
    \subfigure[$p$ chosen within $\tau$ on Face 3]{
    \includegraphics[height=\imgheight]{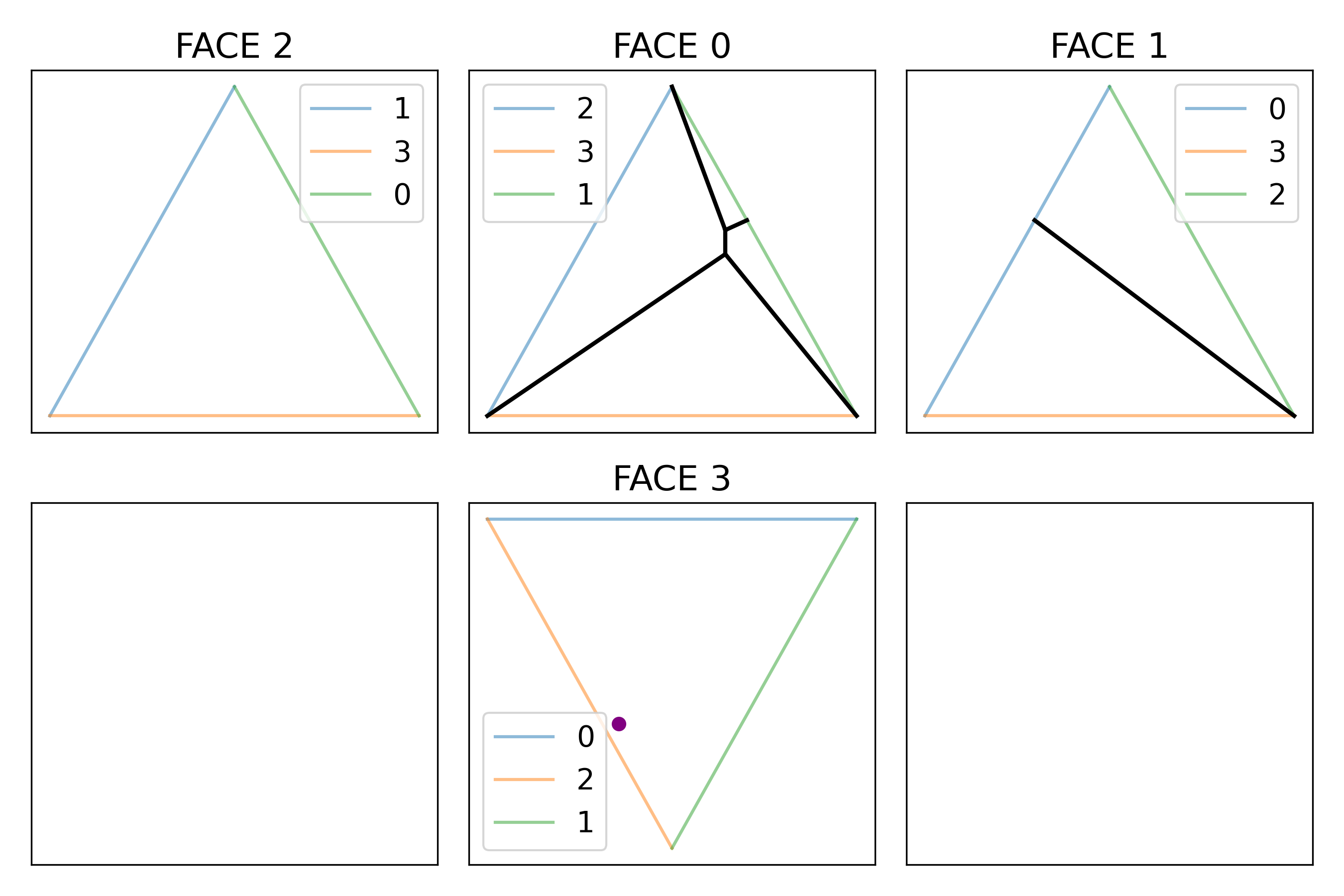}
    \includegraphics[height=\imgheight]{images/tetra/tetra_unfold_three.png}
    }
    \hfill
    \subfigure[$p$ chosen within line $aM$]{
    \includegraphics[height=\imgheight]{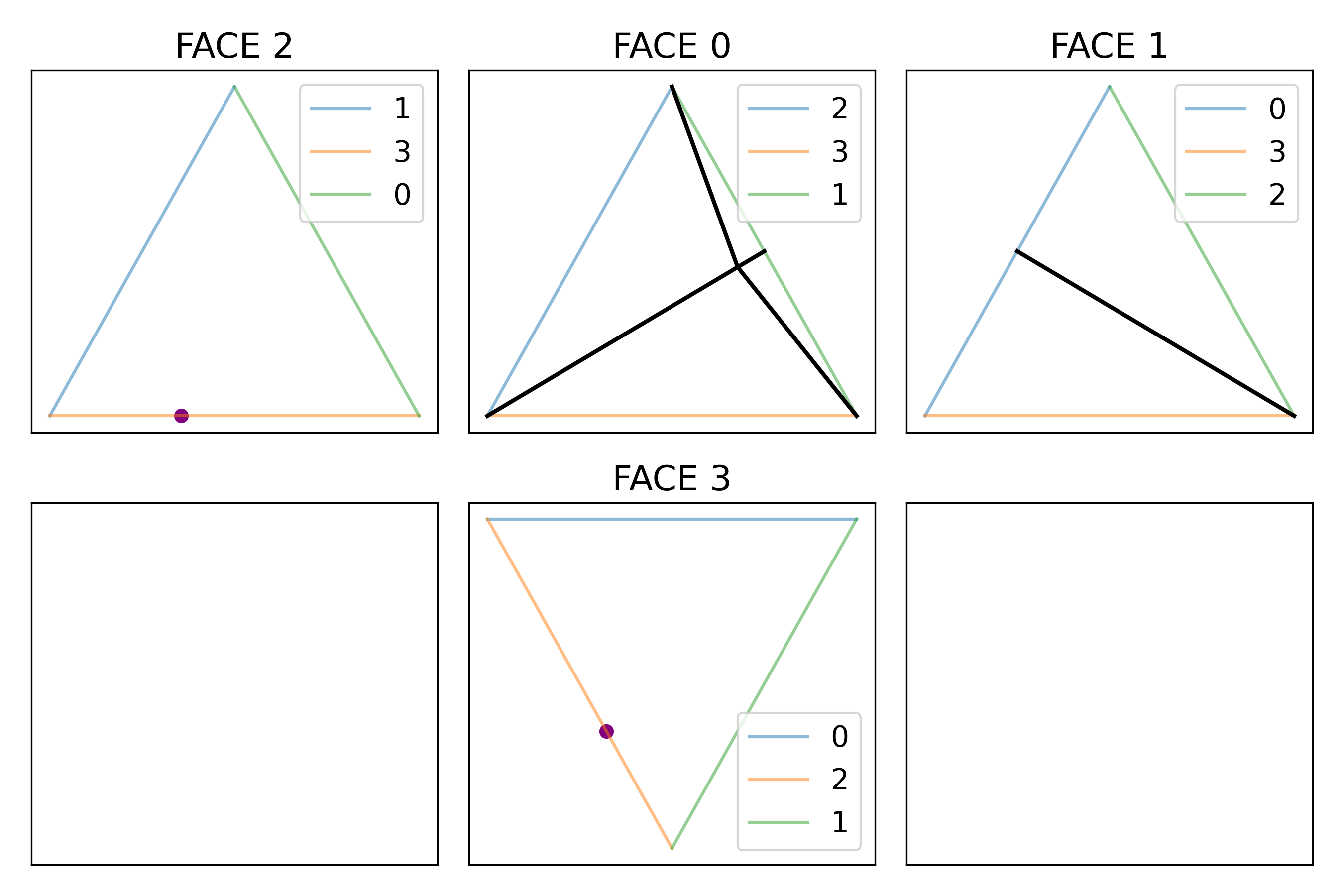}
    \includegraphics[height=\imgheight]{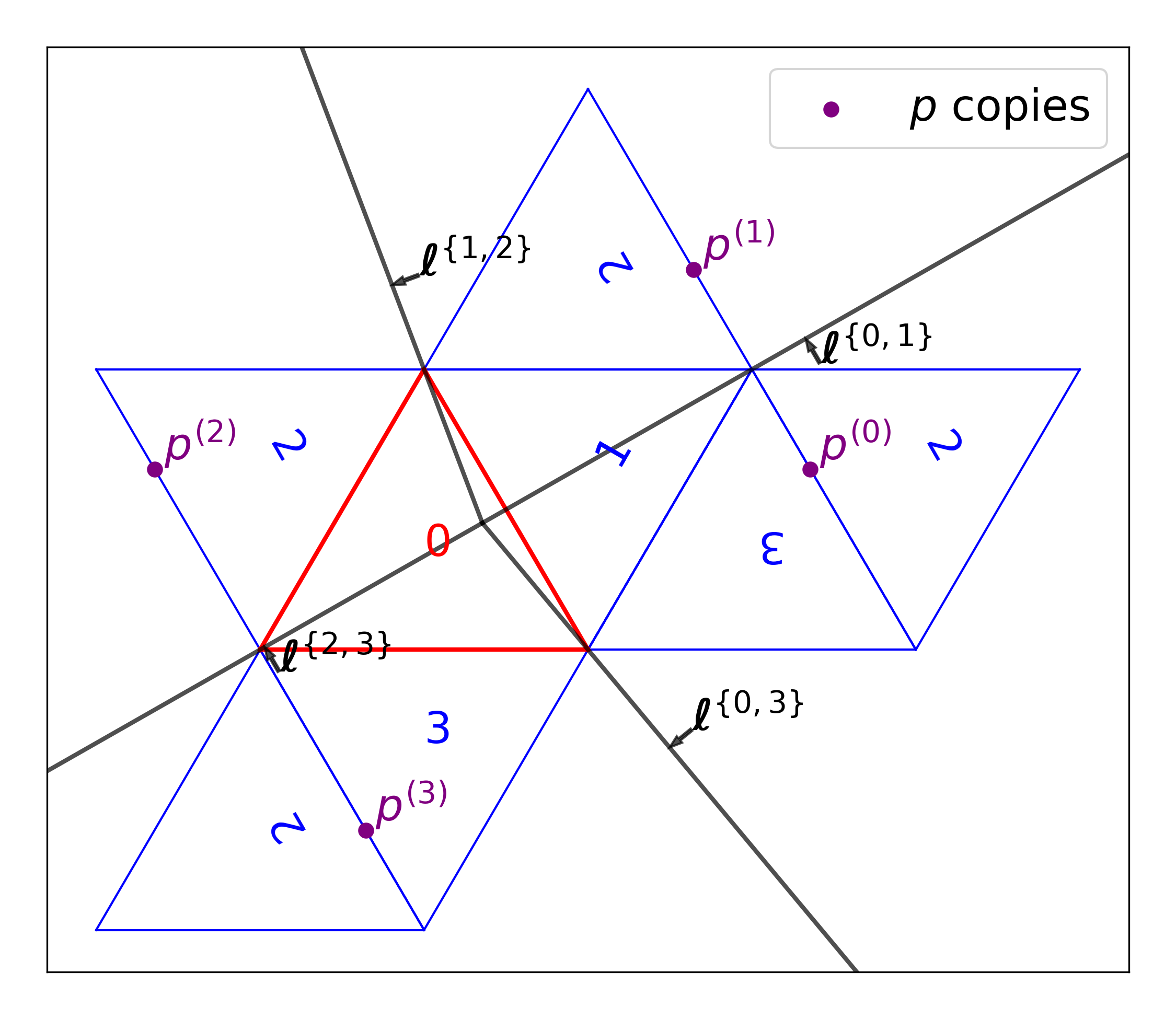}
    }
    \hfill
    \subfigure[$p$ chosen within $\tau'$ on Face 2]{
    \includegraphics[height=\imgheight]{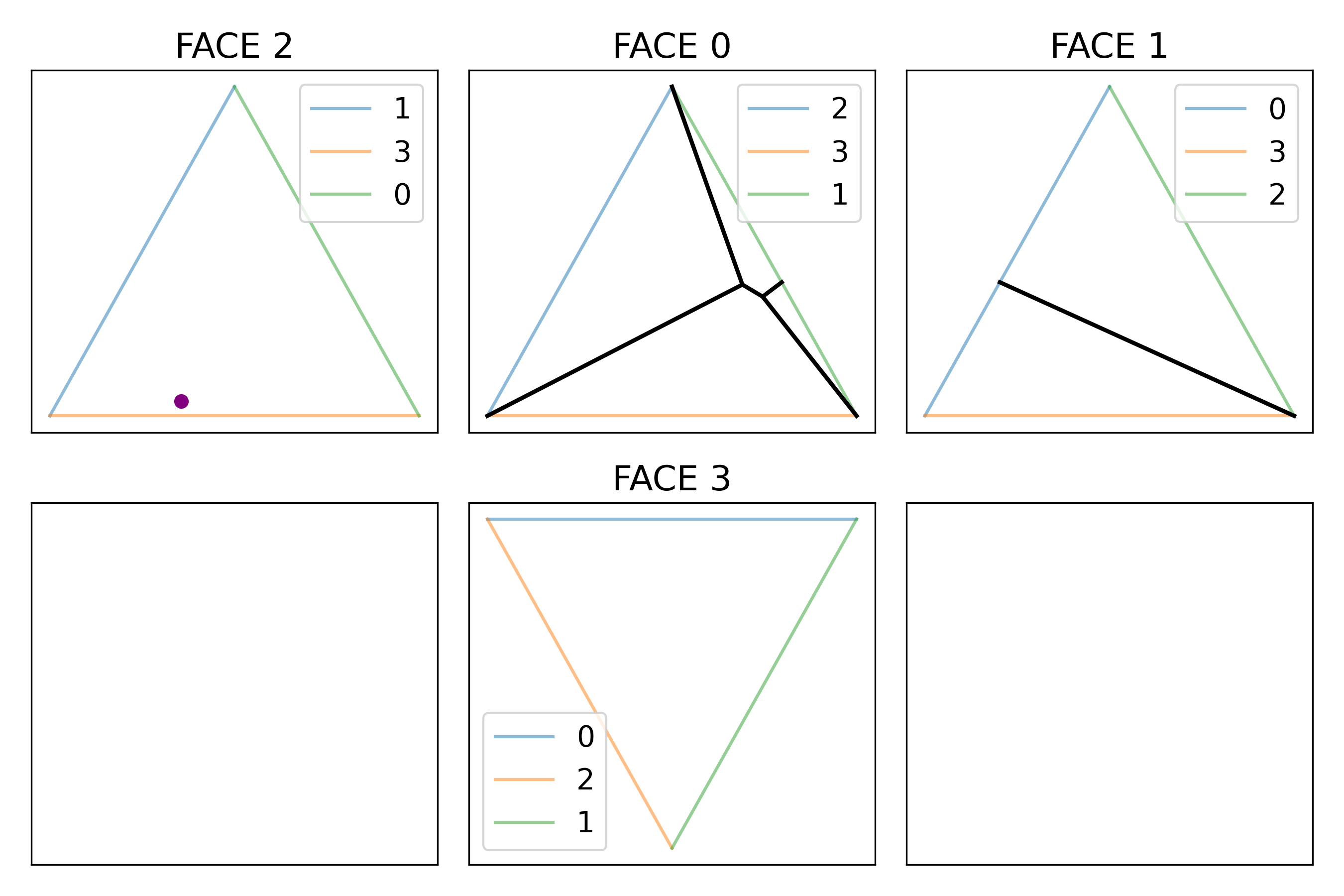}
    \includegraphics[height=\imgheight]{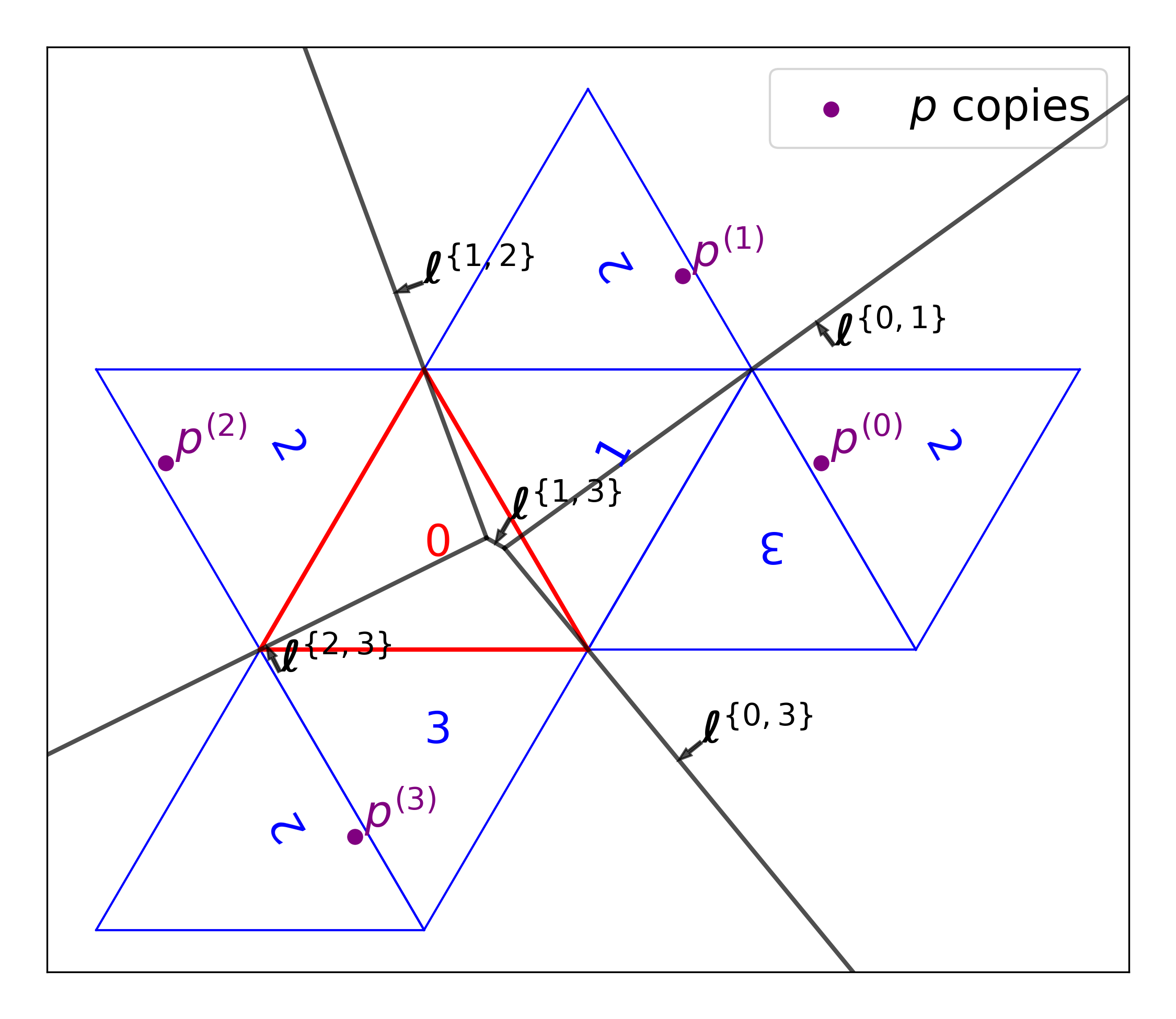}
    }
    \caption{Cut locus on Face 0 while varying $p$}
    \label{fig:tetra_first}
\end{figure}

We will define sets $\mathcal F_i$ to use \autoref{thm:quite_simplex}.
To do this, we will embed simplices into $X\times X$ using \autoref{lem:nice_embeddings}.

First, split Face 2 into six congruent triangles, as in \autoref{fig:dihedral_smmetry}.
Let $\tau'$ be the triangle that shares edge $aM$ with $\tau$.
Let $p^*\in \relint{aM}$, and let $q^*$ be the vertex of the corresponding cut locus (see \autoref{fig:tetra_locus_boundaries}(d)).
Let $p_{3}\colon[0,1]\hookrightarrow X$ be an embedding so that $p_{3}(1)=p^*$ and $p_{3}([0,1))\subseteq \text{int}(\tau)$.
Similarly, let $p_{2}\colon[0,1]\hookrightarrow X$ so that $p_{2}(1)=p^*$ and $p_{2}([0,1))\subseteq \text{int}(\tau')$.
This is possible, a simple example would be linear paths between the points in \autoref{fig:tetra_first}(a) and (b) for $p_3$, and between the points in \autoref{fig:tetra_first}(c) and (b) for $p_2$.

Let $\mathcal J:=\{J\subseteq \{0,1,2,3\}:|J|=3\}$, and $\mathcal I:=\{\{i,j\}\subseteq \{0,1,2,3\}:j\equiv i+1\mod 4\}$ be collections of subsets of $\{0,1,2,3\}$.
For $F\in \{2,3\}$, $I\in \mathcal I$, and $J\in\mathcal J$, with $I\subseteq J$, we will construct map $q^{F,I,J}\colon [0,1]\to X$ to vary linearly between $x^{I}\circ p_F$ and $x^{J}\circ p_F$.
Explicitly, $q^{F,I,J}(t)=(1-t)\cdot x^{I}(p_F(t))+t\cdot x^{J}(p_F(t))$.
Since $x^{I}$ and $x^{J}$ are both on line $\ell^I$, $q^{F,I,J}$ will always be on line $\ell^I$.
Additionally, $q^{F,I,J}(1)=x^{J}(p_F(1))=q^*$.

We will describe one choice of $F$, $I$, and $J$.
We assume $p_3$ and $p_2$ are linear paths suggested by \autoref{fig:tetra_first}.
Let $F=3$, $I=\{0,1\}$, and $J=\{0,1,2\}$.
The line $p_{3}$ varies from the location of $p$ in \autoref{fig:tetra_first}(a) to the location in \autoref{fig:tetra_first}(b).
The point $x^{\{0,1\}}$ is the intersection of $\ell^{\{0,1\}}$ with the edge of Face 0, and $x^{\{0,1,2\}}$ is the top vertex of the cut locus in \autoref{fig:tetra_first}(a).
Thus, $q^{3,\{0,1\},\{0,1,2\}}$ varies between these two points, staying on line $\ell^{\{0,1\}}$ the whole time.
We display a few locations of $q^{3,\{0,1\},\{0,1,2\}}$ in \autoref{fig:tetra_proof}.

We will now describe sets $\mathcal F_i$.
Recall that \autoref{lem:nice_embeddings} takes an embedding $T_k\hookrightarrow Z$ and a map $T_k\hookrightarrow Y^{k'}$ to define a simplex embedding $\Delta_{k+k'-1}\hookrightarrow Z\times Y$.
In this case, we will set $k=1$ and $Z=X$ so that $T_k=[0,1]$, and $p_3$ or $p_2$ may be used as the first embedding.
Since each face of $X$ is a convex subset of $\mathbb R^2$, we let $Y$ be Face 0, resulting in our desired embedding $\Delta_{k+k'-1}\hookrightarrow X\times X$.

\begin{compactitem}
    \item
    $\mathcal F_0=\{\Delta_0\mapsto (p^*,q^*)\}$.
    \item
    $\mathcal F_1$ will be defined by four lines.
    
    Consider $p_{3}$, and the map $x^{J}\circ p_{3}$ for $J\in \{\{0,1,2\},\{0,2,3\}\}$.
    By \autoref{lem:nice_embeddings}, this describes an embedding $f\colon \Delta_{1}\hookrightarrow X\times X$ where $\pi_1(f(t,1-t))=p_{3}(t)$, and $\pi_2(f(t,1-t))= x^{J}(p_{3}(t))$.
    Each map follows $(p,x^J(p))$ as $p$ approaches its endpoint on $aM$.

    We define the other two embeddings similarly using $p_{2}$ and $x^{J}\circ p_{2}$ for $J\in \{\{0,1,3\},\{1,2,3\}\}$.
    
    \item
    $\mathcal F_2$ will be 10 simplex embeddings, corresponding to the 5 lines around the cut locus vertices $x^{J}$ considered when constructing $\mathcal F_1$ (each vertex is adjacent to three lines, with one line counted twice).

    Consider $p_{3}$ and the map $q^{3,I,\{0,1,2\}}\funcprodsymb (x^{\{0,1,2\}}\circ p_{3})$ for $I\in \{\{0,1\},\{1,2\}\}$.
    Since $q^{3,I,\{0,1,2\}}(1)=x^{\{0,1,2\}}(p_{3}(1))=q^*$, and these points are otherwise distinct, we have by \autoref{lem:nice_embeddings} that this defines an embedding $f\colon \Delta_2\hookrightarrow X\times X$ whose image looks like $(p,q)$ for $p\in p_{3}([0,1])$ and $q\in \ell^{I}(p)$.
    
    We can define similar maps using $p_{3}$ and the maps $q^{3,I,\{0,2,3\}}\funcprodsymb (x^{\{0,2,3\}}\circ p_{3})$ for $I\in \{\{2,3\},\{0,3\}\}$.
    This creates four embeddings corresponding to $p_3$.
    
    We get the remaining map by considering just the vertex embeddings $x^{J}\circ p_{3}$.
    Explicitly, consider $p_{3}$ and the map $(x^{\{0,1,2\}}\circ p_{3})\funcprodsymb (x^{\{0,2,3\}}\circ p_{3})$.
    Since $x^{\{0,1,2\}}(p_{3}(1))= x^{\{0,2,3\}}(p_{3}(1))=q^*$, and $x^{\{0,1,2\}}\circ p_{3}\neq x^{\{0,2,3\}}\circ p_{3}$ otherwise, this defines an embedding $\Delta_2\hookrightarrow X\times X$ whose image looks like $(p,q)$ for $p\in p_{3}([0,1])$ and $q\in \ell^{(0,2)}(p)$.

    We create 5 embeddings corresponding to $p_2$ in a similar way to get a total of 10.

    \item
    $\mathcal F_3$ will be a total of 12 simplex embeddings, corresponding to the three Voronoi cells around each of the four cut locus vertices considered when constructing $\mathcal F_1$.

    Consider $p_{3}$ and the map $(x^{\{0,1,2\}}\circ p_{3})\funcprodsymb (x^{\{0,2,3\}}\circ p_{3}) \funcprodsymb q^{3,\{0,1\},\{0,1,2\}}$.
    Note that $x^{\{0,1,2\}}(p_{3}(1)) = x^{\{0,2,3\}}(p_{3}(1)) = q^{3,\{0,1\},\{0,1,2\}}(1)=q^*$, and these maps are otherwise not collinear.
    By \autoref{lem:nice_embeddings}, this defines an embedding $f\colon \Delta_3\hookrightarrow X\times X$ where for each $\textbf{d}\in \Delta_3$, $f(\textbf{d})=(p,q)$ for $p\in p_{3}([0,1])$ and $q$ in the convex hull of $x^{\{0,1,2\}}(p)$, $x^{\{0,2,3\}}(p)$, and $q^{3,\{0,1\},\{0,1,2\}}(\textbf{d}_1)$.

    We will create two similar embeddings by considering $p_{3}$ and the maps
    $(x^{\{0,1,2\}}\circ p_{3})\funcprodsymb q^{3,\{0,1\},\{0,1,2\}}\funcprodsymb q^{3,\{1,2\},\{0,1,2\}}$,
    and $(x^{\{0,1,2\}}\circ p_{3})\funcprodsymb q^{3,\{1,2\},\{0,1,2\}}\funcprodsymb (x^{\{0,2,3\}}\circ p_{3})$.

    These are the three embeddings corresponding to cut locus vertex $x^{\{0,1,2\}}\circ p_{3}$.
    We will obtain the embeddings corresponding to the vertices $(x^{\{0,2,3\}}\circ p_{3})$, $(x^{\{0,1,3\}}\circ p_{2})$, and $(x^{\{1,2,3\}}\circ p_{2})$ in a similar way.
    This accounts for the other nine embeddings.
    
\end{compactitem}

\renewcommand{\imgwidth}{.31\linewidth}
\begin{figure}[ht!]
    \centering
    \subfigure[$t=0$]{
    \includegraphics[width=\imgwidth]{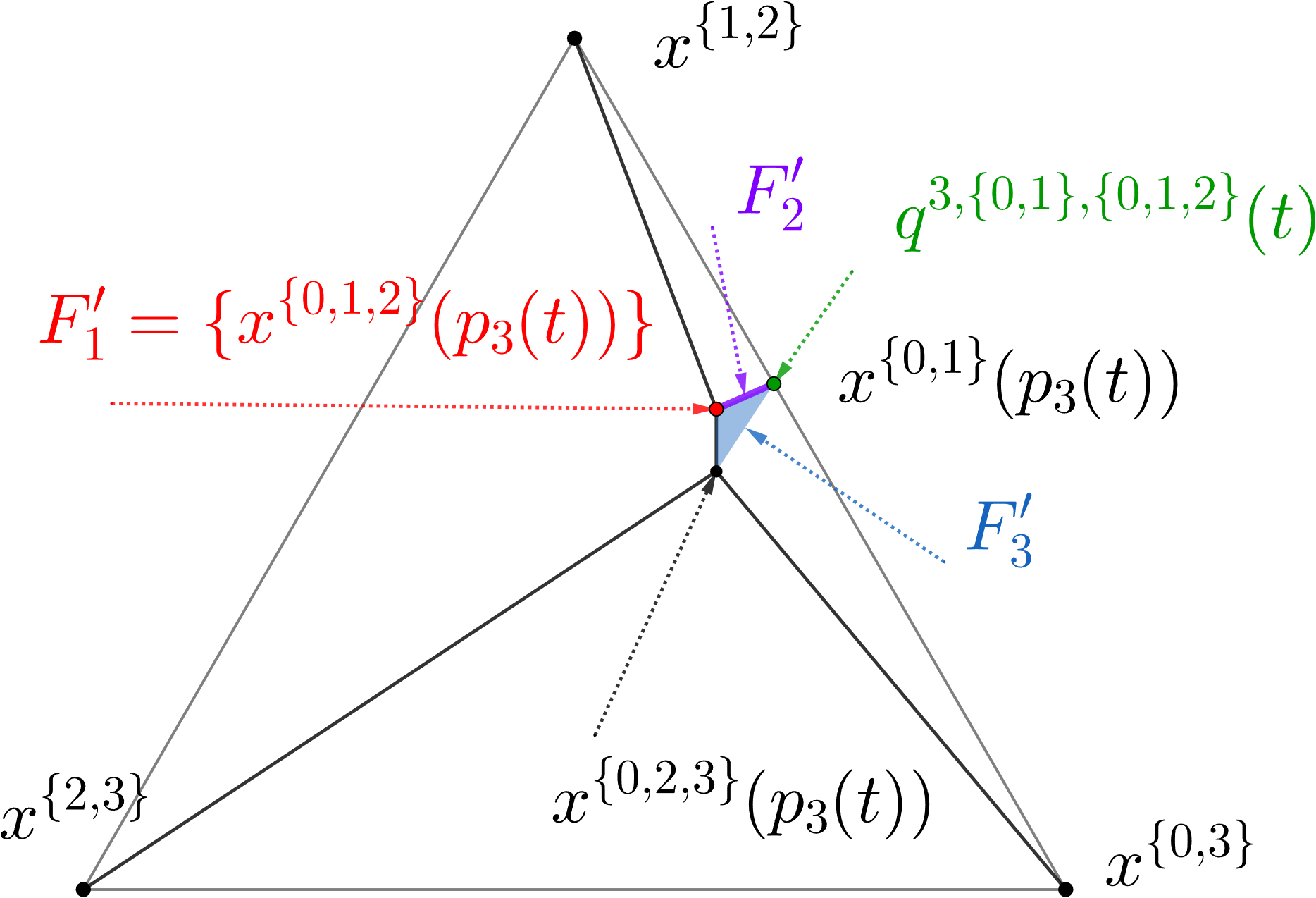}
    }
    \subfigure[$t=0.5$]{
    \includegraphics[width=\imgwidth]{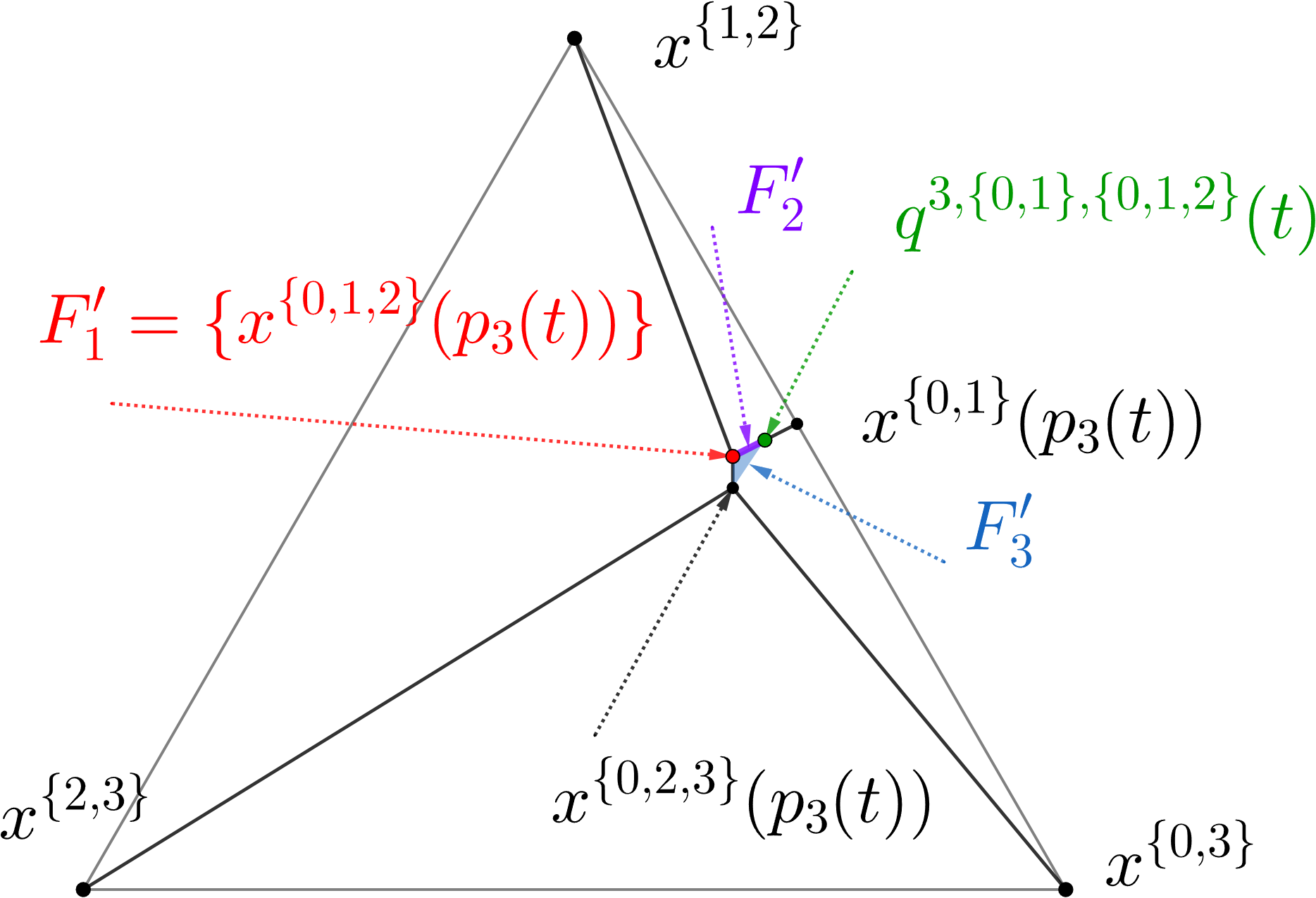}
    }
    \subfigure[$t=1$]{
    \includegraphics[width=\imgwidth]{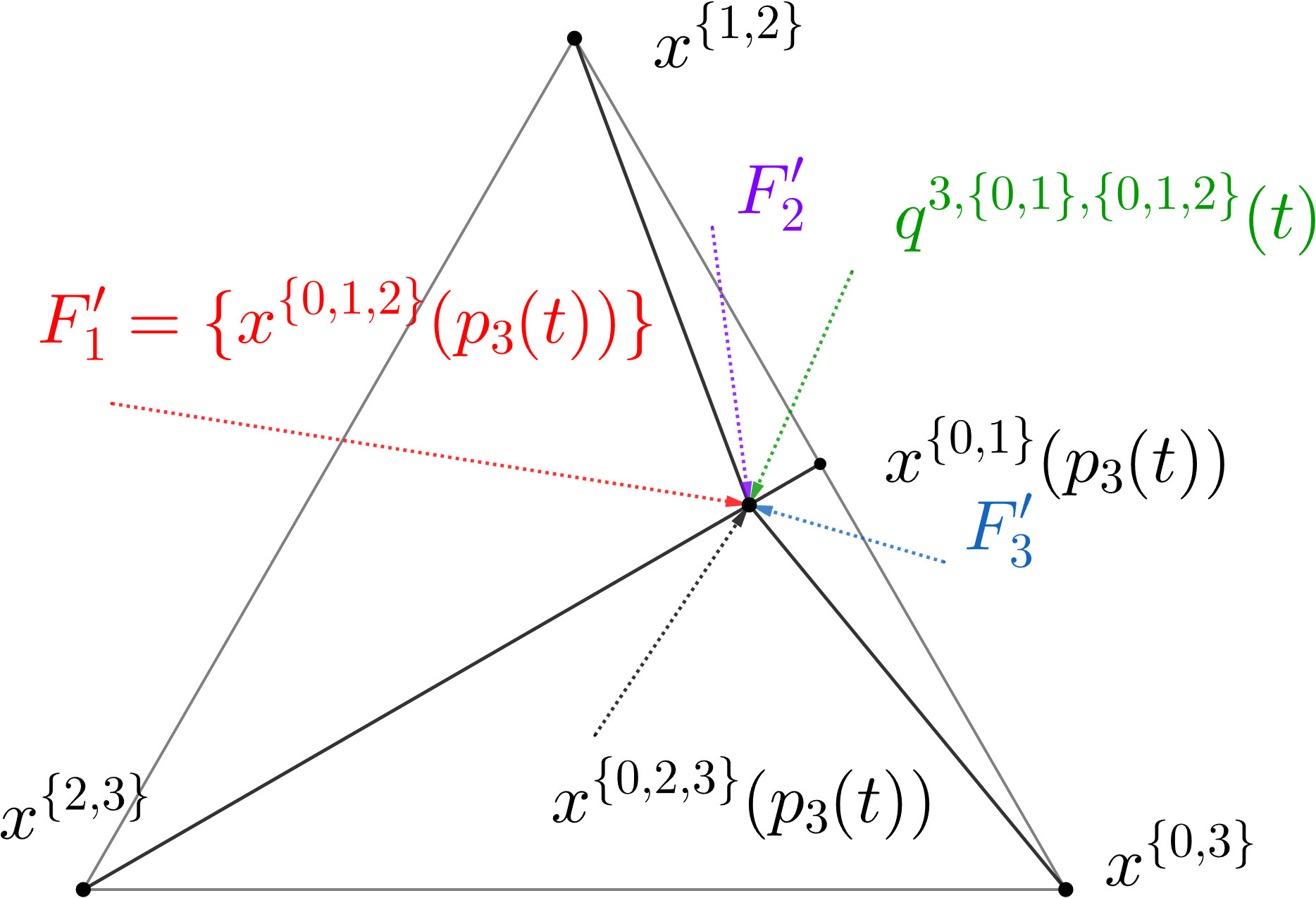}
    }
    \caption{Visualization of selected $F_i\in\mathcal F_i$ as $p$ varies towards $p^*$}
    \label{fig:tetra_proof}
\end{figure}

We will verify the properties of \autoref{thm:quite_simplex} for a few embeddings and geodesic choices, and assert that is is true for all of them by symmetry and by construction.

First, each $\mathcal F_i$ contains the embeddings of $\Delta_i\hookrightarrow X\times X$ by construction, so property (a) always holds.

\begin{compactitem}
    \item
    Consider $F_0\in\mathcal F_0$.
    Property (b) is trivial since $\img{F_0}=\{(p^*,q^*)\}$ is a single point.

    We will consider a \gls{GMPR} $\Gamma_0$ on $\img{F_0}$.
    We have four options, determined by which copy of $p^*$ ($p^{(0)}$, $p^{(1)}$, $p^{(2)}$, or $p^{(3)}$) we connect $q^*$ to.
    Assume our choice connects $q^*$ to $p^{(3)}$.

    Then our choice for $F_1\in \mathcal F_1$ will be the embedding constructed from maps $p_3$ and $x^{\{0,1,2\}}\circ p_{3}$ (e.g. $(p,q)$ where $p$ moves linearly from \autoref{fig:tetra_first}(a) to \autoref{fig:tetra_first}(b), and $q$ is the top vertex of the resulting cut locus).
    (i) is satisfied by construction.
    (ii) is satisfied since for any element  $(p,q)\in \relint{F_1}$, $q$ is a point on the cut locus of $p$ with geodesics only to $p^{(0)}$, $p^{(1)}$, and $p^{(2)}$.
    (iii) is satisfied since there are three \gls{GMPR}s on $\img{F_1}$, each distinguished by which of $p^{(0)}$, $p^{(1)}$, or  $p^{(2)}$ they choose to connect to.
    Thus, property (c) is satisfied.

    \item
    Consider $F_1\in \mathcal F_1$.
    Property (b) is satisfied since any geodesic for $(p,q)\in \relint{F_1}$, defined by connecting $q$ to $p^{(0)}$, $p^{(1)}$, or  $p^{(2)}$, extends to a \gls{GMPR} on $\img{F_1}$ distinguished by the same choice.

    We will consider one of the three \gls{GMPR}s $\Gamma_1$.
    We will choose $\Gamma_1$ to connect to $p^{(2)}$.

    Then our choice for $F_2\in \mathcal F_2$ will be the embedding constructed from maps $p_3$ and $q^{3,\{0,1\},\{0,1,2\}}\funcprodsymb (x^{\{0,1,2\}}\circ p_{3})$.
    This embedding follows a segment of $\ell^{\{0,1\}}(p)$ for each choice of $p$.
    (i) is satisfied by construction.
    (ii) is satisfied since for any element $(p,q)\in \relint{F_2}$, $q$ is a point on the cut locus of $p$ with geodesics only to $p^{(0)}$ and $p^{(1)}$.
    (iii) is satisfied since there are two \gls{GMPR}s on $\img{F_2}$, each distinguished by which of $p^{(0)}$ or $p^{(1)}$ they choose to connect to.
    Thus, property (c) is satisfied.
    
    \item
    Consider $F_2\in\mathcal F_2$.
    Property (b) is satisfied since any geodesic for $(p,q)\in \relint{F_2}$, defined by connecting $q$ to $p^{(0)}$ or $p^{(1)}$, extends to a \gls{GMPR} distinguished by the same choice.
    
    We will consider one of the two \gls{GMPR}s $\Gamma_2$.
    We will choose $\Gamma_2$ to connect to $p^{(1)}$.

    Then our choice for $F_3\in \mathcal F_3$ will be the embedding constructed from maps $p_3$ and $(x^{\{0,1,2\}}\circ p_{3})\funcprodsymb (x^{\{0,2,3\}}\circ p_{3}) \funcprodsymb q^{3,\{0,1\},\{0,1,2\}}$.
    This embedding follows a face bounded by $\ell^{\{0,1\}}(p)$ and $\ell^{\{0,2\}}(p)$ for each choice of $p$.
    (i) is satisfied by construction.
    (ii) is satisfied since for any element $(p,q)\in \relint{F_3}$, $q$ is a point on the cut locus of $p$ with a geodesic only to $p^{(0)}$.
    (iii) is satisfied since there is only one \gls{GMPR} on $\img{F_2}$, as the only choice that extends to the interior is to go to $p^{(0)}$.
    Thus, property (c) is satisfied.

    \item
    Consider $F_3\in \mathcal F_3$.
    Property (b) is satisfied since any geodesic for $(p,q)\in \relint{F_3}$ is defined by connecting $q$ to $p^{(0)}$, and extends to the only \gls{GMPR} over $\img{F_3}$.
    
    We do not need to check (c) as 3 is our maximal dimension.
\end{compactitem}

Then by \autoref{thm:quite_simplex}, $GC(X)\geq 3$ (needing at least four sets in any geodesic motion planner).

We visualize this proof as follows.
For $i>0$, each $F_i$ is a map $\Delta_i\to X\times X$, and is generated using \autoref{lem:nice_embeddings} with $k=1$.
Recall that in this case, $T_1=[0,1]$ and for $t\in T_1$, $D_{t,i}=\{x\in \Delta_{i}:x_1=t\}$.
From construction, we have that the projection on the second factor $\pi_2(F_i(D_{t,i}))$ is the convex hull of $i$ points.
Let this hull be $F_i'$.
In \autoref{fig:tetra_proof}, we plot the $F_i'$ corresponding to the $F_i$ we inspected for a few values of $t$.
We also visualize all of $\mathcal F_i$ in a similar way in \autoref{apx:tetrapendix}.

\subsection{Octahedron}
We will prove \autoref{thm:octahedron_GC}, that the geodesic complexity of an octahedron is exactly four (needing five sets in an efficient geodesic motion planner).

\subsubsection{Cut Locus}
\label{subsubsec:octa_cut_locus_proof}

We will calculate the cut locus of a point $p$ on the octahedron $X$.
Similar to the tetrahedron, we represent all faces as equilateral triangles centered at $(0,0)$ with side length $2\sqrt3$, oriented and labeled as shown in \autoref{fig:octa_zero}(a).
We choose $p$ to be on Face 6, and consider the corresponding cut locus.

We claim that the nineteen regions of Face 6 displayed in \autoref{fig:octa_isomorphic_regions}(a) each have isomorphic cut loci, and the lines of these cut loci vary continuously with the choice of $p$.
These regions consist of four points, nine line segments, and six triangles.
The fact that the cut loci vary continuously with the choice of $p$ follows from the fact that the cut loci are boundaries from a Voronoi diagram, and the points that create this diagram vary continuously.
Thus, we will proceed to find the structure of the cut locus on these regions.

\renewcommand{\imgheight}{110 pt}
\begin{figure}[ht!]
    \centering
\subfigure[All nineteen regions on Face~6]{
\includegraphics[height=\imgheight]{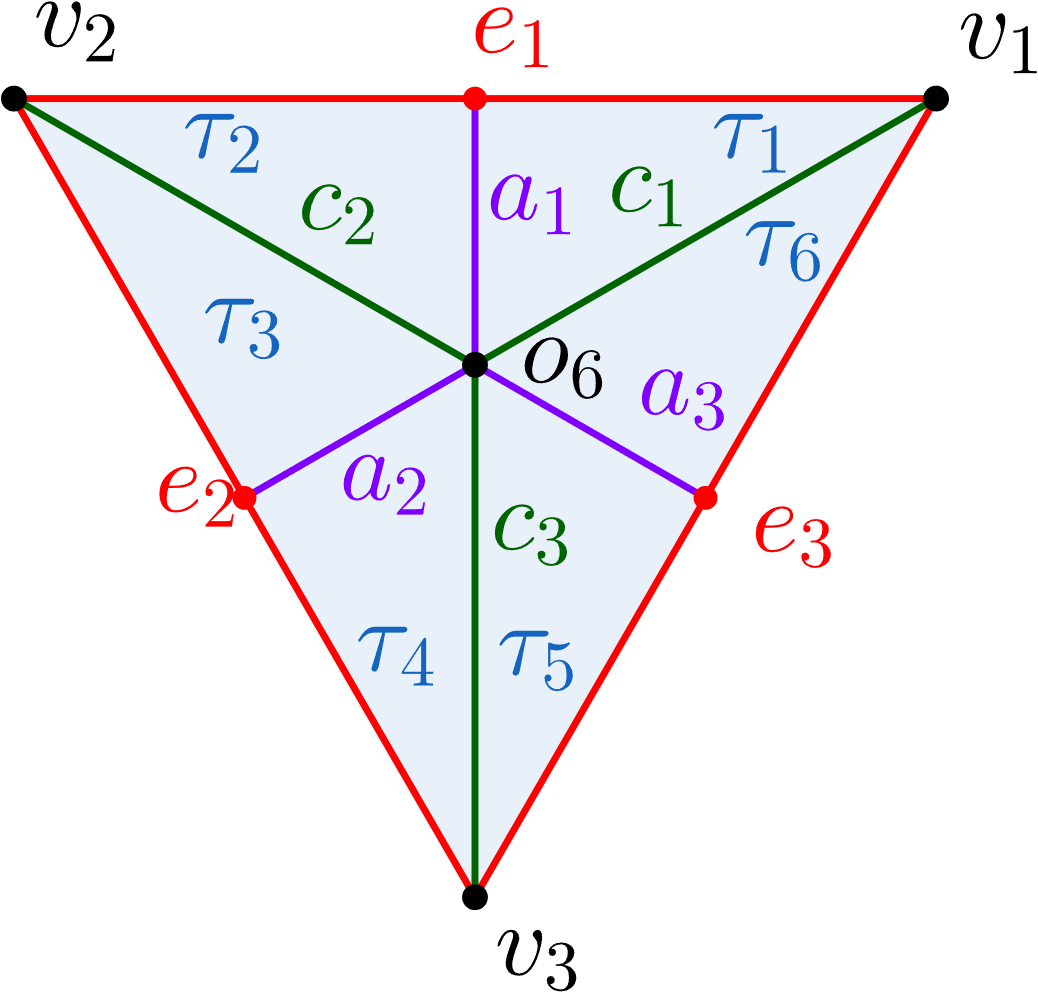}
}
\quad
\subfigure[We focus on $\tau_1$ and its boundary]{
\includegraphics[height=\imgheight]{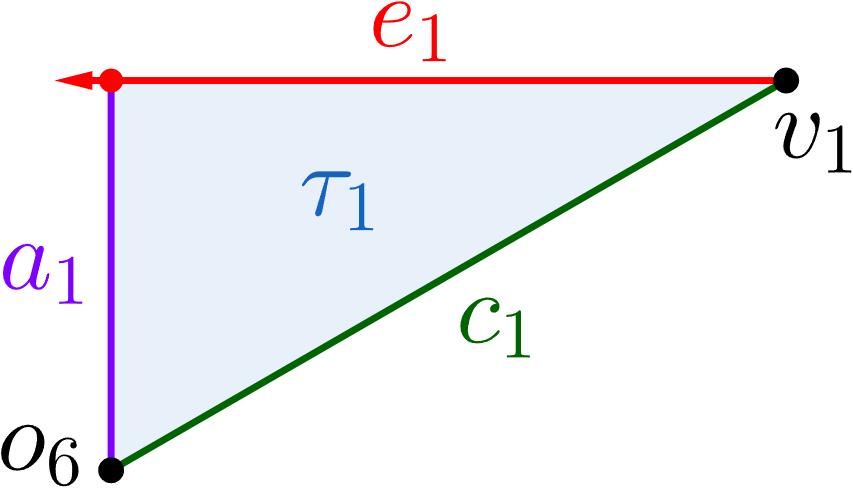}
}
\caption{Partition of Face 6 into regions with isomorphic cut loci}
\label{fig:octa_isomorphic_regions}
\end{figure}

Up to symmetry, a point on one of six congruent triangles on a face (displayed in \autoref{fig:octa_isomorphic_regions}(a)) represents any point on the octahedron surface. 
Thus, we need only consider one of the $\tau_i$, along with its boundary.
We choose $\tau_1$, whose boundary contains $a_1$, $c_1$, $v_1$, $o_6$, and part of edge $e_1$.
Note that \autoref{fig:octa_isomorphic_regions}(a) is a partition, so each $a_i,c_i,e_i$ do not contain their endpoints, and each $\tau_i$ does not contain its boundary.
We will first inspect the sets on the interior of Face 6, $\tau_1$, $c_1$, $a_1$, and $o_6$.
From \autoref{alg:cutlocus}, we find that for $p$ chosen in these regions, the cut locus appears on Faces 0, 1, 3, and 4.

\renewcommand{\imgheight}{135 pt}
\renewcommand{\imgwidth}{.5\linewidth}
\begin{figure}[ht!]
\centering
\subfigure[Cut Locus]{
\includegraphics[width=\imgwidth]{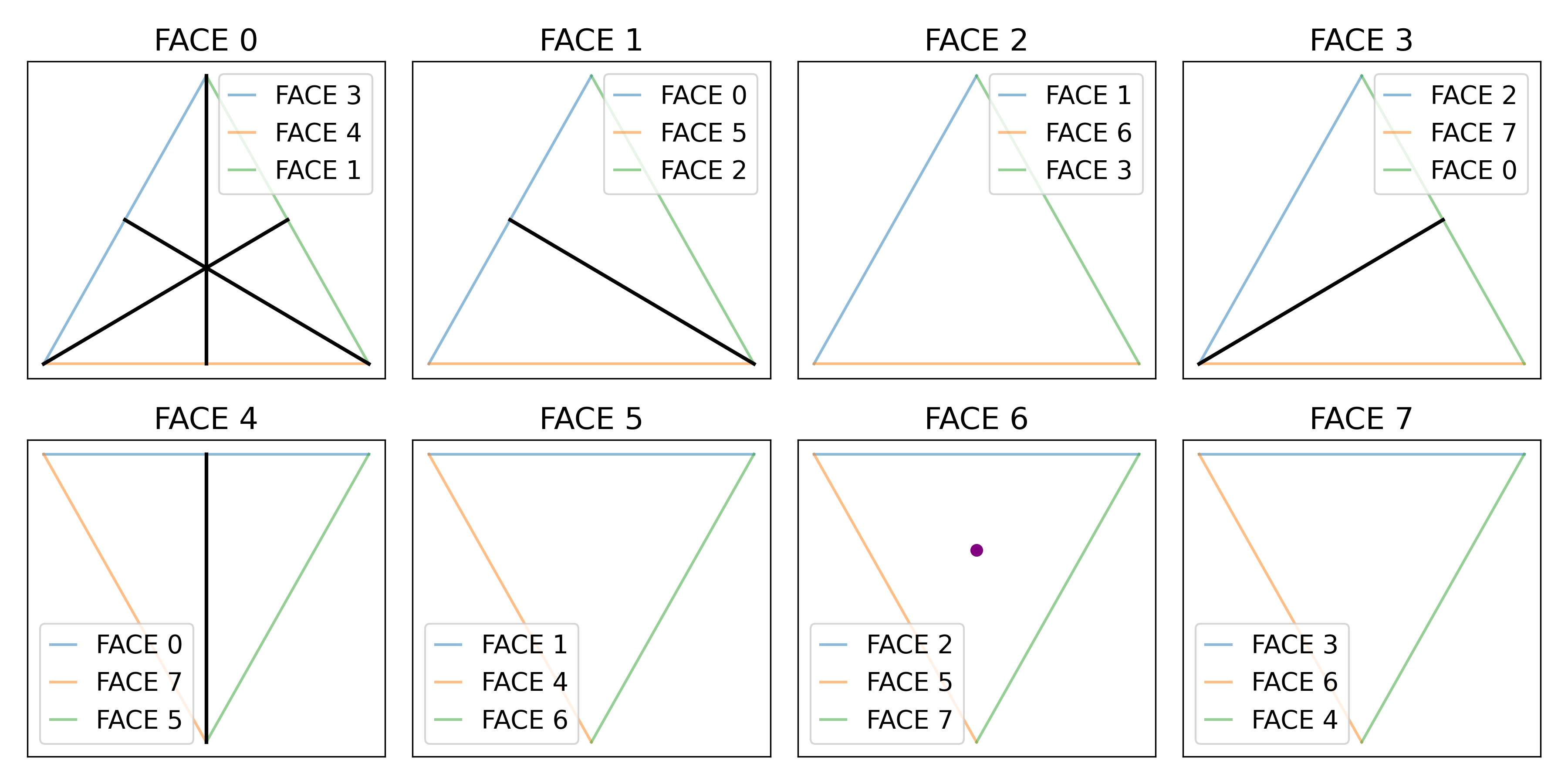}
}
\subfigure[Path unfoldings from Face 0]{
\includegraphics[height=\imgheight]{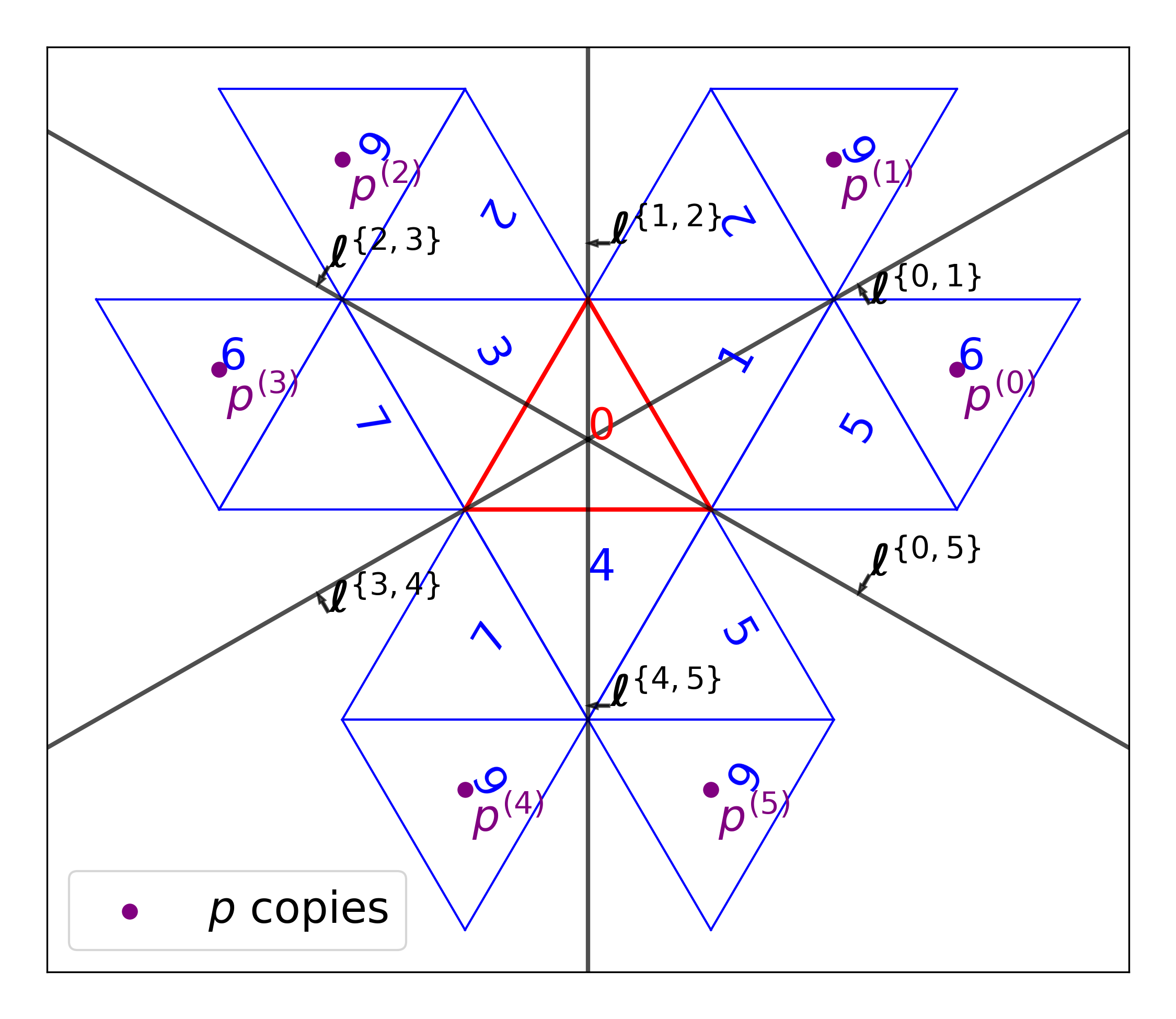}
}
\caption{Cut locus of an octahedron with respect to $p=o_6$}
\label{fig:octa_zero}
\end{figure}

On this region, the cut locus on Face 0 of $p$ arises from the Voronoi diagram of six points, as displayed in \autoref{fig:octa_zero}(b).
We enumerate these $(p^{(0)},\dots,p^{(5)})=\\
\left(
\begin{bmatrix}
    3\sqrt3\\1
\end{bmatrix}+
 p,
\begin{bmatrix}
    2\sqrt3 \\4
\end{bmatrix}+
R_{2\pi/3} p,
\begin{bmatrix}
    -2\sqrt3 \\4
\end{bmatrix}+
R_{-2\pi/3} p,
\begin{bmatrix}
    -3\sqrt3\\1
\end{bmatrix}+
 p,
\begin{bmatrix}
    -\sqrt3\\-5
\end{bmatrix}+
R_{2\pi/3} p,
\begin{bmatrix}
    \sqrt3\\-5
\end{bmatrix}+
R_{-2\pi/3} p,
\right)$.
Analogous to the tetrahedron, we define $\ell^{\{i,j\}}(p)$ as the line of the cut locus that bisects $p^{(i)}$ and $p^{(j)}$, $x^{\{i,j\}}(p)$ as the intersection of selected $\ell^{\{i,j\}}(p)$ with an edge of Face 0, and $x^{\{i,j,k\}}(p)$ as the intersection of $\ell^{\{i,j\}}(p)$ and $\ell^{\{j,k\}}(p)$.
Using a symbolic calculator, we explicitly find their equations in \autoref{apx:octapendix}.

To find the structure of the cut loci, and that they are isomorphic within their respective regions, we use the Wolfram Mathematica symbolic calculator as we did in \autoref{subsubsec:tetra_cut_locus_proof}.
For $p$ chosen within each region, this allows us to verify which $\ell^{\{i,j\}}$ occur on the cut locus as well as the identity and location of each line's endpoints.
To find the cut locus structure on Face 0, it is enough to prove the boundary of each of the six Voronoi cells in this way.
We first inspect $p$ chosen on the interior of Face 6, since we find that for these regions, the cut loci come from copies of $p$ that arise from the same path unfoldings. 

Using \autoref{alg:cutlocus}, we find that when $p$ is on the interior of Face 6, the cut loci on Faces 1, 3, and 4 are as suggested in \autoref{fig:octa_zero}(b).
They each consist of a line arising from two copies of $p$, and are equivalent to extensions of the lines $\ell^{\{0,1\}}$, $\ell^{\{2,3\}}$, and $\ell^{\{4,5\}}$ from the cut locus on Face 0.
We also observe from symmetries of the copies of $p$ that lines $\ell^{\{0,1\}}$, $\ell^{\{2,3\}}$, and $\ell^{\{4,5\}}$ always begin at the indicated vertices of Faces 1, 3, and 4 respectively.
Similarly, $\ell^{\{1,2\}}$, $\ell^{\{3,4\}}$, and $\ell^{\{0,5\}}$ always begin at the vertices of Face 0.

\meminisection{Point $o_6$:} 
First, we consider point $o_6$, the center of Face 6.
We observe that the structure of the cut locus is a star, with 6 lines joined at a point. 
We display this in \autoref{fig:octa_zero}.

\begin{figure}[ht!]
\centering
\subfigure[Cut Locus]{
\includegraphics[width=\imgwidth]{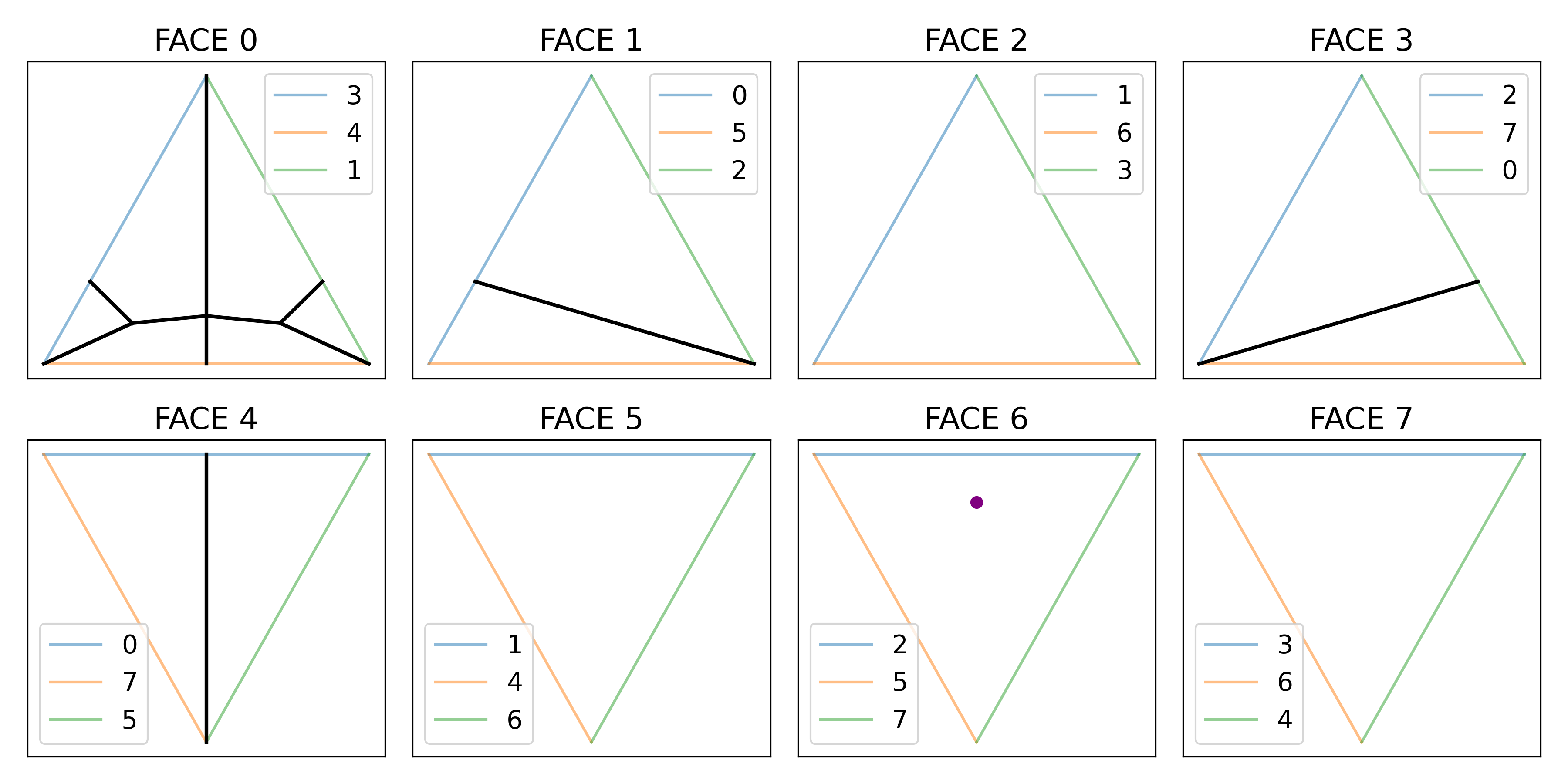}
}
\subfigure[Path unfoldings from Face 0]{
\includegraphics[height=\imgheight]{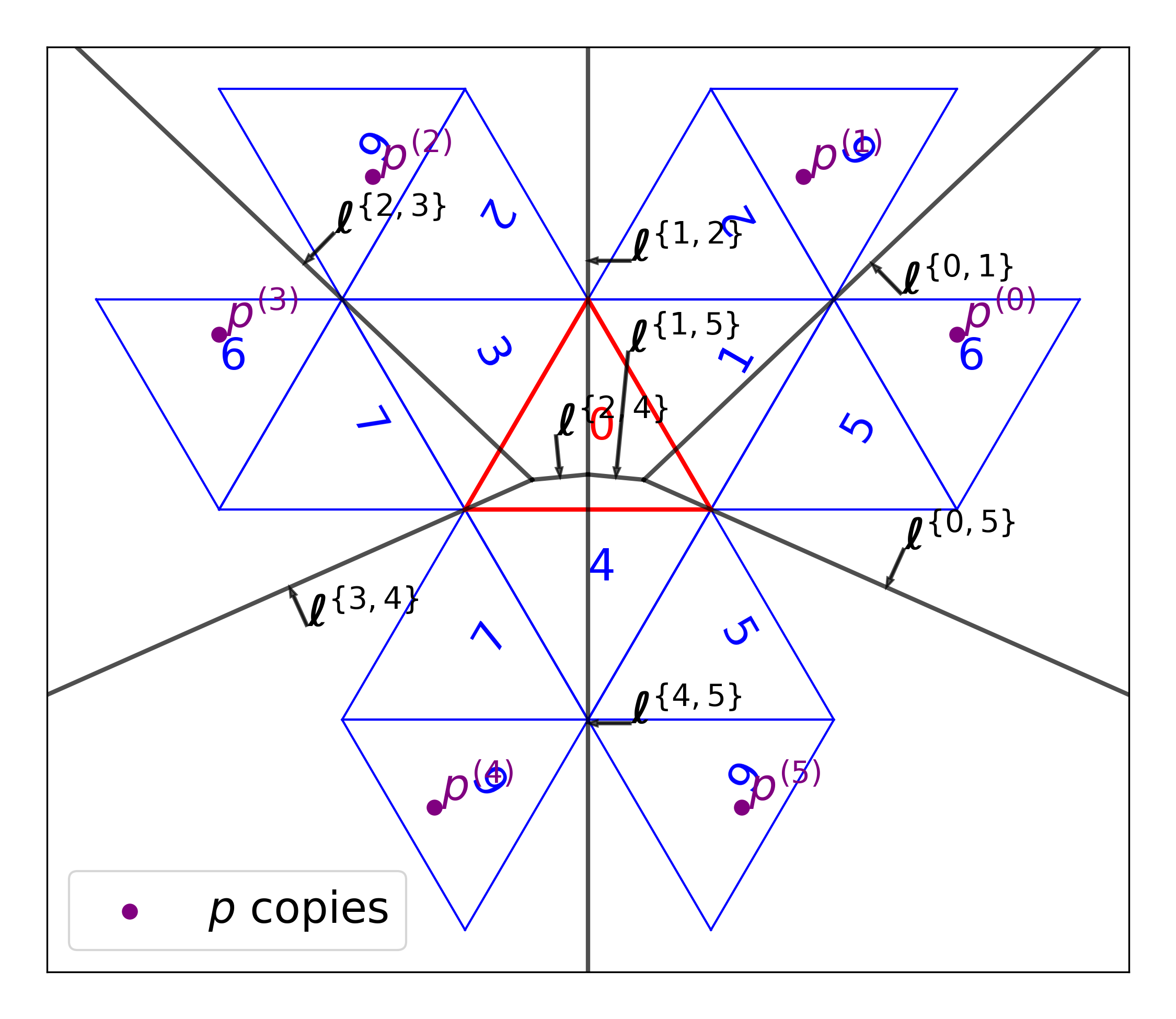}
}
\caption{Cut locus of $p=(0,\frac12)$ chosen on $a_1$}
\label{fig:octaa1locus}
\end{figure}

\meminisection{Line $a_1$:} 
We find that the structure of the cut locus of $p\in a_1$ is isomorphic to \autoref{fig:octaa1locus}, with one vertex $x^{\{1,2,4,5\}}$ incident to four line segments and two vertices $x^{\{0,1,5\}}$ and $x^{\{2,3,4\}}$ incident to three line segments.
All the cut locus vertices are within Face 0.
The cut locus has symmetry with respect to reflection about $x=0$, a result of the path unfoldings that create the six copies of $p$ having the same symmetry.

\begin{figure}[ht!]
\centering
\subfigure[Cut Locus]{
\includegraphics[width=\imgwidth]{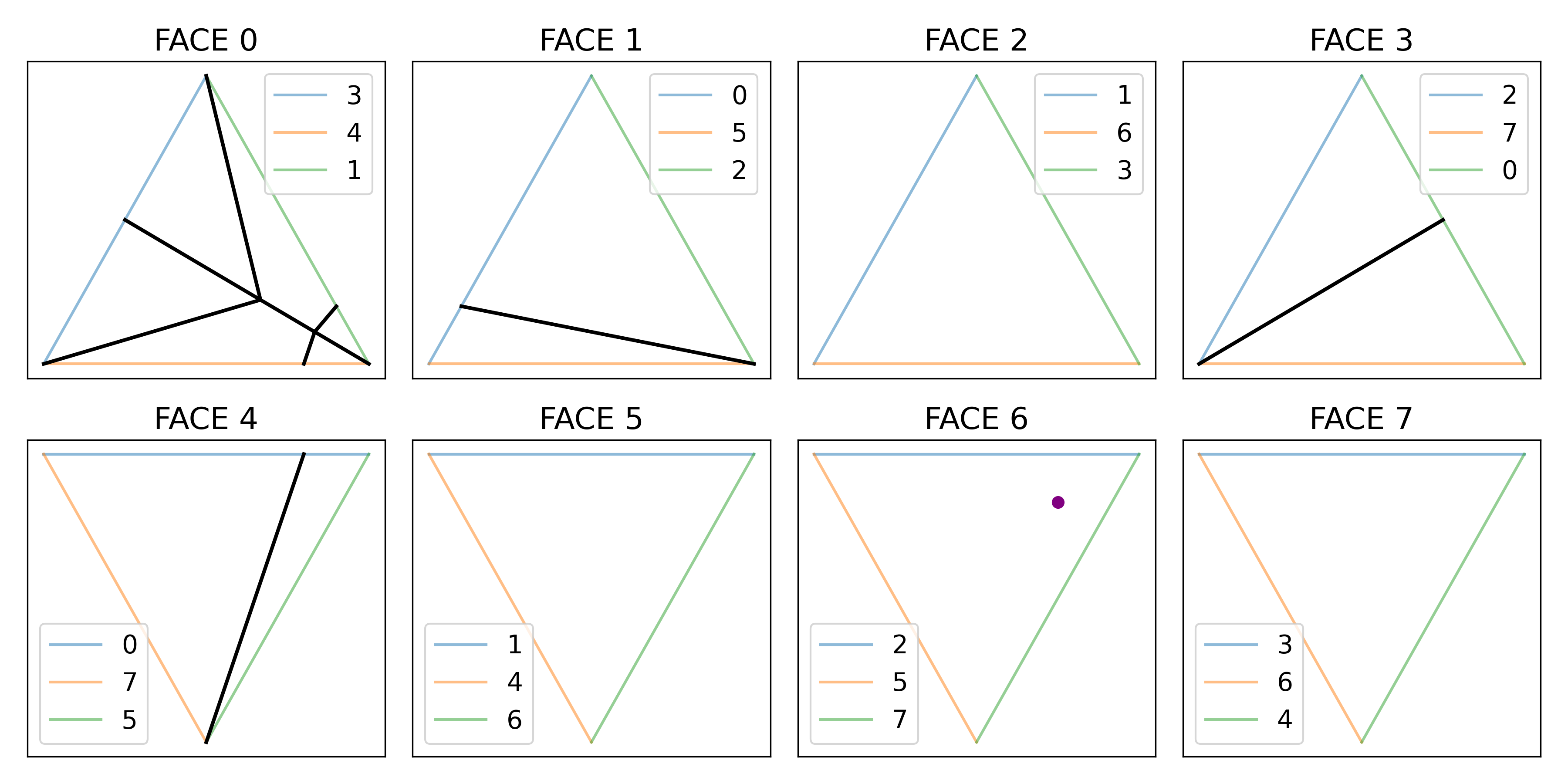}
}
\subfigure[Path unfoldings from Face 0]{
\includegraphics[height=\imgheight]{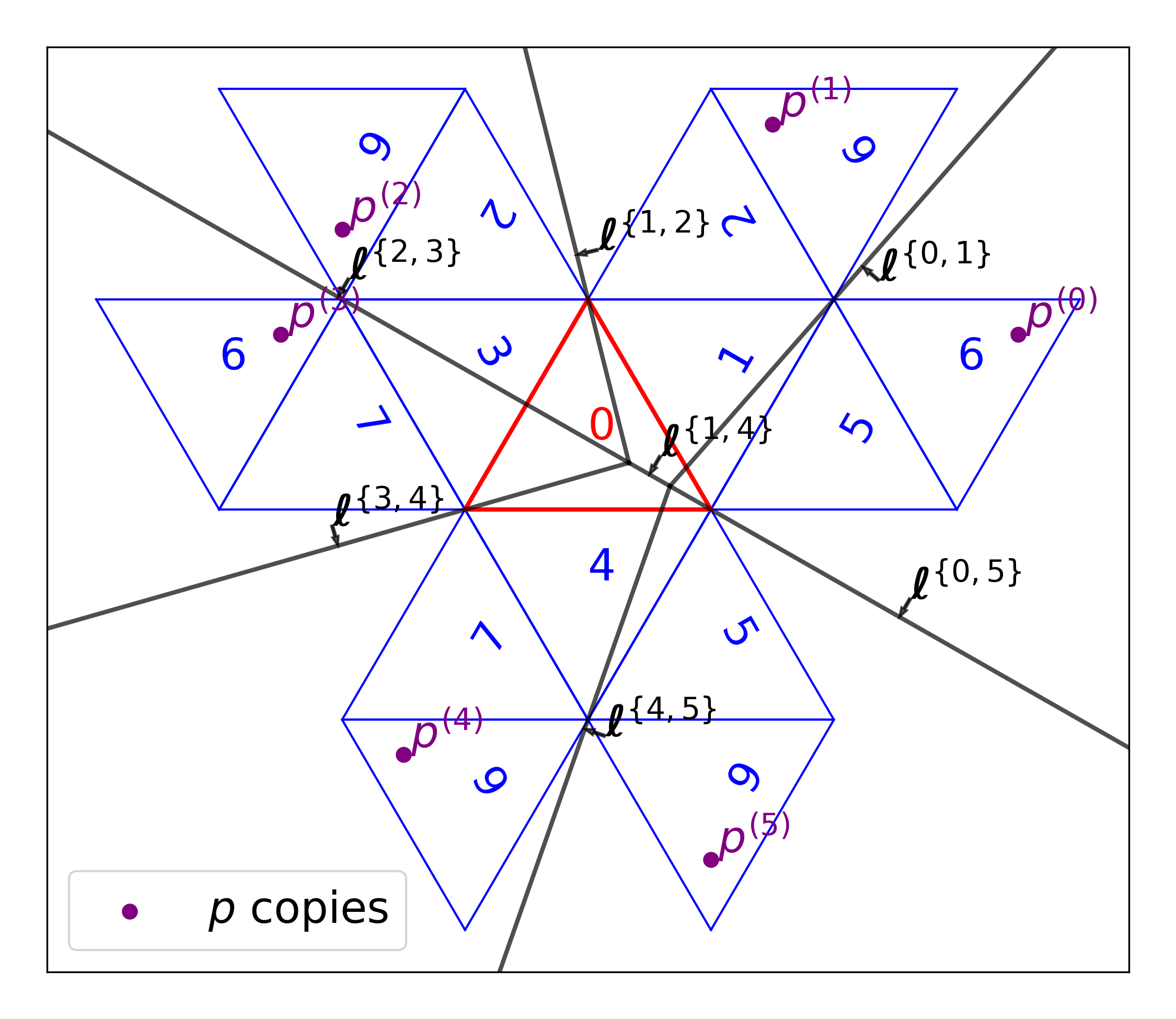}
}
\caption{Cut locus of $p=(\frac{\sqrt3}2,\frac12)$ chosen on $c_1$}
\label{fig:octac1locus}
\end{figure}

\meminisection{Line $c_1$:}
We find that the structure of the cut locus of $p\in c_1$ is isomorphic to \autoref{fig:octac1locus}, with two vertices $x^{\{0,1,4,5\}}$ and $x^{\{1,2,3,4\}}$ both incident to four line segments.
All the cut locus vertices are within Face 0.
The cut locus has symmetry with respect to reflection about the line $y=-\frac{x}{\sqrt3}$, a result of the path unfoldings that create the six copies of $p$ having the same symmetry.

\begin{figure}[ht!]
\centering
\subfigure[Cut Locus]{
\includegraphics[width=\imgwidth]{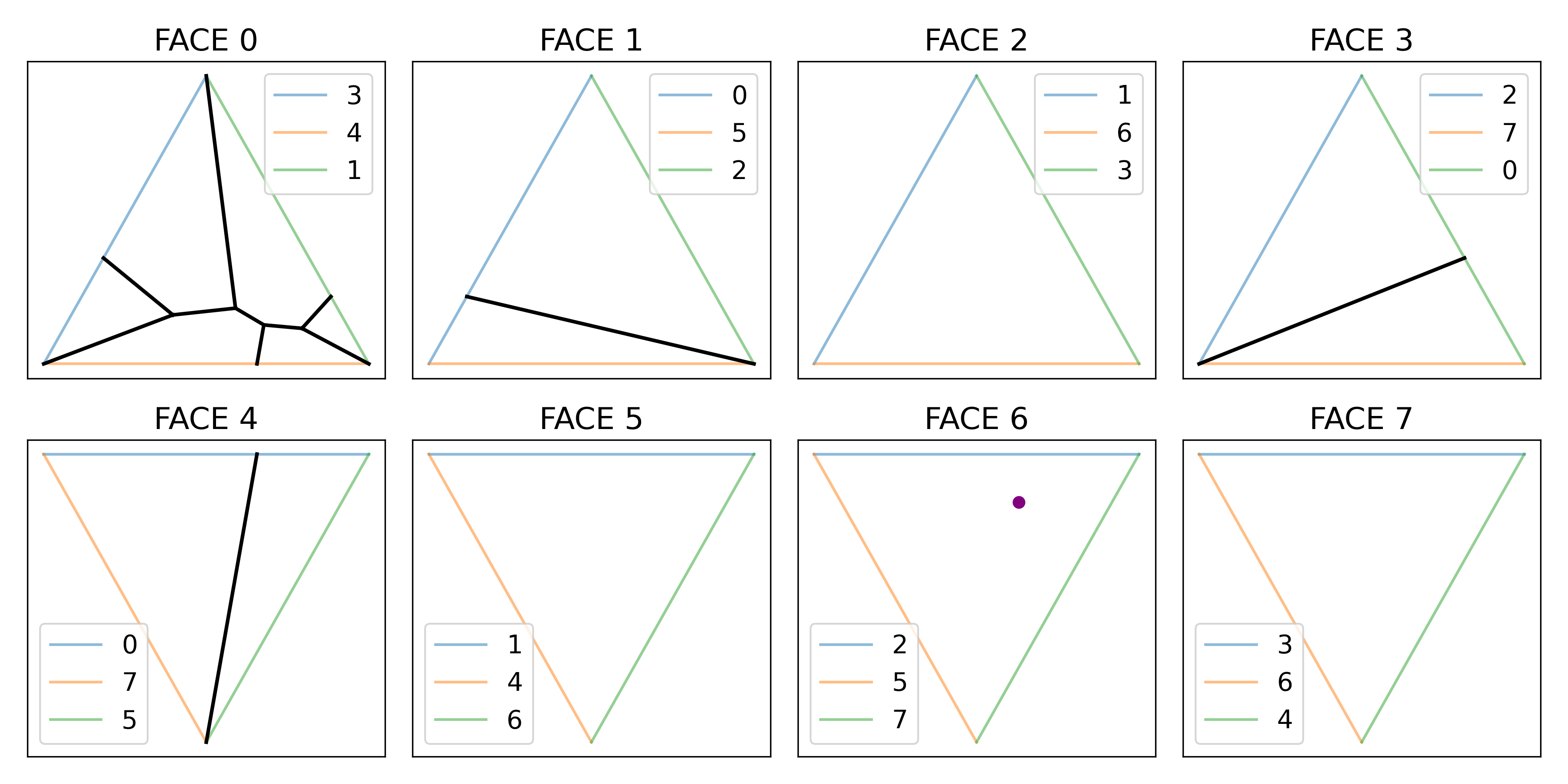}
}
\subfigure[Path unfoldings from Face 0]{
\includegraphics[height=\imgheight]{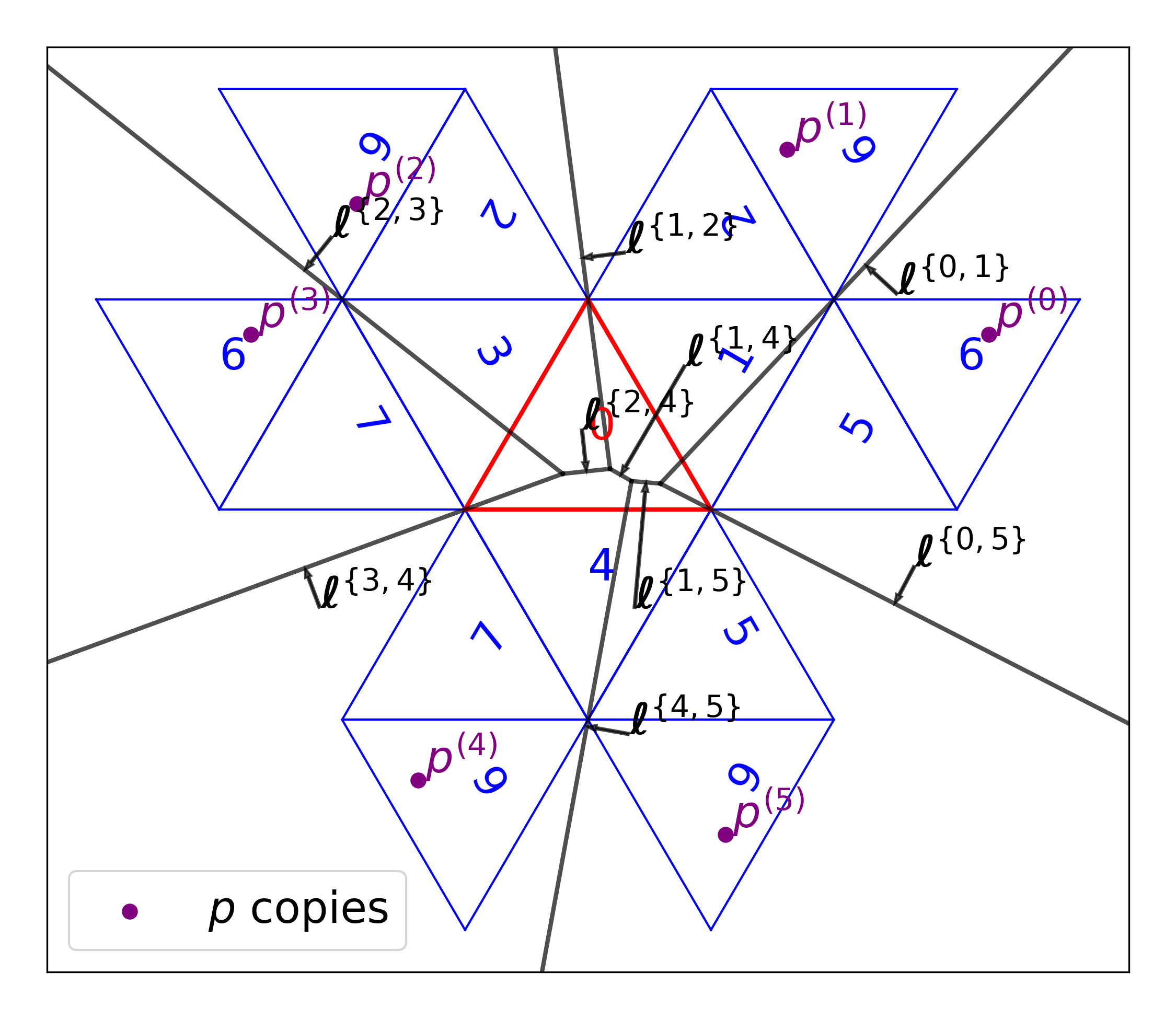}
}
\caption{Cut locus of $p=(0.45,0.5)$ chosen in $\tau_1$}
\label{fig:octatau1locus}
\end{figure}

\meminisection{Region $\tau_1$:}
We find that the structure of the cut locus of $p\in \tau_1$ is isomorphic to \autoref{fig:octatau1locus}, with four vertices $x^{\{0,1,5\}}$, $x^{\{2,3,4\}}$, $x^{\{1,2,4\}}$, and $x^{\{1,4,5\}}$ incident to three line segments. All cut locus vertices are within Face 0.

We will now consider the cut locus of points in $v_1$ and $e_1$, which are not on the interior of Face 6. 
Since $v_1$ is a point, we will simply use \autoref{alg:cutlocus} to calculate it. 
For $e_1$, we use the same method as before, though with a different set of path unfoldings.

\begin{figure}[ht!]
    \centering
\includegraphics[height=175 pt]{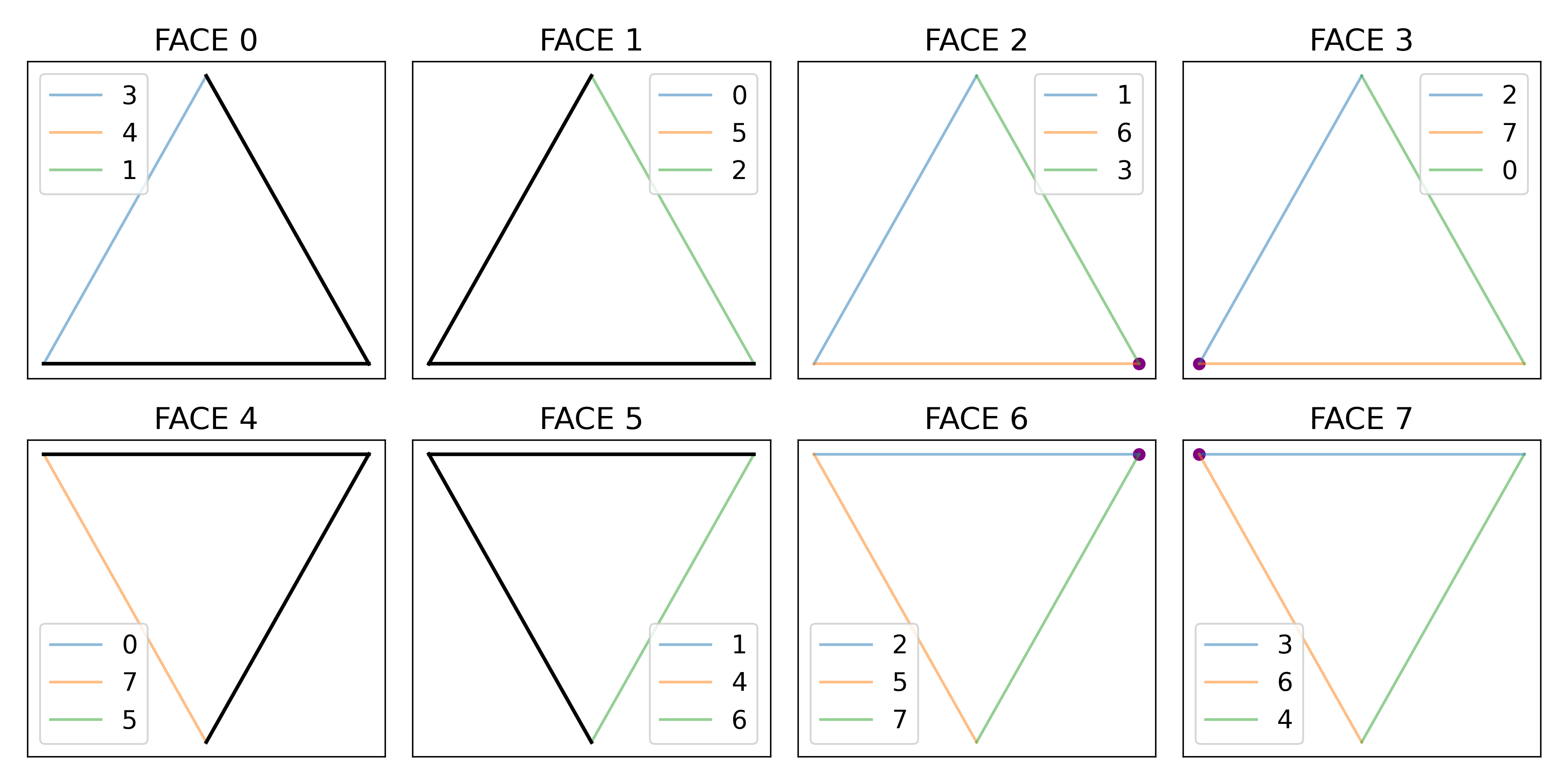}
\caption{Cut locus of $p$ chosen to be $v_1$}
\label{fig:octav1locus}
\end{figure}

\meminisection{Point $v_1$:}
From \autoref{alg:cutlocus}, we find that the cut locus of $p$ chosen on vertex $v_1$ is isomorphic to \autoref{fig:octav1locus}.
When viewed on the octahedron, this structure is a star, with four lines joined at a point.
The point is the antipode of $p$, and the four lines are the four octahedron edges incident to antipode.

\meminisection{Line $e_1$:}
We claim that the cut locus of $p\in e_1$ is isomorphic to \autoref{fig:octae1locus}.
Due to the reflective symmetry about the equator of an octahedron, we have that the cut loci on Faces 0, 1, 2, and 3 are symmetric to the cut loci on Faces 4, 5, 6, and 7 respectively.
We find that the cut loci on Faces 1, 3, 5, and 7 arise from two copies of $p$.
From the symmetry of the two copies, the cut loci is exactly one edge of each of these faces.
We display the case for Face 1 in \autoref{fig:octae1unfold}(b).

\renewcommand{\imgwidth}{.45\linewidth}
\begin{figure}[ht!]
    \centering
    \subfigure[$p=(0,1)$ on Face 6]{
    \includegraphics[width=\imgwidth]{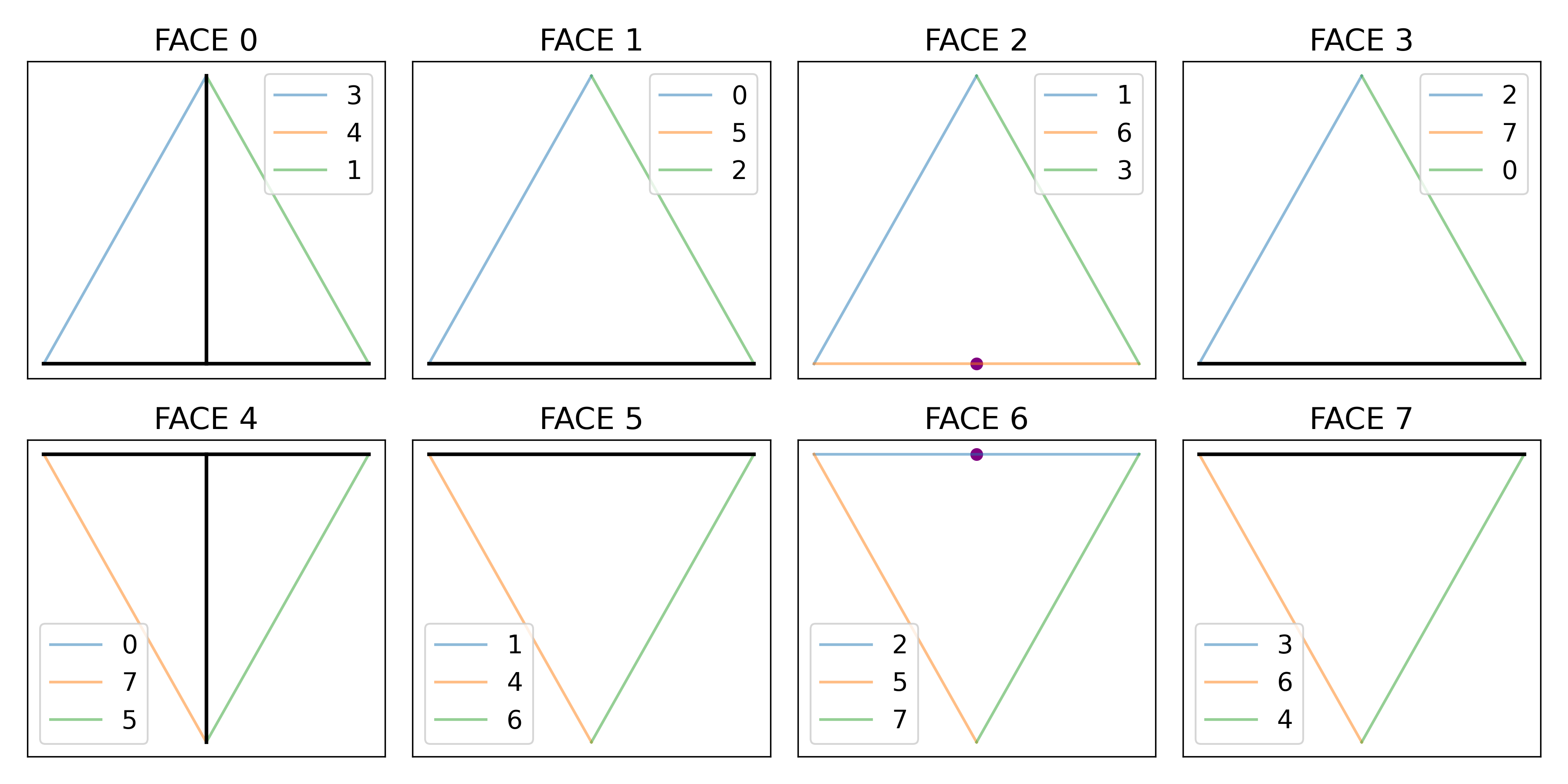}
    }
    \subfigure[$p=(\frac{\sqrt3}2,1)$ on Face 6]{
    \includegraphics[width=\imgwidth]{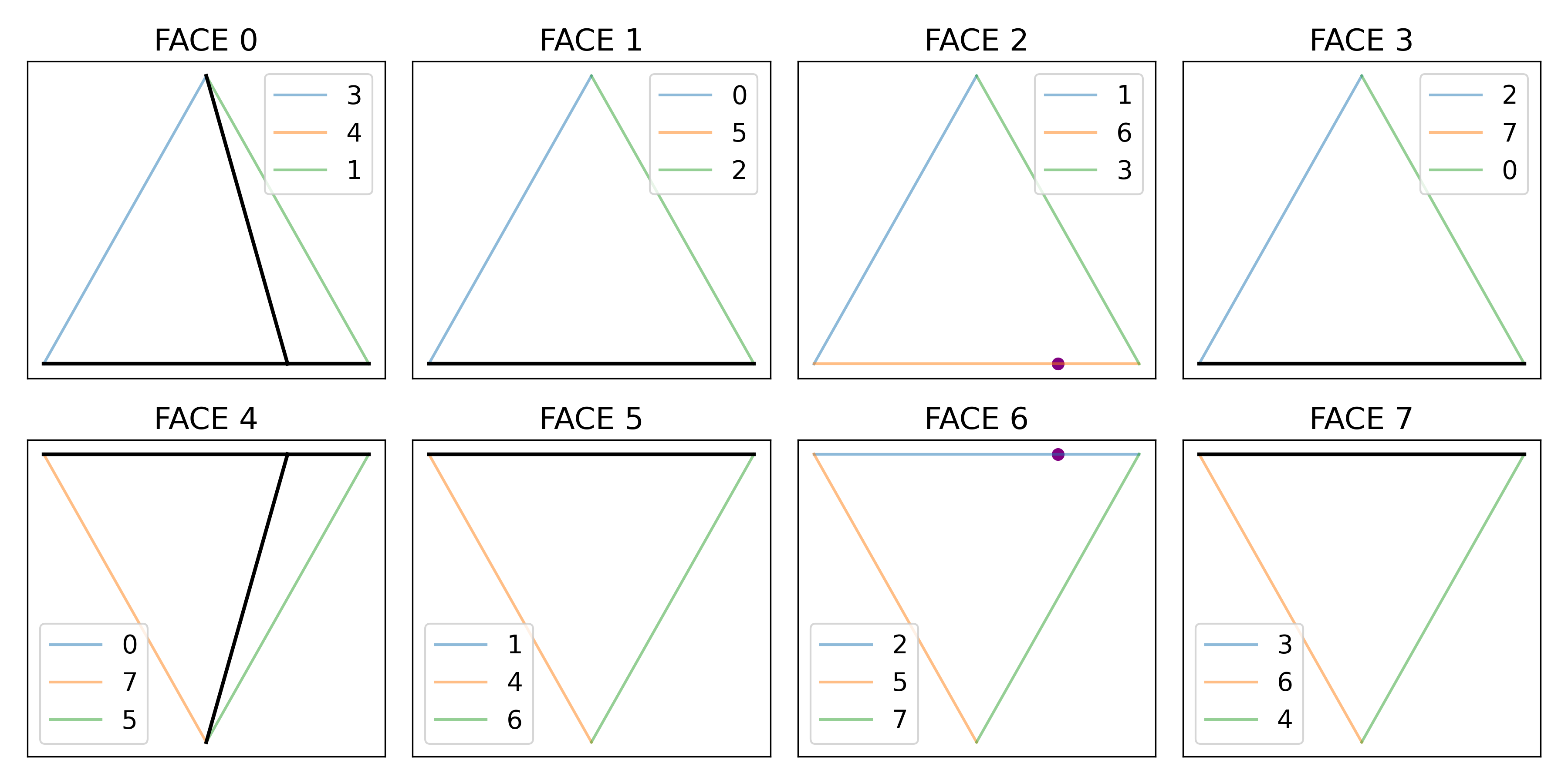}
    }
    \caption{Cut locus of $p$ chosen on $e_1$}
\label{fig:octae1locus}
\end{figure}

\renewcommand{\imgheight}{169 pt}
\begin{figure}[ht!]
    \centering
    \subfigure[Path unfoldings from Face 0]{
    \includegraphics[height=\imgheight]{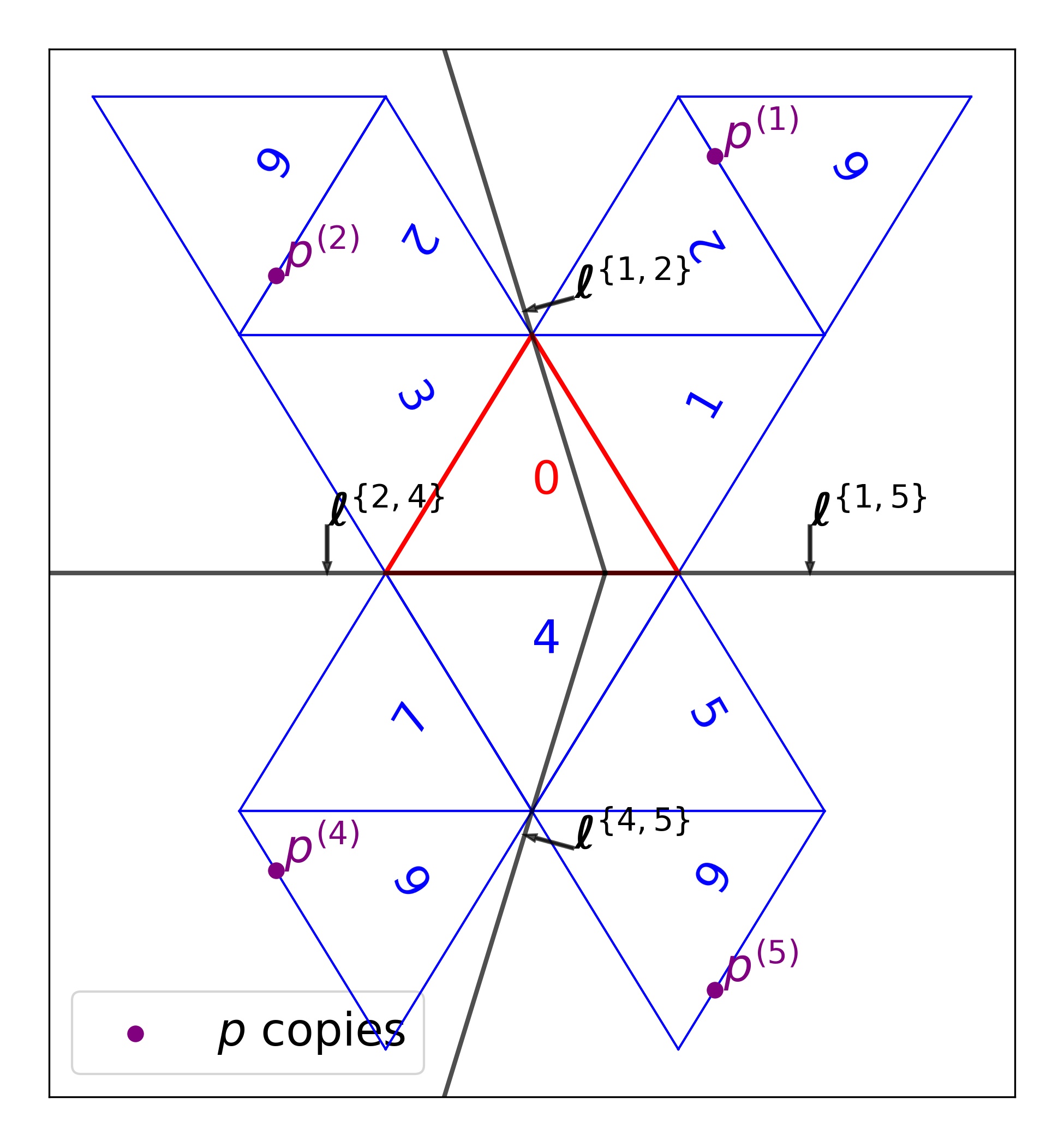}
    }
    \subfigure[Path unfoldings from Face 1]{
    \includegraphics[height=\imgheight]{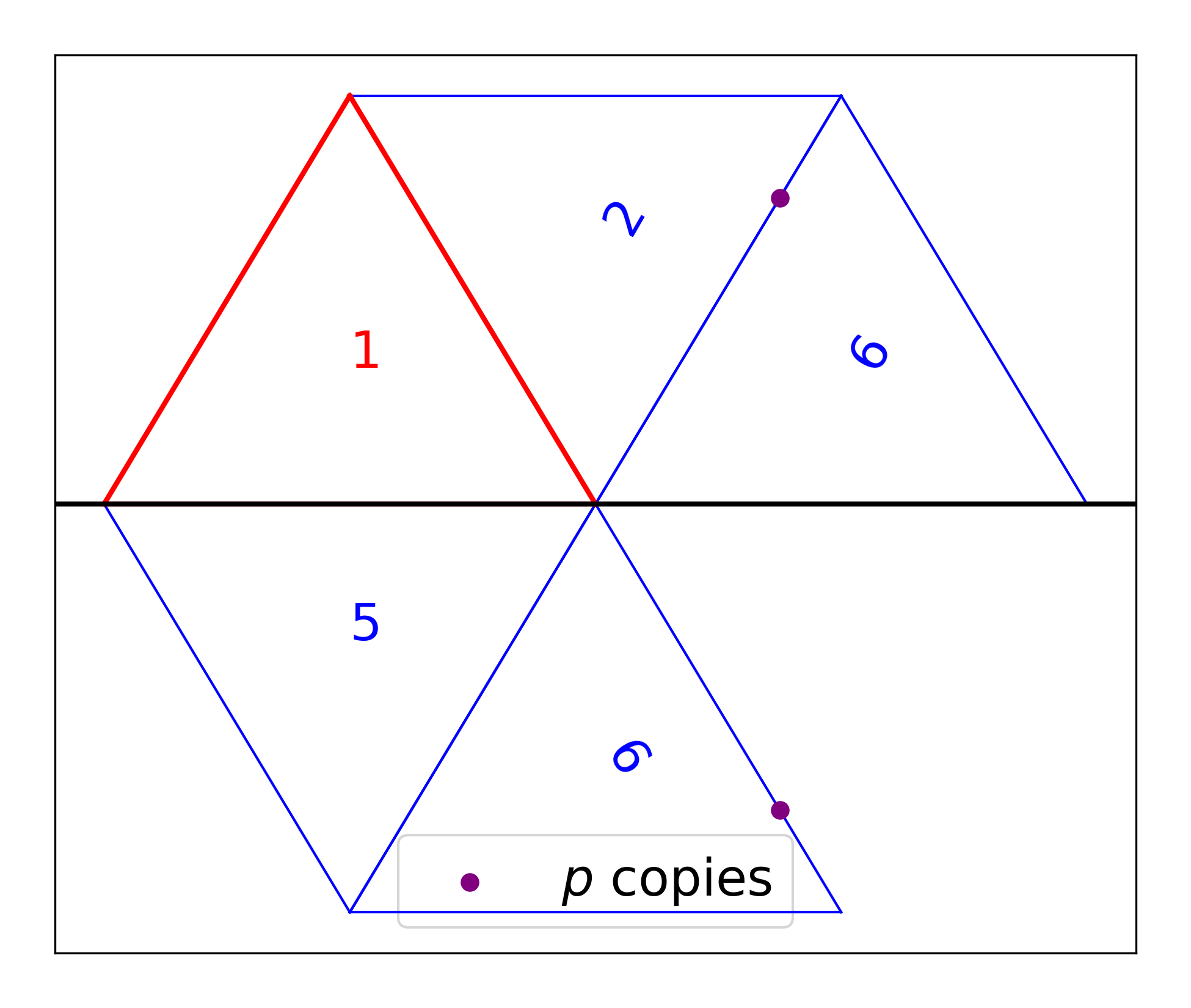}
    }
    \caption{Copies of $p$ and path unfoldings that contribute to the cut locus of $p=(\frac{\sqrt3}2,1)$ chosen on $e_1$}
\label{fig:octae1unfold}
\end{figure}

For Face 0, we find that the cut locus arises from four copies of $p$, labeled\footnote{We choose this so that each point arises from the same path unfolding as the point with the same label earlier in this section.} $p^{(1)},p^{(2)},p^{(4)},p^{(5)}$ as in \autoref{fig:octae1unfold}(a).
The cut locus is isomorphic to a star, with four line segments joined at a vertex $x^{\{1,2,4,5\}}$, the antipode of $p$ on the octahedron.
Three of the line segments occur on Face 0, while $\ell^{\{4,5\}}$ only occurs on Face 4.
The cut locus has symmetry with respect to reflection about the bottom edge of Face 0, a result of the path unfoldings that create the four copies of $p$ having the same symmetry.
Additionally, when $p$ is at the midpoint of edge $e_1$, the cut locus has symmetry with respect to reflection about $x=0$.

\subsubsection{Geodesic Complexity Lower Bound}
\label{subsubsec:octa_lower_bound_proof}

We will use \autoref{thm:quite_simplex} to show that $GX(X)\geq 4$ (needing at least 5 sets in any geodesic motion planner).
To do this, we will embed simplices into $X\times X$ using \autoref{lem:nice_embeddings}.
\begin{figure}[ht!]
    \centering
    \includegraphics[height=150 pt]{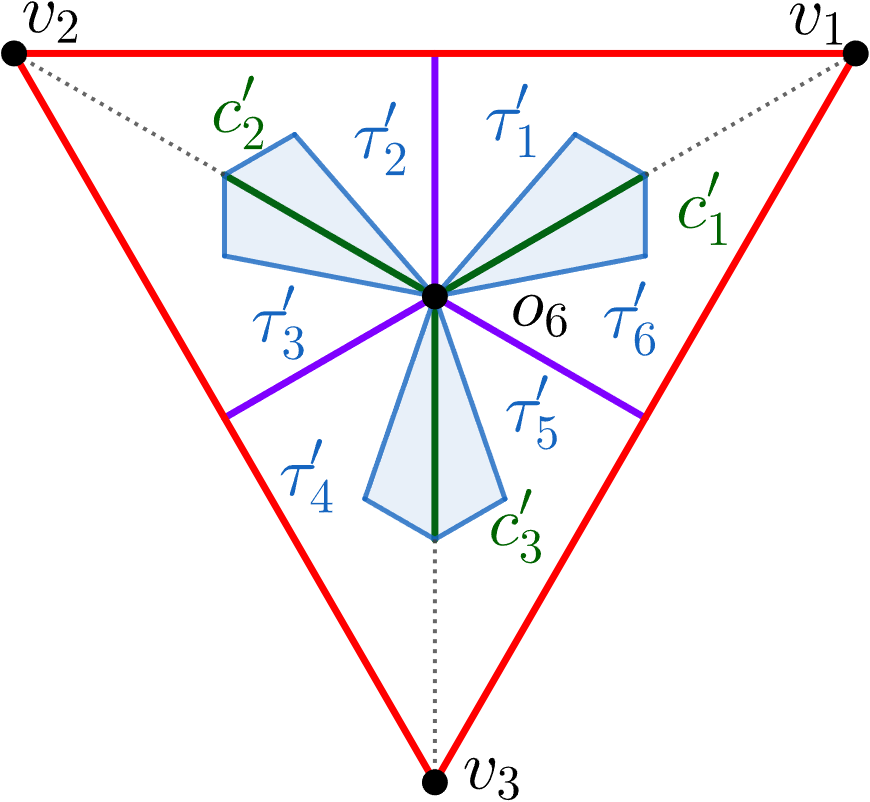}
    \caption{Constructed $\tau_i'$ and $c_i'$.
    Endpoints are the midpoint of each $c_i$ and the centroid of each $\tau_i$.}
    \label{fig:tau_c_prime}
\end{figure}
Let $p^*=o_6$, and $q^*=(0,0)$ on Face 0, and let $\tau_i'\subseteq \overline{\tau_i}$ and $c_i'\subseteq \overline{c_i}$ be closed sets constructed as in \autoref{fig:tau_c_prime}.
Let $p_{c_1}, p_{c_2}, p_{c_3}\colon [0,1]\hookrightarrow X$ be embeddings such that $\img{p_{c_i}}=c_i'$ and $p_{c_i}(1)=o_6=p^*$.
Similarly, let $p_{\tau_1}, \dots,p_{\tau_6}\colon T_2\hookrightarrow X$ be embeddings such that $\img{p_{\tau_i}}=\tau_i'$ and $p_{\tau_i}(S_2)=c_j'$ for $c_j'$ bordering $\tau_i'$.
(Recall from \autoref{def:simplex_shenanigans} that $T_2$ is a triangle in $\mathbb R^2$ and $S_2\subseteq T_2$ is one of its edges.)

Let $\mathcal J:=\{J\subseteq \{0,\dots,5\}:|J|=3\}$, and let $\mathcal I:=\{\{i,j\}\subseteq \{0,\dots,5\}:j\equiv i+1\mod 6\}$ be collections of subsets of $\{0,\dots,5\}$.
We will define $q^{\tau_i,A,J}\colon T_2\hookrightarrow X$ analogously to the construction in \autoref{subsubsec:tetra_lower_bound_proof}.
Explicitly, for $A\in\mathcal I\cup\mathcal J$, $J\in\mathcal J$, and $|A\cap J|=2$, $q^{\tau_i,A,J}$ varies linearly between $x^{A}$ and $x^{J}$ via $q^{\tau_i,A,J}(\textbf{t})=(1-\Sigma \textbf{t})\cdot x^{A}(p_{\tau_i}(\textbf{t}))+(\Sigma \textbf{t})\cdot x^{J}(p_{\tau_i}(\textbf{t}))$, where $\Sigma \textbf{t}$ is the sum $\sum\limits_i \textbf{t}_i$.
Since $x^{A}$ and $x^{J}$ are on line $\ell^{A\cap J}$, $q^{\tau_i,A,J}$ is also on $\ell^{A\cap J}$.
Additionally, for $\textbf{t}\in S_2$, we have $q^{\tau_i,A,J}(\textbf{t})=x^{J}(p_{\tau_i}(\textbf{t}))$.
We display $q^{\tau_1,\{0,5\},\{0,1,5\}}$, $q^{\tau_1,\{1,2,4\},\{1,4,5\}}$, and relevant $x^{A}$ for a few selections of $p\in \tau_1'$ in \autoref{fig:octa_q_vis}.

\renewcommand{\imgwidth}{.48\linewidth}
\begin{figure}[ht!]
    \centering
    \subfigure[$p$ chosen on vertex of $\tau_1'$]{
    \includegraphics[width=\imgwidth]{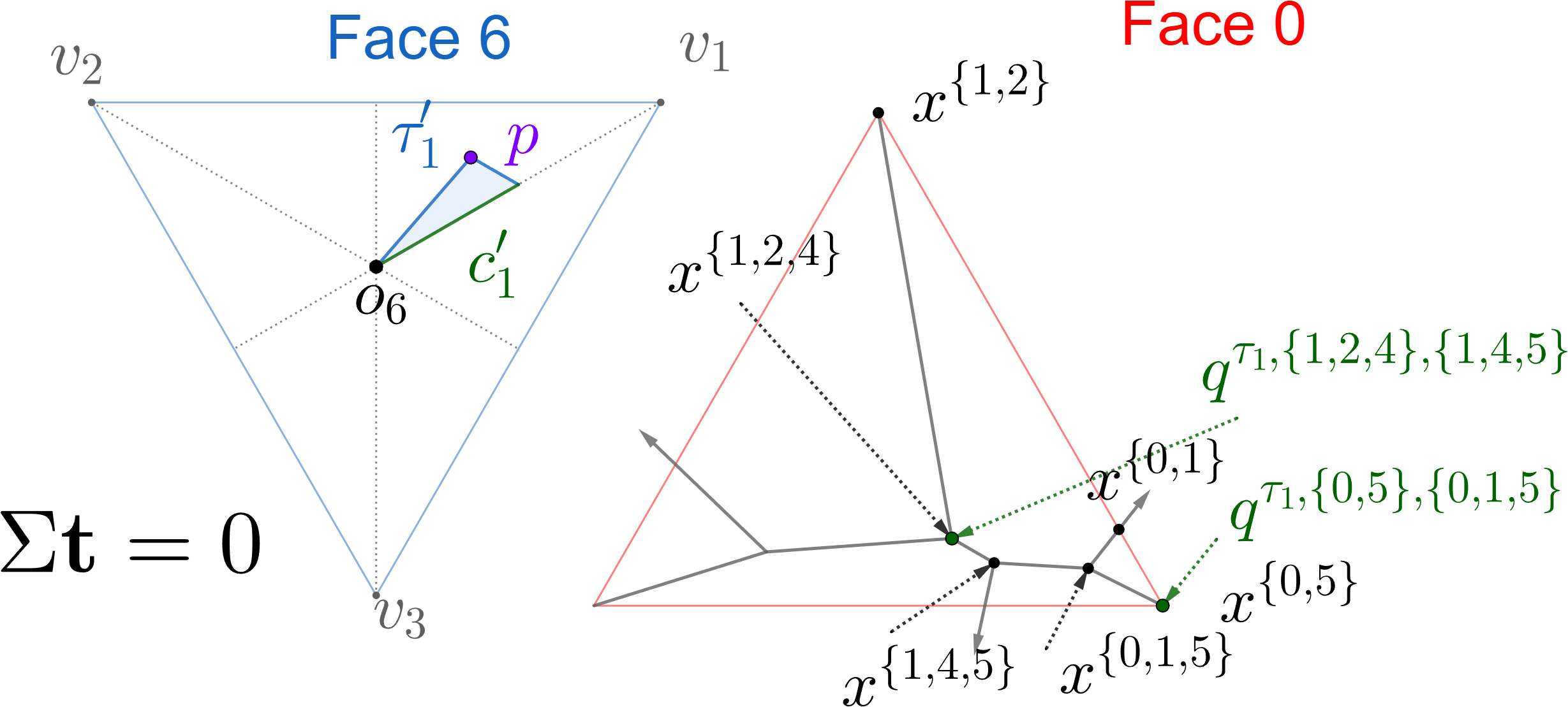}
    }
    \subfigure[$p$ chosen within $\tau_1'$]{
    \includegraphics[width=\imgwidth]{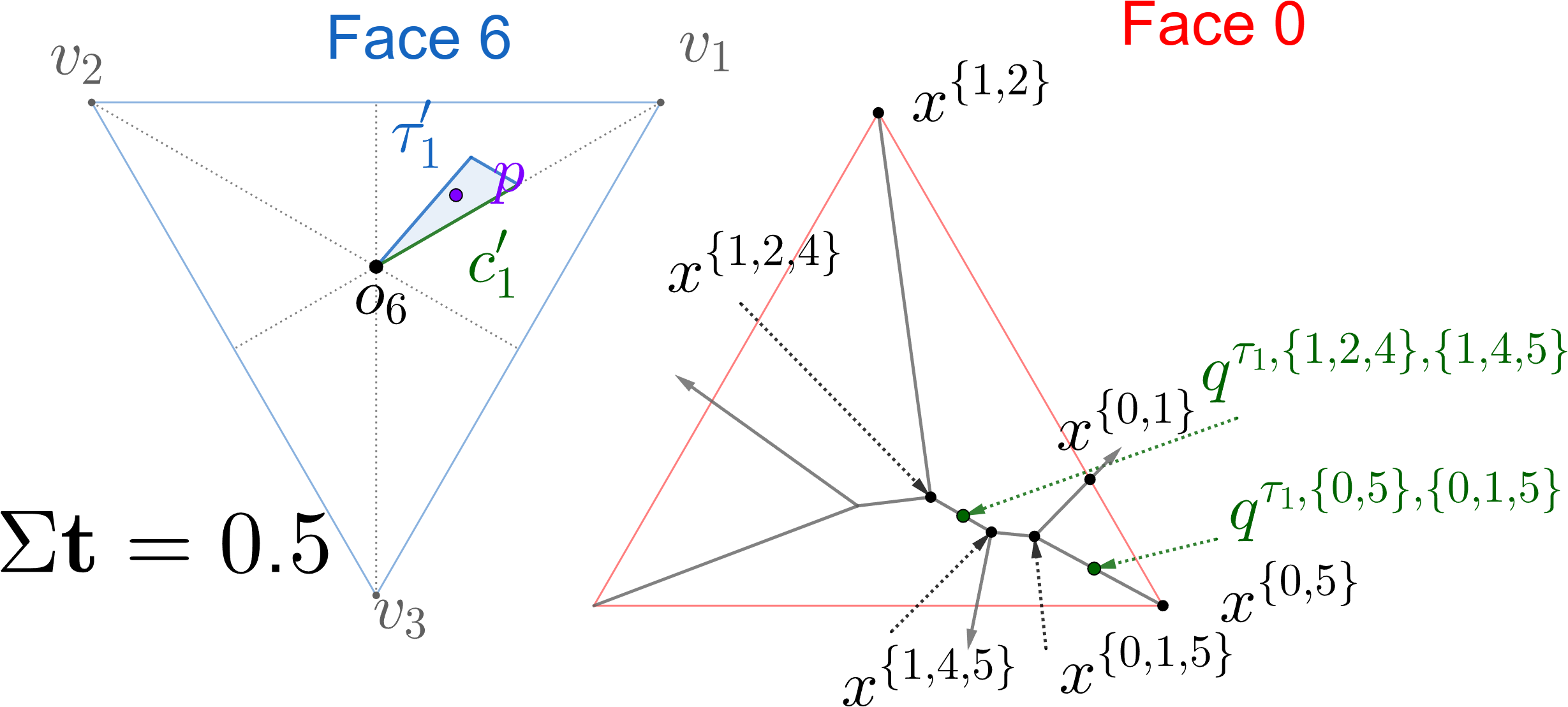}
    }
    \subfigure[$p$ chosen at endpoint of $c_1'$]{
    \includegraphics[width=\imgwidth]{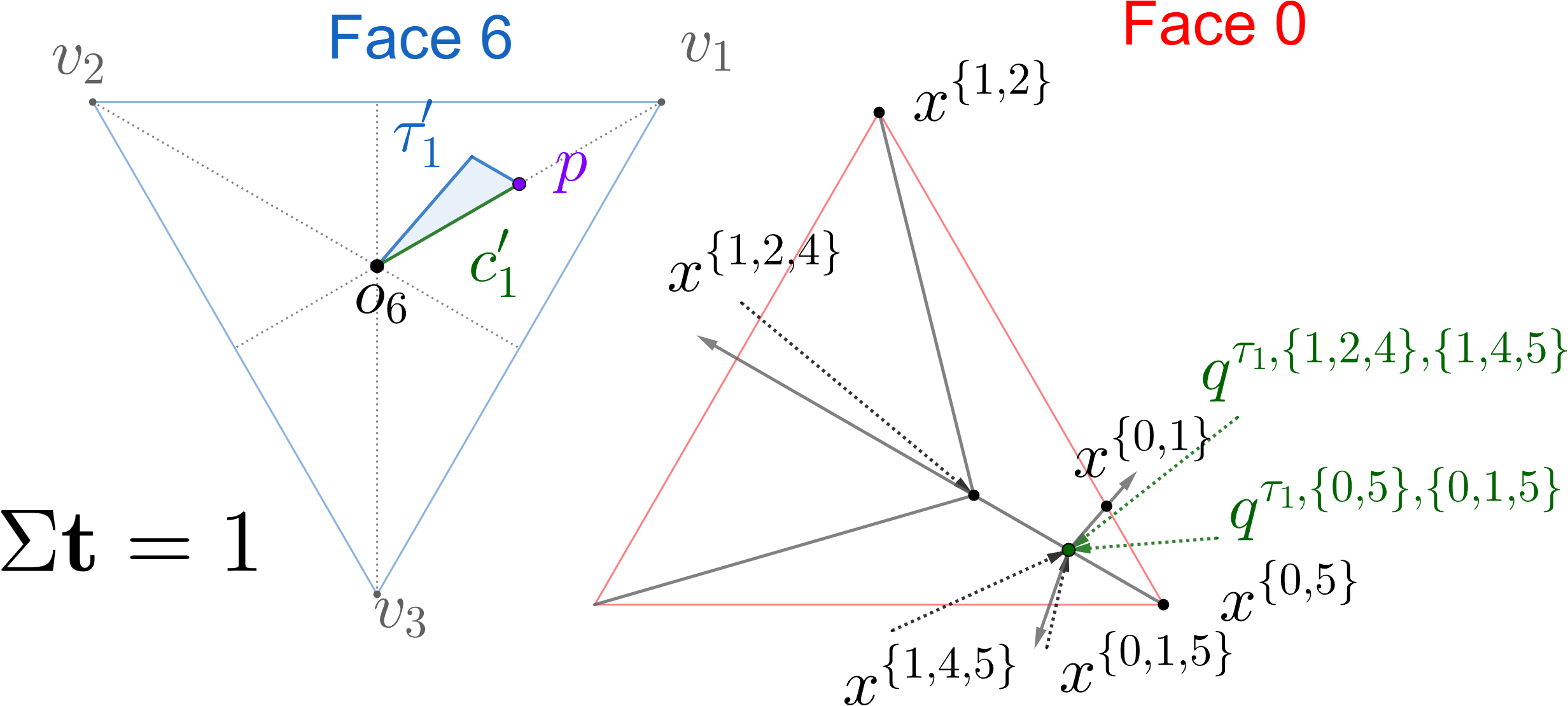}
    }
    \subfigure[$p$ chosen as $o_6$]{
    \includegraphics[width=\imgwidth]{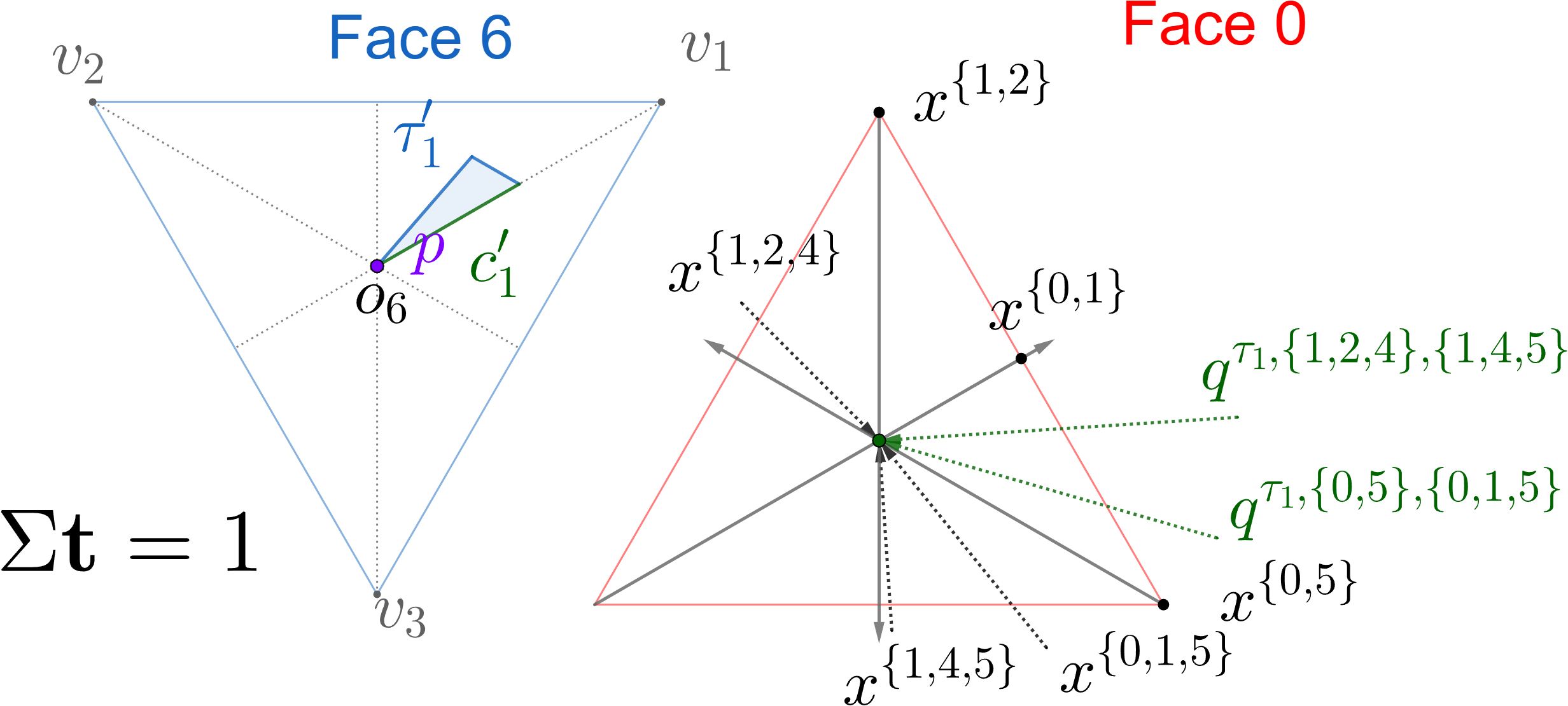}
    }
    \caption{Defined points on Face 0 as $p$ varies about $\tau_1'$}
\label{fig:octa_q_vis}
\end{figure}

We will now describe our $\mathcal F_i$ to use \autoref{thm:quite_simplex} on the octahedron:

\begin{compactitem}
    \item
    $\mathcal F_0=\{\Delta_0\mapsto (p^*,q^*)\}$.
    
    \item
    $\mathcal F_1$ will be defined by three lines.
    Consider $p_{c_1}$ and the map $x^{\{0,1,4,5\}}\circ p_{c_1}$. 
    By \autoref{lem:nice_embeddings}, this describes an embedding $f\colon \Delta_{1}\hookrightarrow X\times X$ where $\pi_1(f(t,1-t))=p_{c_1}(t)$ and $\pi_2(f(t,1-t))= x^{\{0,1,4,5\}}(p_{c_1}(t))$. 
    This map follows $(p,x^{\{0,1,4,5\}}(p))$ as $p$ varies across $c_1'$.
    We will construct symmetric embeddings using maps $p_{c_2}$ and $x^{\{2,3,4,5\}}\circ p_{c_2}$, and maps $p_{c_3}$ and $x^{\{0,1,2,3\}}\circ p_{c_3}$.
    
    \item
    $\mathcal F_2$ will be defined by twelve simplices.
    We will vary $p$ along each of the $\tau_i'$, and choose $q$ to be one of the two vertices of the cut locus whose limit is the vertex chosen in $\mathcal F_1$ as $p$ approaches $c_i'$.
    For the case of $p_{\tau_1}$, these vertices will be $x^{\{0,1,5\}}$ and $x^{\{1,4,5\}}$.
    Explicitly, consider the maps $p_{\tau_1}$ and $x^{\{0,1,5\}}\circ p_{\tau_1}$.
    By \autoref{lem:nice_embeddings}, this describes an embedding $f\colon \Delta_2\hookrightarrow X\times X$ where for $\textbf{t}\in T_2$, $\pi_1(f(\textbf{t},1-\Sigma \textbf{t}))=p_{\tau_1}(\textbf{t})$ and $\pi_2(f(\textbf{t},1-\Sigma \textbf{t}))=
    x^{\{0,1,5\}}(p_{\tau_1}(\textbf{t}))$.
    This map follows $(p,x^{\{0,1,5\}})$ as $p$ varies around $\tau_1'$ (see \autoref{fig:octatau1locus}).
    We will construct an embedding similarly from the maps $p_{\tau_1}$ and $x^{\{1,4,5\}}\circ p_{\tau_1}$.
    We construct the ten remaining embeddings symmetrically from the five other $\tau_i$.
    
    \item
    $\mathcal F_3$ will be defined by 30 simplices, five for each $\tau_i$, corresponding to the five cut locus lines incident to the cut locus vertices used in $\mathcal F_2$.
    Consider $p_{\tau_1}$ and the map $q^{\tau_1,\{0,1\},\{0,1,5\}}\funcprodsymb (x^{\{0,1,5\}}\circ p_{\tau_1})$.
    Since when $\Sigma \textbf{t}=1$, we have $q^{\tau_1,\{0,1\},\{0,1,5\}}(\textbf{t})=x^{\{0,1,5\}}(p_{\tau_1}(\textbf{t}))$, and these points are otherwise distinct, we have by \autoref{lem:nice_embeddings} that this defines an embedding $f\colon \Delta_3\hookrightarrow X\times X$ whose image looks like $(p,q)$ for $p\in \tau_1'$ and $q\in \ell^{\{0,1\}}(p)$.
    We will define two similar embeddings with $p_{\tau_1}$ using the maps $q^{\tau_1,\{0,5\},\{0,1,5\}}\funcprodsymb (x^{\{0,1,5\}}\circ p_{\tau_1})$ and $q^{\tau_1,\{4,5\},\{1,4,5\}}\funcprodsymb (x^{\{1,4,5\}}\circ p_{\tau_1})$.

    For another embedding, consider $p_{\tau_1}$ and the map $(x^{\{0,1,5\}}\circ p_{\tau_1})\funcprodsymb (x^{\{1,4,5\}}\circ p_{\tau_1})$.
    When $\Sigma \textbf{t}=1$, we have $x^{\{0,1,5\}}(p_{\tau_1}(\textbf{t}))=x^{\{1,4,5\}}(p_{\tau_1}(\textbf{t}))$ since $p_{\tau_1}(\textbf{t})\in \overline{c_1}$.
    These points are otherwise distinct, so this defines an embedding $f\colon \Delta_3\hookrightarrow X\times X$ whose image looks like $(p,q)$ for $p\in \tau_1'$ and $q\in \ell^{\{1,5\}}(p)$.

    Finally, consider $p_{\tau_1}$ and the map $(x^{\{1,4,5\}}\circ p_{\tau_1})\funcprodsymb q^{\tau_1,\{1,2,4\},\{1,4,5\}}$.
    Since when $\Sigma \textbf{t}=1$, we have $x^{\{1,4,5\}}(p_{\tau_1}(\textbf{t}))=q^{\tau_1,\{1,2,4\},\{1,4,5\}}(\textbf{t})$, and these points are otherwise distinct, this defines an embedding $f\colon \Delta_3\hookrightarrow X\times X$ whose image looks like $(p,q)$ for $p\in \tau_1'$ and $q\in \ell^{\{1,4\}}(p)$.

    This describes the five embeddings corresponding to $\tau_1$.
    We construct five symmetric embeddings for each of the five remaining $\tau_i$.
    
    \item
    $\mathcal F_4$ will be defined by 36 simplices, six for each $\tau_i$, corresponding to the three cut locus Voronoi cells incident to each each cut locus vertex chosen in $\mathcal F_2$.
    We will inspect $\tau_1$ (with cut locus structure as in \autoref{fig:octatau1locus}) and the vertex $x^{\{0,1,5\}}$.
    Consider $p_{\tau_1}$ and the map $(x^{\{0,1,5\}}\circ p_{\tau_1})\funcprodsymb q^{\tau_1,\{0,1\},\{0,1,5\}}\funcprodsymb (x^{\{1,4,5\}}\circ p_{\tau_1})$.
    When $\Sigma \textbf{t}=1$, we have $x^{\{0,1,5\}}(p_{\tau_1}(\textbf{t}))=
    x^{\{1,4,5\}}(p_{\tau_1}(\textbf{t}))= 
    q^{\tau_1,\{0,1\},\{0,1,5\}}(\textbf{t})$, and these points are otherwise not collinear.
    Thus, by \autoref{lem:nice_embeddings}, this defines an embedding $f\colon \Delta_4\hookrightarrow X\times X$ where for $\textbf{d}\in \Delta_4$, $f(\textbf{d})=(p,q)$ for $p\in \tau_1'$ and $q$ in the convex hull of $x^{\{0,1,5\}}(p)$, $x^{\{1,4,5\}}(p)$, and a point on $\ell^{\{0,1\}}$.
    We will create two similar embeddings by considering the maps $(x^{\{0,1,5\}}\circ p_{\tau_1})\funcprodsymb (x^{\{1,4,5\}}\circ p_{\tau_1})\funcprodsymb q^{\tau_1,\{0,5\},\{0,1,5\}}$ and $(x^{\{0,1,5\}}\circ p_{\tau_1})\funcprodsymb q^{\tau_1,\{0,5\},\{0,1,5\}}\funcprodsymb q^{\tau_1,\{0,1\},\{0,1,5\}}$.

    We will now inspect the other vertex, $x^{\{1,4,5\}}$.
    We will create three similar embeddings by considering the map $(x^{\{1,4,5\}}\circ p_{\tau_1})\funcprodsymb (x^{\{0,1,5\}}\circ p_{\tau_1})
    \funcprodsymb q^{\tau_1,\{1,2,4\},\{1,4,5\}}$,
    the map
    $(x^{\{1,4,5\}}\circ p_{\tau_1})\funcprodsymb q^{\tau_1,\{1,2,4\},\{1,4,5\}}\funcprodsymb q^{\tau_1,\{4,5\},\{1,4,5\}}$, 
    and the map $(x^{\{1,4,5\}}\circ p_{\tau_1})\funcprodsymb q^{\tau_1,\{4,5\},\{1,4,5\}}\funcprodsymb (x^{\{0,1,5\}}\circ p_{\tau_1})$.

    This describes the six embeddings corresponding to $\tau_1$.
    We construct six symmetric embeddings for each of the five remaining $\tau_i$.
\end{compactitem}

We visualize a few examples of sets in $\mathcal F_2$, $\mathcal F_3$, and $\mathcal F_4$ in \autoref{fig:octaproof}.

We will verify the properties of \autoref{thm:quite_simplex} for a few embeddings and geodesic choices.
We assert that similar arguments extend to all embeddings by symmetry in their construction.
First, each $\mathcal F_i$ contains the embeddings of $\Delta_i\hookrightarrow X\times X$ by construction, so property (a) always holds.

\begin{compactitem}
    \item
    Consider $F_0\in\mathcal F_0$.
    Property (b) is trivial since $\img{F_0}=\{(p^*,q^*)\}$ is a single point.

    We will consider a \gls{GMPR} $\Gamma_0$ on $\img{F_0}$.
    We have six options, determined by which copy of $p^*$ ($p^{(0)},\dots,p^{(5)}$) we connect $q^*$ to.
    Assume our choice connects $q^*$ to $p^{(3)}$.

    Then our choice for $F_1\in \mathcal F_1$ will be the embedding constructed from maps $p_{c_1}$ and $x^{\{0,1,4,5\}}\circ p_{c_1}$.

    (i) is trivially satisfied.
    (ii) is satisfied since for any element $(p,q)\in \relint{F_1}$, $q$ is a point on the cut locus of $p$ with geodesics only to $p^{(0)}$, $p^{(1)}$, $p^{(4)}$, and $p^{(5)}$.
    (iii) is satisfied since there are four \gls{GMPR}s on $\img{F_1}$, each distinguished by which of $p^{(0)}$, $p^{(1)}$, $p^{(4)}$, or $p^{(5)}$ they choose to connect to.
    Thus, property (c) is satisfied.

    \item
    Consider $F_1\in \mathcal F_1$.
    Property (b) is satisfied since any geodesic for $(p,q)\in \relint{F_1}$, defined by connecting $q$ to one of $p^{(0)}$, $p^{(1)}$, $p^{(4)}$, or $p^{(5)}$, extends to a \gls{GMPR} distinguished by the same choice.

    We will consider one of the four \gls{GMPR}s $\Gamma_1$.
    We will choose $\Gamma_1$ to connect to $p^{(4)}$.

    Then our choice for $F_2\in \mathcal F_2$ will be the embedding constructed from maps $p_{\tau_1}$ and $x^{\{0,1,5\}}\circ p_{\tau_1}$.
    This embedding follows $x^{\{0,1,5\}}(p)$ as $p$ varies across $\tau_1'$.

    (i) is satisfied by construction.
    (ii) is satisfied since for any element  $(p,q)\in \relint{F_2}$, $q$ is a point on the cut locus of $p$ with geodesics only to $p^{(0)}$, $p^{(1)}$, and $p^{(5)}$.
    (iii) is satisfied since there are three \gls{GMPR}s on $\img{F_2}$, each distinguished by which of $p^{(0)}$, $p^{(1)}$, or $p^{(5)}$ they choose to connect to.
    Thus, property (c) is satisfied.
    \item
    Consider $F_2\in\mathcal F_2$.
    Property (b) is satisfied since any geodesic for $(p,q)\in \relint{F_2}$, defined by connecting $q$ to one of $p^{(0)}$, $p^{(1)}$, or $p^{(5)}$, extends to a \gls{GMPR} distinguished by the same choice.
    
    We will consider one of the three \gls{GMPR}s $\Gamma_2$.
    We will choose $\Gamma_2$ to connect to $p^{(5)}$.

    Then our choice for $F_3\in \mathcal F_3$ will be the embedding constructed from maps $p_{\tau_1}$ and $q^{\tau_1,\{0,1\},\{0,1,5\}}\funcprodsymb (x^{\{0,1,5\}}\circ p_{\tau_1})$.
    This embedding follows a segment of $\ell^{\{0,1\}}(p)$ for each choice of $p$.

    (i) is satisfied by construction.
    (ii) is satisfied since for any element  $(p,q)\in \relint{F_3}$, $q$ is a point on the cut locus of $p$ with geodesics only to $p^{(0)}$ and $p^{(1)}$.
    (iii) is satisfied since there are two \gls{GMPR}s on $\img{F_2}$, each distinguished by which of $p^{(0)}$ or $p^{(1)}$ they choose to connect to.
    Thus, property (c) is satisfied.

    \item
    Consider $F_3\in \mathcal F_3$.
    Property (b) is satisfied since any geodesic for $(p,q)\in \relint{F_3}$, defined by connecting $q$ to one of $p^{(0)}$ or $p^{(1)}$, extends to a \gls{GMPR} distinguished by the same choice.
    
    We will consider one of the two \gls{GMPR}s $\Gamma_3$.
    We will choose $\Gamma_3$ to connect to $p^{(0)}$.

    Then our choice for $F_4\in \mathcal F_4$ will be the embedding constructed from maps $p_{\tau_1}$ and $(x^{\{0,1,5\}}\circ p_{\tau_1})\funcprodsymb q^{\tau_1,\{0,1\},\{0,1,5\}}\funcprodsymb (x^{\{1,4,5\}}\circ p_{\tau_1})$.
    This embedding follows a region bounded partially by $\ell^{\{0,1\}}(p)$ and $\ell^{\{1,5\}}(p)$ for each choice of $p$.
    
    (i) is satisfied by construction.
    (ii) is satisfied since for any element $(p,q)\in \relint{F_4}$, $q$ is a point on the cut locus of $p$ with a geodesic only to $p^{(1)}$.
    (iii) is satisfied since there is only one \gls{GMPR} on $\img{F_4}$, as the only choice that extends to the interior is to go to $p^{(1)}$.
    Thus, property (c) is satisfied.

    \item
    Consider $F_4\in \mathcal F_4$.
    Property (b) is satisfied since any geodesic for $(p,q)\in \relint{F_4}$ is defined by connecting $q$ to $p^{(1)}$, and extends to the only \gls{GMPR} over $\img{F_4}$.

    We do not need to check (c) as 4 is our maximal dimension.
\end{compactitem}

Then by \autoref{thm:quite_simplex}, $GC(X)\geq 4$ (needing at least $5$ sets in any geodesic motion planner).

\renewcommand{\imgwidth}{.48\linewidth}
\begin{figure}[htbp!]
    \centering
    \subfigure[$p$ chosen on vertex of $\tau_1'$]{
    \includegraphics[width=\imgwidth]{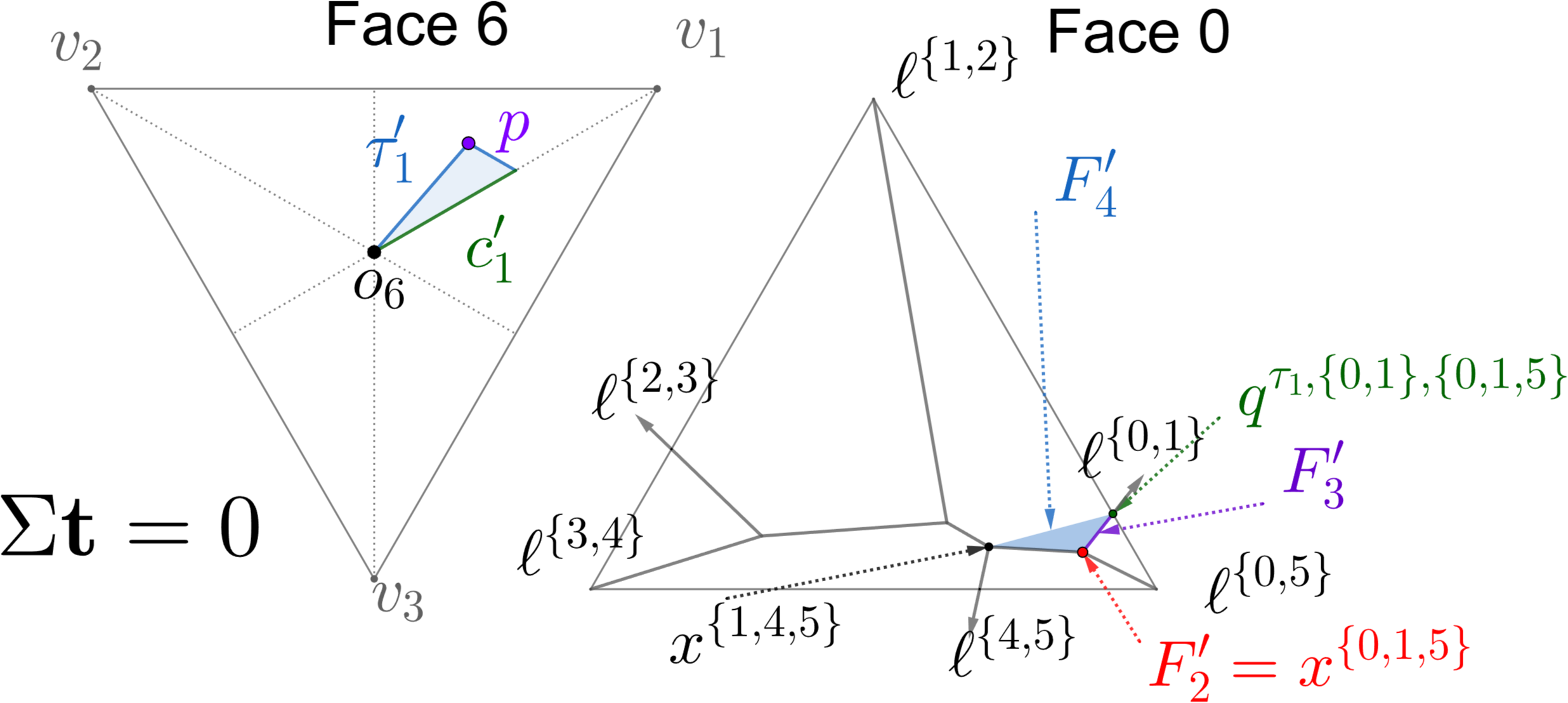}
    }
    \subfigure[$p$ chosen within $\tau_1'$]{
    \includegraphics[width=\imgwidth]{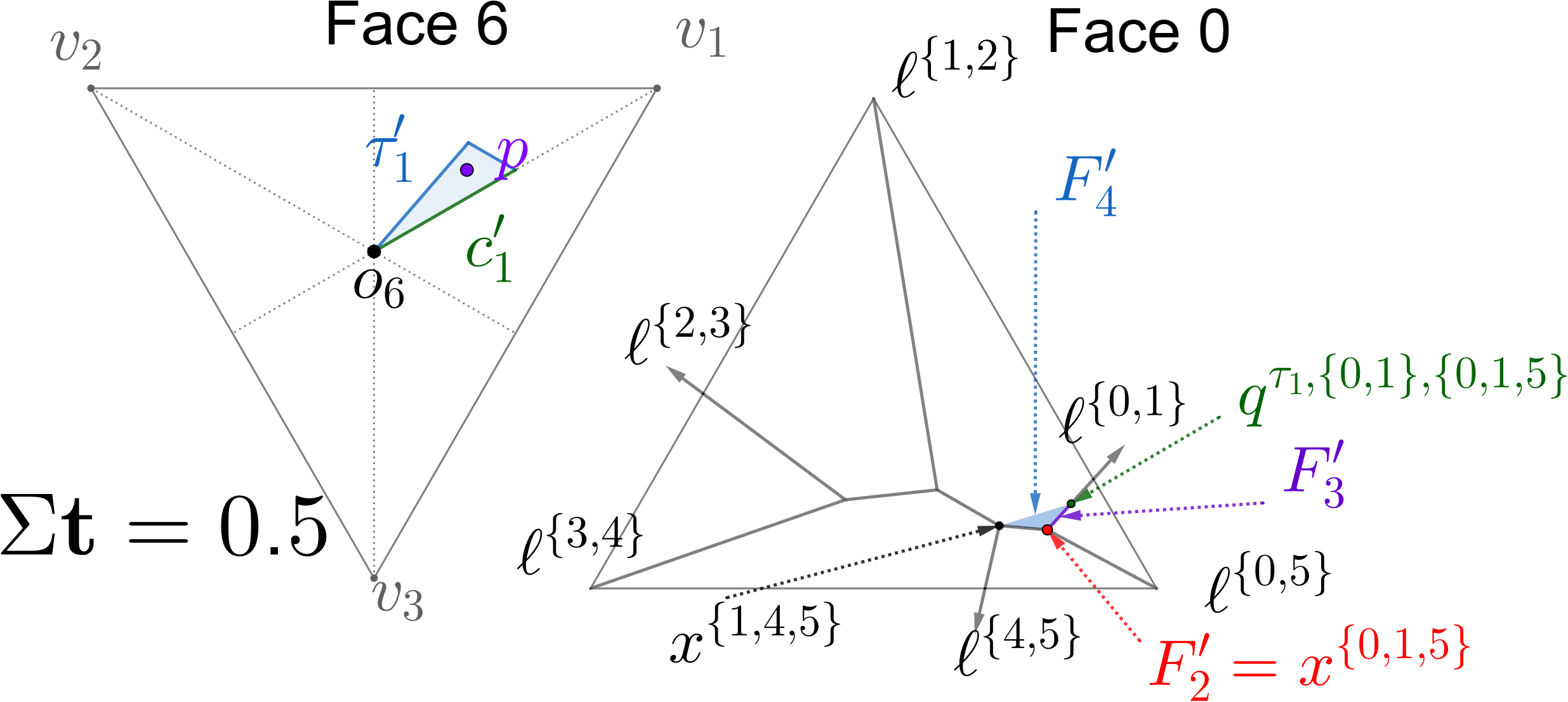}
    }
    \subfigure[$p$ chosen at endpoint of $c_1'$]{
    \includegraphics[width=\imgwidth]{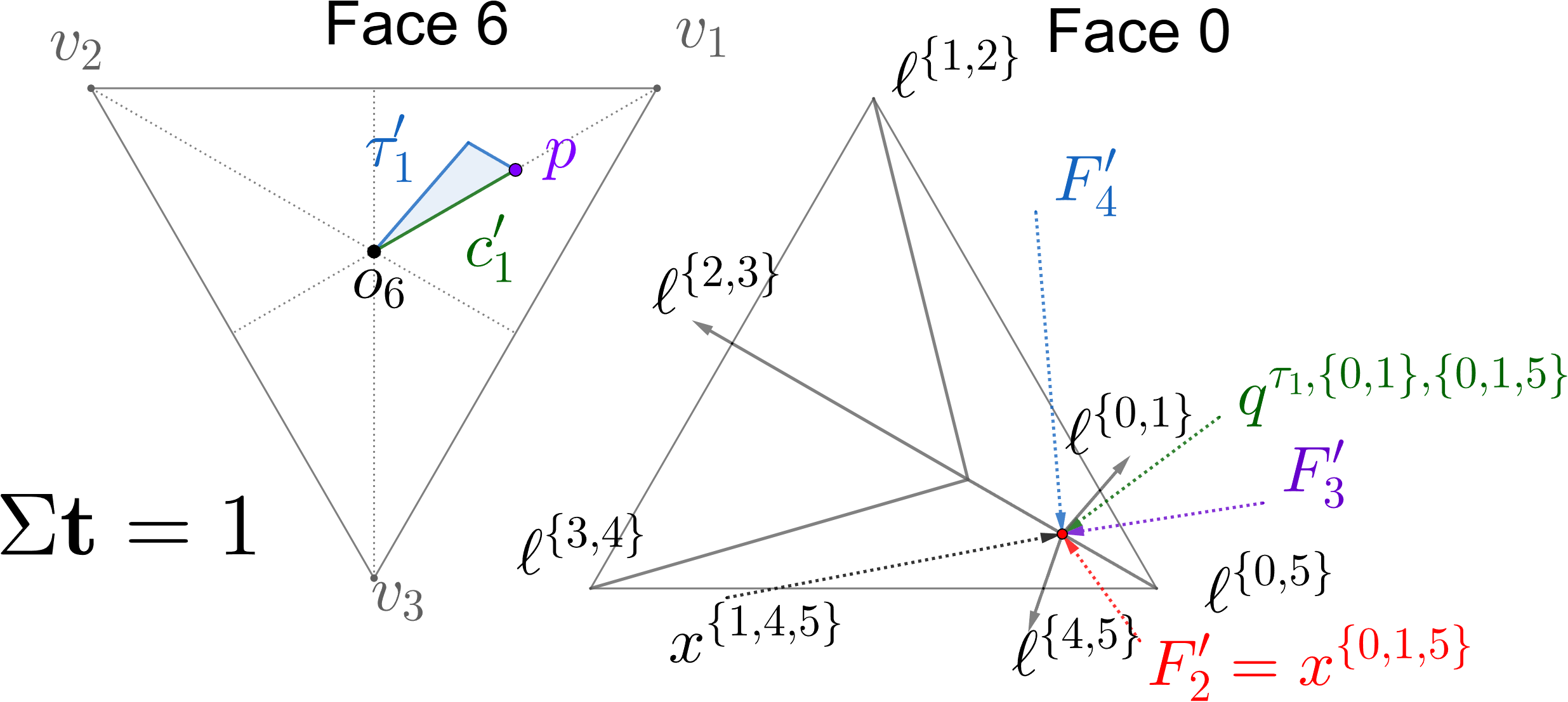}
    }
    \subfigure[$p$ chosen as $o_6$]{
    \includegraphics[width=\imgwidth]{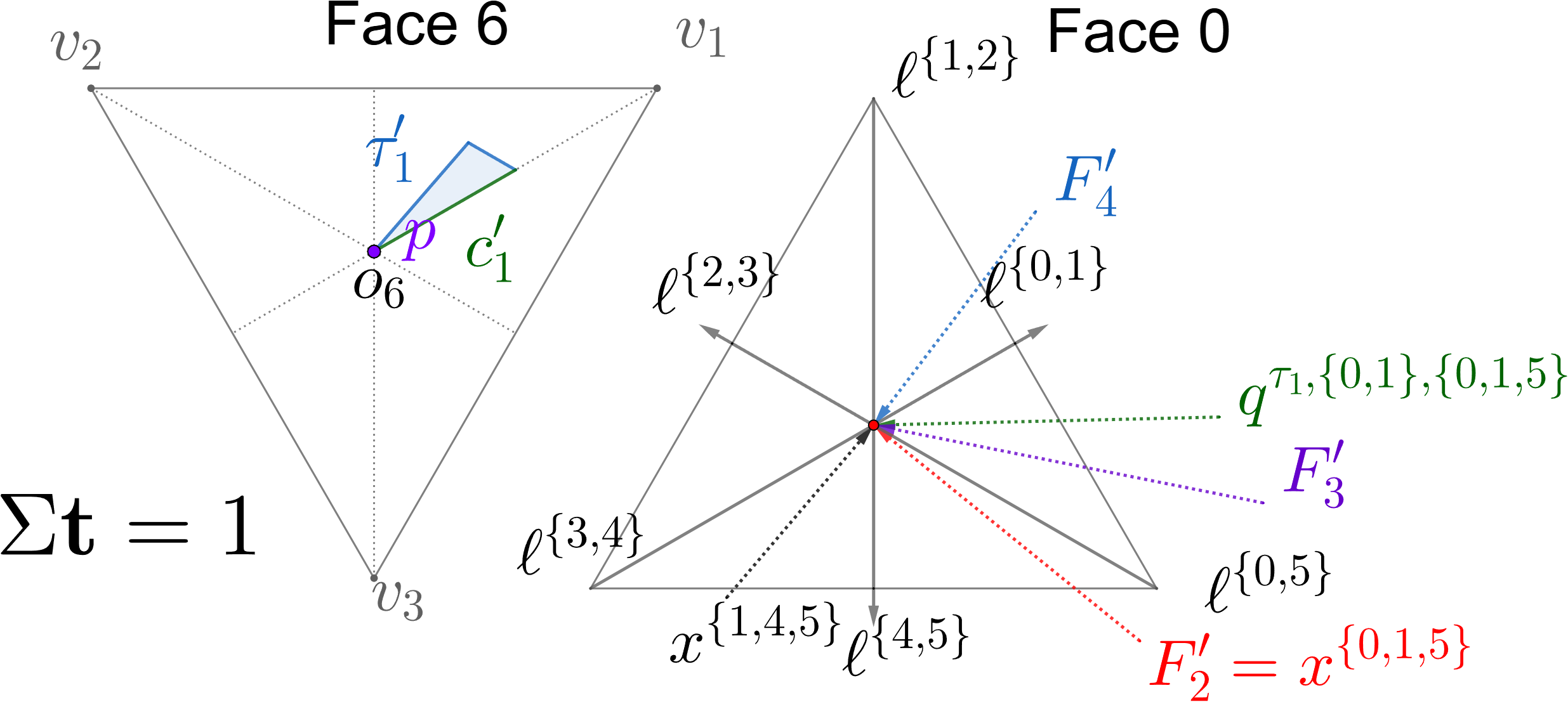}
    }
    \caption{Sets $F'_i$ as $p$ varies about $\tau_1'$}
\label{fig:octaproof}
\end{figure}

We visualize this similar to the method in \autoref{fig:tetra_proof}, with the addition of one dimension.
For $i>1$, each $F_i$ is a map $\Delta_i\hookrightarrow X\times X$, and is generated using \autoref{lem:nice_embeddings} with $k=2$.
Recall that in this case, $T_2$ is a triangle in $\mathbb R^2$, and for $\textbf{t}\in T_2$, $D_{\textbf{t},i-1}=\{\textbf{x}\in \Delta_{i}:(\textbf{x}_{1},\textbf{x}_{2})=\textbf{t}\}$.
From construction, we have that the projection onto the first factor $\pi_1(F_i(D_{\textbf{t},i-1}))$ is a single point, and the image on the second factor $\pi_2(F_i(D_{\textbf{t},i-1}))$ is the convex hull of $i-1$ points.
Let this hull be $F_i'$.
In \autoref{fig:octaproof}, we plot the $F_i'$ corresponding to the $F_i$ we inspected for a few values of $\textbf{t}$.
We also plot the first factor of each $F_i$, which is always a point $p\in\tau_1'$.

\subsubsection{Explicit Geodesic Motion Planner}

We will construct an explicit geodesic motion planner on five sets to show that $GC(X)\leq 4$.
We take inspiration from the work of Davis on the cube \cite{cube} and construct the sets mostly based on \textbf{multiplicity}, the number of distinct geodesics between a pair of points.
From the possible cut loci in \autoref{subsubsec:octa_cut_locus_proof}, the possible multiplicities are 1, 2, 3, 4, and 6, since no pair of points $(p,q)$ has $q$ in exactly five cells of the cut locus of $p$.\footnote{
The multiplicity of $(p,q)$ on a convex polyhedron is the number of Voronoi cells from \autoref{alg:cutlocus} that $q$ is in.
By \autoref{lem:cutlocvoronoi}, we have a geodesic by connecting the image of $q$ to the copy of $p$ associated with the cell, and \autoref{lem:geodesicdecomp} indicates this accounts for all geodesics.}

We extend the enumeration of the regions in \autoref{fig:octa_isomorphic_regions} to the full octahedron.
In total, the octahedron is partitioned into eight points $o_i$, six points $v_i$, twelve edges $e_i$, twenty-four line segments $a_i$, twenty-four line segments $c_i$, and forty-eight cells $\tau_i$.
We partition $E_1\sqcup\dots\sqcup E_5=X\times X$ in the following way:

\begin{compactitem}
    \item
    Let $E_5$ be the set of all pairs in $X\times X$ such that each point is either a vertex, a midpoint of an edge, or the center of a face.
    Since there are six vertices, twelve edges, and eight faces in the octahedron, $E_5$ is a collection of $26^2$ points in $X\times X$.
    All pairs of points with multiplicity 6 are in $E_5$.
    
    \item
    Let $M_1$ be the set of all pairs in $X\times X$ with multiplicity $1$, and let $E_1:=M_1\setminus E_5$.
    
    \item
    Let $M_4$ be the set of all pairs in $X\times X$ with multiplicity $4$, and let $E_4:=M_4\setminus E_5$.
    
    \item
    Let $M_3$ be the set of all pairs in $X\times X$ with multiplicity $3$.

    We additionally consider the following elements with multiplicity 2:
    
    For each $a_i$, the cut locus of $p\in a_i$ has exactly two vertices of degree three.
    For each vertex, two of the three incident lines have an endpoint at octahedron vertices.
    The specified vertices and lines are limits of vertices and lines in the cut locus of $p$ chosen in the neighboring $\tau_j$ or $\tau_k$.
    We will include $(p,q)$ for $p\in a_i\cup \tau_j\cup \tau_k$ and $q$ on one of these four cut locus lines.
    
    In the example of $a_1$ and $\tau_1$ in Figures \ref{fig:octaa1locus} and \ref{fig:octatau1locus}, we consider lines $\ell^{\{0,1\}}$, $\ell^{\{0,5\}}$, $\ell^{\{2,3\}}$, and $\ell^{\{3,4\}}$ on Face 0, and their extensions on Faces 1 and 3. 
    We visualize this in \autoref{fig:octasetannoyance}.

    Formally, we construct set $S$: 
    For each $a_i$ with neighboring $\tau_j$ and $\tau_k$, we include all $(p,q)$ where $p\in a_i\cup\tau_j\cup\tau_k$ and $q$ is on one of the four specified lines of the cut locus.

    We will set $E_3:=(M_3\cup S)\setminus E_5$.
    
    \item
    Let $M_2$ be the set of all pairs in $X\times X$ with multiplicity $2$.
    We will set $E_2:=M_2\setminus (E_5\cup E_3)$.
\end{compactitem}

\renewcommand{\imgwidth}{.48\linewidth}
\begin{figure}[ht!]
    \centering
    \subfigure[$p$ chosen within $\tau_1$]{
    \includegraphics[width=\imgwidth]{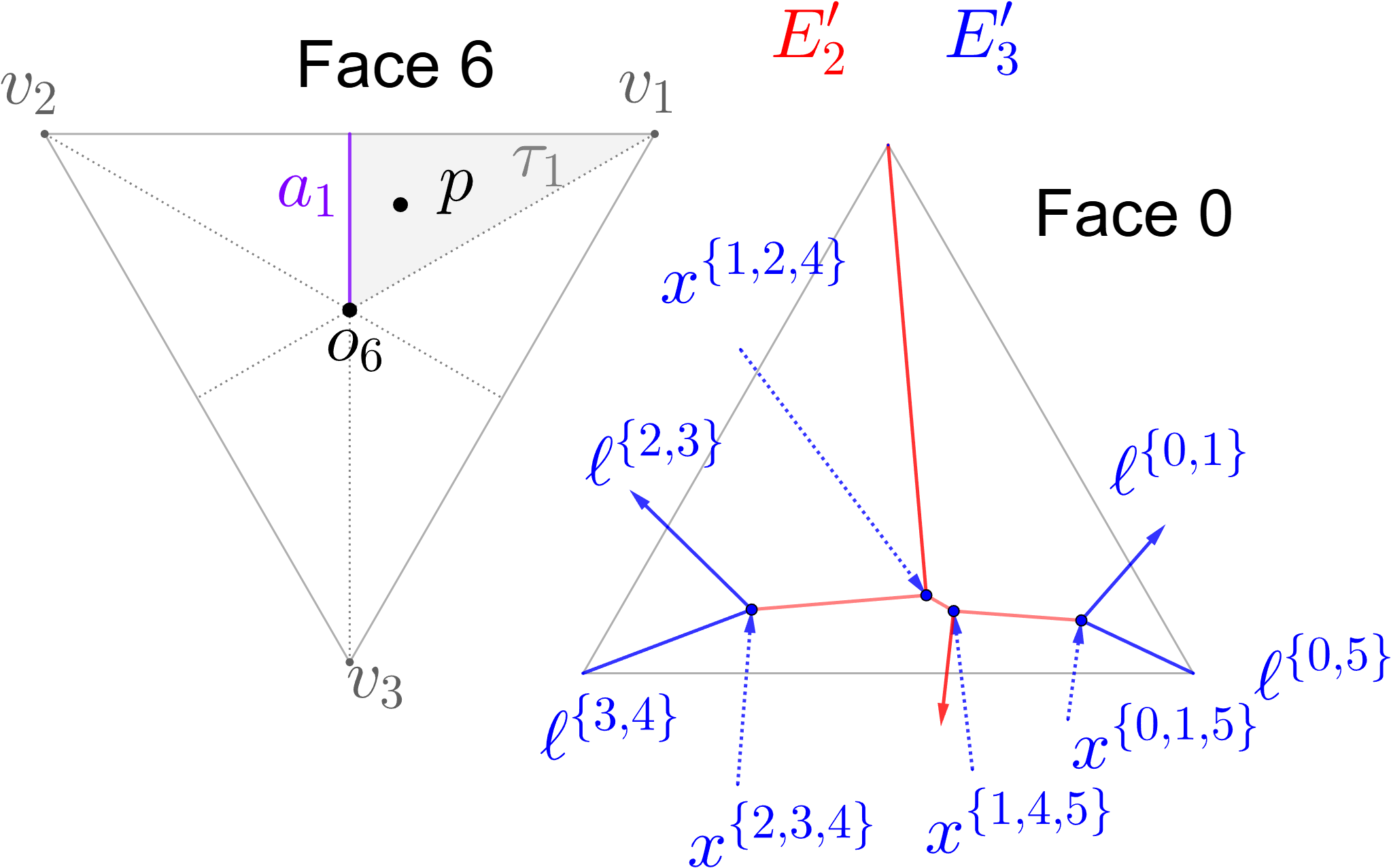}
    }
    \subfigure[$p$ chosen within $a_1$]{
    \includegraphics[width=\imgwidth]{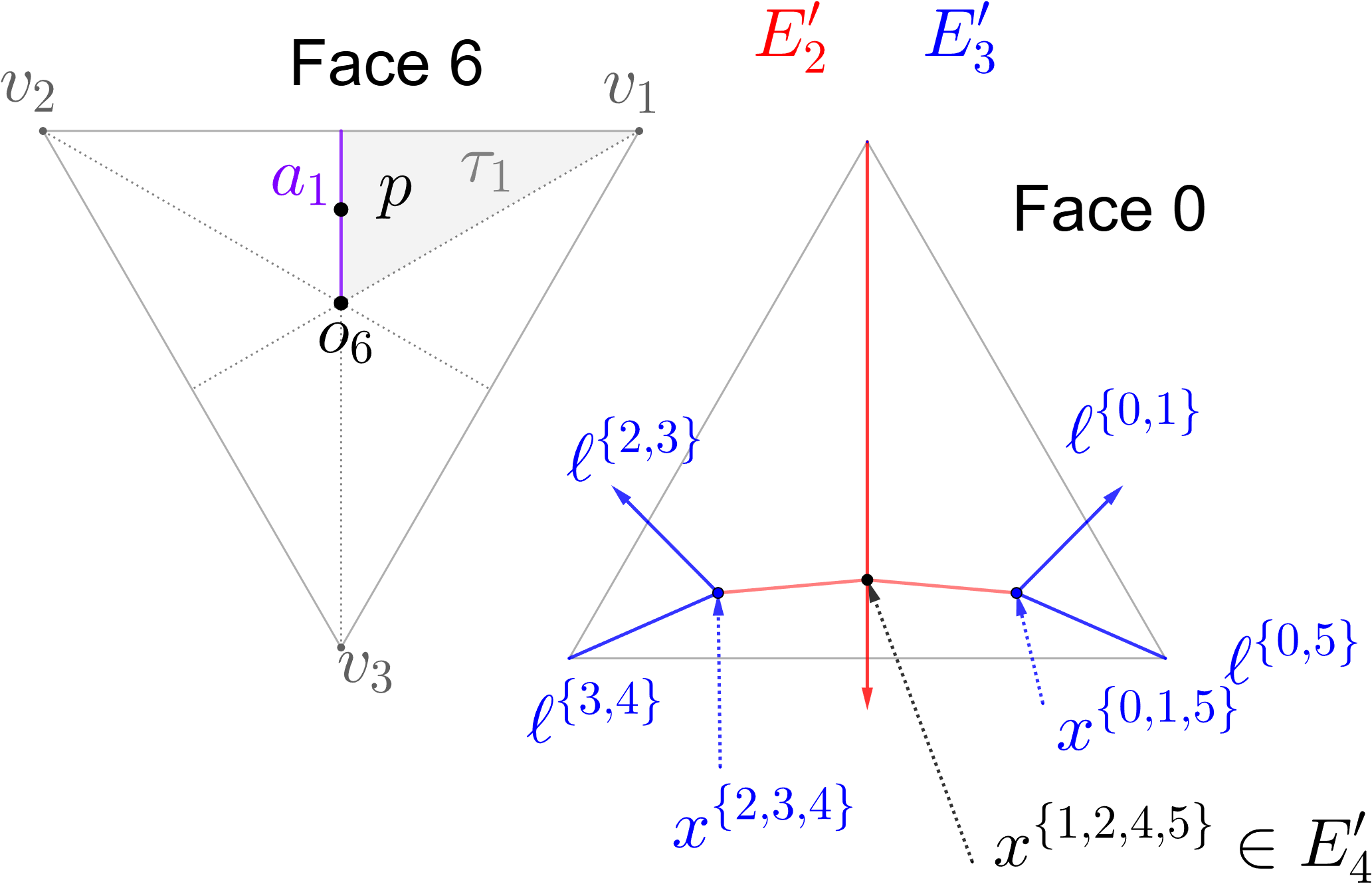}
    }
    \caption{Sets $E'_2$, $E'_3$, and $E'_4$ for $p\in\tau_1\cup a_1$}
\label{fig:octasetannoyance}
\end{figure}
We illustrate the points of multiplicity 2 in set $E_3$.
For $p\in X$, we define $E_i':=\{q\in X:(p,q)\in E_i\}$.
In \autoref{fig:octasetannoyance}, we depict the sets $E'_2$, $E'_3$, and $E'_4$ on Face 0 for $p$ chosen in $\tau_1$ and $a_1$.
In particular, points on cut locus lines $\ell^{\{0,1\}}$, $\ell^{\{0,5\}}$, $\ell^{\{2,3\}}$, and $\ell^{\{3,4\}}$ have multiplicity 2 with $p$, but have been included in set $E_3'$.

\begin{lemma}
\label{lem:point_separation}
    For $X$ compact Hausdorff and $A,B\subseteq X$, let $f\colon A\to X$ and $g\colon \overline{B}\to X$ be maps such that for all $x\in A\cap \overline{B},\;f(x)\neq g(x)$.
    Then an open set separates $S_A:=\{(a,f(a)):a\in A\}$ from $S_B:=\{(b,g(b)):b\in \overline B\}$.
\end{lemma}
Since $X$ is compact, $\overline B$ is compact, and since $X$ is Hausdorff, $X\times X$ is Hausdorff.
The map $b\mapsto (b,g(b))$ is continuous from continuity of $g$.
Since $S_B$ is the image of $\overline{B}$ under this map, it is compact and thus closed.
From our hypotheses, $S_A$ is disjoint from $S_B$.
Then $(S_B)^C$ is open and separates $S_A$ from $S_B$.
\proofclose{\autoref{lem:point_separation}}

We will now construct a \gls{GMPR} $\sigma_i\colon E_i\to GX$ for each $1\leq i\leq 5$.
To construct some of these \gls{GMPR}s, we use \autoref{lem:point_separation} to split $E_i$ into a \textbf{clopen partition}.\footnote{Where each set is clopen - both closed and open.}

\meminisection{Sets $E_1$ and $E_5$:}
There is only one choice for $\sigma_1\colon E_1\to GX$, as there is a unique geodesic between each $(p,q)\in E_1$.
Continuity follows from Lemma 3.12 of Chapter 1 in \cite{metic_non_pos}.
For $\sigma_5\colon E_5\to GX$, we arbitrarily select a geodesic for each element of $E_5$.
This is continuous since $E_5$ is discrete.

\meminisection{Set $E_4$:}
The points that have multiplicity 4 are $(p,q)$ such that $p$ is on some $a_i$, $c_i$, $v_i$, or $e_i$.
Given $p$, the cut locus is isomorphic to those in Figures \ref{fig:octaa1locus}, \ref{fig:octac1locus}, \ref{fig:octav1locus}, or \ref{fig:octae1locus} respectively, and $q$ must be a cut locus vertex with degree four.
For $p\in c_i$, there are two choices for $q$, as both cut locus vertices have degree four.
For $p$ chosen on some $e_i$, $v_i$, or $a_i$, there is only one choice for $q$.
We may ignore $(p,q)$ with multiplicity 4 where $p$ is some $v_i$ or on the midpoint of some $e_i$, as in this case $(p,q)\in E_5$.
Thus, we enumerate line segments $e_j'$, where removing the midpoint from each $e_i$ creates two $e_j'$.
Then for $(p,q)\in E_4$, $p$ is on some $a_i$, $c_i$, or $e_j'$, and these line segments are pairwise separated from each other.

\meminisubsection{Clopen Partition:}
To create a \gls{GMPR}, we partition $E_4$ into the following types of sets:
\begin{compactitem}
    \item
    Let $L$ be either some $a_i$ or some $e_i'$, and define $f\colon\overline{L}\to X$ in the following way: 
    For any point $p\in L$ there is a unique vertex of the cut locus with degree four.
    Let $f(p)$ be this point, and note that $f(p)$ varies continuously for $p\in L$.
    An endpoint of $L$ is either on some edge $e_i$, some $v_i$ or some $o_i$.
    As $p$ approaches an endpoint, the limit $f(p)$ is well defined.
    In all cases, the limit will be the (only) vertex of the cut locus of $p$.
    We will define $f(p)$ as this limit.

    We construct a set $\{(p,f(p)):p\in L\}$.
    We additionally note that $f$ is well defined and continuous on $\overline L$.
    There are forty-eight sets of this type, arising from the twenty-four choices of $a_i$ and twelve choices of $e_j$ (each edge $e_j$ creates two $e_i'$).
    \item
    
    Let $L$ be some $c_i$, and define $f\colon\overline{L}\to X$ in the following way:
    For $p\in L$, there are two vertices of the cut locus with degree four, and both vary continuously with $p$.
    Pick one and let this be $f(p)$ for $p\in L$.
    An endpoint of $L$ is either some $v_i$ or some $o_i$.
    As $p$ approaches an endpoint, the limit $f(p)$ is well defined as the (only) vertex of the cut locus of $p$.
    We will define $f(p)$ as this limit for $p$ an endpoint of $L$.
    
    We construct sets $\{(p,f(p)):p\in L\}$.
    We additionally note that $f$ is well defined and continuous on $\overline L$.
    There are forty-eight sets of this type, as there are two choices of $f$ for each of the twenty-four choices of $c_i$.
\end{compactitem}
Each pair of sets satisfies the hypotheses for \autoref{lem:point_separation}.
For most sets this is trivial, as for $L_1,L_2$ each some $a_i$, $c_i$, $v_i$, or $e_j'$, if $L_1\neq L_2$ then $L_1\cap \overline{L_2}=\emptyset$.
The only case we must consider is if $L_1=L_2$, which can only yield two different sets for $L_1=c_i$.
This case satisfies the hypotheses since $c_i\cap \overline{c_i}=c_i$ and for any $p\in c_i$, the two vertices of the cut locus of $p$ are distinct.
Thus, this collection of sets is a clopen partition of $E_4$.
To define a \gls{GMPR} on all of $E_4$, it is enough to define a \gls{GMPR} on each of these ninety-six sets.

\meminisubsection{\acrlong{GMPR}:}
Consider any of these sets $S$.
We claim that there are exactly four \gls{GMPR}s on $S$.
For $(p,q)\in S$, $q$ is adjacent to exactly four cells in the cut locus of $p$ by construction.
These cells vary continuously wrt. $p$, as they arise from a Voronoi diagram generated by copies of $p$.
The four \gls{GMPR}s are defined by which copy of $p$ we choose to create a path to.

Consider the case of $a_1$, displayed in \autoref{fig:octaa1locus}.
In this case, $q$ is the point equidistant from copies $p^{(1)}$, $p^{(2)}$, $p^{(4)}$, and $p^{(5)}$.
We will define a \gls{GMPR} by choosing one of these copies of $p$ arbitrarily, and letting $\gamma\colon[0,1]\to \mathbb R^2$ be the path from this copy to the image of $q$.
With $f$ the path unfolding of the chosen copy of $p$, $f^{-1}\circ \gamma$ defines a path from $p$ to $q$ in the octahedron.\footnote{While $f$ itself does not necessarily have an inverse, we may decompose $\gamma$ so that each segment is only on one face of the path unfolding. 
Since $f$ restricted to each face is invertible, we may use it to pull back each segment of $\gamma$ into a path on $X$.
}
Continuity of this path follows from continuity of the locations of the copies of $p$.
For each of the sets composing $E_4$, we arbitrarily choose one of the four \gls{GMPR}s.
Our \gls{GMPR} $\sigma_4$ is the piecewise function with our ninety-six chosen \gls{GMPR}s as components.

\meminisection{Set $E_3$:}
The points that have multiplicity 3 are $(p,q)$ such that $p$ is either on some $\tau_i$ or some $a_i$. 
The cut locus of $p$ will be isomorphic to those in Figures \ref{fig:octatau1locus} and \ref{fig:octaa1locus} respectively.
If $p$ is on some $\tau_i$, $q$ must be one of the four vertices of the resulting cut locus.
If $p$ is on some $a_i$, $q$ must be one of the two vertices of the resulting cut locus with degree three.
We additionally must consider the points of multiplicity 2 that we added to $E_3$.
These points are $(p,q)$ such that $p$ is on some connected $a_i\cup\tau_j\cup \tau_k$, and $q$ is on one of the four specified cut locus lines.

\meminisubsection{Clopen Partition:}
We partition $E_3$ into the following types of sets:
\begin{compactitem}
\item
    Consider some $a_i$ and its adjacent sets $\tau_j,\tau_k$, and define $A:=a_i\cup \tau_j\cup \tau_k$.
    For $p\in A$, there are two vertices of the cut locus that always have degree three.
    These vertices vary continuously with $p$, as they always arise from the same copies of $p$ throughout $A$.
    Two of the cut locus lines connected to each vertex are also included in set $E_3$.
    In the case of $a_1$ and $\tau_1$ (displayed in \autoref{fig:octasetannoyance}), we inspect vertex $x^{\{0,1,5\}}$ with lines $\ell^{\{0,1\}}$ and $\ell^{\{0,5\}}$, and vertex $x^{\{2,3,4\}}$ with lines $\ell^{\{2,3\}}$ and $\ell^{\{3,4\}}$.
        
    We pick one of these two and construct a set $S$ of all $(p,q)$ with $p\in A$ and $q$ either the vertex itself, or on one of the two specified cut locus lines connected to the vertex.
        There are forty-eight of these sets, as we construct two for each of the twenty-four $a_i$.
    \item
    Consider some $\tau_i$, and note that for $p$ chosen on $\tau_i$, there are a total of four vertices in the cut locus of $p$.
    We have already included two of these vertices in sets of the first type, so we will consider the other two.
    In the case of \autoref{fig:octasetannoyance}(a), we consider the vertices $x^{\{1,2,4\}}$ and $x^{\{1,4,5\}}$.
    
    As before, we will choose one vertex and define map $f\colon \overline{\tau_i}\to X$.
    Within $\tau_i$, $f$ will be the specified vertex.
    As $p$ approaches the boundary of $\tau_i$, each vertex has a well defined limit.
    As $p$ approaches some $c_i$, $f(p)$ approaches a vertex of the resulting cut locus.
    In this case, which vertex $f(p)$ approaches depends on our initial choice of cut locus vertex to follow.
    In all other cases, as $p$ approaches the boundary, $f(p)$ approaches the unique vertex of the cut locus with maximal degree.
    We construct set $\{(p,f(p)):p\in \tau_i\}$ and note that $f$ is well defined and continuous on $\overline{\tau_i}$.
    There are ninety-six of these sets, as we construct two for each of the forty-eight $\tau_i$.
\end{compactitem}

We will show that these sets form a clopen partition of $E_3$.
Let $S_1$ and $S_2$ be constructed sets.
In most cases, $S_1$ and $S_2$ are easily separated on their first factors, so we will only consider when this is not true.
We will also use the observation that in a metric space, if $S_1$ is not separable from $S_2$ by an open set, then $S_1\cap \overline{S_2}$ is nonempty, and there is a sequence $(p_n,q_n)_n\subseteq S_2$ such that $\lim\limits_{n\to\infty}(p_n,q_n)$ exists and is in $S_1$.

\begin{compactitem}
    \item
    Assume $S_1$ and $S_2$ are both sets of the second type.
    In this case, they satisfy the hypotheses for \autoref{lem:point_separation}, and thus are separable by an open set.
    \item
    Assume $S_1$ and $S_2$ are both sets of the first type.
    We may assume they arise from the same $a_i$ and adjacent $\tau_j$, $\tau_k$.
    
    \AFSOC that $S_1$ is not separable from $S_2$.
    Then find sequence $(p_n,q_n)_n\subseteq S_2$ whose limit is some $(p,q)\in S_1$.    
    Then $p$ and all $p_n$ must be in $a_i\cup \tau_j\cup\tau_k$.
    However, the cut locus of a point $p'$ chosen in this region arises from a Voronoi diagram of points that vary continuously with $p'$.
    In particular, since each $q_n$ is on particular lines of the cut locus of $p_n$, $q$ must be on the closure of those same cut locus lines with respect to $p$.
    However, the cut locus lines that $q_n$ can be on must be different from the cut locus lines that $q$ can be on, from construction of $S_1$ and $S_2$.
    The closures of these lines never intersect, which gives us the contradiction that $\lim\limits_{n\to\infty}(p_n,q_n)_n\neq (p,q)$.
    In an example using $\tau_1$ and $a_1$, assume $S_1$ is created with vertex $x^{\{2,3,4\}}$ and lines $\ell^{\{2,3\}}$ and $\ell^{\{3,4\}}$, and $S_2$ is created with vertex $x^{\{0,1,5\}}$ and lines $\ell^{\{0,1\}}$ and $\ell^{\{0,5\}}$ (see \autoref{fig:octasetannoyance}).
    We observe from our cut locus calculations that for no cut locus of a point $p'\in a_1\cup \tau_1\cup \tau_2$ does the set $\overline{\ell^{\{0,1\}}(p')\cup x^{\{0,1,5\}}(p')\cup\ell^{\{0,5\}}(p')}$ intersect $\ell^{\{2,3\}}(p')\cup x^{\{2,3,4\}}(p')\cup\ell^{\{3,4\}}(p')$.
    
    \item
    Assume $S_1$ is of the first type and $S_2$ is of the second type.
    We may assume $S_1$ arises from $a_i$ and adjacent $\tau_j$, $\tau_k$, and $S_2$ arises from $\tau_j$.
    
    First, take $(p_n,q_n)_n\subseteq S_1$ such that the limit $(p,q)\in S_2$.
    This implies $p\in \tau_j$.
    Similar to the case where both sets were of the first type, we reach a contradiction since the closure of the cut locus lines that each $q_n$ may be on do not intersect the cut locus vertices that $q$ may be on. 
    In an example using $\tau_1$ and $a_1$, for no cut locus of a point $p'\in\tau_1$ does the set $\overline{\ell^{\{0,1\}}(p')\cup x^{\{0,1,5\}}(p')\cup\ell^{\{0,5\}}(p')}$ intersect $x^{\{1,4,5\}}(p')$ or $x^{\{1,2,4\}}(p')$.

    Now take $(p_n,q_n)_n\subseteq S_2$ such that the limit $(p,q)\in S_1$. 
    From construction, $p_n\in \tau_j$ and $p\in a_i\cup\tau_j\cup\tau_k$.
    Since $p$ must be the limit of the $p_n$, $p\in \overline{\tau_j}\cap (a_i\cup\tau_j\cup\tau_k)$, which implies either $p\in\tau_j$ or $p\in a_i$.
    By construction, $q_n$ is a particular vertex of the cut locus of $p_n$. 
    Since this vertex varies continuously for points on $\tau_j\cup a_i$, $q$ must be the limit of this vertex as $p_n$ approaches $p$.
    This results in a similar contradiction, as these vertices do not lie on any of the cut locus lines that $q$ must be on due to $(p,q)\in S_1$.
    In an example using $\tau_1$ and $a_1$, for no cut locus of a point $p'\in\tau_1\cup a_1$ do vertices $x^{\{1,4,5\}}(p')$, $x^{\{1,2,4\}}(p')$, or $x^{\{1,2,4,5\}}(p')$ intersect $\ell^{\{0,1\}}(p')\cup x^{\{0,1,5\}}(p')\cup\ell^{\{0,5\}}(p')$.
    
\end{compactitem}

Thus, these 144 sets form a clopen partition of $E_3$, and to define a \gls{GMPR} over $E_3$, it is enough to define a \gls{GMPR} for each set.

\meminisubsection{\acrlong{GMPR}:}
Each set of the second type contains $(p,q)$, where $q$ is adjacent to three Voronoi cells of the cut locus of $p$. 
Similar to the case for $E_4$, the set has three possible \gls{GMPR}s, each corresponding to a choice of cell.
We choose one arbitrarily for each of these 96 sets.
Each set $S$ of the first type contains $(p,q)$ for $p\in a_i\cup\tau_j\cup\tau_k$ for neighboring $a_i,\tau_j,\tau_k$, and $q$ chosen as one of the resulting cut locus vertices, or on specified cut locus lines connected to the vertex.
From construction, there is a single Voronoi cell of the cut locus of $p$ that every choice of $q$ is adjacent to. 
We will choose our \gls{GMPR} to always approach the copy of $p$ that creates this Voronoi cell.
Our overall \gls{GMPR} $\sigma_3$ is the piecewise function with the \gls{GMPR} for each set as components.
The continuity of each constituent \gls{GMPR} is evident from the same argument as $\sigma_4$, where we first construct a \gls{GMPR} within the appropriate path unfolding on $\mathbb R^2$, then use $f^{-1}$ to pull it into the octahedron.

\meminisection{Set $E_2$:}
The points that have multiplicity 2 are a choice of $p\in X$ and $q$ on a line of the cut locus of $p$.

\meminisubsection{Partition:}
We will partition $E_2$ into sets of two types.
\begin{compactitem}

    \item
    Consider some $c_i$, and its adjacent sets $\tau_j$, $\tau_k$, and define $A:=c_i\cup \tau_j\cup \tau_k$.
    For $p\in c_i$, there is one cut locus line that does not touch any edge of the octahedron (for example, $\ell^{\{1,4\}}$ in \autoref{fig:octac1locus}).
    When $p\in\tau_i$ or $p\in\tau_j$, there is one cut locus line that approaches our specified cut locus line as $p$ approaches $c_i$. 
    This line is always the middle of three cut locus lines that do not touch any edge of the octahedron, and in the example of \autoref{fig:octatau1locus}, we would consider $\ell^{\{1,4\}}$.
    We construct a set $S$ of all $(p,q)$ with $p\in A$ and $q$ on the relative interior of this cut locus line. 
    
    There are twenty-four of these sets, corresponding to one for each $c_i$.
    \item
    The remaining set is all other points in $E_2$.
\end{compactitem}

\meminisubsection{\acrlong{GMPR}:}
Every point in $E_2$ is $(p,q)$ with $q$ on the relative interior of a line of the cut locus of $p$.
Thus, to describe a \gls{GMPR} $\sigma_2$, it is enough to describe which side of the line $\sigma_2$ approaches $q$ from.
This is equivalent to constructing $\sigma_2$ to approach one of the two copies of $p$ that generate the cut locus line.
For each set of the first type, we arbitrarily choose which of the two copies of $p$ the \gls{GMPR} approaches from.
For the set $S$ of the second type, we will use the fact that the octahedron is orientable. 
The allows us to have a consistent notion of \lrq{clockwise} for each point on the octahedron, which we use to choose a side of each line.

\renewcommand{\imgwidth}{.48\linewidth}
\begin{figure}[ht!]
    \centering
    \subfigure[$p$ chosen within $e_1$]{
    \includegraphics[width=\imgwidth]{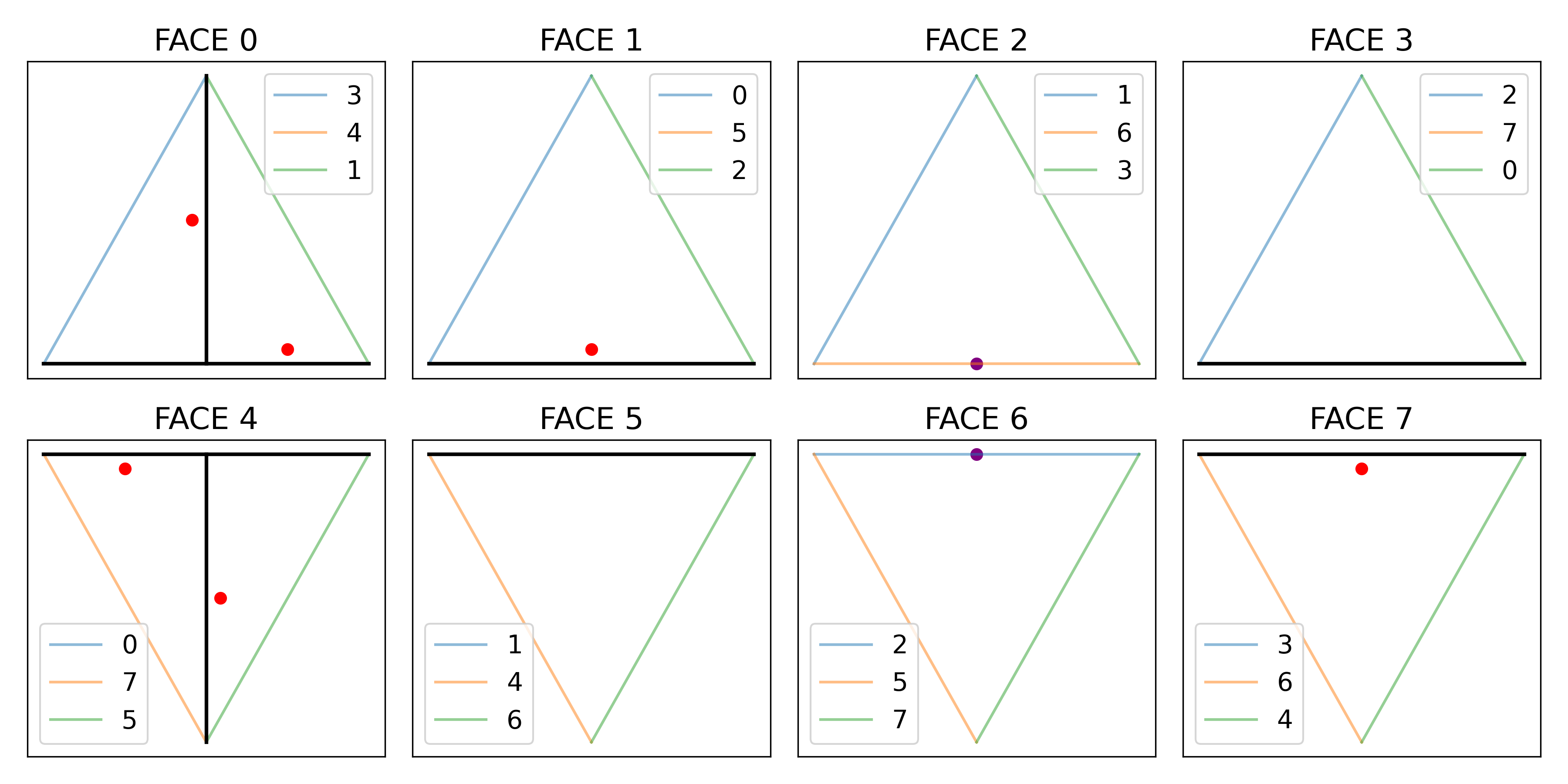}
    }
    \subfigure[$p$ chosen within $c_1$]{
    \includegraphics[width=\imgwidth]{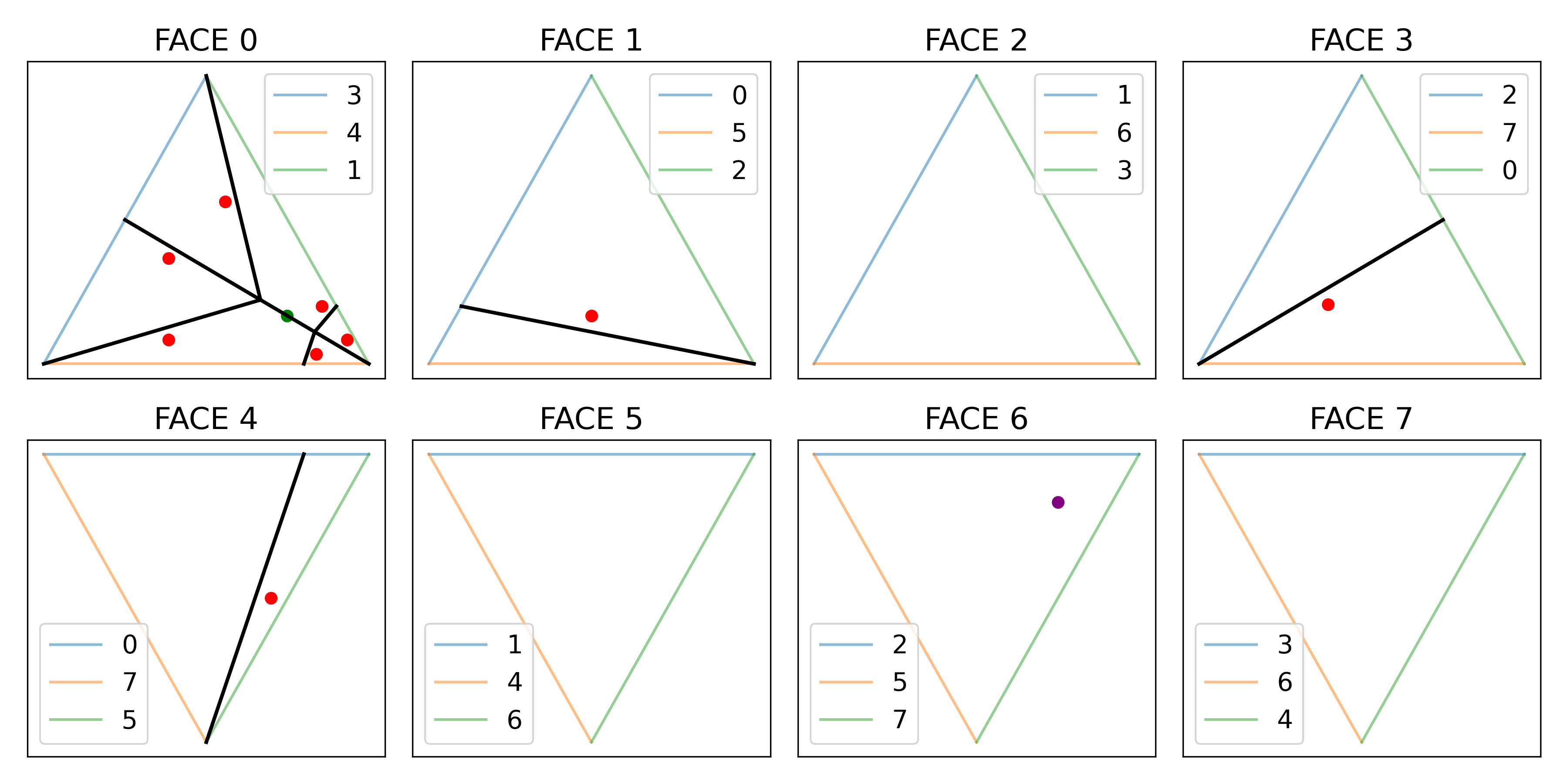}
    }
    \subfigure[$p$ chosen as $o_6$]{
    \includegraphics[width=\imgwidth]{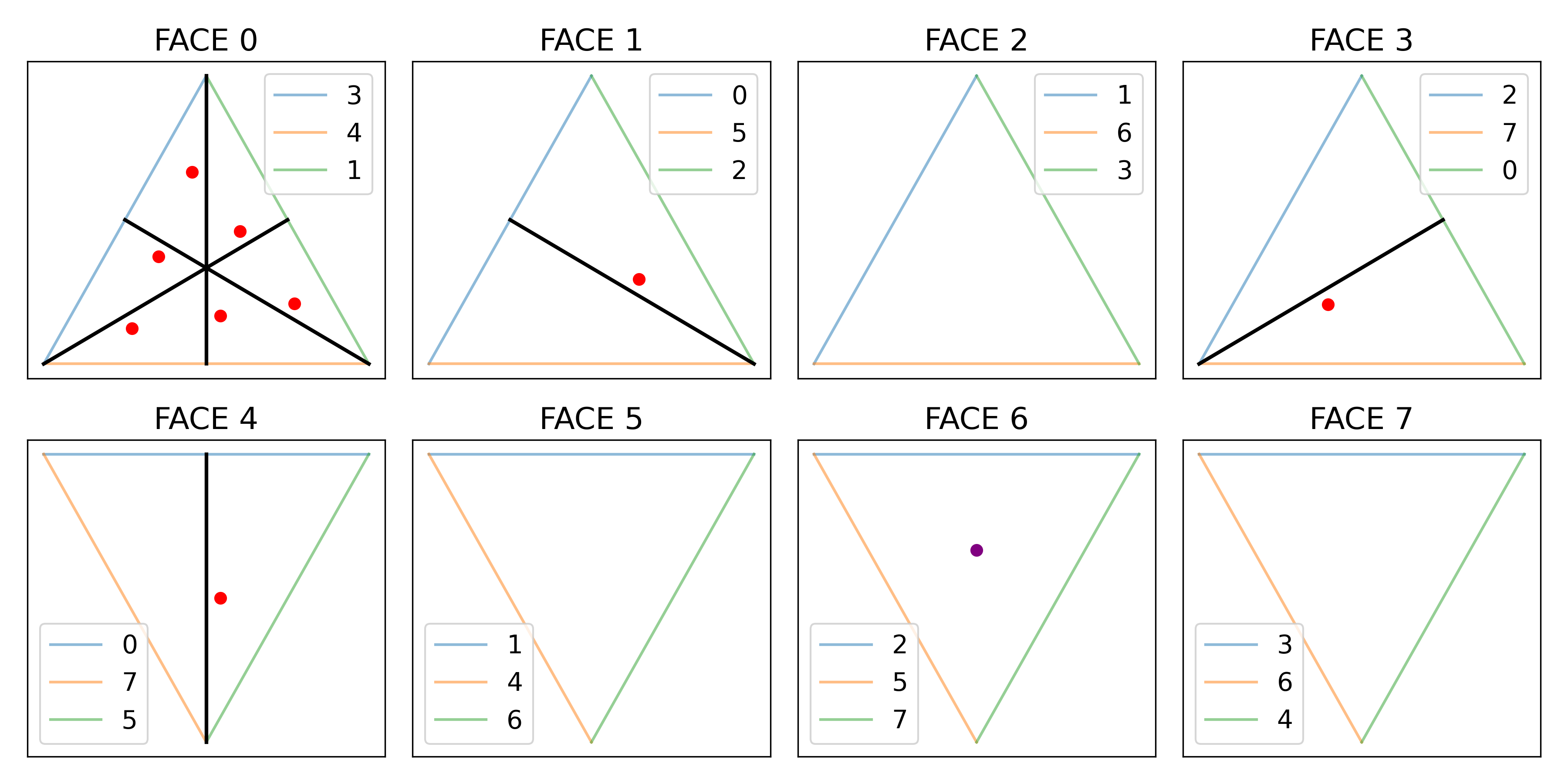}
    }
    \subfigure[$p$ chosen as $v_1$]{
    \includegraphics[width=\imgwidth]{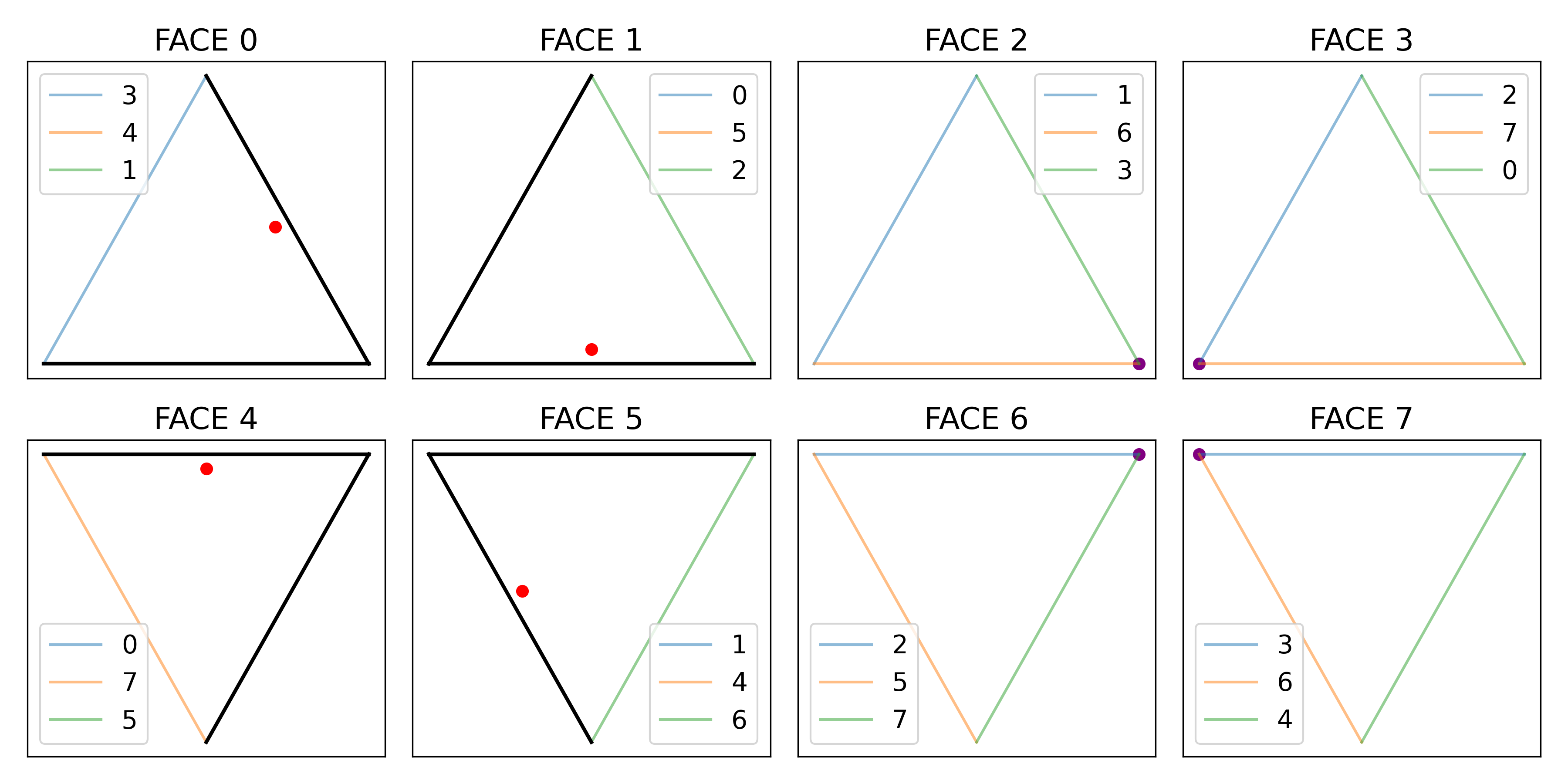}
    }
    \subfigure[$p$ chosen within $a_1$]{
    \includegraphics[width=\imgwidth]{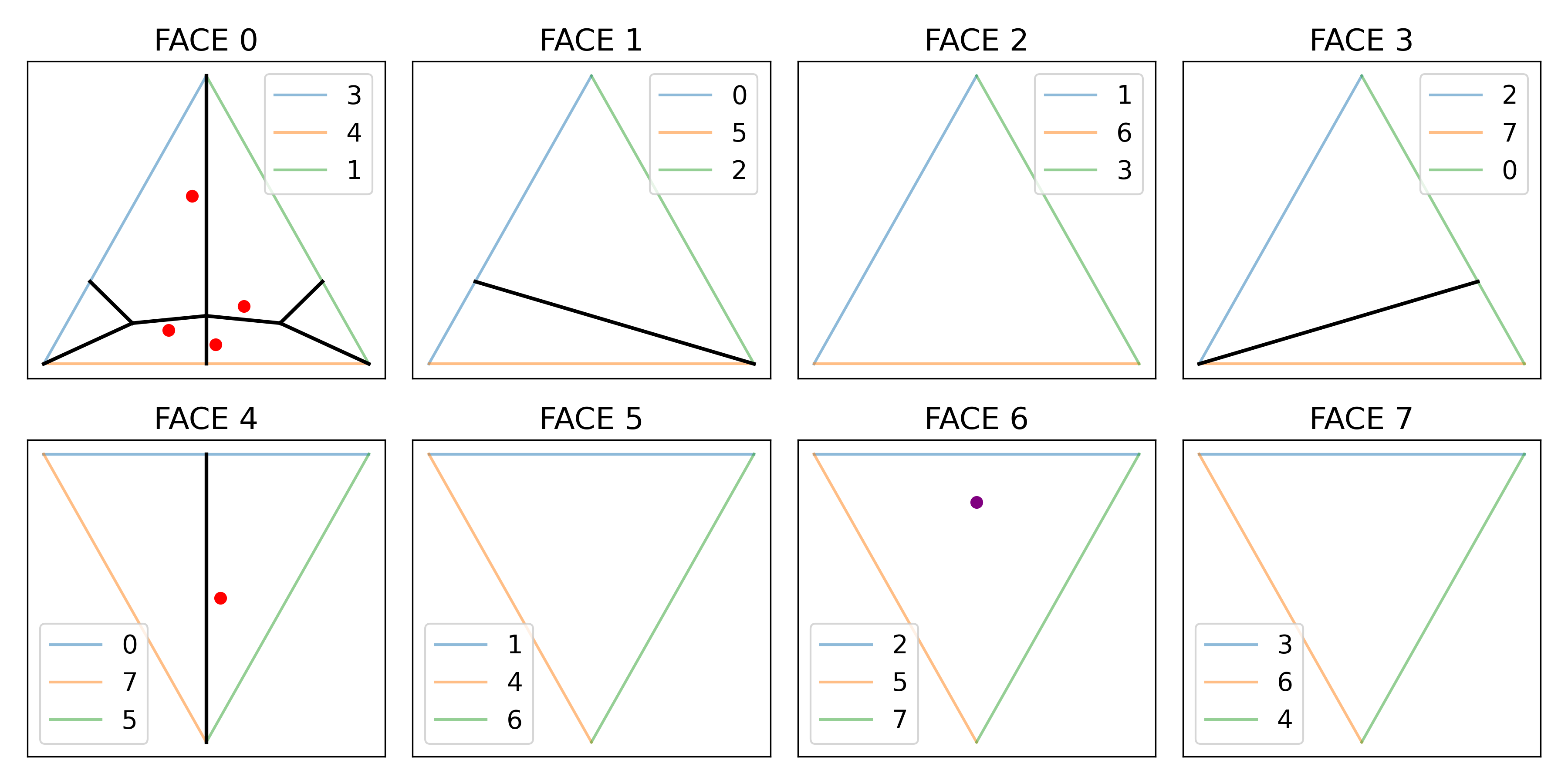}
    }
    \subfigure[$p$ chosen as $\tau_1$]{
    \includegraphics[width=\imgwidth]{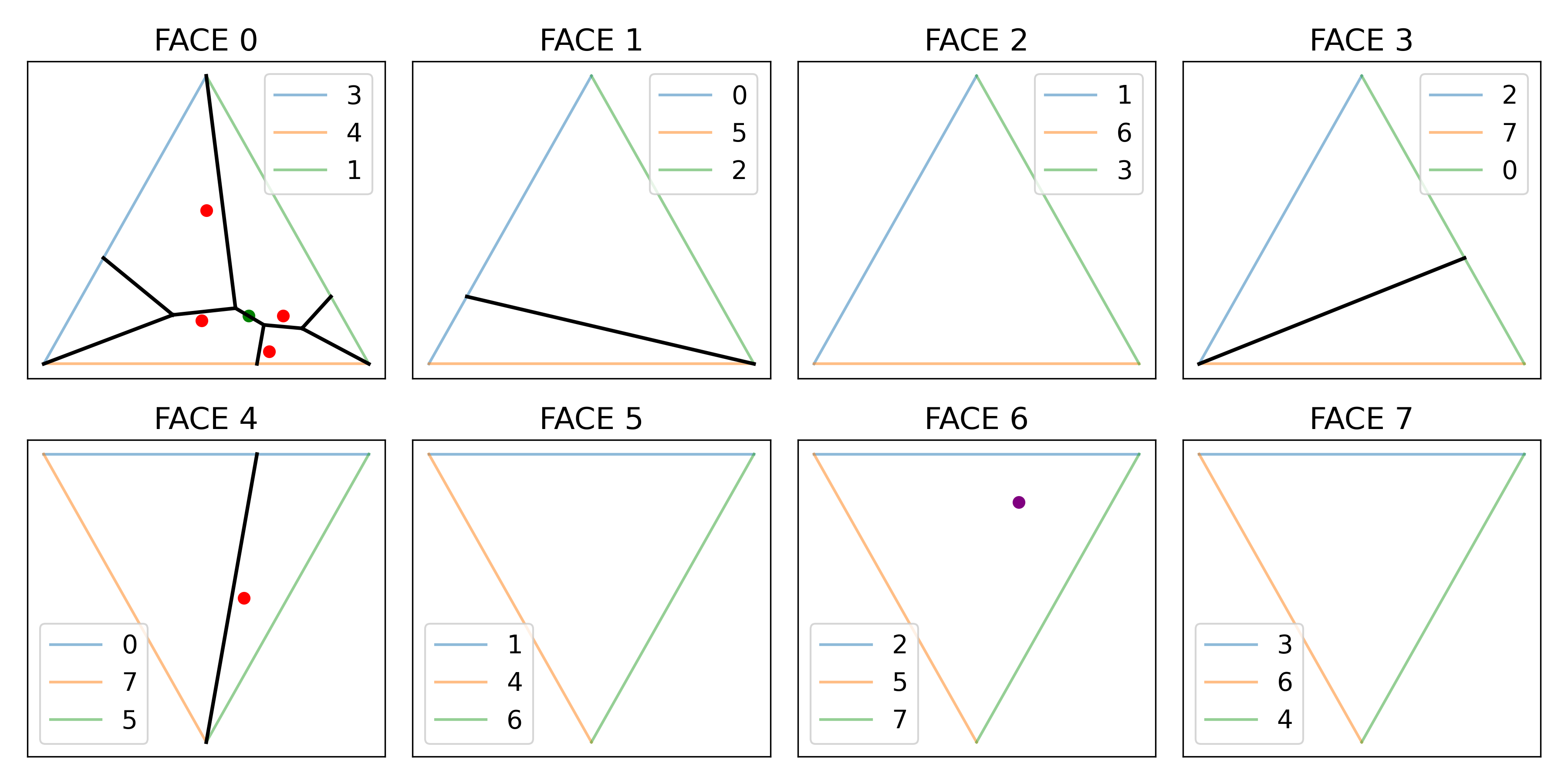}
    }
    \caption{Geodesics chosen by our \gls{GMPR} for a few start points (displayed as purple dots).
    Red dots indicate the direction from which the \gls{GMPR} on $E_2$ approaches each cut locus line.
    Green dots indicate an arbitrary choice.
    }
\label{fig:octa_gmpr_E2_examples}
\end{figure}
\begin{compactitem}
    \item
    For $p$ chosen on some $e_i$, $c_i$, or $o_i$ (e.g. \autoref{fig:octa_gmpr_E2_examples}(a, b, c)), all lines of the cut locus that we must account for have a unique endpoint at a vertex of the octahedron.
    With respect to this endpoint, we have a notion of the clockwise side of each line.
    We construct $\sigma_2$ to approach each point $q$ on a cut locus line from the clockwise side.
    \item 
    For $p$ chosen on some $v_i$ (e.g. \autoref{fig:octa_gmpr_E2_examples}(d)), the lines of the resulting cut locus have both of their endpoints at vertices of the octahedron.
    We can distinguish one vertex as the antipode of $p$, and we will use the vertex that is \textit{not} the antipode to determine the clockwise side of each line.
    We construct $\sigma_2$ to approach each point on the cut locus line from the clockwise side.
    \item
    For $p$ in some $A:=a_i\cup \tau_j\cup \tau_k$ (with $\tau_i$ and $\tau_k$ adjacent to $a_i$, and edge $e_\ell$ being the edge adjacent to $A$), a few cut locus lines have an endpoint on vertices of the octahedron.
    We construct $\sigma_2$ to approach each point on these lines from the clockwise side (with respect to the octahedron vertex endpoint).
    For Figures \ref{fig:octaa1locus} and \ref{fig:octatau1locus}, these lines are $\ell^{\{1,2\}}$ and $\ell^{\{4,5\}}$.
        There remain two lines not considered in other cases (one line is considered in the sets of the first kind, and four lines are considered in set $E_3$).
    As $p$ approaches the (unique) adjacent edge $e_\ell$, the limit of each of these lines is an octahedron edge.
    As before, we have a notion of the clockwise side of each cut locus line for $p\in e_\ell$, which allows us to define a clockwise side for the lines that approach these edges for $p\in A$.
    We construct $\sigma_2$ to approach each point on these cut locus lines from the clockwise side.
    In an example with $\tau_1$, $a_1$, and $e_1$, lines $\ell^{\{1,5\}}$ and $\ell^{\{2,4\}}$ in Figures \ref{fig:octaa1locus} and \ref{fig:octatau1locus} approach the lines of the same name in \autoref{fig:octae1unfold}(a).
    We display the chosen directions for $a_1$ and $\tau_1$ in \autoref{fig:octa_gmpr_E2_examples}(e, f)
\end{compactitem}

\meminisubsection{Continuity:}
We will show the continuity of $\sigma_2$.
Since we construct $\sigma_2$ similar to the \gls{GMPR} of $E_2$ in \cite{cube}, we will prove continuity in a similar way.
First, let $A$ be a region of the octahedron with isomorphic cut locus.
For a point on $A$, the cut locus vertices and line segments vary continuously.
For the relative interior of a line segment, $\sigma_2$ is constructed to always approach from a consistent side, and thus always projects to a linear path in the same continuously varying Voronoi cell.
As before, this implies that the \gls{GMPR} is continuous, as we can use the path embedding associated with this Voronoi cell to pull each path in $\mathbb R^2$ back into $X$.
Thus, $\sigma_2$ must be continuous with its first factor restricted to any region $A$ of isomorphic cut locus.

\AFSOC that there exists a discontinuity at point $(p,q)$.
We claim it must be witnessed by two distinct regions of isomorphic cut locus $A_1$ and $A_2$, and sequence $(p_n,q_n)_n$ that converges to $(p,q)$ such that $(p_n)_n\subseteq A_1$ and $\lim\limits_{n\to\infty}p_n\in A_2$.
It is certainly true that the discontinuity is witnessed by a sequence, as both $X\times X$ and $GX$ are metric spaces.
Then we may produce a sequence $(p_n,q_n)_n\subseteq E_2$ that approaches $(p,q)\in E_2$ such that each $\sigma_2(p_n,q_n)$ is $\varepsilon$-far from $\sigma_2(p,q)$ for fixed $\varepsilon>0$.
Since we partitioned the octahedron into finitely many regions of isomorphic cut locus, there is a region $A_1$ of isomorphic cut locus such that an infinite subsequence $(p_n)_{n\in M}\subseteq A_1$.
\WLOG we may assume this is the whole sequence.
Finally, $p\notin A_1$, since $\sigma_2$ is continuous when $p$ is restricted to $A_1$.
We may additionally assume that for each $(p_n,q_n)_n$, $q_n$ is always on the same cut locus line of $p_n$.
We may do this since there are finitely many of these for $p_n$ chosen in $A_1$, and $q_n$ must be on the relative interior of one of them.
This indicates (in rough terms) that any discontinuity arises from a chosen direction for a line in the cut locus of points in $A_1$ not matching the chosen direction of the line's limit in $A_2$.
We will show that this cannot be true, as our construction ensures consistency of the choice of side for each line.
Since the limit $p\in A_2$, we need only check regions $A_1, A_2$ such that $\overline{A_1}\cap A_2\neq \emptyset$.

Assume $(p_n)_n\subseteq a_i$.
There are four cut locus lines we must consider, and the limit $p$ is either some $o_j$ or the midpoint of an edge $e_j$.
Since other $a_i$ are symmetric, we assume $a_i=a_1$, $e_j=e_1$, and $o_j=o_6$.
This case is visualized in Figures \ref{fig:octaa1locus} and \ref{fig:octa_gmpr_E2_examples}(e).
\begin{compactitem}
    \item
    We find from \autoref{alg:cutlocus} that if $(p_n)_n$ approaches $o_6$, cut locus lines $\ell^{\{1,5\}}$ and $\ell^{\{2,4\}}$ both approach $o_0$, the center of Face 0.
    Thus, if the $(q_n)_n$ are on either of these lines, then $q=o_0$.
    In this case, $(p,q)=(o_6,o_0)\in E_5$, and $(p_n,q_n)_n$ cannot witness a discontinuity in $\sigma_2$.
    
    For cut locus lines $\ell^{\{1,2\}}$ and $\ell^{\{4,5\}}$, we find that their limits are the lines of the same name as $(p_n)_n$ approaches $o_6$.
    Thus, if $(q_n)$ is on any of these lines, $q$ must be on the closure of the limiting cut locus line.
    Additionally, as $(p_n)_n$ approaches $o_6$, each line's endpoint on an octahedron vertex is fixed.
    From construction of $\sigma_2$, this implies that the clockwise side of these cut locus lines is consistent with the clockwise side of their limits.
    Then if $q$ is on the interior of a cut locus line, $(\sigma_2(p_n,q_n))_n$ must approach $\sigma_2(p,q)$, as the geodesics chosen are on a consistent side of a continuously varying cut locus line.
    If $q$ is instead on an endpoint, $(p,q)\in E_5$, and $(p_n,q_n)_n$ cannot witness a discontinuity in $\sigma_2$.
    
    Thus, $(p_n,q_n)_n$ cannot witness a discontinuity if $p=o_6$.
    In the analysis of other cases, we will not mention the case where $q$ is the endpoint of a cut locus line, as this always results in $(p,q)\notin E_2$.
    \item
    Now assume $(p_n)_n$ approaches the midpoint of $e_1$.
    We find cut locus lines $\ell^{\{1,2\}}$ and $\ell^{\{4,5\}}$ approach lines of the same name (in \autoref{fig:octae1unfold}) as $(p_n)_n$ approaches $p$, and their respective endpoints on octahedron vertices are fixed.
    Then as before, the clockwise side of these lines are consistent and $(p_n,q_n)_n$ cannot witness a discontinuity if $(q_n)_n$ is chosen on these lines.
    
    The cut locus lines $\ell^{\{1,5\}}$ and $\ell^{\{2,4\}}$ approach the lines of the same name as $(p_n)_n$ approaches $e_1$. 
    Our construction of $\sigma_2$ ensures that our chosen side for these cut locus lines is consistent with their limits at $p\in e_1$.
    Then, $(p_n,q_n)_n$ cannot witness a discontinuity if $(q_n)_n$ is chosen on these lines.
    \end{compactitem}
Therefore, there cannot be a discontinuity witnessed by $(p_n)_n\subseteq a_i$.

Assume $(p_n)_n\subseteq c_i$.
There are seven cut locus lines we must consider, and the limit $p$ is some $v_j$ or $o_j$. 
Since other $c_i$ are symmetric, we assume $c_i=c_1$, $v_j=v_1$, and $o_j=o_6$.
This case is visualized in Figures \ref{fig:octac1locus} and \ref{fig:octa_gmpr_E2_examples}(b).
\begin{compactitem}
    \item
    If $p=o_6$, we find that cut locus line $\ell^{\{1,4\}}$ approaches $o_0$, the center point of Face 0.
    If $(q_n)_n$ is chosen on this line, $(p,q)=(o_6,o_0)\in E_5$, and $(p_n,q_n)_n$ cannot witness a discontinuity.
    
    The six other cut locus lines approach lines of the same name as $(p_n)_n$ approaches $o_6$.
    Their respective endpoints on octahedron vertices are fixed, so the clockwise side of each cut locus line is consistent with the clockwise side of their limits.
    Then, as before, $(p_n,q_n)_n$ cannot witness a discontinuity if $(q_n)_n$ is chosen on any of these lines.
    \item
    If $p=v_1$, we find as $(p_n)_n$ approaches $v_1$, cut locus lines $\ell^{\{1,4\}}$ and $\ell^{\{0,5\}}$ both approach the antipode of $v_1$ (let this be $v_4$).
    Then if $(q_n)_n$ is chosen on these lines, $q=v_4$ and $(p,q)=(v_1,v_4)\in E_5$.
    
    The cut locus lines $\ell^{\{0,1\}}$, $\ell^{\{1,2\}}$, $\ell^{\{3,4\}}$, and $\ell^{\{4,5\}}$ approach the four octahedron edges that are incident to $v_4$ as $(p_n)_n$ approaches $v_1$.
    Their respective endpoints on octahedron vertices are fixed, and their other endpoint becomes $v_4$.
    Since we ignore $v_4$ when determining orientation on these cut locus lines, the clockwise side of each of these cut locus lines is consistent with the clockwise side of their limit.
    Thus, $(p_n,q_n)_n$ cannot witness a discontinuity if $(q_n)_n$ is chosen on any of these lines.

    Finally, as $(p_n)_n$ approaches $v_1$, the cut locus line $\ell^{\{2,3\}}$ approaches a line across Faces 3 and 0 whose endpoints are $v_1$ and $v_4$.
    However, this line does not exist on the cut locus of $v_1$.
    Then if $(q_n)_n$ is on this line, $(p,q)\in E_5$ (if $q$ is $v_1$, $v_4$, or the midpoint of an octahedron edge) or $(p,q)\in E_1$ (otherwise).
\end{compactitem}
Therefore, there cannot be a discontinuity witnessed by $(p_n)_n\subseteq c_i$.

Assume $(p_n)_n\subseteq e_i$.
There are four cut locus lines we must consider, and the limit $p$ is one of two octahedron vertices $v_j$, $v_k$. 
Since other $e_i$ are symmetric, we assume $e_i=e_1$ and $p=v_1$ (we need not check both endpoints of $e_1$).
This case is visualized in Figures \ref{fig:octae1locus}, \ref{fig:octae1unfold}(a), and \ref{fig:octa_gmpr_E2_examples}(a).
\\
As $(p_n)_n$ approaches $v_1$, cut locus lines $\ell^{\{1,2\}}$, $\ell^{\{1,5\}}$, and $\ell^{\{4,5\}}$ approach three of the cut locus lines of $v_1$.
The endpoints of these lines on their respective octahedron vertices are fixed, and their other endpoint becomes $v_4$ (the antipode of $v_1$).
Then as before, the clockwise side of each of these lines is consistent with their limit, and $(p_n,q_n)_n$ cannot witness a discontinuity if $(q_n)_n$ is chosen on any of these lines.
\\
The remaining line $\ell^{\{2,4\}}$ approaches a line that starts at $v_1$, ends at $v_4$, consists of two octahedron edges, and includes the remaining cut locus line of $v_1$ (equivalently, this line extends the remaining cut locus line of $v_1$ by one more octahedron edge).
For $p_n\in e_1$ we use the endpoint that does \textit{not} approach $v_4$ to determine the clockwise side of this line.
When $p=v_1$, we use the endpoint that is not $v_4$ to determine the clockwise side, and thus the part of this line that approaches the cut locus line of $v_1$ has a clockwise side consistent with its limit.
Thus, if the limit $q$ is on this cut locus line, $(p_n,q_n)_n$ cannot witness a discontinuity.
Otherwise, $(p,q)$ must either be in $E_5$ (if $q$ is an octahedron vertex or a midpoint of any octahedron edge) or in $E_1$ (in any other case).
Thus, there cannot be a discontinuity witnessed by $(p_n,q_n)_n\subseteq e_i$.

Assume $(p_n)_n\subseteq \tau_i$.
There are five cut locus lines to consider, and $p$ is on some $o_j$, $v_j$, $e_j$, $a_j$, or $c_j$.
Since other $\tau_i$ are symmetric, we assume $\tau_i=\tau_1$, $v_j=v_1$, $e_j=e_1$, $a_j=a_1$, $c_j=c_1$, and $o_j=o_6$.
This case is visualized in Figures \ref{fig:octatau1locus} and \ref{fig:octa_gmpr_E2_examples}(f).
\begin{compactitem}
    \item
    If $p=o_6$, cut locus lines $\ell^{\{1,4\}}$, $\ell^{\{1,5\}}$, and $\ell^{\{2,4\}}$ approach the antipode $o_0$ as $(p_n)_n$ approaches $o_6$.
    Thus, $(q_n)_n$ chosen on any of these lines results in $q=o_0$ and $(p,q)\in E_5$.
    
    The two other cut locus lines approach lines of the same name as $(p_n)_n$ approaches $o_6$.
    Their respective endpoints on octahedron vertices are fixed, so the clockwise side of each cut locus line is consistent with the clockwise side of their limits.
    Then, as before, $(p_n,q_n)_n$ cannot witness a discontinuity if $(q_n)_n$ is chosen on any of these lines.
    \item
    If $p=v_1$, cut locus lines $\ell^{\{1,4\}}$ and $\ell^{\{1,5\}}$ approach the antipode $v_4$ as $(p_n)_n$ approaches $v_1$.
    If $(q_n)_n$ were chosen on any of these lines, $(p,q)=(v_1,v_4)\in E_5$.
    
    As $(p_n)_n$ approaches $v_1$, the lines $\ell^{\{1,2\}}$ and $\ell^{\{4,5\}}$ approach cut locus lines of $v_1$.
    The endpoints of these lines on their respective octahedron vertices are fixed, and each line's other endpoint approaches $v_4$.
    Thus, the clockwise side of each cut locus line is consistent with the clockwise side of their limits, and $(p_n,q_n)_n$ cannot witness a discontinuity if $(q_n)_n$ is chosen on these lines.
    
    Finally, any $(q_n)_n$ chosen on cut locus line $\ell^{\{2,4\}}$ approaches the bottom edge of Face 0, which is a cut locus line of $v_1$.
    From construction of $\sigma_2$, the chosen side of $\ell^{\{2,4\}}$ is consistent with the clockwise side of the bottom edge of Face 0.
    Then $(p_n,q_n)_n$ cannot witness a discontinuity if $(q_n)_n$ is chosen on this line.
    \item
    If $p\in e_1$, cut locus line $\ell^{\{1,4\}}$ approaches the vertex of the cut locus as $(p_n)_n$ approaches $p$.
    Then if $(q_n)_n$ is chosen on this line, $q$ is this vertex, and $(p,q)\notin E_2$.
    
    The lines $\ell^{\{1,2\}}$ and $\ell^{\{4,5\}}$ continuously approach the lines of the same name as $(p_n)_n$ approaches $p$, keeping their respective endpoints at an octahedron vertex constant.
    Then, as before, the clockwise sides of these lines are consistent with their limits, and $(p_n,q_n)_n$ cannot witness a discontinuity if $(q_n)_n$ is chosen on these lines.
    
    For $(q_n)_n$ chosen on the two remaining lines, $q$ must be on the bottom edge of Face 0, which is a cut locus line of $p$.
    Specifically, $q$ must be on the cut locus line with the same name as those that $(q_n)_n$ are on.
    As before, from construction of $\sigma_2$, the chosen side of each cut locus line is consistent with the clockwise side of the line $q$ is on, and $(p_n,q_n)_n$ cannot witness a discontinuity if $(q_n)_n$ is chosen on these lines.
    
    \item
    If $p\in a_1$, the cut locus line $\ell^{\{1,4\}}$ approaches a cut locus vertex as $(p_n)_n$ approaches $p$.
    Thus, if $(q_n)_n$ were on this line, $q$ is this vertex and $(p,q)\notin E_2$.
    
    The lines $\ell^{\{1,2\}}$ and $\ell^{\{4,5\}}$ approach the lines of the same name in \autoref{fig:octaa1locus}, with their respective endpoints at octahedron vertices fixed.
    Then $(p_n,q_n)_n$ cannot witness a discontinuity with $(q_n)_n$ chosen on these lines.

    The lines $\ell^{\{1,5\}}$ and $\ell^{\{2,4\}}$ approach the lines of the same name in \autoref{fig:octaa1locus}.
    From construction of $\sigma_2$, our choice of direction to approach this line is consistent with its limit, so $(p_n,q_n)_n$ cannot witness a discontinuity with $(q_n)_n$ chosen on this line.
    
    \item
    If $p\in c_1$, the cut locus lines $\ell^{\{1,5\}}$ and $\ell^{\{2,4\}}$ approach vertices of the cut locus as $(p_n)_n$ approaches $p$.
    If $(q_n)_n$ were on either line, then $q$ is a cut locus vertex and $(p,q)\notin E_2$.
    
    The line $\ell^{\{1,4\}}$ continuously approaches the line of the same name as $(p_n)_n$ approaches $p$, and by construction of $\sigma_2$, our choice of direction to approach this line is consistent with its limit, so $(p_n,q_n)_n$ cannot witness a discontinuity with $(q_n)_n$ chosen on this line.
    
    The cut locus lines $\ell^{\{1,2\}}$ and $\ell^{\{4,5\}}$ approach the lines of the same name as $(p_n)_n$ approaches $p$, with their respective endpoints at octahedron vertices fixed.
    Thus, the clockwise side of these lines is consistent with their limits, and $(p_n,q_n)_n$ cannot witness a discontinuity with $(q_n)_n$ is chosen on these lines.
\end{compactitem}
Therefore, there cannot be a discontinuity witnessed by $(p_n)_n\subseteq \tau_i$.

Finally, if $(p_n)_n\subseteq\{v_i\}$, then $p=v_i$, so $(p_n,q_n)_n$ cannot witness a discontinuity.
Similarly, $(p_n)_n\subseteq\{o_i\}$ results in $p=o_i$.

Thus, $GC(X)\leq 4$, needing at most five sets in a minimal geodesic motion planner.
With the lower bound in \autoref{subsubsec:octa_lower_bound_proof}, this implies that the geodesic complexity of an octahedron is four, and a minimal geodesic motion planner on this space requires a partition of exactly five sets.

\section*{Acknowledgments}
The authors are grateful for feedback from Donald Davis on early versions of this manuscript.

\bibliographystyle{plain} 
\bibliography{refs}
\clearpage

\appendix
\section{Tetrahedron}
\label{apx:tetrapendix}

\renewcommand{\imgwidth}{.31\linewidth}
\begin{figure}[hb!]
    \centering
    \subfigure[$t=0$]{
    \includegraphics[width=\imgwidth]{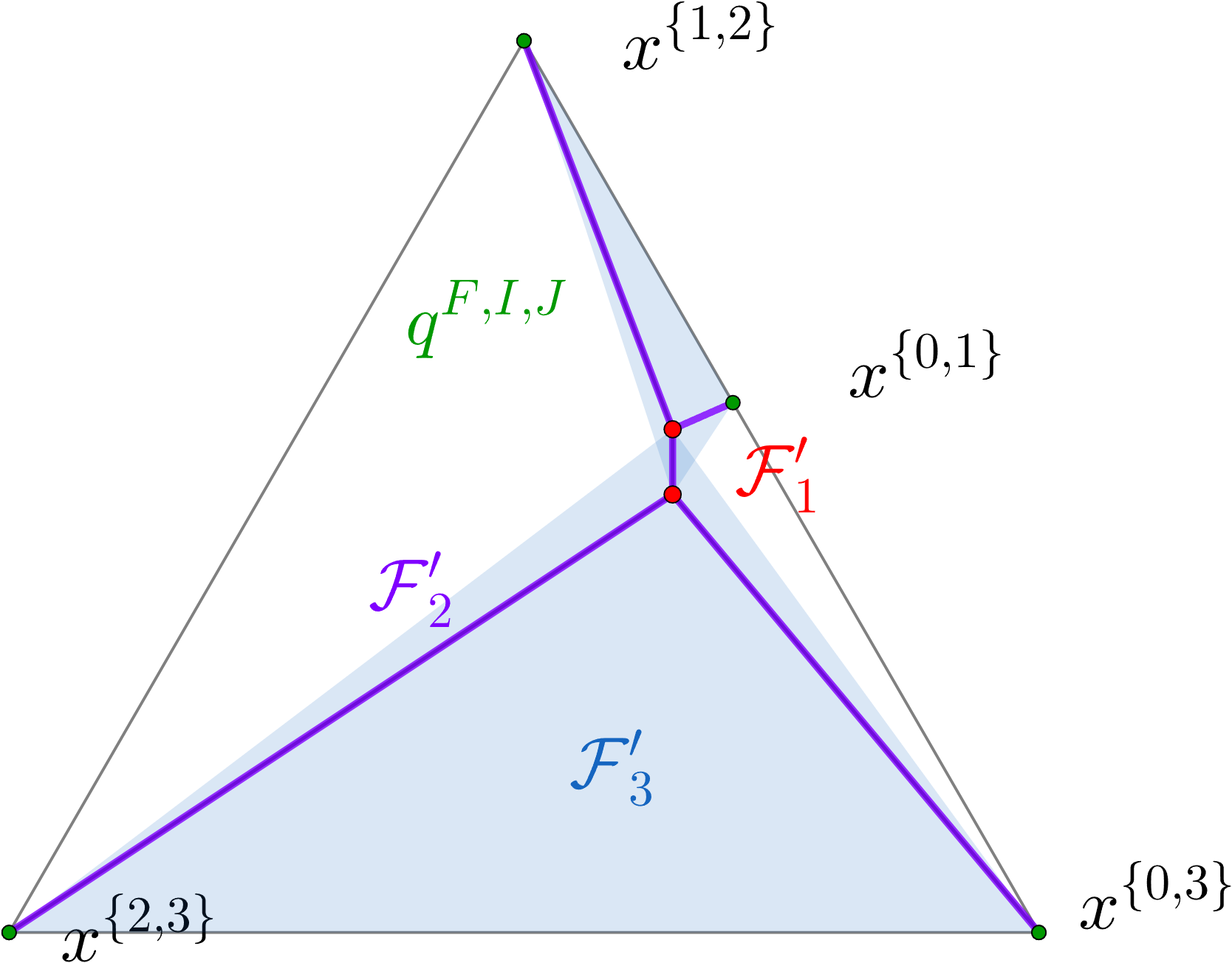}
    }
    \subfigure[$t=0.5$]{
    \includegraphics[width=\imgwidth]{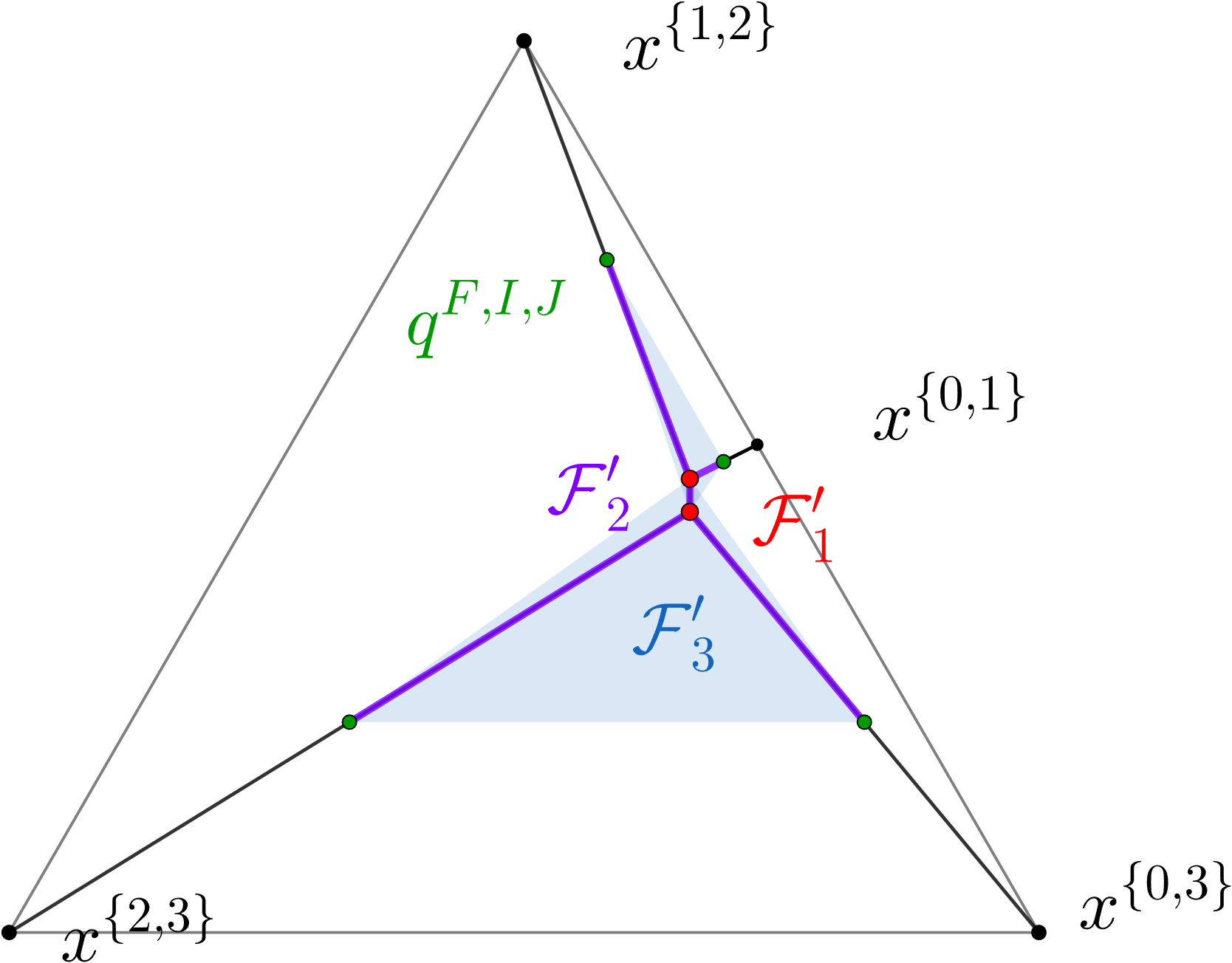}
    }
    \subfigure[$t=1$]{
    \includegraphics[width=\imgwidth]{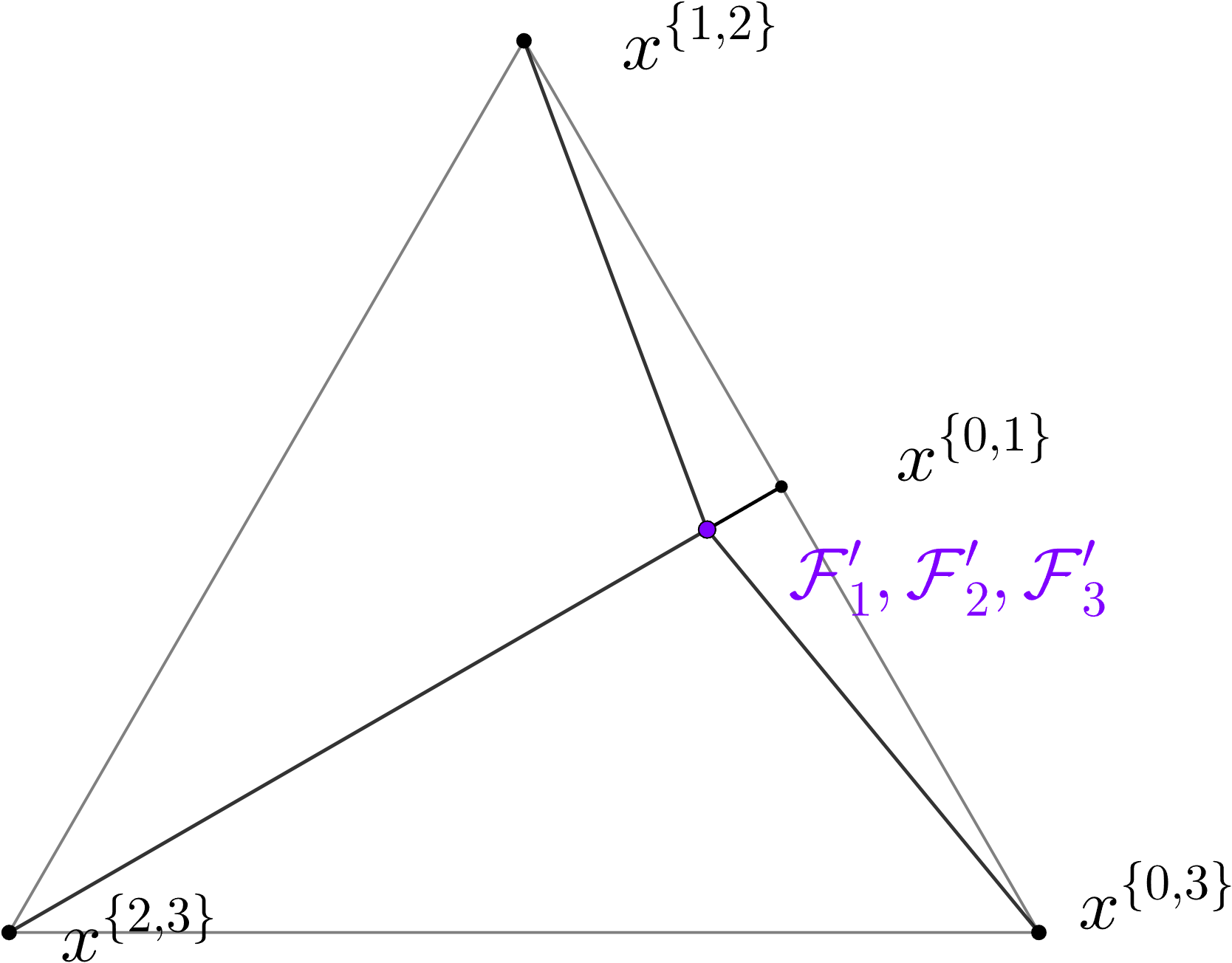}
    }
    \caption{Visualization of $\mathcal F_i$ as $p$ varies towards $p^*$ along Face 3}
    \label{fig:tetra_fullproof}
\end{figure}

We would like to visualize our method's application to the tetrahedron by displaying our sets $\mathcal F_i$ in a way that imitates \autoref{fig:torus_proof}(a).
Since some embedded simplices are 3-dimensional, we use the same method as in \autoref{fig:tetra_proof} and plot a few slices of one dimension.
For $i>0$, each $F_i\in \mathcal F_i$ is a map $\Delta_i\to X\times X$.
Recall from \autoref{def:simplex_shenanigans} that $T_1=[0,1]$ and for $t\in T_1$, $D_{t,i}=\{x\in \Delta_{i}:x_0=t\}$.
By \autoref{lem:nice_embeddings}, we have that the image $\pi_2(D_{t,i})$ is the convex hull of $i$ points.
Let this hull be $F_i'$, and let $\mathcal F_i'$ be the union $\bigcup\limits_{F_i\in\mathcal F_i}F_i'$.
In \autoref{fig:tetra_fullproof}, we plot $\mathcal F_i'$ for a few values of $t$ for cases where the first factor embedding is $p_3$.

\renewcommand{\imgwidth}{.23\linewidth}\begin{figure}[hb!]
\centering
\subfigure[Figure 2.1 of \cite{tetra}]{
\includegraphics[width=\imgwidth]{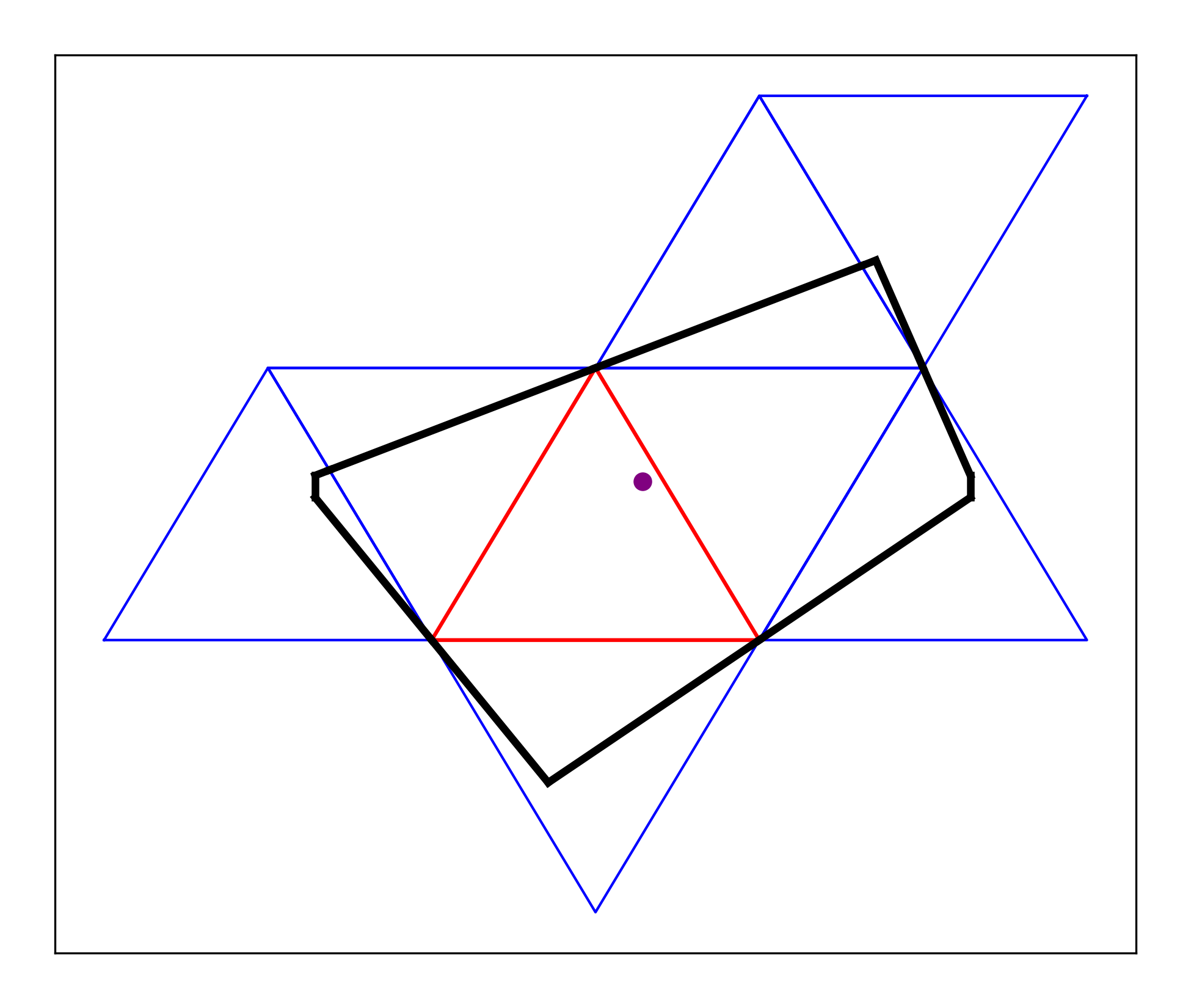}}
\subfigure[Figure 2.2 of \cite{tetra}]{
\includegraphics[width=\imgwidth]{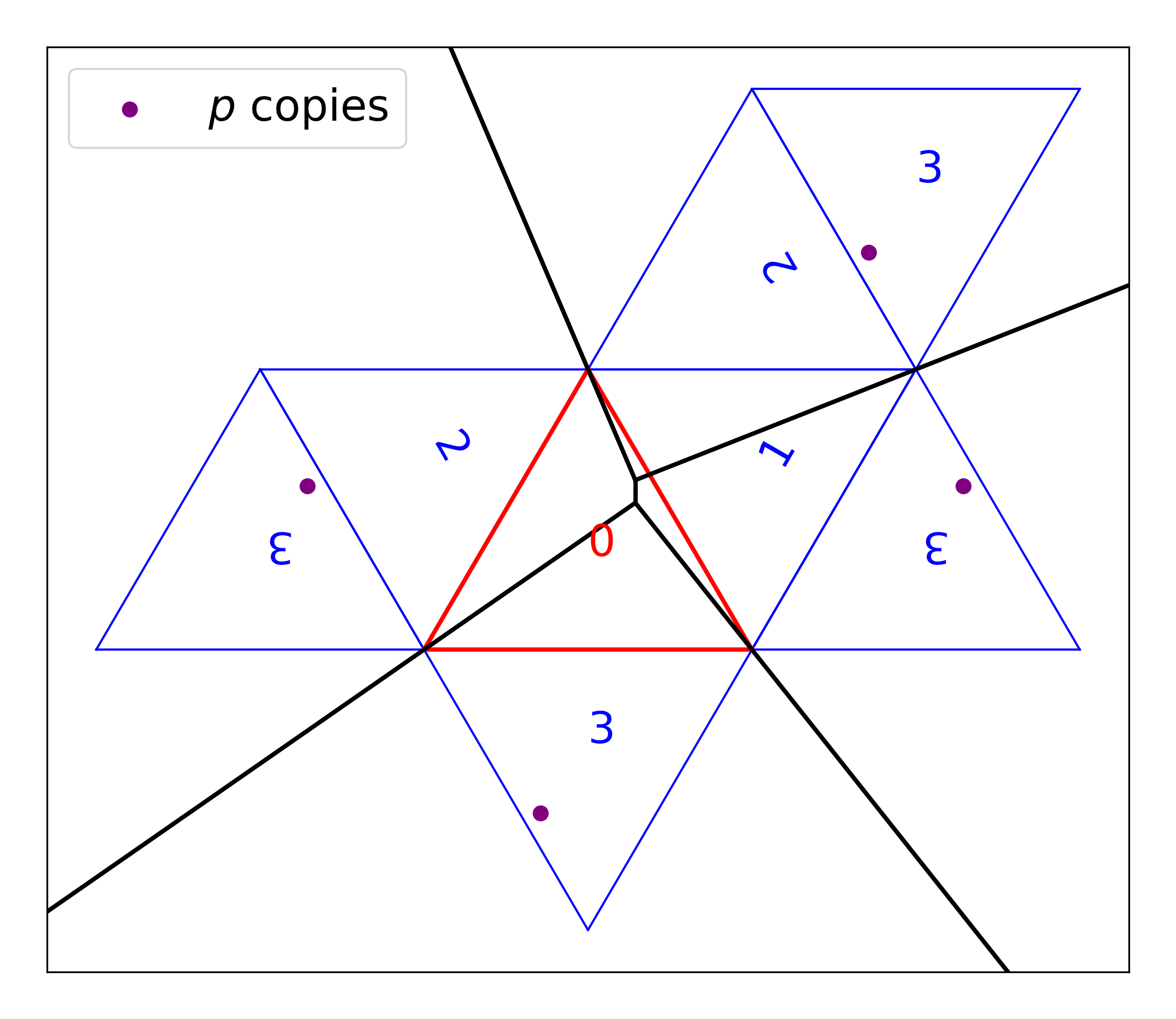}
}
\subfigure[Figure 2.5 (left) of \cite{tetra}]{
\includegraphics[width=\imgwidth]{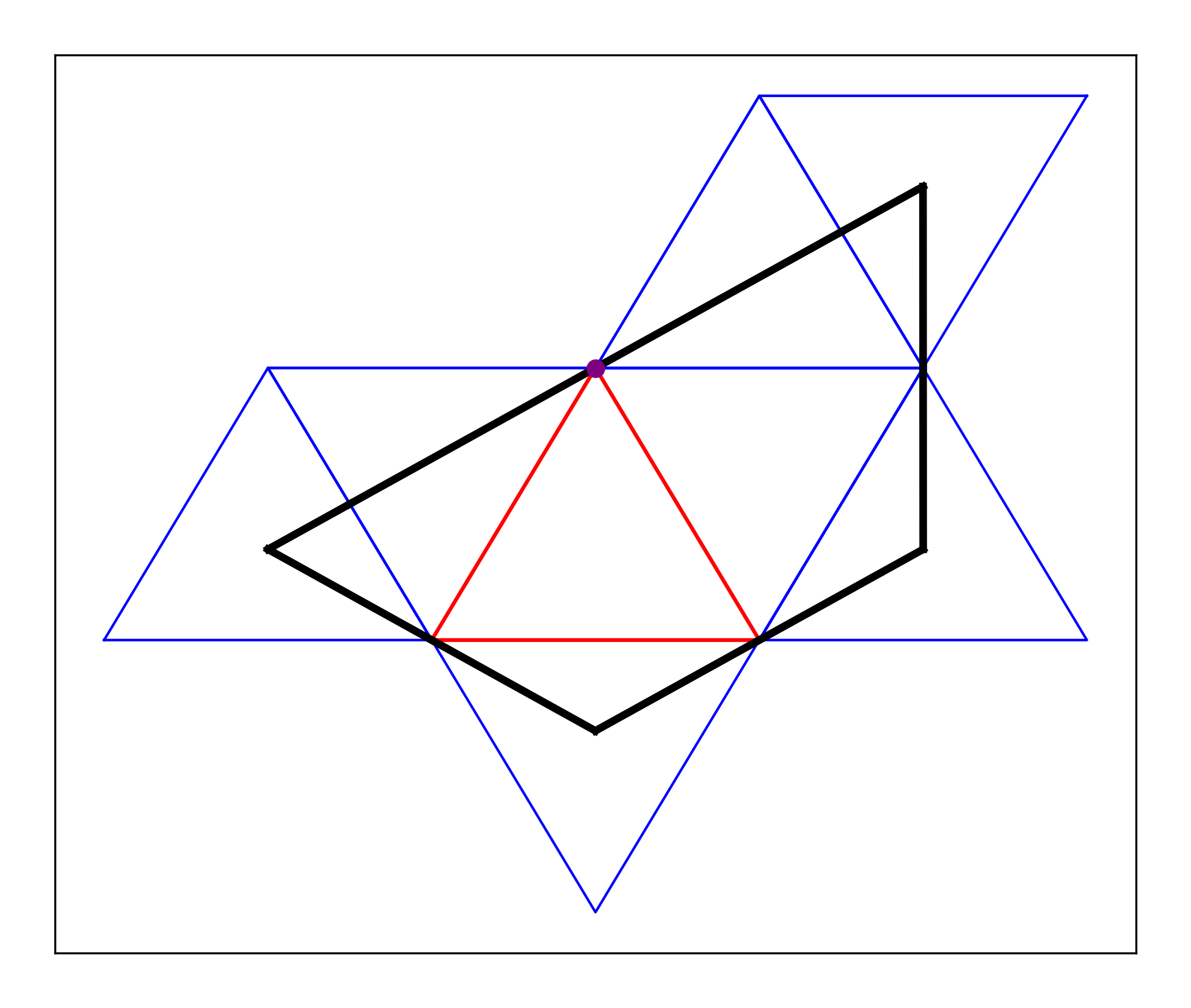}
}
\subfigure[Figure 2.5 (right) of \cite{tetra}]{
\includegraphics[width=\imgwidth]{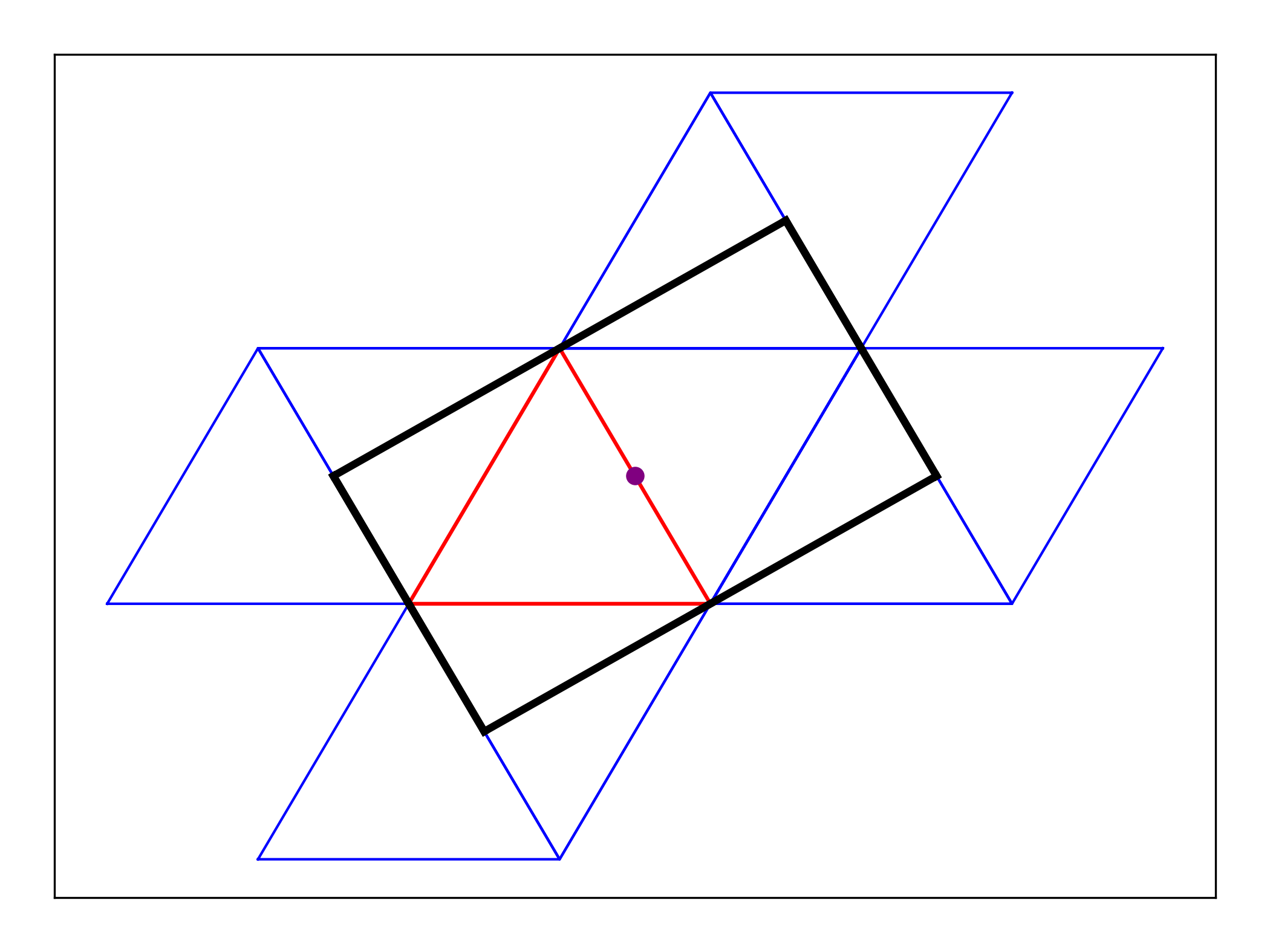}
}
\subfigure[Figure 2.6 of \cite{tetra}]{
\includegraphics[width=\imgwidth]{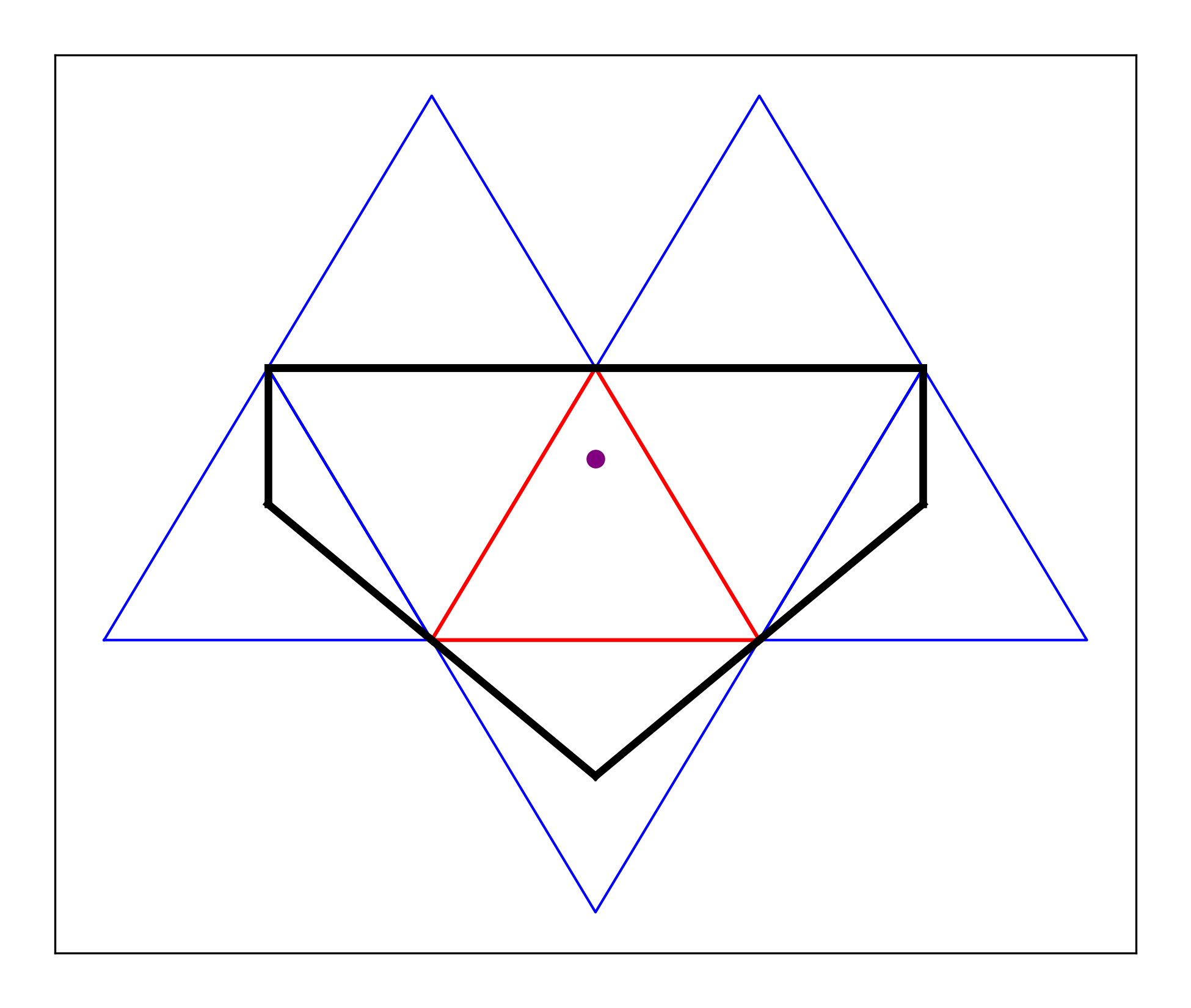}
}
\subfigure[Figure 2.7 of \cite{tetra}]{
\includegraphics[width=\imgwidth]{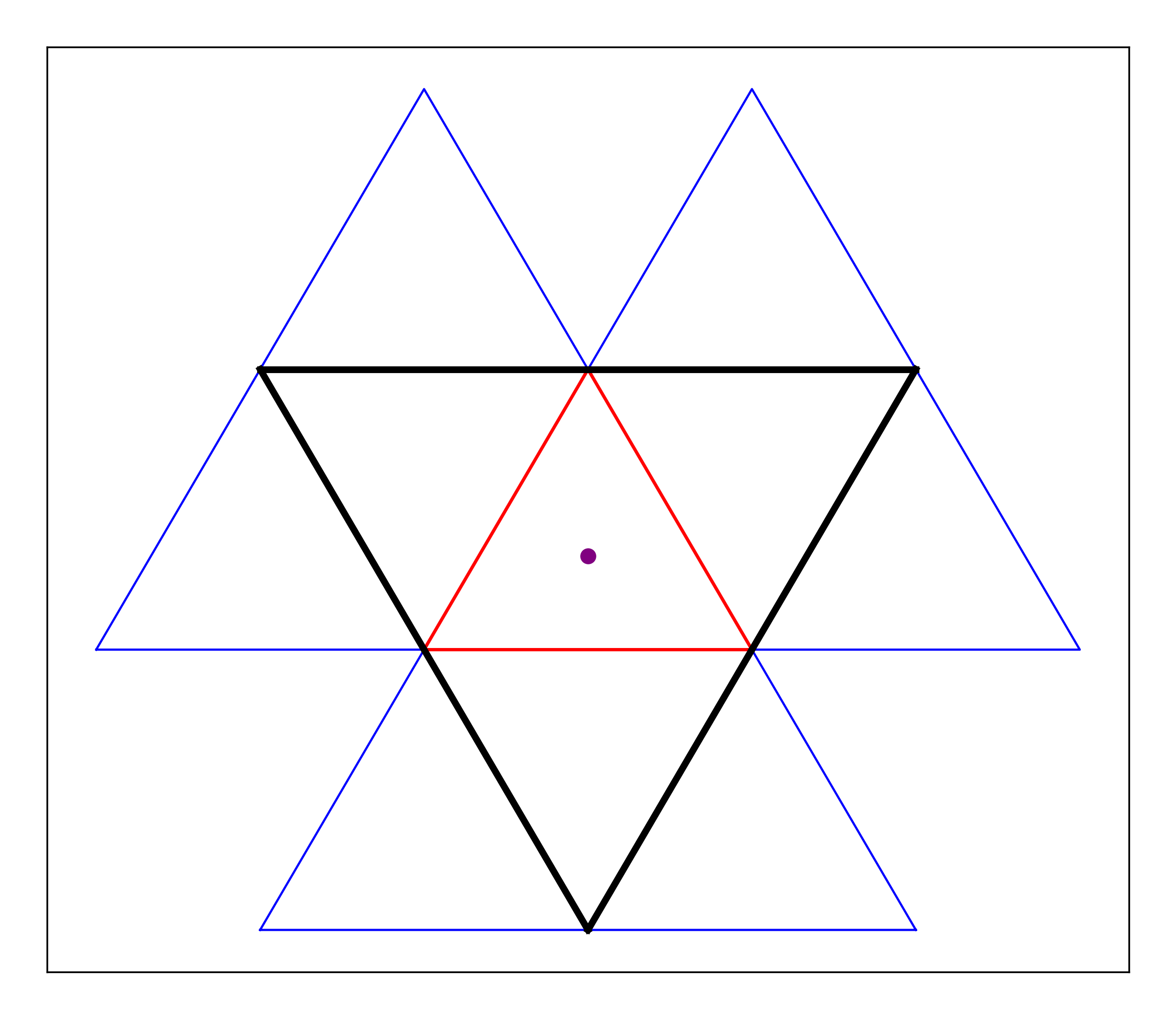}
}
\subfigure[Figure 2.8 of \cite{tetra}]{
\includegraphics[width=\imgwidth]{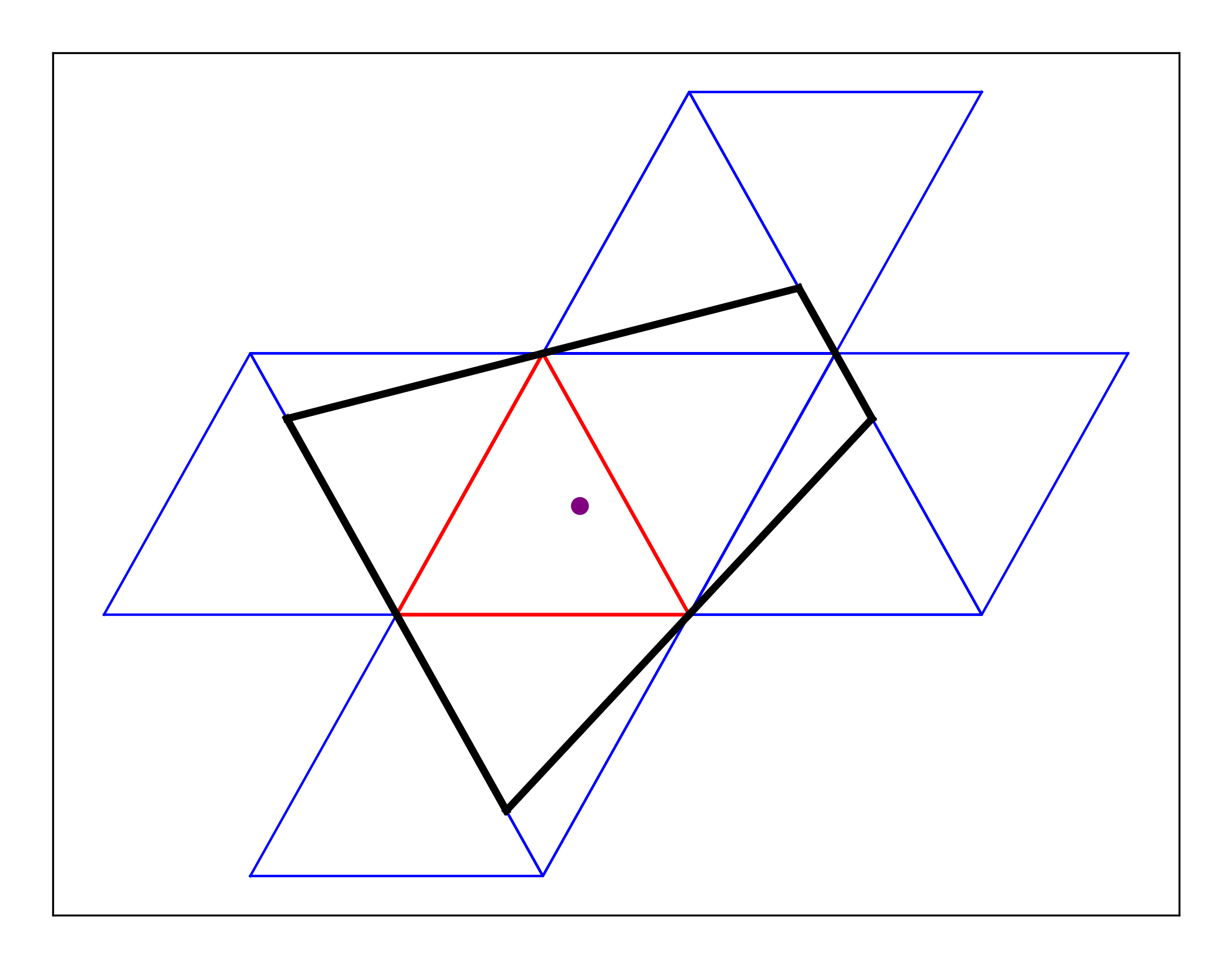}
}
\caption{Recreation of figures in~\cite[Sec.~2]{tetra} using \autoref{alg:cutlocus}}
\end{figure}

\section{Octahedron}

We will find the equations of the lines $\ell^{\{i,j\}}$ for choices of $p$ on Face 6 of the octahedron.
These equations are well defined for $p$ on all of Face 6 other than its vertices.

\begin{compactitem}
    \item
    $\ell^{\{ 0 , 1 \}}(p) := t \mapsto  \left(\frac{p_1}{4} - \sqrt{3} \frac{p_2}{4} + 5 \frac{\sqrt{3}}{2}, \sqrt{3} \frac{p_1}{4} + \frac{p_2}{4} + \frac{5}{2}\right) + \left(-\sqrt{3} \frac{p_1}{4} + 3 \frac{p_2}{4} - \frac{3}{2}, -3 \frac{p_1}{4} - \sqrt{3} \frac{p_2}{4} - \frac{\sqrt{3}}{2}\right) t$.
    \item
    $\ell^{\{ 0 , 2 \}}(p) := t \mapsto  \left(\frac{p_1}{4} + \sqrt{3} \frac{p_2}{4} + \frac{\sqrt{3}}{2}, -\sqrt{3} \frac{p_1}{4} + \frac{p_2}{4} + \frac{5}{2}\right) + \left(\sqrt{3} \frac{p_1}{4} + 3 \frac{p_2}{4} - \frac{3}{2}, -3 \frac{p_1}{4} + \sqrt{3} \frac{p_2}{4} - 5 \frac{\sqrt{3}}{2}\right) t$.
    \item
    $\ell^{\{ 0 , 3 \}}(p) := t \mapsto  \left(p_1, p_2 + 1\right) + \left(0, -3 \sqrt{3}\right) t$.
    \item
    $\ell^{\{ 0 , 4 \}}(p) := t \mapsto  \left(\frac{p_1}{4} - \sqrt{3} \frac{p_2}{4} + \sqrt{3}, \sqrt{3} \frac{p_1}{4} + \frac{p_2}{4} - 2\right) + \left(-\sqrt{3} \frac{p_1}{4} + 3 \frac{p_2}{4} + 3, -3 \frac{p_1}{4} - \sqrt{3} \frac{p_2}{4} - 2 \sqrt{3}\right) t$.
    \item
    $\ell^{\{ 0 , 5 \}}(p) := t \mapsto  \left(\frac{p_1}{4} + \sqrt{3} \frac{p_2}{4} + 2 \sqrt{3}, -\sqrt{3} \frac{p_1}{4} + \frac{p_2}{4} - 2\right) + \left(\sqrt{3} \frac{p_1}{4} + 3 \frac{p_2}{4} + 3, -3 \frac{p_1}{4} + \sqrt{3} \frac{p_2}{4} - \sqrt{3}\right) t$.
    \item
    $\ell^{\{ 1 , 2 \}}(p) := t \mapsto  \left(\frac{-p_1}{2}, 4 - \frac{p_2}{2}\right) + \left(\sqrt{3} \frac{p_1}{2}, \sqrt{3} \frac{p_2 - 4}{2}\right) t$.
    \item
    $\ell^{\{ 1 , 3 \}}(p) := t \mapsto  \left(\frac{p_1}{4} - \sqrt{3} \frac{p_2}{4} - \frac{\sqrt{3}}{2}, \sqrt{3} \frac{p_1}{4} + \frac{p_2}{4} + \frac{5}{2}\right) + \left(\sqrt{3} \frac{p_1}{4} - 3 \frac{p_2}{4} + \frac{3}{2}, 3 \frac{p_1}{4} + \sqrt{3} \frac{p_2}{4} - 5 \frac{\sqrt{3}}{2}\right) t$.
    \item
    $\ell^{\{ 1 , 4 \}}(p) := t \mapsto  \left(\frac{-p_1}{2} - \sqrt{3} \frac{p_2}{2} + \frac{\sqrt{3}}{2}, \sqrt{3} \frac{p_1}{2} - \frac{p_2}{2} - \frac{1}{2}\right) + \left(\frac{9}{2}, -3 \frac{\sqrt{3}}{2}\right) t$.
    \item
    $\ell^{\{ 1 , 5 \}}(p) := t \mapsto  \left(\frac{-p_1}{2} + 3 \frac{\sqrt{3}}{2}, \frac{-p_2}{2} - \frac{1}{2}\right) + \left(\sqrt{3} \frac{p_1}{2} + \frac{9}{2}, \sqrt{3} \frac{p_2 - 1}{2}\right) t$.
    \item
    $\ell^{\{ 2 , 3 \}}(p) := t \mapsto  \left(\frac{p_1}{4} + \sqrt{3} \frac{p_2}{4} - 5 \frac{\sqrt{3}}{2}, -\sqrt{3} \frac{p_1}{4} + \frac{p_2}{4} + \frac{5}{2}\right) + \left(-\sqrt{3} \frac{p_1}{4} - 3 \frac{p_2}{4} + \frac{3}{2}, 3 \frac{p_1}{4} - \sqrt{3} \frac{p_2}{4} - \frac{\sqrt{3}}{2}\right) t$.
    \item
    $\ell^{\{ 2 , 4 \}}(p) := t \mapsto  \left(\frac{-p_1}{2} - 3 \frac{\sqrt{3}}{2}, \frac{-p_2}{2} - \frac{1}{2}\right) + \left(-\sqrt{3} \frac{p_1}{2} + \frac{9}{2}, \sqrt{3} \frac{1 - p_2}{2}\right) t$.
    \item
    $\ell^{\{ 2 , 5 \}}(p) := t \mapsto  \left(\frac{-p_1}{2} + \sqrt{3} \frac{p_2}{2} - \frac{\sqrt{3}}{2}, -\sqrt{3} \frac{p_1}{2} - \frac{p_2}{2} - \frac{1}{2}\right) + \left(\frac{9}{2}, 3 \frac{\sqrt{3}}{2}\right) t$.
    \item
    $\ell^{\{ 3 , 4 \}}(p) := t \mapsto  \left(\frac{p_1}{4} - \sqrt{3} \frac{p_2}{4} - 2 \sqrt{3}, \sqrt{3} \frac{p_1}{4} + \frac{p_2}{4} - 2\right) + \left(-\sqrt{3} \frac{p_1}{4} + 3 \frac{p_2}{4} + 3, -3 \frac{p_1}{4} - \sqrt{3} \frac{p_2}{4} + \sqrt{3}\right) t$.
    \item
    $\ell^{\{ 3 , 5 \}}(p) := t \mapsto  \left(\frac{p_1}{4} + \sqrt{3} \frac{p_2}{4} - \sqrt{3}, -\sqrt{3} \frac{p_1}{4} + \frac{p_2}{4} - 2\right) + \left(\sqrt{3} \frac{p_1}{4} + 3 \frac{p_2}{4} + 3, -3 \frac{p_1}{4} + \sqrt{3} \frac{p_2}{4} + 2 \sqrt{3}\right) t$.
    \item
    $\ell^{\{ 4 , 5 \}}(p) := t \mapsto  \left(\frac{-p_1}{2}, \frac{-p_2}{2} - 5\right) + \left(\sqrt{3} \frac{p_1}{2}, \sqrt{3} \frac{p_2 + 2}{2}\right) t$.

\end{compactitem}

For $j\equiv i+1\mod 6$, let $x^{\{i,j\}}(p)$ be the intersection of $\ell^{\{i,j\}}(p)$ with the boundary of Face 0. We find $x^{\{i,j\}}(p)$ for arbitrary $p$:

\begin{compactitem}
    \item     
    $x^{\{0,1\}}(p)=
    \left(3\frac{\sqrt{3}p_1 + p_2 + 2}{3p_1 - \sqrt{3}p_2 + 4\sqrt{3}}, \frac{-3p_1 - 5\sqrt{3}p_2 + 2\sqrt{3}}{3p_1 - \sqrt{3}p_2 + 4\sqrt{3}}\right)$.
    \item     
    $x^{\{1,2\}}(p)=(0,2)$.
    \item 
    $x^{\{2,3\}}(p)=
    \left(3\frac{-\sqrt{3}p_1 + p_2 + 2}{3p_1 + \sqrt{3}p_2 - 4\sqrt{3}}, \frac{-3p_1 + 5\sqrt{3}p_2 - 2\sqrt{3}}{3p_1 + \sqrt{3}p_2 - 4\sqrt{3}}\right)$.
    \item
    $x^{\{3,4\}}(p)=(-\sqrt3,-1)$.
    \item $x^{\{4.5\}}(p)=
    \left(\frac{3p_1}{p_2 + 2}, -1\right)$.
    \item $x^{\{0,5\}}(p)=(\sqrt3,-1)$.
\end{compactitem}

From the intersections of lines $\ell^{\{i,j\}}$ and $\ell^{\{j,k\}}$, we calculate the $x^{(i,j,k)}$:

\begin{compactitem}

\item $x^{\{0, 1, 2\}}(p)=
        \left(2 p_1 \frac{\sqrt{3} p_1 + p_2 + 2}{p_1^{2} + 2 \sqrt{3} p_1 + p_2^{2} - 6 p_2 + 8}, 2 \frac{p_1^{2} + \sqrt{3} p_1 p_2 - 2 \sqrt{3} p_1 + 2 p_2^{2} - 8 p_2}{p_1^{2} + 2 \sqrt{3} p_1 + p_2^{2} - 6 p_2 + 8}\right)$.

\item $x^{\{0, 1, 3\}}(p)=
        \left(p_1, \frac{3 p_1^{2} + \sqrt{3} p_1 p_2 - 2 \sqrt{3} p_1 - 12 p_2}{\sqrt{3} p_1 - 3 p_2 + 6}\right)$.

\item $x^{\{0, 1, 4\}}(p)=
        \left(\frac{p_1^{2} - 2 \sqrt{3} p_1 p_2 + 6 \sqrt{3} p_1 + 3 p_2^{2} + 6 p_2}{4 (p_1 + \sqrt{3})}, \frac{\sqrt{3} p_1^{2} - 2 p_1 p_2 - 2 p_1 - \sqrt{3} p_2^{2} - 6 \sqrt{3} p_2}{4 (p_1 + \sqrt{3})}\right)$.

\item $x^{\{0, 1, 5\}}(p)=
        \left(\frac{3 p_1^{2} + 2 \sqrt{3} p_1 p_2 + 10 \sqrt{3} p_1 + 9 p_2^{2} + 18 p_2}{\sqrt{3} p_1^{2} + 12 p_1 + \sqrt{3} p_2^{2} + 8 \sqrt{3}}, \frac{-\sqrt{3} p_1^{2} - 6 p_1 p_2 - 6 p_1 + \sqrt{3} p_2^{2} - 10 \sqrt{3} p_2}{\sqrt{3} p_1^{2} + 12 p_1 + \sqrt{3} p_2^{2} + 8 \sqrt{3}}\right)$.

\item $x^{\{0, 2, 3\}}(p)=
        \left(p_1, \frac{-3 p_1^{2} + \sqrt{3} p_1 p_2 - 2 \sqrt{3} p_1 + 12 p_2}{\sqrt{3} p_1 + 3 p_2 - 6}\right)$.

\item $x^{\{0, 2, 4\}}(p)=
        \left(\frac{3 \sqrt{3} p_1^{2} + 2 p_1 p_2 - 26 p_1 - 3 \sqrt{3} p_2^{2} - 6 \sqrt{3} p_2}{p_1^{2} + p_2^{2} - 28}, \frac{-p_1^{2} + 6 \sqrt{3} p_1 p_2 - 6 \sqrt{3} p_1 + p_2^{2} + 26 p_2}{p_1^{2} + p_2^{2} - 28}\right)$.

\item $x^{\{0, 2, 5\}}(p)=
        \left(\frac{p_1^{2} + 2 \sqrt{3} p_1 p_2 + 6 \sqrt{3} p_1 + 3 p_2^{2} + 6 p_2}{4 (p_1 + 2 \sqrt{3})}, \frac{-\sqrt{3} p_1^{2} - 2 p_1 p_2 - 2 p_1 + \sqrt{3} p_2^{2} - 6 \sqrt{3} p_2}{4 (p_1 + 2 \sqrt{3})}\right)$.

\item $x^{\{0, 3, 4\}}(p)=
        \left(p_1, \frac{-3 p_1^{2} - \sqrt{3} p_1 p_2 + 2 \sqrt{3} p_1 - 6 p_2}{-\sqrt{3} p_1 + 3 p_2 + 12}\right)$.

\item $x^{\{0, 3, 5\}}(p)=
        \left(p_1, \frac{-3 p_1^{2} + \sqrt{3} p_1 p_2 - 2 \sqrt{3} p_1 - 6 p_2}{\sqrt{3} p_1 + 3 p_2 + 12}\right)$.

\item $x^{\{0, 4, 5\}}(p)=
        \left(2 p_1 \frac{\sqrt{3} p_1 + p_2 + 8}{p_1^{2} + 2 \sqrt{3} p_1 + p_2^{2} + 6 p_2 + 8}, 2 \frac{-2 p_1^{2} + \sqrt{3} p_1 p_2 - 2 \sqrt{3} p_1 - p_2^{2} - 2 p_2}{p_1^{2} + 2 \sqrt{3} p_1 + p_2^{2} + 6 p_2 + 8}\right)$.

\item $x^{\{1, 2, 3\}}(p)=
        \left(2 p_1 \frac{-\sqrt{3} p_1 + p_2 + 2}{p_1^{2} - 2 \sqrt{3} p_1 + p_2^{2} - 6 p_2 + 8}, 2 \frac{p_1^{2} - \sqrt{3} p_1 p_2 + 2 \sqrt{3} p_1 + 2 p_2^{2} - 8 p_2}{p_1^{2} - 2 \sqrt{3} p_1 + p_2^{2} - 6 p_2 + 8}\right)$.

\item $x^{\{1, 2, 4\}}(p)=
        \left(p_1 \frac{p_1 - \sqrt{3} p_2 - 2 \sqrt{3}}{p_1 + \sqrt{3} p_2 - 4 \sqrt{3}}, \frac{p_1 p_2 - 2 p_1 - \sqrt{3} p_2^{2} + 4 \sqrt{3} p_2}{p_1 + \sqrt{3} p_2 - 4 \sqrt{3}}\right)$.

\item $x^{\{1, 2, 5\}}(p)=
        \left(p_1 \frac{p_1 + \sqrt{3} p_2 + 2 \sqrt{3}}{p_1 - \sqrt{3} p_2 + 4 \sqrt{3}}, \frac{p_1 p_2 - 2 p_1 + \sqrt{3} p_2^{2} - 4 \sqrt{3} p_2}{p_1 - \sqrt{3} p_2 + 4 \sqrt{3}}\right)$.

\item $x^{\{1, 3, 4\}}(p)=
        \left(\frac{p_1^{2} - 2 \sqrt{3} p_1 p_2 - 6 \sqrt{3} p_1 + 3 p_2^{2} + 6 p_2}{4 (p_1 - 2 \sqrt{3})}, \frac{\sqrt{3} p_1^{2} - 2 p_1 p_2 - 2 p_1 - \sqrt{3} p_2^{2} + 6 \sqrt{3} p_2}{4 (p_1 - 2 \sqrt{3})}\right)$.

\item $x^{\{1, 3, 5\}}(p)=
        \left(\frac{-3 \sqrt{3} p_1^{2} + 2 p_1 p_2 - 26 p_1 + 3 \sqrt{3} p_2^{2} + 6 \sqrt{3} p_2}{p_1^{2} + p_2^{2} - 28}, \frac{-p_1^{2} - 6 \sqrt{3} p_1 p_2 + 6 \sqrt{3} p_1 + p_2^{2} + 26 p_2}{p_1^{2} + p_2^{2} - 28}\right)$.

\item $x^{\{1, 4, 5\}}(p)=
        \left(p_1 \frac{p_1 - \sqrt{3} p_2 + 4 \sqrt{3}}{p_1 + \sqrt{3} p_2 + 2 \sqrt{3}}, \frac{p_1 p_2 - 2 p_1 - \sqrt{3} p_2^{2} - 2 \sqrt{3} p_2}{p_1 + \sqrt{3} p_2 + 2 \sqrt{3}}\right)$.

\item $x^{\{2, 3, 4\}}(p)=
        \left(\frac{-3 p_1^{2} + 2 \sqrt{3} p_1 p_2 + 10 \sqrt{3} p_1 - 9 p_2^{2} - 18 p_2}{\sqrt{3} p_1^{2} - 12 p_1 + \sqrt{3} p_2^{2} + 8 \sqrt{3}}, \frac{-\sqrt{3} p_1^{2} + 6 p_1 p_2 + 6 p_1 + \sqrt{3} p_2^{2} - 10 \sqrt{3} p_2}{\sqrt{3} p_1^{2} - 12 p_1 + \sqrt{3} p_2^{2} + 8 \sqrt{3}}\right)$.

\item $x^{\{2, 3, 5\}}(p)=
        \left(\frac{p_1^{2} + 2 \sqrt{3} p_1 p_2 - 6 \sqrt{3} p_1 + 3 p_2^{2} + 6 p_2}{4 (p_1 - \sqrt{3})}, \frac{-\sqrt{3} p_1^{2} - 2 p_1 p_2 - 2 p_1 + \sqrt{3} p_2^{2} + 6 \sqrt{3} p_2}{4 (p_1 - \sqrt{3})}\right)$.

\item $x^{\{2, 4, 5\}}(p)=
        \left(p_1 \frac{-p_1 - \sqrt{3} p_2 + 4 \sqrt{3}}{-p_1 + \sqrt{3} p_2 + 2 \sqrt{3}}, \frac{-p_1 p_2 + 2 p_1 - \sqrt{3} p_2^{2} - 2 \sqrt{3} p_2}{-p_1 + \sqrt{3} p_2 + 2 \sqrt{3}}\right)$.       

\item $x^{\{3, 4, 5\}}(p)=
        \left(2 p_1 \frac{-\sqrt{3} p_1 + p_2 + 8}{p_1^{2} - 2 \sqrt{3} p_1 + p_2^{2} + 6 p_2 + 8}, 2 \frac{-2 p_1^{2} - \sqrt{3} p_1 p_2 + 2 \sqrt{3} p_1 - p_2^{2} - 2 p_2}{p_1^{2} - 2 \sqrt{3} p_1 + p_2^{2} + 6 p_2 + 8}\right)$.
\end{compactitem}

\label{apx:octapendix}
\section{Future Work}

We believe the following are true, in order of expected difficulty:
\begin{conj}
    The geodesic complexity of an icosahedron is at least four (needing at least five sets).
    \label{conj:icosa_lower}
\end{conj}
\begin{conj}
    The geodesic complexity of a dodecahedron is at least four (needing at least five sets).
    \label{conj:dodeca_lower}
\end{conj}
\begin{conj}
    The geodesic complexity of an icosahedron is at most four (needing at most five sets).
    \label{conj:icosa_upper}
\end{conj}
\begin{conj}
    The geodesic complexity of a dodecahedron is at most four (needing at most five sets).
    \label{conj:dodeca_upper}
\end{conj}

\renewcommand{\imgwidth}{.45\linewidth}
\begin{figure}[ht!]
\centering
\subfigure[Cut locus of point on icosahedron]{
\includegraphics[width=\imgwidth]{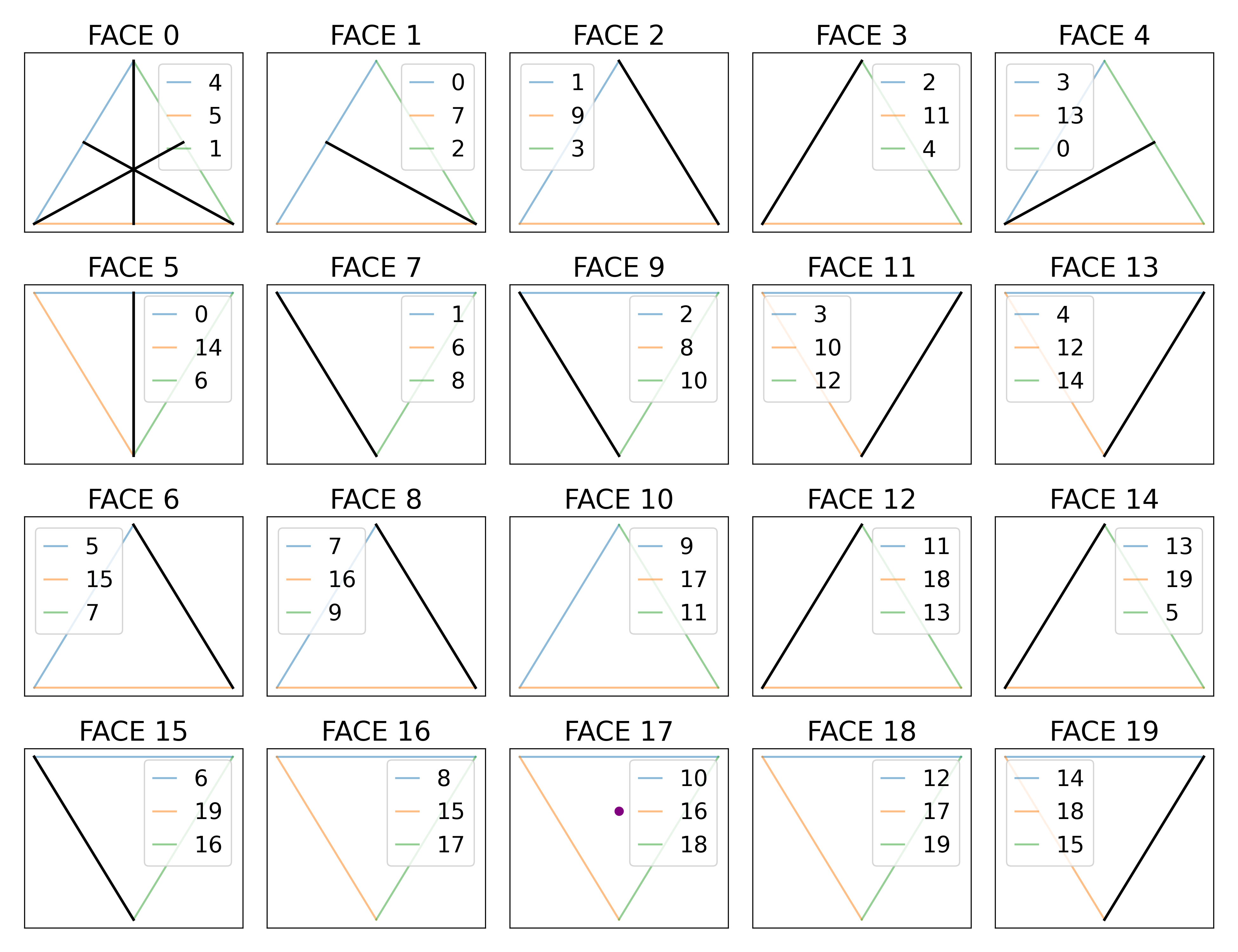}
}
\subfigure[Cut locus of point on dodecahedron]{
\includegraphics[width=\imgwidth]{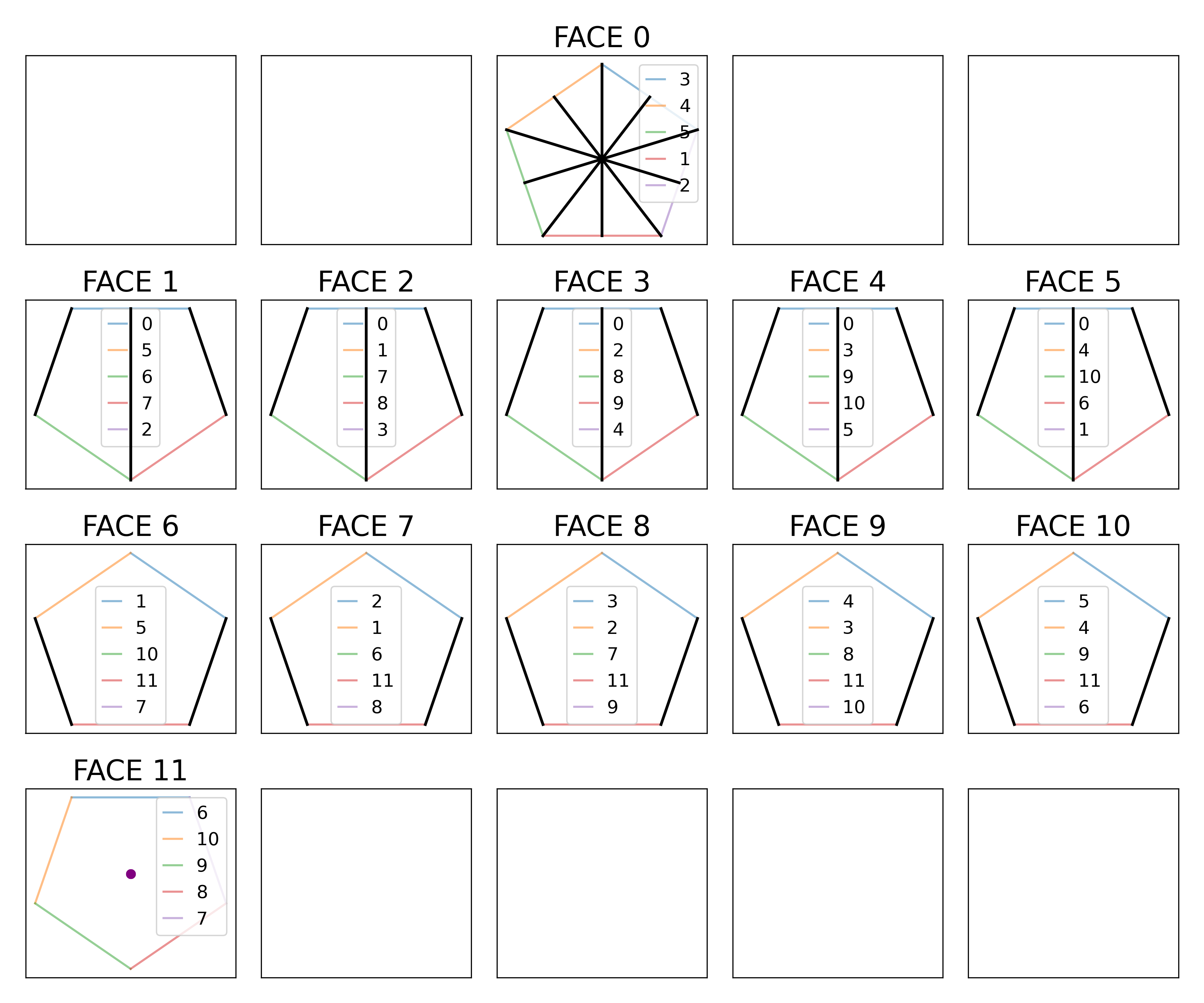}
}
\subfigure[Voronoi star diagram of point on icosahedron]{
\includegraphics[width=\imgwidth]{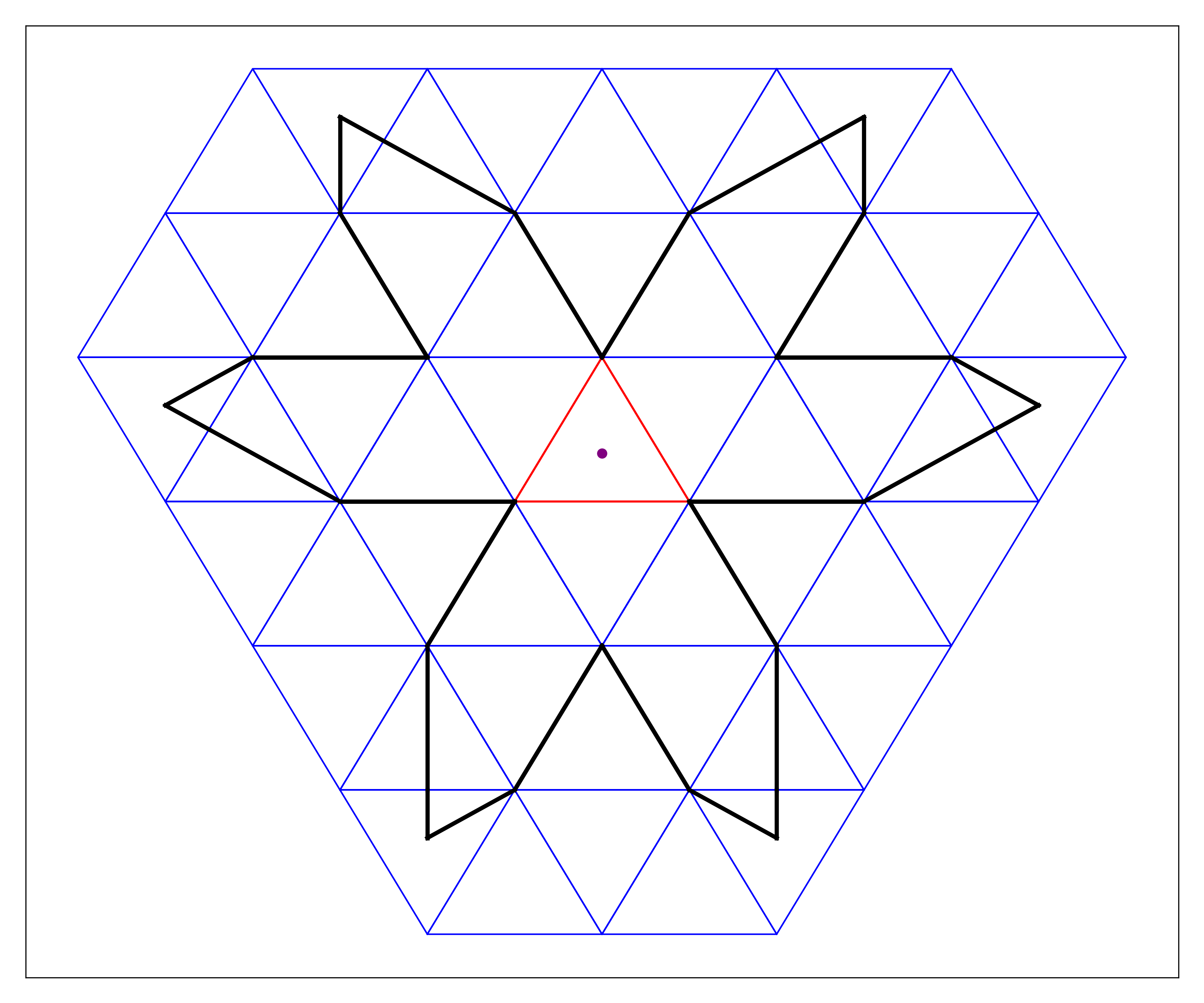}
}
\subfigure[Voronoi star diagram of point on dodecahedron]{
\includegraphics[width=\imgwidth]{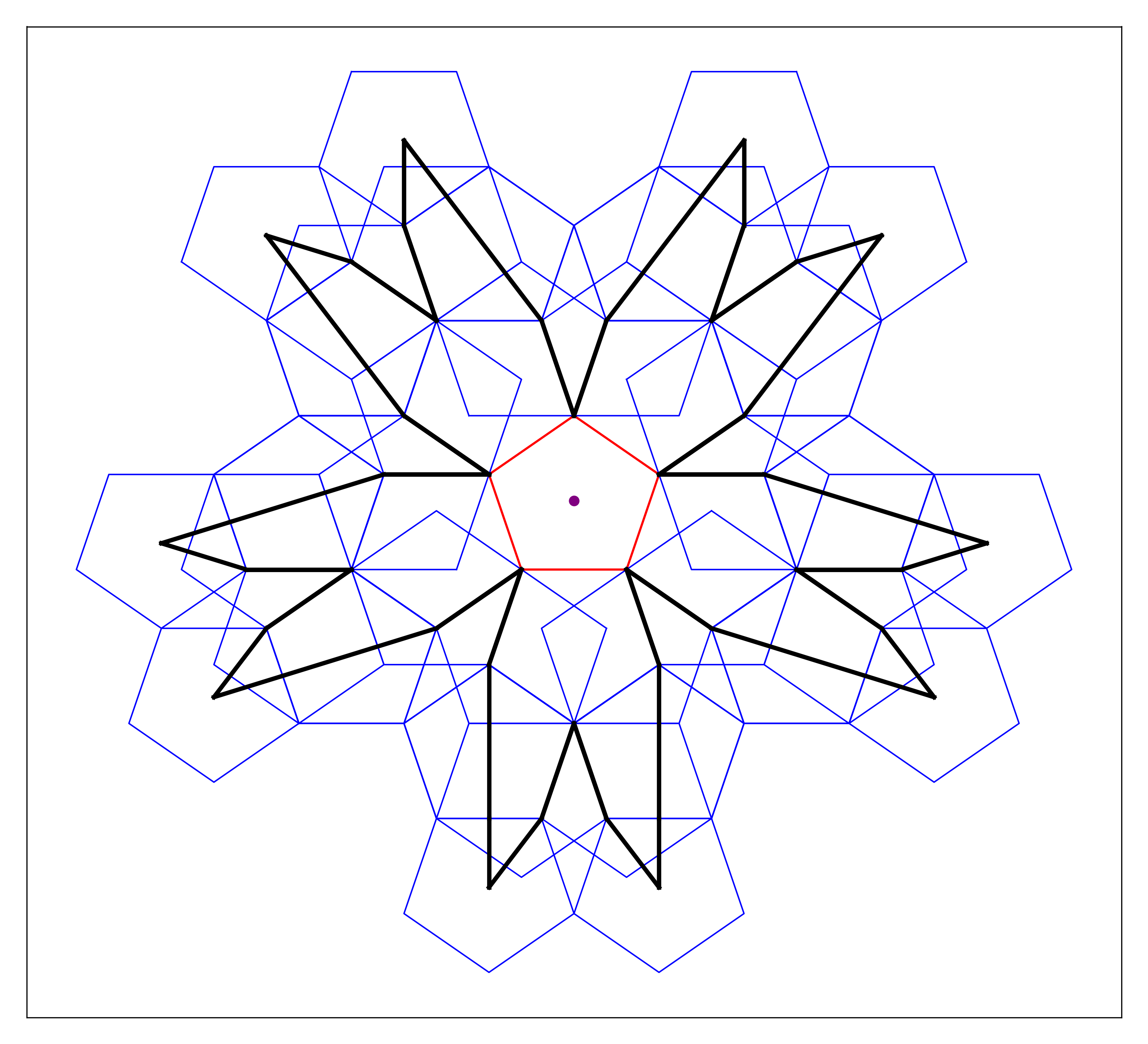}
}
\caption{Output of \autoref{alg:cutlocus} on the center of an icosahedron and dodecahedron face}
\label{fig:icosadodecazero}
\end{figure}

Our motivation for believing Conjectures \ref{conj:icosa_lower} and \ref{conj:dodeca_lower} is that the cut locus of the center of an icosahedron (resp. dodecahedron) face is a star with six (resp. ten) lines (\autoref{fig:icosadodecazero}).
We believe \autoref{thm:quite_simplex} can be used with $\mathcal F_i$ constructed similar to \autoref{subsubsec:octa_lower_bound_proof}, and we expect \autoref{conj:icosa_lower} to be easier since the cut locus structure is more similar to the octahedron.

Our motivation for Conjectures \ref{conj:icosa_upper} and \ref{conj:dodeca_upper} is that the tetrahedron \cite{tetra}, cube \cite{cube}, and now the octahedron have geodesic complexity at most four.
Increasing complexities of the cut loci on these spaces have required explicit geodesic motion planners to be made with more care.
However, we do not notice differences with cut loci on the icosahedron and dodecahedron that would suggest this construction becomes impossible.

\end{document}